\documentclass[reqno,10pt]{amsart}
\usepackage{bm}
\usepackage{amsmath}
\usepackage{}
\usepackage{a4wide}
\usepackage{mathbbol}
\usepackage{mathrsfs}
\usepackage{amssymb,mathrsfs,graphicx,extpfeil}
\usepackage{epsfig}
\usepackage{indentfirst, latexsym, amssymb, enumerate,amsmath,graphicx}
\usepackage{float}
\usepackage{colortbl}
\usepackage{epsfig,subfigure}

\usepackage{color}

\usepackage{hyperref}
\allowdisplaybreaks[4]
\usepackage{refcheck}
\usepackage{amsmath}
\title[Classical limit of the relativistic Vlasov-Maxwell-Landau system]{Classical limit of the relativistic Vlasov-Maxwell-Landau system }

\author[C. Q. Cao]{Chuqi Cao}
\address[C. Q. Cao]{{\newline Department of  Applied Mathematics, The Hong Kong Polytechnic University, Hong Kong, China}}
\email{chuqicao@gmail.com}

\author[L.-B. He]{Ling-Bing He}
\address[L.-B. He]{{\newline Department of Mathematical Sciences, Tsinghua University, Beijing 100084, China}}
\email{hlb@tsinghua.edu.cn}

\author[Y. J. Lei]{Yuanjie Lei}
\address[Y. J. Lei]{{\newline School of Mathematics and Statistics, and Hubei Key Laboratory of Engineering Modeling and Scientific Computing, Huazhong University of Science and Technology, Wuhan, 430074, China}}
\email{leiyuanjie@hust.edu.cn}

\author[Q. H. Xiao]{Qinghua Xiao}
\address[Q. H. Xiao]{\newline Innovation Academy for Precision Measurement Science and Technology, Chinese Academy of Sciences, Wuhan 430071, China}
\email{xiaoqh@apm.ac.cn}

\newtheorem{theorem}{Theorem}[section]
\newtheorem{lemma}{Lemma}[section]
\newtheorem{corollary}{Corollary}[section]
\newtheorem{proposition}{Proposition}[section]

\newtheorem{remark}{Remark}[section]

\def\charf {\mbox{{\text 1}\kern-.30em {\text l}}}

\newcommand{\R}{\mathbb{R}}

\newcommand{\pa}{\partial}

\newcommand{\lag}{\langle}
\newcommand{\rag}{\rangle}

\numberwithin{equation}{section}

\begin{document}

\date{\today}




\begin{abstract} The physical essence of the non-relativistic limit, from the relativistic Vlasov-Maxwell-Landau system to the Vlasov-Poisson-Landau system, lies in the transition from finite-speed electromagnetic waves to instantaneous Coulomb interactions, and from relativistic to Newtonian particle dynamics. We rigorously justify this limit (mathematically corresponding to the light speed $c \to \infty$) in a periodic box via three key technical advances: establishing a uniform-in-$c$ coercivity estimate for the relativistic Landau collision operator, constructing a novel weighted energy functional to overcome the weakening dissipation of the electromagnetic field at large $c$, and proving a corresponding global well-posedness result.
\end{abstract}
\maketitle


 \setcounter{tocdepth}{1}
\tableofcontents

\section{Introduction}

\setcounter{equation}{0}

A plasma, often described as the fourth state of matter, is a gaseous medium of ions and electrons that constitutes more than 90\% of the visible universe. Its ubiquity and unique properties make it a cornerstone of modern physics. The study of plasmas is not only central to groundbreaking technological endeavors, such as controlled thermonuclear fusion, but also fundamental to our understanding of diverse phenomena in astrophysics, space science, and materials processing. To describe the collective behavior of such a complex, many-body system, the framework of kinetic theory provides the most detailed description, governing the evolution of the particle distribution function in phase space.

Among kinetic models, the relativistic Vlasov-Maxwell-Landau (RVML) system stands as a complete and self-consistent description for a dilute, collisional plasma. It masterfully integrates several physical layers: the Vlasov equation captures the advection of particles in self-consistent electromagnetic fields, the Maxwell equations govern the dynamics of these fields at finite speed $c$ (the speed of light), and the Landau collision operator accounts for the binary interactions between charged particles. This synthesis makes the RVML system the model of choice for high-energy scenarios where particle speeds approach $c$, such as in astrophysical jets or advanced accelerator designs.

Formally, in the regime where particle velocities are much slower than the speed of light, the full relativistic theory should converge to a classical description. This is embodied in the non-relativistic limit $c \to \infty$, where the RVML system is expected to reduce to the Vlasov-Poisson-Landau (VPL) system. The physical essence of this singular limit is profound: it marks a transition from the relativistic paradigm of finite-speed electromagnetic waves and Lorentz-invariant dynamics, to the classical world of instantaneous Coulomb interactions and Newtonian particle motion. The ``Maxwell" part of the model simplifies to the``Poisson" equation, signifying the demise of magnetic induction and electromagnetic radiation, leaving only an electrostatic potential determined instantaneously by the charge density.

Despite the clear physical intuition, a rigorous mathematical justification of this limit presents formidable challenges. The primary difficulties stem from the singular perturbation in the Maxwell equations, leading to weakly dissipative electromagnetic fields for large $c$, and the need for analytic results that hold uniformly in $c$. This work bridges this gap. Our goal is to provide a rigorous derivation of the Vlasov-Poisson-Landau system from the relativistic Vlasov-Maxwell-Landau system in a periodic box. This achievement rests on three key technical advances: 

$\bullet$  establishing a uniform-in-$c$ coercivity estimate for the relativistic Landau operator; 

$\bullet$  constructing a novel weighted energy functional capable of taming the weak field dissipation; 

$\bullet$  proving a corresponding global well-posedness result, independent of the large light speed, for the relativistic system by a delicated energy functional.

 This work thus not only secures a fundamental limit in kinetic theory but also demonstrates the intricate interplay between collisional dissipation and long-range forces in a singular asymptotic regime.

\subsection{Related Works}

This subsection provides a brief review of previous works closely related to our main results.

\subsubsection{Non-relativistic Equations}
In the near-Maxwellian regime, Guo pioneered the global existence of classical solutions to the Landau equation using a weighted energy method \cite{Guo-CMP-2002}, later extending this result to the Vlasov-Poisson-Landau system \cite{Guo-JAMS-2012}. For the more complex Vlasov-Maxwell-Landau system, global classical solutions were constructed in \cite{ Yu-JMP-2004, Duan-AIHP-2014, Lei-Zhao-JFA-2014}, with a rigorous justification of the Hilbert expansion given in \cite{LLXZ-JLMS-2024}. In the near-vacuum regime, Luk established the global well-posedness for the Landau equation with moderately soft potentials \cite{Luk-ANPDE-2019}.

\subsubsection{Relativistic Equations}
For the relativistic Landau equation, Lemou provided a detailed analysis of the linearized collision operator \cite{Lemou-MMAS-2000}. Strain and Guo then established the global well-posedness for the RVML  system near equilibrium in a periodic box \cite{Strain-Guo-CMP-2004}, a result later extended to the whole space with studies on smoothing effects and optimal time decay \cite{Yu-JMP-2004, Yu-JDE-2009, Yang-Yu-JMPA-2012, Liu-Zhao-JDE-2014, Xiao-JDE-2015}. Parallel results for relativistic Boltzmann-type equations can be found in \cite{Glassey-Strauss-PRIMS-1993, Glassey-Strauss-TTSP-1993, Guo-Strain-CMP-2012, Jang-Strian-ANPDE-2022}, while the Hilbert expansion for the RVML system was treated in \cite{Ouyang-Wu-Xiao-QAM-2024}. For the relativistic Vlasov-Maxwell system, the nonlinear stability of vacuum was first proven by Glassey and Strauss \cite{Glassey-Strauss-CMP-1987}, with further developments in \cite{Glassey-Schaeffer-CMP-1988, Schaeffer-IMJ-2004, Bigorgne-CMP-2020, Wei-Yang-CMP-2021, Wang-CMP-2022}.

\subsubsection{The Classical Limit}
The classical (Newtonian) limit of the relativistic Boltzmann equation was first addressed by Calogero, who constructed local-in-time solutions uniform in the light speed $c$ \cite{Calogero-JMP-2004}. Strain later justified this limit over arbitrary finite time intervals \cite{Strain-SIAM-2010}. More recently, Wang-Xiao \cite{Wang-Xiao-Arxiv-2023} established a uniform-in-$c$ Hilbert expansion, and Cao-Ouyang-Wang-Xiao \cite{Cao-Ouyang-Wang-Xiao-JFA-2026} proved the global-in-time validity of the Newtonian limit. For the relativistic Vlasov-Maxwell system, results are given in \cite{Schaeffer-CMP-1986}, with the Newtonian limit of the relativistic Vlasov-Maxwell-Boltzmann system studied in \cite{Jiang-Lei-Zhao-JFA-2024}. To the best of our knowledge, the Newtonian limit for the relativistic Landau and Vlasov-Maxwell-Landau equations remains an open problem.

\subsection{Relativistic Vlasov-Maxwell-Landau System and Vlasov-Poisson-Landau system}
\subsubsection{Relativistic Vlasov-Maxwell-Landau System} For light speed $c\geq1$, consider the Cauchy problem of the RVML system\cite{Balescu-1988, Strain-Guo-CMP-2004}:
\begin{align}\label{main1-00}
\begin{aligned}
&\partial_t F^c_+ + \frac{cp}{p^0_+}\cdot \nabla_x F^c_+ + Z\mathrm{e}\big(E^c+\frac{p}{p^0_+}\times  B^c\big)\cdot\nabla_p F^c_+ =\mathcal{C}\left(F^c_-,F^c_+\right)+\mathcal{C}\left(F^c_+,F^c_+\right),\\
&\partial_t F^c_- + \frac{cp}{p^0_-}\cdot \nabla_x F^c_- - \mathrm{e}\big(E^c+\frac{p}{p^0_-}\times  B^c\big)\cdot\nabla_p F^c_- =\mathcal{C}\left(F^c_-,F^c_-\right)+\mathcal{C}\left(F^c_+,F^c_-\right),\\
&\partial_t E^c-  c\nabla_x \times B^c =-4\pi \int_{\mathbb R^3}\Big(Z\mathrm{e}\frac{cp}{p^0_+} F^c_+-\mathrm{e}\frac{cp}{p^0_-} F^c_-\Big)\, \mathrm{d} p, \\\rule{0ex}{1.0em}
&\partial_t B^c+ c\nabla_x \times E^c=0,\\\rule{0ex}{1.5em}
& \nabla_x\cdot E^c=4\pi \int_{\mathbb R^3}  \left(Z\mathrm{e}F^c_+-\mathrm{e}F^c_-\right)\, \mathrm{d}p, \qquad \nabla_x\cdot B^c=0,\\\rule{0ex}{1.0em}
&F^c_{\pm}(0,x,p)=F^c_{\pm,0}(x,p),\qquad E^c(0,x)=E^c_{0}(x),\qquad B^c(0,x)=B^c_{0}(x).
\end{aligned}
\end{align}
Here $F^c_{+}(t,x,p)$ is the number density function for ions at time $t\geq0$, position $x=(x_1,x_2,x_3) \in \mathbb T^3=[-\pi, \pi]^3$, and momentum $p=(p_1,p_2,p_3)\in \mathbb R^3$. $p^0_+=\sqrt{m^2_+c^2 + |p|^2}$ is the energy of an ion, constants $Z$ and $m_+$ are the charge number and rest mass for the ion species, respectively. Correspondingly, $F^c_{-}(t,x,p)$, $p^0_-$, $m_-$ are the number density function, energy, rest mass for the electro species. $e$ is the magnitude of an electro's charge, and $\big(E^c(t,x), B^c(t,x)\big)$ are the electromagnetic fields.

Denote the four-momentums $p^{\mu}_{\pm} = \left(p^0_{\pm}, p\right)$ and $q^{\mu}_{\pm} = \left(q^0_{\pm}, q\right)$.  We use the Einstein convention that repeated up-down indices be
summed and we raise and lower indices using the Minkowski metric $g^{\mu\nu}
:=
\text{diag}(-1, 1, 1, 1)$. The Lorentz inner product is then given by
\begin{align*}
p^{\mu}_{\pm}q_{\pm,\mu}:=-p^0_{\pm}q^0_{\pm} +
\sum_{i=1}^3 p_iq_i .
\end{align*}
The relativistic Landau collision operator $\mathcal{C}\left(G_- ,H_+\right)$ on the R.H.S. of \eqref{main1-00}, which registers binary collisions between particles, takes the following form:
\begin{align} \label{coll}
\mathcal{C}\left(G_- ,H_+\right) :=  2\pi Z^2\mathrm{e}^4\ln(\wedge )\nabla_{p}\cdot\left\{\int_{\mathbb{R}^3} \Phi^c_{-,+}(p,q) \Big[G_-(q) \nabla_pH_+(p) -\nabla_qG_-(q) H_+(p)\Big]\,\mathrm{d} q \right\},
\end{align}
where the collision kernel $\Phi^c_{-,+}(p,q)$ is a $3\times3$ non-negative matrix
\begin{align}\label{cker}
    \Phi^c_{-,+}(p,q):=\left(\frac{m_-c}{q^0_-}\right)\left(\frac{m_+c}{p^0_+}\right)\Lambda^c_{-,+}(p,q)\mathcal{S}^c_{-,+}(p,q)
\end{align}
with
\begin{align}
\Lambda^c_{-,+}(p,q):=&\;\frac{1}{m_-^2m_+^2c^4}\left(q^{\mu}_-p_{+,\mu}\right)^2
\left(\frac{1}{m_-^2m_+^2c^2}\left(q^{\mu}_-p_{+,\mu}\right)^2-c^2\right)^{-\frac{3}{2}},\label{definition Phi c p q}\\
\mathcal{S}^c_{-,+}(p,q):=&\;\left(\frac{1}{m_-^2m_+^2c^2}\left(q^{\mu}_-p_{+,\mu}\right)^2-c^2\right)I_3-\Big(\frac{q}{m_-}
-\frac{p}{m_+}\Big)\otimes\Big(\frac{q}{m_-}-\frac{p}{m_+}\Big)\nonumber\\
&\;-\left(\frac{1}{m_-m_+c^2}\left(q^{\mu}_-p_{+,\mu}\right)+1\right)\Big(\frac{q}{m_-}\otimes \frac{p}{m_+}+\frac{p}{m_+}\otimes \frac{q}{m_-}\Big).\label{definition S c p q}
\end{align}
Here $\ln(\wedge)$ denotes the Coulomb logarithm for the interaction of ions and electros.
 From \cite{Strain-Guo-CMP-2004}, for our purposes, we may simply regard $\ln(\wedge)$ as a fixed constant independent of other parameters.
 Corresponding to \eqref{coll}, other collision operators on the R.H.S. of \eqref{main1-00} can be written as
 \begin{align*} 
&\mathcal{C}\left(G_+ ,H_+\right) :=  2\pi Z^2\mathrm{e}^4\ln(\wedge )\nabla_{p}\cdot\left\{\int_{\mathbb{R}^3} \Phi^c_{+,+}(p,q) \Big[G_+(q) \nabla_pH_+(p) -\nabla_qG_+(q) H_+(p)\Big]\,\mathrm{d} q \right\},\\
&\mathcal{C}\left(G_- ,H_-\right) :=  2\pi Z^2\mathrm{e}^4\ln(\wedge )\nabla_{p}\cdot\left\{\int_{\mathbb{R}^3} \Phi^c_{-,-}(p,q) \Big[G_-(q) \nabla_pH_-(p) -\nabla_qG_-(q) H_-(p)\Big]\,\mathrm{d} q \right\},\\
&\mathcal{C}\left(G_+ ,H_-\right) :=  2\pi Z^2\mathrm{e}^4\ln(\wedge )\nabla_{p}\cdot\left\{\int_{\mathbb{R}^3} \Phi^c_{+,-}(p,q) \Big[G_+(q) \nabla_pH_-(p) -\nabla_qG_+(q) H_-(p)\Big]\,\mathrm{d} q \right\}.
\end{align*}

It is well-known that $\Phi^c_{\mp,\pm}(p,q)$ satisfies
\begin{align}\label{Phi}
    \sum_{i=1}^3\Phi^{c,ij}_{\mp,\pm}(p,q)\left(\frac{p_i}{p^0_{\pm}}-\frac{q_i}{q^0_{\mp}}\right)=\sum_{j=1}^3\Phi^{c,ij}_{\mp,\pm}(p,q)\left(\frac{p_j}{p^0_{\pm}}-\frac{q_j}{q^0_{\mp}}\right)=0.
\end{align}
With the fact \eqref{Phi}, it is straightforward to verify that the collision operator $\mathcal{C}$ satisfies the orthogonality property:
\begin{align*}
    \int_{\mathbb{R}^3}\begin{pmatrix}1\\m_+p\\m_+p^0_+\end{pmatrix}\mathcal{C}\left(G_- ,H_+\right)(p)\,\mathrm{d} p+\int_{\mathbb{R}^3}\begin{pmatrix}1\\m_-p\\m_-p^0_-\end{pmatrix}\mathcal{C}\left(H_+ ,G_-\right)(p)\,\mathrm{d} p=\mathbf{0}.
\end{align*}
Then, the classical solution to $\eqref{main1-00}$ satisfies conservation laws
\begin{align}\label{cons-RVML}
    &\frac{\mathrm{d}}{\mathrm{d} t}\int_{\mathbb{T}^3\times\mathbb{R}^3} m_+F^c_+\,\mathrm{d} p\mathrm{d} x=\frac{\mathrm{d}}{\mathrm{d} t}\int_{\mathbb{T}^3\times\mathbb{R}^3} m_-F^c_-\,\mathrm{d} p\mathrm{d} x=0,\nonumber\\
    &\frac{\mathrm{d}}{\mathrm{d} t}\bigg\{\int_{\mathbb{T}^3\times\mathbb{R}^3}p \left(m_+F_+^c+m_-F^c_-\right)\,\mathrm{d} p\mathrm{d} x +\frac{1}{4\pi c}\int_{\mathbb{T}^3}\Big(E^{c}\times B^{c}\Big)\,\mathrm{d} x\bigg\}=\mathbf{0},\\
    &\frac{\mathrm{d}}{\mathrm{d} t}\bigg\{\int_{\mathbb{T}^3\times\mathbb{R}^3}\left(m_+p^0_+F^c_+ + m_-p^0_-F^c_-\right)\,\mathrm{d} p\mathrm{d} x+\frac{1}{8\pi c} \int_{\mathbb{T}^3} \Big(|E^{c}|^2+|B^{c}|^2\Big)\,\mathrm{d} x\bigg\}=0.\nonumber
\end{align}

To construct solutions to the Cauchy problem \eqref{main1-00} around an equilibrium, we introduce global relativistic Maxwellians:
\begin{align}\label{Maxwellian-c}
\mu_{+}^{c}(p)= \frac{n_+\exp\left\{-\gamma_{+}p^0_{+} /(m_{+}c) \right\} }{4 \pi Z\mathrm{e}m^2_{+}
c  k_BT K_2(\gamma_{+})} ,\qquad \mu_{-}^{c}(p)= \frac{n_-\exp\left\{-\gamma_{-}p^0_{-} /(m_{-}c) \right\} }{4 \pi  \mathrm{e} m^2_{-}
c k_BT K_2(\gamma_{-})} ,
\end{align}
where constants $n_{\pm}$ are number densities for ions and electros,  $ \gamma_{\pm}=\frac{m_{\pm}c^2}{k_B T}$ with the temperature $T$ and the Boltzmann constant  $k_B$, and $K_2(\cdot)$ is the modified second order Bessel function defined in \eqref{defini-kj}.

 \subsubsection{Vlasov-Poisson-Landau  system} Formally, we let $c\rightarrow \infty$ in \eqref{main1-00} to obtain the Cauchy problem of the Vlasov-Poisson-Landau(VPL) system:
\begin{align}\label{main2-00}
\begin{aligned}
&\partial_t F^{\infty}_+ + p\cdot \nabla_x F^{\infty}_+ + Z\mathrm{e} E^{\infty} \cdot\nabla_p F^{\infty}_+ =\mathbb{C}\left(F^{\infty}_-,F^{\infty}_+\right)+\mathbb{C}\left(F^{\infty}_+,F^{\infty}_+\right),\\
&\partial_t F^{\infty}_- + p\cdot \nabla_x F^{\infty}_- - \mathrm{e}E^{\infty}\cdot\nabla_p F^{\infty}_- =\mathbb{C}\left(F^{\infty}_-,F^{\infty}_-\right)+\mathbb{C}\left(F^{\infty}_+,F^{\infty}_-\right),\\
& \nabla_x\cdot E^{\infty}=4\pi \int_{\mathbb R^3}  \left(Z\mathrm{e}F^{\infty}_+-\mathrm{e}F^{\infty}_-\right)\, \mathrm{d}p, \\\rule{0ex}{1.0em}
&F^{\infty}_{\pm}(0,x,p)=F^{\infty}_{\pm,0}(x,p),
\end{aligned}
\end{align}
and
\begin{align}\label{B0}
\begin{aligned}
\partial_t E^{\infty} =&-4\pi \int_{\mathbb R^3} p\big(Z\mathrm{e} F^{\infty}_+-\mathrm{e} F^{\infty}_-\big)\, \mathrm{d} p,\qquad \nabla_x \times E^{\infty} =0,\\\rule{0ex}{1.0em}
\partial_t B^{\infty}=&\nabla_x \times B^{\infty}=\nabla_x\cdot B^{\infty}=0.
\end{aligned}
\end{align}
The  Landau collision operators  on the R.H.S. of \eqref{main2-00} are given by
\begin{align*}
&\mathbb{C}\left(G_- ,H_+\right) := 2\pi Z^2\mathrm{e}^4\ln(\wedge )\nabla_{p}\cdot\left\{\int_{\mathbb{R}^3} \Phi^{\infty}_{-,+}(p,q) \Big[G_-(q) \nabla_pH_+(p) -\nabla_qG_-(q) H_+(p)\Big]\,\mathrm{d} q \right\},\\
&\mathbb{C}\left(G_+ ,H_+\right) := 2\pi Z^2\mathrm{e}^4\ln(\wedge )\nabla_{p}\cdot\left\{\int_{\mathbb{R}^3} \Phi^{\infty}_{+,+}(p,q) \Big[G_+(q) \nabla_pH_+(p) -\nabla_qG_+(q) H_+(p)\Big]\,\mathrm{d} q \right\},\\
&\mathbb{C}\left(G_- ,H_-\right) :=  2\pi Z^2\mathrm{e}^4\ln(\wedge )\nabla_{p}\cdot\left\{\int_{\mathbb{R}^3} \Phi^{\infty}_{-,-}(p,q) \Big[G_-(q) \nabla_pH_-(p) -\nabla_qG_-(q) H_-(p)\Big]\,\mathrm{d} q \right\},\\
&\mathbb{C}\left(G_+ ,H_-\right) :=  2\pi Z^2\mathrm{e}^4\ln(\wedge )\nabla_{p}\cdot\left\{\int_{\mathbb{R}^3} \Phi^{\infty}_{+,-}(p,q) \Big[G_+(q) \nabla_pH_-(p) -\nabla_qG_+(q) H_-(p)\Big]\,\mathrm{d} q \right\}.
\end{align*}
Here the collision kernel $\Phi^{\infty}_{-,+}(p,q)$ is a $3\times3$ non-negative matrix
\begin{align*}
    \Phi^{\infty}_{-,+}(p,q):=\Lambda^{\infty}_{-,+}(p,q)\mathcal{S}^{\infty}_{-,+}(p,q)
\end{align*}
with
\begin{align*}
\Lambda^{\infty}_{-,+}(p,q):=&\;\Big|\frac{q}{m_-}-\frac{p}{m_+}\Big|^{-3},\\
\mathcal{S}^{\infty}_{-,+}(p,q):=&\;\Big|\frac{q}{m_-}-\frac{p}{m_+}\Big|^2I_3
-\Big(\frac{q}{m_-}-\frac{p}{m_+}\Big)\otimes\Big(\frac{q}{m_-}-\frac{p}{m_+}\Big),\nonumber
\end{align*}
and other collision kernels $\Phi^{\infty}_{+,+}(p,q),\Phi^{\infty}_{-,-}(p,q), \Phi^{\infty}_{+,-}(p,q)$ can be given in the same way.

$\eqref{B0}_1$ implies that there exists a scalar function $\phi^{\infty}$ such that $E^{\infty}(t,x) \equiv-\nabla_x\phi^{\infty}(t,x)$,  while $\eqref{B0}_2$ implies that $B^{\infty}$ should be a constant vector.  Since $B^c(t,x)$ should converges to $B^{\infty}$ as $c\rightarrow \infty$, one can further obtain that $B^{\infty}=\bar{B}^0$. 

As is known to all, for a classical solution to $\eqref{main2-00}$, its mass, total momentum, and total energy are conservative
\begin{align}\label{cons-VPL}
    &\frac{\mathrm{d}}{\mathrm{d} t}\int_{\mathbb{T}^3\times\mathbb{R}^3} m_+F^{\infty}_+\,\mathrm{d} p\mathrm{d} x=\frac{\mathrm{d}}{\mathrm{d} t}\int_{\mathbb{T}^3\times\mathbb{R}^3} m_-F^{\infty}_-\,\mathrm{d} p\mathrm{d} x=0,\nonumber\\
    &\frac{\mathrm{d}}{\mathrm{d} t}\bigg\{\int_{\mathbb{T}^3\times\mathbb{R}^3}p \left(m_+F_+^{\infty}+m_-F^{\infty}_-\right)\,\mathrm{d} p\mathrm{d} x \bigg\}=\mathbf{0},\\
    &\frac{\mathrm{d}}{\mathrm{d} t}\bigg\{\int_{\mathbb{T}^3\times\mathbb{R}^3}\frac{|p|^2}{2}\left(m_+F^{\infty}_+ + m_-F^{\infty}_-\right)\,\mathrm{d} p\mathrm{d} x+\frac{1}{8\pi} \int_{\mathbb{T}^3} |E^{\infty}|^2\,\mathrm{d} x\bigg\}=0.\nonumber
\end{align}

Corresponds to \eqref{Maxwellian-c}, the equilibrium of $F^{\infty}_{\pm}$ in \eqref{main2-00} are
\begin{align}\label{Maxwellian-inf}
\mu_{+}^{\infty}(p)= \frac{n_+}{Z\mathrm{e}}\Big(\frac{ m_{+}}{2\pi T }\Big)^{\frac32} \exp\left\{-\frac{m_{+}|p|^2}{2T} \right\},\qquad \mu_{-}^{\infty}(p)= \frac{n_-}{\mathrm{e}}\Big(\frac{ m_{-}}{2\pi T }\Big)^{\frac32} \exp\left\{-\frac{m_{-}|p|^2}{2T} \right\}.
\end{align}

\subsubsection{Series-expansion-in-$c^{-1}$ for the non-relativistic limit} In this subsection, we address the following mathematical challenges inherent in the non-relativistic limit $c\rightarrow \infty$: (i). Characterizing the asymptotic behavior of the collision operator; (ii). Justifying the transition of the equilibrium state from \eqref{Maxwellian-c} to \eqref{Maxwellian-inf}. (iii).  Controlling the electromagnetic field in the limit. 

Let us explain clearly for $(iii)$. From \eqref{Maxwellian-c}, one may easily verify that $[E^c, B^c]$ satisfy the wave equation 
\begin{align}\label{toymodelwaveequation} \pa_{tt} u-c^2\triangle u=\mathcal{S}. \end{align}
This suggests that the pointwise decay of $u$  behaves like $(1+t/c)^{-\mathsf{s}}$ over the whole space, which implies that the dissipation of the electromagnetic field $[E^c, B^c]$ weakens as $c\rightarrow \infty$. Consequently, controlling nonlinear terms such as $\big(E^c+\frac{p}{p^0_\pm}\times  B^c\big)\cdot\nabla_p F^c_\pm$ becomes a major challenge in establishing global well-posedness.
 
\medskip

Since the limiting magnetic field $B^\infty(t,x)$ is a constant vector,   it is natural to expand the initial data $B^\infty(0,x)$ in powers of 
 $c^{-1}$: $B^\infty(0,x)=\sum_{j=1}^N c^{-i} B^{i,c}(0,x)$. To fix the number $N$, we seek a  solution  $[F^c_{\pm}(t,x,p), E^c(t,x), B^c(t,x)]$ to the RVML system \eqref{main1-00} that admits the expansion:
\begin{align*} 
F^c_{\pm}(t,x,p)=\sum_{i=0}^N\frac{1}{c^i}F^{i,c}_{\pm}(t,x,p),\quad E^c(t,x)=\sum_{i=0}^N\frac{1}{c^i}E^{i,c}(t,x),\quad B^c(t,x)=\sum_{i=0}^N\frac{1}{c^i}B^{i,c}(t,x).
\end{align*}
Under this ansatz, the typical nonlinear term can be rewritten as 
\begin{align}\label{typicalnonlinearterm}
(E^c+\frac{p}{p^0_\pm}\times  B^c\big)\cdot\nabla_p F^c_\pm=\sum_{i=0}^N\sum_{j=0}^N\frac{1}{c^{i+j}} (E^{i,c}+\frac{p}{p^0_\pm}\times  B^{i,c}\big)\cdot\nabla_p F^{j,c}_\pm.
\end{align}
Comparing orders of $c^{-1}$ reveals the following advantages:

\noindent$\bullet$ At leading order ($0$-th order),   the associated system for $(F^{0,c}_{\pm}(t,x),E^{0,c}(t,x),B^{0,c}(t,x))$ closely resembles the VPL system \eqref{main2-00}. The   nonlinear term \eqref{typicalnonlinearterm} reduces to a standard nonlinear expression involving only these leading-order components $(F^{0,c}_{\pm}(t,x),E^{0,c}(t,x),B^{0,c}(t,x))$.

\noindent$\bullet$ For $m$-th order with $1\le m\le N$, the   term \eqref{typicalnonlinearterm} includes both linear contributions in $(F^{m,c}_{\pm}(t,x),E^{m,c}(t,x),B^{m,c}(t,x))$ and a specific nonlinear term of the form $c^{-m}(E^{m,c}+\frac{p}{p^0_\pm}\times  B^{m,c}\big)\cdot\nabla_p F^{m,c}_\pm$. This structure is promising for analysis, as the factor $c^{-m}$ can counterbalance the weakening dissipation from the wave equation \eqref{toymodelwaveequation}.
\medskip

 By further computation, we finally fix $N=3$ which will yield the following expansion: 
 assume that the solution $\big[\begin{pmatrix}\begin{smallmatrix} F^c_{+}(t,x,p)\\F^c_{-}(t,x,p)\\\end{smallmatrix}\end{pmatrix}, E^c(t,x), B^c(t,x)\big]$ to Cauchy problem \eqref{main1-00} and its initial datum $\big[\begin{pmatrix}\begin{smallmatrix} F^c_{+,0}(x,p)\\F^c_{-,0}(x,p)\\\end{smallmatrix}\end{pmatrix}, E^c_{0}(x), B^c_0(x)\big]$ satisfy the following expansions:
\begin{align}\label{expanion0t}
F^c_{\pm}(t,x,p)=\sum_{i=0}^3\frac{1}{c^i}F^{i,c}_{\pm}(t,x,p),\qquad E^c(t,x)=\sum_{i=0}^3\frac{1}{c^i}E^{i,c}(t,x),\qquad B^c(t,x)=\sum_{i=0}^3\frac{1}{c^i}B^{i,c}(t,x),
\end{align}
and
\begin{align}\label{expanion00}
F^c_{\pm,0}(x,p)=\sum_{i=0}^3\frac{1}{c^i}F^{i,c}_{0}(x,p),\qquad E^c_{0}(x)=\sum_{i=0}^3\frac{1}{c^i}E^{i,c}_{0}(x),\qquad B^c_{0}(x)=\sum_{i=0}^3\frac{1}{c^i}B^{i,c}_{0}(x).
\end{align}
Now we plug the above expansions into \eqref{main1-00} and compare the order of ${c}^{-1}$ to obtain that
\begin{itemize}
          \item [i)] For $i=0, 1, 2$, $B^{i,c}(t,x)$ satisfy
          \begin{align}\label{mag-0cnd}
          \partial_t B^{i,c}(t,x)=\nabla_x \times B^{i,c}(t,x)=\nabla_x\cdot B^{i,c}(t,x)=0,
          \end{align}
          which implies that $B^{i,c}(t,x)$ are nothing but constant vectors $\bar{B}^{i}$.
          
          \item [ii)] For given constant magnetic fields $\bar{B}^{i}$ with $i=0, 1, 2$, $\big[\begin{pmatrix}\begin{smallmatrix} F^{0,c}_{+}(t,x,p)\\F^{0,c}_{-}(t,x,p)\\\end{smallmatrix}\end{pmatrix}, E^{0,c}(t,x)\big]$ solves the following Cauchy problem:
\begin{align}\label{mainF0-00}
&\partial_t F^{0,c}_+ + \frac{cp}{p^0_+}\cdot \nabla_x F^{0,c}_+ + Z\mathrm{e}\Big(E^{0,c}+\frac{p}{p^0_+}\times  \sum_{j=0}^2\frac{1}{c^j}\bar{B}^{j}\Big)\cdot\nabla_p F^{0,c}_+ \nonumber\\ &\qquad=\mathcal{C}\left(F^{0,c}_-,F^{0,c}_+\right)+\mathcal{C}\left(F^{0,c}_+,F^{0,c}_+\right),\nonumber\\
&\partial_t F^{0,c}_- + \frac{cp}{p^0_-}\cdot \nabla_x F^{0,c}_- - \mathrm{e}\Big(E^{0,c}+\frac{p}{p^0_-}\times  \sum_{j=0}^2\frac{1}{c^j}\bar{B}^{j}\Big)\cdot\nabla_p F^{0,c}_- \nonumber\\ &\qquad=\mathcal{C}\left(F^{0,c}_-,F^{0,c}_-\right)+\mathcal{C}\left(F^{0,c}_+,F^{0,c}_-\right),\\
&\partial_t E^{0,c}=-4\pi \int_{\mathbb R^3}\Big(Z\mathrm{e}\frac{cp}{p^0_+} F^{0,c}_+-\mathrm{e}\frac{cp}{p^0_-} F^{0,c}_-\Big)\, \mathrm{d} p, \nonumber\\\rule{0ex}{1.5em}
& \nabla_x\cdot E^{0,c}=4\pi \int_{\mathbb R^3}  \left(Z\mathrm{e}F^{0,c}_+-\mathrm{e}F^{0,c}_-\right)\, \mathrm{d}p, \qquad \nabla_x \times E^{0,c}=0,\nonumber\\\rule{0ex}{1.0em}
&F^{0,c}_{\pm}(0,x,p)=F^{0,c}_{\pm,0}(x,p),\qquad E^{0,c}(0,x)=E^{0,c}_{0}(x).\nonumber
\end{align}
	\item [iii)] For $i=1, 2$, $\big[\begin{pmatrix}\begin{smallmatrix} F^{i,c}_{+}(t,x,p)\\F^{i,c}_{-}(t,x,p)\\\end{smallmatrix}\end{pmatrix}, E^{i,c}(t,x)\big]$ is a solution to the following Cauchy problem of linear system:
		\begin{align}\label{mainFi-00}
&\partial_t F^{i,c}_+ + \frac{cp}{p^0_+}\cdot \nabla_x F^{i,c}_+ + Z\mathrm{e}\sum_{\substack{j_1+j_2=i\\0\leq j_1,j_2\leq i}}E^{j_1,c}\cdot\nabla_p F^{j_2,c}_+ + Z\mathrm{e}\Big(\frac{p}{p^0_+}\times  \sum_{j=0}^2\frac{1}{c^j}\bar{B}^{j}\Big)\cdot\nabla_p F^{i,c}_+ \nonumber\\
&\qquad=\sum_{\substack{j_1+j_2=i\\0\leq j_1,j_2\leq i}}\Big[\mathcal{C}\left(F^{j_1,c}_-,F^{j_2,c}_+\right)+\mathcal{C}\left(F^{j_1,c}_+,F^{j_2,c}_+\right)\Big],\nonumber\\
&\partial_t F^{i,c}_- + \frac{cp}{p^0_-}\cdot \nabla_x F^{i,c}_- - \mathrm{e}\sum_{\substack{j_1+j_2=i\\0\leq j_1,j_2\leq i}}E^{j_1,c}\cdot\nabla_p F^{j_2,c}_- - \mathrm{e}\Big(\frac{p}{p^0_-}\times  \sum_{j=0}^2\frac{1}{c^j}\bar{B}^{j}\Big)\cdot\nabla_p F^{i,c}_-  \nonumber\\
&\qquad=\sum_{\substack{j_1+j_2=i\\0\leq j_1,j_2\leq i}}\Big[\mathcal{C}\left(F^{j_1,c}_-,F^{j_2,c}_-\right)+\mathcal{C}\left(F^{j_1,c}_+,F^{j_2,c}_-\right)\Big],\\
&\partial_t E^{i,c}=-4\pi \int_{\mathbb R^3}\Big(Z\mathrm{e}\frac{cp}{p^0_+} F^{i,c}_+-\mathrm{e}\frac{cp}{p^0_-} F^{i,c}_-\Big)\, \mathrm{d} p, \nonumber\\\rule{0ex}{1.5em}
& \nabla_x\cdot E^{i,c}=4\pi \int_{\mathbb R^3}  \left(Z\mathrm{e}F^{i,c}_+-\mathrm{e}F^{i,c}_-\right)\, \mathrm{d}p, \qquad \nabla_x \times E^{i,c}=0,\nonumber\\\rule{0ex}{1.0em}
&F^{i,c}_{\pm}(0,x,p)=F^{i,c}_{\pm,0}(x,p),\qquad E^{i,c}(0,x)=E^{i,c}_{0}(x).\nonumber
\end{align}
	\item [iv)] $\big[\begin{pmatrix}\begin{smallmatrix} F^{3,c}_{+}(t,x,p)\\F^{3,c}_{-}(t,x,p)\\\end{smallmatrix}\end{pmatrix}, E^{3,c}(t,x), B^{3,c}(t,x)\big]$ is a solution to the Cauchy problem of RVML system:
		\begin{align}\label{mainF3-00}
&\partial_t F^{3,c}_+ + \frac{cp}{p^0_+}\cdot \nabla_x F^{3,c}_+ + Z\mathrm{e}\sum_{\substack{j_1+j_2\geq3\\0\leq j_1,j_2\leq 3}}\frac{1}{c^{j_1+j_2-3}}E^{j_1,c}\cdot\nabla_p F^{j_2,c}_+\nonumber\\
 &\qquad+ Z\mathrm{e}\frac{p}{p^0_+}\times  \Big(\sum_{j=0}^2\frac{1}{c^j}B^{j,c}\cdot\nabla_p F^{3,c}_++B^{3,c}\cdot \sum_{j=0}^2\frac{1}{c^j}\nabla_pF^{j,c}_++\frac{1}{c^3}B^{3,c}\cdot\nabla_p F^{3,c}_+\Big) \nonumber\\
&\qquad=\sum_{\substack{j_1+j_2\geq3\\0\leq j_1,j_2\leq 3}}\frac{1}{c^{j_1+j_2-3}}\Big[\mathcal{C}\left(F^{j_1,c}_-,F^{j_2,c}_+\right)
+\mathcal{C}\left(F^{j_1,c}_+,F^{j_2,c}_+\right)\Big],\nonumber\\
&\partial_t F^{3,c}_- + \frac{cp}{p^0_-}\cdot \nabla_x F^{3,c}_- - \mathrm{e}\sum_{\substack{j_1+j_2\geq3\\0\leq j_1,j_2\leq 3}}\frac{1}{c^{j_1+j_2-3}}E^{j_1,c}\cdot\nabla_p F^{j_2,c}_-\nonumber\\
 &\qquad- \mathrm{e}\frac{p}{p^0_-}\times  \Big(\sum_{j=0}^2\frac{1}{c^j}\bar{B}^{j}\cdot\nabla_p F^{3,c}_++B^{3,c}\cdot \sum_{j=0}^2\frac{1}{c^j}\nabla_pF^{j,c}_-+\frac{1}{c^3}B^{3,c}\cdot\nabla_p F^{3,c}_-\Big)  \nonumber\\
&\qquad=\sum_{\substack{j_1+j_2\geq3\\0\leq j_1,j_2\leq 3}}\frac{1}{c^{j_1+j_2-3}}\Big[\mathcal{C}\left(F^{j_1,c}_-,F^{j_2,c}_+\right)
+\mathcal{C}\left(F^{j_1,c}_+,F^{j_2,c}_+\right)\Big],\\
&\partial_t E^{3,c}-  c\nabla_x \times B^{3,c}=-4\pi \int_{\mathbb R^3}\Big(Z\mathrm{e}\frac{cp}{p^0_+} F^{3,c}_+-\mathrm{e}\frac{cp}{p^0_-} F^{3,c}_-\Big)\, \mathrm{d} p, \nonumber\\\rule{0ex}{1.0em}
&\partial_t B^{3,c}+ c\nabla_x \times E^{3,c}=0,\nonumber\\\rule{0ex}{1.5em}
& \nabla_x\cdot E^{3,c}=4\pi \int_{\mathbb R^3}  \left(Z\mathrm{e}F^{3,c}_+-\mathrm{e}F^{3,c}_-\right)\, \mathrm{d}p, \qquad \nabla_x\cdot B^{3,c}=0,\nonumber\\\rule{0ex}{1.0em}
&F^{3,c}_{\pm}(0,x,p)=F^{3,c}_{\pm,0}(x,p),\qquad E^{3,c}(0,x)=E^{3,c}_{0}(x),\qquad B^{3,c}(0,x)=E^{3,c}_{0}(x).\nonumber
\end{align}
	\end{itemize}

\begin{remark} Assuming a constant magnetic field initially, the local well-posedness result implies that our series expansion in $c^{-1}$ captures the essential structure of the solution $[F^c_{\pm}(t,x,p), E^c(t,x),\\ B^c(t,x)]$ to the RVML system \eqref{main1-00}. Therefore, the expansion is a faithful representation of the solution dynamics. Consequently, truncating at order $N=3$ is both sufficient and necessary for our purpose of analyzing the limit, while also being technically expedient.
 \end{remark}

\subsection{Reformulation}
In this subsection, we reformulate the systems \eqref{mainF0-00}-\eqref{mainF3-00} and the VPL system \eqref{main2-00} in the perturbative framework.

\subsubsection{Normalization the equilibrium} From now on, for simplicity, we normalize all constants except for the light speed $c$ in \eqref{main1-00} and \eqref{main2-00} to be 1 since they don't lead to essential difficulties in the proof of main results in this article. Accordingly, the normalized Maxwellians in \eqref{Maxwellian-c} and \eqref{Maxwellian-inf} are
\begin{align*}
\mu^{c}(p):=&\mu_{+}^{c}(p)=\mu_{-}^{c}(p)= \frac{\exp\left\{-cp^0 \right\} }{4 \pi
c  K_2(c^2)},\qquad \mbox{and}\\\rule{0ex}{1.0em}
\mu^{\infty}(p):=&\mu_{+}^{\infty}(p)=\mu_{-}^{\infty}(p)= \frac{ 1}{(2\pi )^{\frac32} } \exp\left\{-\frac{|p|^2}{2} \right\}.
\end{align*}
In order to establish the global existence of a solution to \eqref{main1-00} for given $c\geq1$ and derive the classical limit from the RVML system \eqref{main1-00} to the VPL system \eqref{main2-00},
we will perturb the solution $\big[\begin{pmatrix}\begin{smallmatrix} F^c_{+}(t,x,p)\\F^c_{-}(t,x,p)\\\end{smallmatrix}\end{pmatrix}, E^c(t,x), B^c(t,x)\big]$ to Cauchy problem \eqref{main1-00} around the equilibrium $\big[\begin{pmatrix}\begin{smallmatrix} \mu^c\\\mu^c\\\end{smallmatrix}\end{pmatrix}, \mathrm{0}, B^{\infty}\big]$ and the solution $\big[\begin{pmatrix}\begin{smallmatrix} F^{\infty}_{+}(t,x,p)\\F^{\infty}_{-}(t,x,p)\\\end{smallmatrix}\end{pmatrix}, E^{\infty}(t,x)\big]$ to Cauchy problem \eqref{main2-00} around the equilibrium $\big[\begin{pmatrix}\begin{smallmatrix} \mu^{\infty}\\\mu^{\infty}\\\end{smallmatrix}\end{pmatrix}, \mathrm{0}\big]$  separately.

\subsubsection{Reformulation of the RVML system}
For $i=0, 1, 2, 3$, we introduce the standard perturbation $f^{i,c}_{\pm}(t,x,p)$ as
\begin{align*}
f^{0,c}_{\pm}:=\left(\mu^c\right)^{-\frac12}\left(F^{0,c}_{\pm}-\mu^c\right),\qquad f^{i,c}_{\pm}:=\left(\mu^c\right)^{-\frac12}F^{i,c}_{\pm},\qquad i=1, 2, 3.
\end{align*}
Denote $f^{i,c}(t,x,p)=[f^{i,c}_+(t,x,p), f^{i,c}_-(t,x,p)]^{\textit{t}}$. Here and below $[\cdot,\cdot]^{\textit{t}}$ is used as the transpose of a row vector $[\cdot,\cdot]$. Then for the Cauchy problems \eqref{mainF0-00}, \eqref{mainFi-00}, and \eqref{mainF3-00}, we can obtain that
\begin{itemize}
          \item [i)] $[f^{0,c}(t,x,p), E^{0,c}(t,x)]$ satisfies
\begin{align}\label{mainF0}
&\partial_t f^{0,c} + \frac{cp}{p^0}\cdot \nabla_x f^{0,c}+ E^{0,c}\cdot\frac{cp}{p^0}\sqrt{\mu^c}\zeta_0-\frac{1}{2}\zeta_1\frac{cp}{p^0}\cdot E^{0,c} f^{0,c}\nonumber\\
&\hspace{0.2cm}+\zeta_1 \big(E^{0,c}+\frac{p}{p^0}\times  \sum_{j=0}^2\frac{\bar{B}^{j}}{c^j}\big)\cdot\nabla_p f^{0,c}+\mathcal{L}f^{0,c} =\varGamma\left(f^{0,c}, f^{0,c}\right),\nonumber\\
&\partial_t E^{0,c} =- \int_{\mathbb R^3}\frac{cp}{p^0}\sqrt{\mu^c}\left(f^{0,c}_{+}-f^{0,c}_-\right)\, \mathrm{d} p, \\\rule{0ex}{1.5em}
& \nabla_x\cdot E^{0,c}= \int_{\mathbb R^3}  \sqrt{\mu^c}\left(f^{0,c}_+-f^{0,c}_-\right)\, \mathrm{d}p, \qquad
 \nabla_x \times E^{0,c}=0,\nonumber\\\rule{0ex}{1.0em}
&f^{0,c}(0,x,p)=f^{0,c}_{0}(x,p),\qquad E^{0,c}(0,x)=E^{0,c}_{0}(x), \nonumber
\end{align}
where $\zeta_0=[-1,1]^{\textit{t}}$, $\zeta_1=\mathrm{diag}(1,-1)$. For $h=[h_+,h_-]^{\textit{t}}$, the linear relativistic Landau collision operator $\mathcal{L}h=[\mathcal{L}_+h, \mathcal{L}_-h]^{\textit{t}}$  is given by
\begin{align}\label{linearAK-c}
\begin{aligned}
&\mathcal{L}h=[\mathcal{L}_+h, \mathcal{L}_-h]^{\textit{t}}, \qquad \mathcal{L}_{\pm}h:=\mathcal{A}_{\pm}h+\mathcal{K}_{\pm}h\\
&\mathcal{A}_{\pm}h:=-2\left(\mu^c\right)^{-\frac12}\mathcal{C}\left(\mu^c,\sqrt{\mu^c}h_{\pm}\right),\\
&\mathcal{K}_{+}h=\mathcal{K}_{-}h:=-\left(\mu^c\right)^{-\frac12}\left[\mathcal{C}\left(\sqrt{\mu^c}h_{+},\mu^c\right)+\mathcal{C}\left(\sqrt{\mu^c}h_-,\mu^c\right)\right].
\end{aligned}
\end{align}
And for $\tilde{h}=[\tilde{h}_+, \tilde{h}_-]^{\textit{t}}$, the nonlinear relativistic Landau collision operator $\varGamma(h,\tilde{h})$ is defined as
\begin{align}\label{Nonlin-c}
\begin{aligned}
\varGamma(h,\tilde{h})=&[\varGamma_+(h,\tilde{h}), \varGamma_-(h,\tilde{h})]^{\textit{t}},\\
\varGamma_+(h,\tilde{h}):=&\left(\mu^c\right)^{-\frac12}\left[\mathcal{C}\left(\sqrt{\mu^c}h_+,\sqrt{\mu^c}\tilde{h}_{+}\right)+\mathcal{C}\left(\sqrt{\mu^c}h_-,\sqrt{\mu^c}\tilde{h}_{+}\right)\right],\\
\varGamma_-(h,\tilde{h}):=&\left(\mu^c\right)^{-\frac12}\left[\mathcal{C}\left(\sqrt{\mu^c}h_-,\sqrt{\mu^c}\tilde{h}_{-}\right)+\mathcal{C}\left(\sqrt{\mu^c}h_+,\sqrt{\mu^c}\tilde{h}_{-}\right)\right].
\end{aligned}
\end{align}
\item [ii)] $[f^{i,c}(t,x,p), E^{i,c}(t,x)]$ with $i=1, 2$ satisfies
\begin{align}\label{mainFi}
\begin{aligned}
&\partial_t f^{i,c} + \frac{cp}{p^0}\cdot \nabla_x f^{i,c}+ E^{i,c}\cdot\frac{cp}{p^0}\sqrt{\mu^c}\zeta_0-\frac{1}{2}\zeta_1\frac{cp}{p^0}\cdot \sum_{\substack{j_1+j_2=i\\0\leq j_1,j_2\leq i}}E^{j_1,c} f^{j_2,c}\\
&\hspace{0.2cm}+\zeta_1 \sum_{\substack{j_1+j_2=i\\0\leq j_1,j_2\leq i}}E^{j_1,c}\cdot\nabla_p f^{j_2,c}+\zeta_1 \frac{p}{p^0}\times  \sum_{j=0}^2\frac{1}{c^j}\bar{B}^{j}\cdot\nabla_p f^{i,c}+\mathcal{L}f^{i,c} \\
&\hspace{0.2cm}=\sum_{\substack{j_1+j_2=i\\0\leq j_1,j_2\leq i}}\varGamma\left(f^{j_1,c}, f^{j_2,c}\right),\\
&\partial_t E^{i,c} =- \int_{\mathbb R^3}\frac{cp}{p^0}\sqrt{\mu^c}\left(f^{i,c}_{+}-f^{i,c}_-\right)\, \mathrm{d} p, \\\rule{0ex}{1.5em}
& \nabla_x\cdot E^{i,c}= \int_{\mathbb R^3}  \sqrt{\mu^c}\left(f^{i,c}_+-f^{i,c}_-\right)\, \mathrm{d}p, \qquad
 \nabla_x \times E^{i,c}=0,\\\rule{0ex}{1.0em}
&f^{i,c}(0,x,p)=f^{i,c}_{0}(x,p),\qquad E^{i,c}(0,x)=E^{i,c}_{0}(x).
\end{aligned}
\end{align}
\item [iii)] $[f^{3,c}(t,x,p), E^{3,c}(t,x), B^{3,c}(t,x)]$ satisfies
\begin{align}\label{mainF3}
\begin{aligned}
&\partial_t f^{3,c} + \frac{cp}{p^0}\cdot \nabla_x f^{3,c}+ E^{3,c}\cdot\frac{cp}{p^0}\sqrt{\mu^c}\zeta_0-\frac{1}{2}\zeta_1\frac{cp}{p^0}\cdot E^{3,c} \sum_{j=0}^2\frac{f^{j,c}}{c^j}-\frac{1}{2}\zeta_1\frac{cp}{p^0}\cdot\sum_{j=0}^2\frac{E^{j,c} }{c^j} f^{3,c}\\
&\hspace{0.2cm}+\zeta_1 E^{3,c}\cdot\sum_{j=0}^2\frac{\nabla_p f^{j,c}}{c^j}+\zeta_1 \sum_{j=0}^2\frac{E^{j,c}}{c^j}\cdot\nabla_p f^{3,c}+\Big(\zeta_1 \frac{p}{p^0}\times  \sum_{j=0}^2\frac{\bar{B}^{j}}{c^j}\Big)\cdot\nabla_p f^{3,c} \\
&\hspace{0.2cm}-\frac{1}{2}\zeta_1\frac{cp}{p^0}\cdot \frac{E^{3,c}}{c^3} f^{3,c}+\zeta_1 \big(E^{3,c}+\frac{p}{p^0}\times B^{3,c}\big)\cdot\frac{\nabla_p f^{3,c}}{c^3}+\mathcal{L}f^{3,c}\\
&\hspace{0.2cm}=\sum_{j=0}^2\frac{1}{c^j}\left[\varGamma\left(f^{j,c}, f^{3,c}\right)+\varGamma\left(f^{3,c}, f^{j,c}\right)\right]+\frac{1}{c^3}\varGamma\left(f^{3,c}, f^{3,c}\right)+\mathcal{Q},\\
&\partial_t E^{3,c}-  c\nabla_x \times B^{3,c} =- \int_{\mathbb R^3}\frac{cp}{p^0}\sqrt{\mu^c}\left(f^{3,c}_{+}-f^{3,c}_-\right)\, \mathrm{d} p, \\\rule{0ex}{1.5em}
&\partial_t B^{3,c}+ c\nabla_x \times E^{3,c}=0,\\\rule{0ex}{1.5em}
& \nabla_x\cdot E^{3,c}= \int_{\mathbb R^3}  \sqrt{\mu^c}\left(f^{3,c}_+-f^{3,c}_-\right)\, \mathrm{d}p, \\\rule{0ex}{1.0em}
&f^{3,c}(0,x,p)=f^{3,c}_{0}(x,p),\qquad E^{3,c}(0,x)=E^{3,c}_{0}(x),\qquad B^{3,c}(0,x)=B^{3,c}_{0}(x),
\end{aligned}
\end{align}
where
\begin{align*}
 \mathcal{ Q}= & \frac{1}{2}\zeta_1\frac{cp}{p^0}\cdot \sum_{\substack{j_1+j_2\geq3\\1\leq j_1,j_2< 3}}\frac{E^{j_1,c} f^{j_2,c}}{c^{j_1+j_2-3}}-\zeta_1 \sum_{\substack{j_1+j_2\geq3\\ 1\leq j_1,j_2< 3}}\frac{E^{j_1,c}\cdot\nabla_p f^{j_2,c}}{c^{j_1+j_2-3}}\\
   & +\sum_{\substack{j_1+j_2\geq3   \\1\leq j_1,j_2<3}}\frac{\varGamma\left(f^{j_1,c}, f^{j_2,c}\right)}{c^{j_1+j_2-3}}.
\end{align*}
The equation of $B^{3,c}(t,x)$ in \eqref{mainF3} implies that $\frac{\mathrm{d}}{\mathrm{d}t}\int_{\mathbb{T}^3} B^{3,c}(t,x)\,\mathrm{d} x=0$. Then there exists a constant $\bar{B}^{3}$ such that
\begin{align}\label{B3-infty}
\bar{B}^{3}:=\frac{\int_{\mathbb{T}^3} B^{3,c}(t,x)\,\mathrm{d} x}{|\mathbb{T}^3|}.
\end{align}
\end{itemize}

\subsubsection{Properties of linearized operator $\mathcal{L}$} It is well known that the null space $\mathcal {N}$ of the linearized operator $\mathcal{L}$ is given by \cite{Strain-Guo-CMP-2004}:
\begin{align*}
    \mathcal {N}=\mbox{span}\left\{\begin{pmatrix} 1\\0\\\end{pmatrix}\mu^c , \begin{pmatrix} 0\\1\\\end{pmatrix}\mu^c ,   \begin{pmatrix} p_i\\p_i\\\end{pmatrix}\mu^c  (1\leq i\leq3), \begin{pmatrix} p^0\\p^0\\\end{pmatrix}\mu^c \right\}.
\end{align*}
For $h=[h_+,h_-]^{\textit{t}}$, we introduce the well-known micro-macroscopic decomposition
$$h=\mathcal{P}h+\{I-\mathcal{P}\}h.$$
Here  $\mathcal{P}h=[\mathcal{P}_+h, \mathcal{P}_-h]^{\textit{t}}$ is the orthogonal projection from $L^2_p$ onto $\mathcal{N}$:
\begin{align}\label{macfe}
\mathcal{P}_{\pm}h=\mathbf{a}^{h}_{\pm}\sqrt{\mu^c}+  \mathbf{b}^{h}\cdot \frac{p}{C_{\mathbf{b}}} \sqrt{\mu^c}+ \mathbf{c}^{h}\frac{p^{0}-C_0}{C_{\mathbf{c}}}\sqrt{\mu^c},
\end{align}
where $\mathbf{a}^{h}_{\pm}$, $\mathbf{b}^{h}$, and $\mathbf{c}^{h}$ are coefficients, explicit form of constants are as follows
\begin{align*}
   C_{0}=&\int_{{\mathbb R}^3} p^0\mu^c \,\mathrm{d} p=c\frac{K_3(c^2)}{K_2(c^2)} -\frac{1}{c} ,\quad
   C_{\mathbf{b}}=\Big(\int_{{\mathbb R}^3}\frac{1}{3}|p|^2\mu^c \,\mathrm{d} p\Big)^{\frac{1}{2}} =\Big(\frac{K_3(c^2)}{K_2(c^2)}\Big)^{\frac{1}{2}},\\
   C_{\mathbf{c}}=&\Big(\int_{{\mathbb R}^3}\left(p^0\right)^2\mu^c \,\mathrm{d} p-C_{0}^2\Big)^{\frac{1}{2}} =\Big(-c^2\frac{K^2_3(c^2)}{K^2_2(c^2)}+5\frac{K_3(c^2)}{K_2(c^2)}+c^2-\frac{1}{c^2}\Big)^{\frac{1}{2}}.
\end{align*}
\begin{remark}\label{macro-r}
 Note that by \eqref{transform} and \eqref{acurate} in Lemma \ref{K01p},
\begin{align*}
  \frac{K_3(c^2)}{K_2(c^2)}= &\frac{4}{c^2}+\frac{K_1(c^2)}{K_2(c^2)}=\frac{4}{c^2}+\frac{1}{2/c^2
  +K_0(c^2)/K_1(c^2)}=1+\frac{5}{2c^2}
  +\frac{15}{8c^4}-\frac{15}{8c^6}+\frac{O(1)}{c^{8}}, \\
  C_{\mathbf{b}}=&1+\frac{5}{2c^2}+\frac{O(1)}{c^{4}},\quad C_{0}=c+\frac{3}{2c}
  +\frac{15}{8c^3}-\frac{15}{8c^5}+\frac{O(1)}{c^{7}},\quad C_{\mathbf{c}}^2=\frac{3}{2c^2}+\frac{O(1)}{c^{4}}.
\end{align*}
    Then we have
   \begin{align*}
     \frac{p^0-C_{0}}{C_{\mathbf{c}}}&=\frac{p^0-c-3/(2c)+O(1)/c^3}
  {\sqrt{3/(2c^2)+O(1)/c^4}}=\frac{|p|^2-3}{\sqrt{6}}+\frac{O(1)}{c^2}.
   \end{align*}
  This together with \eqref{diff-mu1} implies that we can write \eqref{macfe} as
  \begin{align*}
  \mathcal{P}_{\pm}h=&\mathbf{a}^{h}_{\pm}\sqrt{\mu^c}+  \left(1+\frac{O(1)}{c^2}\right)\mathbf{b}^{h}\cdot p \sqrt{\mu^c}+ \mathbf{c}^{h}\left(\frac{|p|^2-3}{\sqrt{6}}+\frac{O(1)}{c^2}\right)\sqrt{\mu^c}\\
  =&\mathbf{a}^{h}_{\pm}\sqrt{\mu^{\infty}}+\mathbf{b}^{h}\cdot p \sqrt{\mu^{\infty}}+ \mathbf{c}^{h}\frac{|p|^2-3}{\sqrt{6}}\sqrt{\mu^{\infty}}+\frac{O(1)}{c^2}.
  \end{align*}
  Moreover,  we can further obtain
  \begin{align}\label{rp-cp}
    |\mathcal{P}h|_{L^2}= \left(1+\frac{O(1)}{c^2}\right)|{\bf P}h|_{L^2},
  \end{align}
  where ${\bf P}h$ is the orthogonal projection for the classical case introduced in the next subsection.
\end{remark}
For convenience, we denote the macroscopic part of $f^{i,c}$ with $i=0, 1, 2, 3$ as
\begin{align*}
\mathcal{P}_{\pm}f^{i,c}=\mathbf{a}^{i}_{\pm}\sqrt{\mu^c}+  \mathbf{b}^{i}\cdot \frac{p}{C_{\mathbf{b}}} \sqrt{\mu^c}+ \mathbf{c}^{i}\frac{p^{0}-C_0}{C_{\mathbf{c}}}\sqrt{\mu^c}.
\end{align*}

\subsubsection{Reformulation of the VPL system}
For the VPL system, we  introduce the standard perturbation $f^{\infty}_{\pm}(t,x,p)$ to $\mu^{\infty}$ as
\begin{align*}
f^{\infty}_{\pm}:=\left(\mu^{\infty}\right)^{-\frac12}\left(F^{\infty}_{\pm}-\mu^{\infty}\right).
\end{align*}
Denote $f^{\infty}(t,x,p)=[f^{\infty}_+(t,x,p), f^{\infty}_-(t,x,p)]^{\textit{t}}$. We can rewrite the Cauchy problem of the VPL system \eqref{main2-00} as
\begin{align}\label{main2}
\begin{aligned}
&\partial_t f^{\infty} + p\cdot \nabla_x f^{\infty}+E^{\infty}\cdot p\sqrt{\mu^{\infty}}\zeta_0+\zeta_1  E^{\infty}\cdot\nabla_p f^{\infty}\\
&\hspace{0.2cm}-\frac{1}{2}\zeta_1p\cdot E^{\infty} f^{\infty}+\mathbf{L} f^{\infty} =\Gamma\left(f^{\infty}, f^{\infty}\right),\\
&\partial_t E^{\infty} =- \int_{\mathbb R^3} p\big(Z\mathrm{e} F^{\infty}_+-\mathrm{e} F^{\infty}_-\big)\, \mathrm{d} p,\qquad \nabla_x \times E^{\infty} =0,\\\rule{0ex}{1.0em}
& \nabla_x\cdot E^{\infty}= \int_{\mathbb R^3}  \sqrt{\mu^{\infty}}\left(f^{\infty}_+-f^{\infty}_-\right)\, \mathrm{d}p,\\\rule{0ex}{1.0em}
&f^{\infty}_{\pm}(0,x,p)=f^{\infty}_{\pm,0}(x,p).
\end{aligned}
\end{align}
Here the linear Landau collision operator $\mathbf{L}h$  is given by
\begin{align}\label{linearAK-inf}
\begin{aligned}
&\mathbf{L}h=[\mathbf{L}_+h, \mathbf{L}_-h]^{\textit{t}}, \qquad \mathbf{L}_{\pm}h:=\mathbf{A}_{\pm}h+\mathbf{K}_{\pm}h,\\
&\mathbf{A}_{\pm}h:=-2\left(\mu^{\infty}\right)^{-\frac12}\mathbb{C}\left(\mu^{\infty},\sqrt{\mu^{\infty}}h_{\pm}\right),\\
&\mathbf{K}_{+}h=\mathbf{K}_{-}h:=-\left(\mu^{\infty}\right)^{-\frac12}\left[\mathbb{C}\left(\sqrt{\mu^{\infty}}h_{+},\mu^{\infty}\right)+\mathbb{C}\left(\sqrt{\mu^{\infty}}h_-,\mu^{\infty}\right)\right].
\end{aligned}
\end{align}
And for $\tilde{h}=[\tilde{h}_+,\tilde{h}_-]^{\textit{t}}$, the nonlinear Landau collision operator $\Gamma(h,\tilde{h})$ is defined as
\begin{align}\label{Nonlin-inf}
\begin{aligned}
\Gamma(h,\tilde{h})=&[\Gamma_+(h,\tilde{h}), \Gamma_-(h,\tilde{h})]^{\textit{t}},\\
\Gamma_+(h,\tilde{h}):=&\left(\mu^{\infty}\right)^{-\frac12}\left[\mathbb{C}\left(\sqrt{\mu^{\infty}}h_+,\sqrt{\mu^{\infty}}\tilde{h}_{+}\right)+\mathbb{C}\left(\sqrt{\mu^{\infty}}h_-,\sqrt{\mu^{\infty}}\tilde{h}_{+}\right)\right],\\
\Gamma_-(h,\tilde{h}):=&\left(\mu^{\infty}\right)^{-\frac12}\left[\mathbb{C}\left(\sqrt{\mu^{\infty}}h_-,\sqrt{\mu^{\infty}}\tilde{h}_{-}\right)+\mathbb{C}\left(\sqrt{\mu^{\infty}}h_+,\sqrt{\mu^{\infty}}\tilde{h}_{-}\right)\right].
\end{aligned}
\end{align}
\subsubsection{Properties of linearized operator $\mathbf{L}$}It is well known that the null space of the linearized operator $\mathbf{L}$ is given by \cite{Guo-JAMS-2012}
\begin{align*}
    \mathbb {N}=\mbox{span}\left\{\begin{pmatrix} 1\\0\\\end{pmatrix}\mu^{\infty} , \begin{pmatrix} 0\\1\\\end{pmatrix}\mu^{\infty} ,   \begin{pmatrix} p_i\\p_i\\\end{pmatrix}\mu^{\infty}  (1\leq i\leq3),  \begin{pmatrix}|p|^{2}\\|p|^{2}\\\end{pmatrix}\mu^{\infty} \right\}.
\end{align*}
Correspondingly, for $h=[h_+,h_-]^{\textit{t}}$, we can also introduce the micro-macroscopic decomposition
$$h=\mathbf{P}h+\{I-\mathbf{P}\}h,$$
where $\mathbf{P}h=[\mathbf{P}_+h, \mathbf{P}_-h]^{\textit{t}}$ is the orthogonal projection from $L^2_p$ onto $\mathbb{N}$.

\subsubsection{Conservation of the system} 

Suppose $[F^c(t, x, p), E^c(t, x), B^c(t, x)]$ and $[F^{\infty}(t, x, p), E^{\infty}(t,x)]$ have the same mass, total momentum and total energy
as their equilibriums $[\begin{pmatrix}\begin{smallmatrix} \mu^c\\\mu^c\\\end{smallmatrix}\end{pmatrix}, \mathbf{0}, \sum_{i=0}^3\frac{\bar{B}^{i}}{c^i}]$ and  $[\begin{pmatrix}\begin{smallmatrix} \mu^{\infty}\\\mu^{\infty}\\\end{smallmatrix}\end{pmatrix}, \mathbf{0}]$, respectively.
Then, corresponding to  conservation laws \eqref{cons-RVML} and \eqref{cons-VPL}, we assume that
\begin{align}\label{cons-RVML0}
    &\int_{\mathbb{T}^3\times\mathbb{R}^3} \sqrt{\mu^{c}}\sum_{i=0}^3\frac{f^{i,c}_+(t)}{c^i}\,\mathrm{d} p\mathrm{d} x=\int_{\mathbb{T}^3\times\mathbb{R}^3} \sqrt{\mu^{c}}\sum_{i=0}^3\frac{f^{i,c}_-(t)}{c^i}\,\mathrm{d} p\mathrm{d} x\equiv0,\nonumber\\
    &\int_{\mathbb{T}^3\times\mathbb{R}^3}p \sqrt{\mu^{c}}\sum_{i=0}^3\frac{f^{i,c}_{+}(t)+f^{i,c}_{-}(t)}{c^i}\,\mathrm{d} p\mathrm{d} x +\frac{1}{ c}\int_{\mathbb{T}^3}\sum_{i,j=0}^3\frac{E^{i,c}(t)\times B^{j,c}(t)}{c^{i+j}}\,\mathrm{d} x\equiv\mathbf{0},\\
    &\int_{\mathbb{T}^3\times\mathbb{R}^3}p^0\sqrt{\mu^{c}}\sum_{i=0}^3\frac{f^{i,c}_{+}(t)+f^{i,c}_{-}(t)}{c^i}\,\mathrm{d} p\mathrm{d} x+\frac{1}{2c} \int_{\mathbb{T}^3} \Big(\Big|\sum_{i=0}^3\frac{E^{i,c}(t)}{c^i}\Big|^2+\Big|\frac{B^{3,c}(t)-\bar{B}^{3}}{c^3}\Big|^2\Big)\,\mathrm{d} x\equiv0,\nonumber
\end{align}
and
\begin{align}\label{cons-VPL0}
&\int_{\mathbb{T}^3\times\mathbb{R}^3} \sqrt{\mu^{\infty}}f^{\infty}_{+}(t)\,\mathrm{d} p\mathrm{d} x=\int_{\mathbb{T}^3\times\mathbb{R}^3} \sqrt{\mu^{\infty}}f^{\infty}_{-}(t)\,\mathrm{d} p\mathrm{d} x\equiv0,\nonumber\\
    &\int_{\mathbb{T}^3\times\mathbb{R}^3}p \sqrt{\mu^{\infty}}\left(f^{\infty}_{+}(t)+f^{\infty}_{-}(t)\right)\,\mathrm{d} p\mathrm{d} x \equiv\mathbf{0},\\
    &\int_{\mathbb{T}^3\times\mathbb{R}^3}\frac{|p|^2}{2}\sqrt{\mu^{\infty}}\left(f^{\infty}_{+}(t) + f^{\infty}_{-}(t)\right)\,\mathrm{d} p\mathrm{d} x+\frac{1}{2} \int_{\mathbb{T}^3} |E^{\infty}(t)|^2\,\mathrm{d} x\equiv0.\nonumber
\end{align}

From \eqref{mainF0}, \eqref{mainFi}, and \eqref{mainF3}, we can obtain
\begin{align}\label{cons-RVML00}
    &\int_{\mathbb{T}^3\times\mathbb{R}^3} \sqrt{\mu^{c}}f^{0,c}_+(t)\,\mathrm{d} p\mathrm{d} x=\int_{\mathbb{T}^3\times\mathbb{R}^3} \sqrt{\mu^{c}}f^{0,c}_-(t)\,\mathrm{d} p\mathrm{d} x\equiv0,\nonumber\\
    &\int_{\mathbb{T}^3\times\mathbb{R}^3}p \sqrt{\mu^{c}}\left[f^{0,c}_{+}(t)+f^{0,c}_{-}(t)\right]\,\mathrm{d} p\mathrm{d} x +\frac{1}{ c}\int_{\mathbb{T}^3}E^{0,c}(t)\times \sum_{j=0}^2 \frac{\overline{B}^{j}}{c^j}\,\mathrm{d} x\equiv\mathbf{0},\\
    &\int_{\mathbb{T}^3\times\mathbb{R}^3}p^0\sqrt{\mu^{c}}\left[f^{0,c}_{+}(t)+f^{0,c}_{-}(t)\right]\,\mathrm{d} p\mathrm{d} x+\frac{1}{2c} \int_{\mathbb{T}^3} \big|E^{0,c}(t)\big|^2\,\mathrm{d} x\equiv0,\nonumber
\end{align}
\begin{align}\label{cons-RVML0i}
    &\int_{\mathbb{T}^3\times\mathbb{R}^3} \sqrt{\mu^{c}}f^{i,c}_+(t)\,\mathrm{d} p\mathrm{d} x=\int_{\mathbb{T}^3\times\mathbb{R}^3} \sqrt{\mu^{c}}f^{i,c}_-(t)\,\mathrm{d} p\mathrm{d} x\equiv0,\qquad i=1, 2,\nonumber\\
    &\int_{\mathbb{T}^3\times\mathbb{R}^3}p \sqrt{\mu^{c}}\left[f^{i,c}_{+}(t)+f^{i,c}_{-}(t)\right]\,\mathrm{d} p\mathrm{d} x +\frac{1}{ c}\int_{\mathbb{T}^3}E^{i,c}(t)\times \sum_{j=0}^2 \frac{\overline{B}^{j}}{c^j}\,\mathrm{d} x\equiv\mathbf{0},\\
    &\int_{\mathbb{T}^3\times\mathbb{R}^3}p^0\sqrt{\mu^{c}}\left[f^{i,c}_{+}(t)+f^{i,c}_{-}(t)\right]\,\mathrm{d} p\mathrm{d} x+\frac{1}{2c} \int_{\mathbb{T}^3} \sum_{j_1+j_2=i}\big[E^{j_1,c}(t)\cdot E^{j_2,c}(t)\big]\,\mathrm{d} x\equiv0,\nonumber
\end{align}
and 
\begin{align}\label{cons-RVML03}
    &\int_{\mathbb{T}^3\times\mathbb{R}^3} \sqrt{\mu^{c}}f^{3,c}_+(t)\,\mathrm{d} p\mathrm{d} x=\int_{\mathbb{T}^3\times\mathbb{R}^3} \sqrt{\mu^{c}}f^{3,c}_-(t)\,\mathrm{d} p\mathrm{d} x\equiv0,\nonumber\\
    &\int_{\mathbb{T}^3\times\mathbb{R}^3}p \sqrt{\mu^{c}}\left[f^{3,c}_{+}(t)+f^{3,c}_{-}(t)\right]\,\mathrm{d} p\mathrm{d} x +\frac{1}{ c}\int_{\mathbb{T}^3}\sum_{j_1+j_2\geq3}\frac{E^{j_1,c}(t)\times B^{j_2,c}(t)}{c^{j_1+j_2}}\,\mathrm{d} x\nonumber\\
    &\hspace{5.7cm}+\frac{1}{ c}\int_{\mathbb{T}^3}\sum_{j=0}^3 \frac{E^{j,c}(t)}{c^j}\times B^{3,c}(t)\,\mathrm{d} x\equiv\mathbf{0},\\
    &\int_{\mathbb{T}^3\times\mathbb{R}^3}p^0\sqrt{\mu^{c}}\left[f^{3,c}_{+}(t)+f^{3,c}_{-}(t)\right]\,\mathrm{d} p\mathrm{d} x+\frac{1}{2c} \int_{\mathbb{T}^3} \Big[\sum_{j_1+j_2\geq3}\frac{E^{j_1,c}(t)\cdot E^{j_2,c}(t)}{c^{j_1+j_2-3}}+\frac{\big|B^{3,c}(t)-\bar{B}^{3}\big|^2}{c^3}\Big]\,\mathrm{d} x\equiv0.\nonumber
\end{align}

As in \eqref{macfe}, we denote $[\mathbf{a}^i_{\pm}, \mathbf{b}^i, \mathbf{c}^i]$ as the corresponding coefficients of $\mathcal{P}_{\pm} f^{i,c}$ for $i=0, 1, 2, 3$. By direct computations, we can obtain from \eqref{cons-RVML00}, \eqref{cons-RVML0i}, and \eqref{cons-RVML03}  that
\begin{align}\label{cons-RVMLi}
    &\int_{\mathbb{T}^3} \mathbf{a}^{i}_+(t)\,\mathrm{d} x=\int_{\mathbb{T}^3} \mathbf{a}^{i}_-(t)\,\mathrm{d} x=0,\qquad i=0, 1, 2, 3,\nonumber\\
    &\Big|\int_{\mathbb{T}^3}\mathbf{b}^{i}(t)\,\mathrm{d} x\Big|\leq \frac{C}{c} \Big( \sum_{j=0}^i\|E^{j,c}(t)\|+\delta_{i3}\sum_{j=0}^3\frac{\|E^{j,c}(t)\|}{c^j} \big(\big\|B^{3,c}(t)-\bar{B}^{3}\big\|+|\bar{B}^{3}|\big)\Big),\\
    &\Big|\int_{\mathbb{T}^3}\mathbf{c}^{i}(t)\,\mathrm{d} x\Big|\leq \frac{C}{c} \Big( \sum_{j=0}^i\big\|E^{j,c}(t)\big\|^2+\delta_{i3} \big\|B^{3,c}(t)-\bar{B}^{3}\big\|^2\Big),\nonumber
\end{align}
where $C$ is a constant independent of $c$.

\subsection{Notations}
Throughout this paper, we adopt the following notations:

\begin{itemize}
    \item \textbf{Constants and Inequalities:} The letter $C$ denotes a generic positive constant independent of the light speed $c$, whose value may change from line to line. We write $a \lesssim b$ if $a \leq C b$ for such a constant $C$, and $a \approx b$ if both $a \lesssim b$ and $b \lesssim a$ hold.  Finally the symbol $\backsimeq$ is used to denote the equivalent relationship between two functionals.

    \item \textbf{Derivatives and Multi-indices:} We use the multi-indices $\alpha = (\alpha_1, \alpha_2, \alpha_3)$ and $\beta = (\beta_1, \beta_2, \beta_3)$ to denote spatial and momentum derivatives, respectively. Define
    \[
    \partial^\alpha_\beta = \partial_{x_1}^{\alpha_1} \partial_{x_2}^{\alpha_2} \partial_{x_3}^{\alpha_3} \partial_{p_1}^{\beta_1} \partial_{p_2}^{\beta_2} \partial_{p_3}^{\beta_3}.
    \]
    We write $\partial^\alpha$ when $\beta = 0$, and $\partial_\beta$ when $\alpha = 0$. The length of $\alpha$ is $|\alpha| = \alpha_1 + \alpha_2 + \alpha_3$. For two multi-indices, $\alpha' \leq \alpha$ means that each component of $\alpha'$ is no greater than the corresponding component of $\alpha$, and $\alpha' < \alpha$ means that $\alpha' \leq \alpha$ and $|\alpha'| < |\alpha|$. For brevity, we set $\partial_{x_i} = \partial_i$ for $i = 1, 2, 3$.

    \item \textbf{Function Spaces:} The inner product $(\cdot, \cdot)$ denotes the $L^2$ inner product over $\mathbb{R}^3_p$, with the corresponding norm $|\cdot|_{L^2}$. The notation $\langle \cdot, \cdot \rangle$ and the norm $\| \cdot \|$ are used for the $L^2$ inner product and norm over either $\mathbb{T}^3_x \times \mathbb{R}^3_p$ or $\mathbb{T}^3_x$, with the specific domain being clear from the context. For any integer $k \geq 0$ and $1 \leq r \leq \infty$, $W^{k,r}$ denotes the standard Sobolev space over $\mathbb{T}^3 \times \mathbb{R}^3$ or $\mathbb{T}^3$, and we set $H^k = W^{k,2}$, with the norm
    \[
    \| f \|_{H^k}^2 = \sum_{|\alpha| = 0}^k \| \partial^\alpha f \|^2.
    \]
    Finally, $\nabla_x = (\partial_{x_1}, \partial_{x_2}, \partial_{x_3})$ and $\nabla_p = (\partial_{p_1}, \partial_{p_2}, \partial_{p_3})$.
\end{itemize}

\subsection{Function spaces for the RVML system}
We first define two important things:

\noindent$\bullet$ We recall the collision kernel $\Phi^{c,ij}(p,q)=\Phi^{c,ij}_{\mp,\pm}(p,q)$ from \eqref{cker} and define the collision frequency
$$
\sigma^{c,ij}(p):=\int_{{\R}^{3}}\Phi^{c,ij}(p,q)\mu^c(q)\, \mathrm{d} q.
$$

\noindent$\bullet$ For constants $0\leq\vartheta\leq 1$ and $\ell\geq0$, we introduce the weight function
\begin{align}\label{tt01}
w_{\ell}(p):=\langle p\rangle^{\frac{5}{2}\ell}\exp\left(\frac{\vartheta\langle p\rangle}{\ln(\mathrm{e}+t)}\right), \quad \langle p\rangle:=\sqrt{1+|p|^2}.
\end{align}
Here and throughout, we suppress the $t$-dependence in $w_{\ell_i}$ for notational simplicity.

\subsubsection{Sobolev norms}
For   $l\in \mathbb{R}$, $m\in \mathbb{N}_0$, $\ell_i\geq m$ with $i=0, 1, 2, 3$, and function $h(t,x,p)=[h_+(t,x,p), h_-(t,x,p)]^{\textit{t}}$, we  define weighted Sobolev norms
\begin{align*}
   |h(t,x)|^2_{H^m_{l}}:&=\sum_{i=0}^m\sum_{\pm}\left|\lag p\rag^{l}\nabla_x^i h_{\pm}(t,x)\right|^2_{L^2},\qquad  \|h(t)\|^2_{H^m_{l}}:=\sum_{i=0}^m\sum_{\pm}\left\|\lag p\rag^{l}\nabla_x^i h_{\pm}(t)\right\|^2,\\
 |h(t,x)|^2_{H^{m,i}_{w}}:&=\sum_{j=0}^m\sum_{\pm}\left|w_{\ell_i-j}\nabla_x^j h_{\pm}(t,x)\right|^2_{L^2},\qquad  \|h(t)\|^2_{H^{m,i}_{w}}:=\sum_{j=0}^m\sum_{\pm}\left\|w_{\ell_i-j}\nabla_x^j h_{\pm}(t)\right\|^2.  
\end{align*}
When $m=0$, we denote these norms by $|\cdot|_{L^2_{l}}, \|\cdot\|_{L^2_{l}}, |\cdot|_{L^{2,i}_{w}}$, and $\|\cdot\|_{L^{2,i}_{w}}$, respectively.

\subsubsection{Dissipation norms}
We define weighted dissipation norms in momentum space: 
\begin{align}
|h(t,x)|_{\sigma, l}^2&:=\sum_{\pm}\int_{\mathbb{R}^3}
\sigma^{c,ij}(p)\lag p\rag^{2l}\Big[2\big(\partial_{p_i}h_{\pm}\partial_{p_j}h_{\pm}\big)(t, x, p)
+\frac{1}{2}\frac{cp_i}{p^0}\frac{cp_j}{p^0}|h_{\pm}|^2(t, x, p)\Big]\,\mathrm{d} p, \label{l-norm}\\
|h(t,x)|_{\sigma, w_{\ell_i}}^2&:=\sum_{\pm}\int_{\mathbb{R}^3}
\sigma^{c,ij}(p)w^2_{\ell_i}(p)\Big[2\big(\partial_{p_i}h_{\pm}\partial_{p_j}h_{\pm}\big)(t, x, p)
+\frac{1}{2}\frac{cp_i}{p^0}\frac{cp_j}{p^0}|h_{\pm}|^2(t, x, p)\Big]\,\mathrm{d} p, \label{ell-norm}
\end{align}
%
The corresponding dissipation norms on $\mathbb{T}^3\times\mathbb{R}^3$ are:
\begin{eqnarray*}
    \|h(t)\|_{\sigma, l}^2:=\int_{{\mathbb T}^3}|h(t,x)|_{\sigma,l}^2\, \mathrm{d}x,\qquad \|h(t)\|_{\sigma, w_{\ell_i}}^2:=\int_{{\mathbb T}^3}|h(t,x)|_{\sigma,w_{\ell_i}}^2\, \mathrm{d}x.
\end{eqnarray*}
We also define the Sobolev dissipation norms:
\begin{align*}
    \|h(t)\|^2_{H^m_{\sigma,l}}:=\sum_{|\alpha|=0}^m\left\|\partial^{\alpha}h(t)\right\|^2_{\sigma,l},\qquad  \|h(t)\|^2_{H^{m,i}_{\sigma,w}}:=
    \sum_{|\alpha|=0}^m\left\|\partial^{\alpha}h(t)\right\|^2_{\sigma,w_{\ell_i-|\alpha|}},\qquad m\in \mathbb{N}_0.
\end{align*}
When $l=0$, we use the simplified notation: $\|\cdot\|_{H^m}$, $|\cdot|_{\sigma}$,  $\|\cdot\|_{\sigma}$, and $\|\cdot\|_{H^m_{\sigma}}$.   In particular, the unweighted dissipation norm is  
\begin{align}\label{Sigma-norm}
|h(t,x)|_{\sigma}^2&:=\sum_{\pm}\int_{\mathbb{R}^3}
\sigma^{c,ij}(p)\Big[2\big(\partial_{p_i}h_{\pm}\partial_{p_j}h_{\pm}\big)(t, x, p)
+\frac{1}{2}\frac{cp_i}{p^0}\frac{cp_j}{p^0}|h_{\pm}|^2(t, x, p)\Big]\,\mathrm{d} p.
\end{align}

\subsubsection{Energy functional and dissipation rate} 
We fix parameters $0<\vartheta\leq \frac{1}{32}$ and $\ell=\ell_i=\ell_3+2(3-i)\geq8$ for $i=0, 1, 2, 3$ in \eqref{tt01}.

\noindent $\bullet$($i$-th order energy functional and dissipation rate with $i=0,1,2$) For the solutions $[f^{0,c}(t,x,p), E^{0,c}(t,x)]$ to \eqref{mainF0} and $[f^{i,c}(t,x,p), E^{i,c}(t,x)]$ ($i=1,2$) to \eqref{mainFi}, we define:
\begin{align} \label{enerfun-ic}
\mathcal{E}^{i,c}(t):=&\|f^{i,c}(t)\|^2_{H^{5,i}_{w}}+\|E^{i,c}(t)\|^2_{H^5},\\
\label{dissirat-ic}
\mathcal{D}^{i,c}(t)\backsimeq&\|f^{i,c}(t)\|^2_{H^{5,i}_{\sigma,w}}
+Y(t)\|\langle p\rangle^{\frac{1}{2}}f^{i,c}(t)\|^2_{H^{5,i}_{w}}+\|\mathbf{a}^i_+-\mathbf{a}^i_-\|^2+\|E^{i,c}(t)\|^2_{H^5},
\end{align}
where $Y(t):=(1+t)^{-1}\ln ^{-2}(\mathrm{e}+t)$.

 For $2\leq m\leq 5$, we also define unweighted versions: 
\begin{align}\label{nowei-edi}
\begin{aligned}
\mathcal{E}^{i,c}_{m}(t):=&\|f^{i,c}(t)\|^2_{H^{m}}+\|E^{i,c}(t)\|^2_{H^m}, \\
\mathcal{D}^{i,c}_{m}(t)\backsimeq&\|f^{i,c}(t)\|^2_{H^{m}_{\sigma}}
+\|f^{i,c}(t)\|^2_{H^{m}}+\|\mathbf{a}^i_+-\mathbf{a}^i_-\|^2+\|E^{i,c}(t)\|^2_{H^m}.
\end{aligned}
\end{align}

\noindent$\bullet$($3$-th order energy functional and dissipation rate) For the solution $[f^{3,c}(t,x,p), E^{3,c}(t,x), B^{3,c}(t,x)]$ to \eqref{mainF3},  we define
\begin{align} \label{enerfun-3c} 
\mathcal{E}^{3,c}(t):=&\mathcal{E}^{3,c}_{5}(t)+\|f^{3,c}(t)\|^2_{H^{4,3}_{w}}
+(1+t)^{-\frac{1+\epsilon_{0}}{2}}\|f^{3,c}(t)\|^2_{H^{5,3}_{w}},\\ \label{dissirat-3c}
\mathcal{D}^{3,c}(t)\backsimeq& \mathcal{D}^{3,c}_{5}(t)+\|f^{3,c}(t)\|^2_{H^{4,3}_{\sigma,w}}+Y(t)\|\langle p\rangle^{\frac{1}{2}}f^{3,c}(t)\|^2_{H^{4,3}_{w}}\\
&+(1+t)^{-\frac{1+\epsilon_{0}}{2}}\|f^{3,c}(t)\|^2_{H^{5,3}_{\sigma,w}}+(1+t)^{-\frac{1+\epsilon_{0}}{2}}Y(t)\|\langle p\rangle^{\frac{1}{2}}f^{3,c}(t)\|^2_{H^{5,3}_{w}},\nonumber
\end{align}
with the unweighted norms:
\begin{align}\label{nowei-ed3}
\begin{aligned}
   \mathcal{E}^{3,c}_{m}(t):=&\big\|f^{3,c}(t)\big\|^2_{H^m}+\big\|E^{3,c}(t)\big\|^2_{H^m}
   +\big\|B^{3,c}(t)-\bar{B}^{3}\big\|^2_{H^m},\qquad 2\leq m\leq5,\\
  \mathcal{D}^{3,c}_{m}(t)\backsimeq& \big\|f^{3,c}(t)\big\|^2_{H^{m}_{\sigma}}+\big\|\big[\mathbf{a}^3_+-\mathbf{a}^3_-\big](t)\big\|^2
  +\frac{1}{c^2}\Big(\big\|E^{3,c}(t)\big\|^2_{H^{m-1}}+\big\|B^{3,c}(t)-\bar{B}^{3}\big\|^2_{H^{m-1}}\Big).
\end{aligned}
\end{align}
Here $\epsilon_0 > 0$ is a small constant and $\bar{B}^3$ represents the background magnetic field.

\subsection{Function spaces for the VPL system}

We now introduce the function spaces and norms specific to the Vlasov-Poisson-Landau (VPL) system. For a function $h(t,x,p)=[h_+(t,x,p), h_-(t,x,p)]^{\text{t}}$, we define the following weighted Sobolev norms:

\begin{align*}
|h(t,x)|^2_{\tilde{H}^m_{l}}:=&\sum_{|\alpha|+|\beta|=0}^m\sum_{\pm}\left|\langle p\rangle^{l-|\alpha|-|\beta|}\partial^{\alpha}_{\beta} h_{\pm}(t,x)\right|^2_{L^2}, \\ 
  \|h(t)\|^2_{\tilde{H}^m_{l}}:=&\sum_{|\alpha|+|\beta|=0}^m\sum_{\pm}\left\|\langle p\rangle^{l-|\alpha|-|\beta|}\partial^{\alpha}_{\beta} h_{\pm}(t)\right\|^2,\end{align*}
  where $\partial^{\alpha}_{\beta}$ denotes derivatives in both spatial and momentum variables.

\subsubsection{Dissipation norms for VPL system}

We define the collision frequency for the VPL system:
$$
\sigma^{\infty,ij}(p):=\int_{{\R}^{3}}\Phi^{\infty,ij}(p,q)\mu^{\infty}(q)\, \mathrm{d} q.
$$
The weighted dissipation norms are given by:
  \begin{align*} |h(t,x)|_{\bm{\sigma}, l}^2:=&\sum_{\pm}\int_{\mathbb{R}^3}
\sigma^{\infty,ij}(p)\langle p\rangle^{2l}\Big[2\big(\partial_{p_i}h_{\pm}\partial_{p_j}h_{\pm}\big)(t, x, p)
+\frac{1}{2}p_ip_j|h_{\pm}(t, x, p)|^2(p)\Big]\,\mathrm{d} p, \\ \|h(t)\|^2_{\tilde{H}^m_{\bm{\sigma},l}}:=&\sum_{|\alpha|+|\beta|=0}^m\sum_{\pm}\int_{{\mathbb T}^3}\left|\partial^{\alpha}_{\beta} h_{\pm}(t,x)\right|_{\bm{\sigma},l-|\alpha|-|\beta|}^2\, \mathrm{d}x,\qquad l\geq m,
\end{align*}
 
\subsubsection{Unweighted dissipation norm}

For the case $l=0$, we have the equivalent characterization:
\begin{eqnarray*}
|h(t,x)|_{\bm{\sigma}}^2&:=&\sum_{\pm}\int_{\mathbb{R}^3}
\sigma^{\infty,ij}(p)\Big[2\big(\partial_{p_i}h_{\pm}\partial_{p_j}h_{\pm}\big)(t, x, p)
+\frac{1}{2}p_ip_j|h_{\pm}(t, x, p)|^2(p)\Big]\,\mathrm{d} p\\
&\backsimeq&\sum_{\pm} \Big(\Big|\langle p\rangle^{-\frac{1}{2}} h_{\pm} (t,x)\Big|^2_{L^2}+\Big|\langle p\rangle^{-\frac{3}{2}}\frac{p}{|p|}\cdot\nabla_{p} h_{\pm} (t,x)\Big|^2_{L^2}\nonumber\\
    &&\qquad+ \Big|\langle p\rangle^{-\frac{1}{2}}\frac{p}{|p|}\times\nabla_{p} h_{\pm} (t,x)\Big|^2_{L^2}\Big).
\end{eqnarray*}

\subsubsection{Energy and dissipation functionals}

For the solution $[f^{\infty}(t,x,p), E^{\infty}(t,x)]$ to the Cauchy problem \eqref{main2}, we define the energy functional and dissipation rate:
\begin{align*}
\mathcal{E}^{\infty}_{m,l}(t):=&\|f^{\infty}(t)\|^2_{\tilde{H}^m_{l}}+\|E^{\infty}(t)\|^2_{H^m},\\
\mathcal{D}^{\infty}_{m,l}(t)\backsimeq&\|f^{\infty}(t)\|^2_{\tilde{H}^m_{\bm{\sigma},l}}+\|E^{\infty}(t)\|^2_{H^m}.
\end{align*}

These norms and functionals are tailored to capture the specific analytical structure of the Vlasov-Poisson-Landau system and will be used in establishing the uniform energy estimates.
\subsection{Main results}
In this part, we present two main results in our paper: the global existence of unique solution $[F^c(t,x,p), E^c(t,x), B^c(t,x)]$ to the Cauchy problem of the RVML system \eqref{main1-00} with uniform energy estimates independent of light speed $c\geq1$, and the classical limit from the RVML system \eqref{main1-00}  to the VPL system \eqref{main2-00}.

Let us  recall a well-known existence result of the VPL system \eqref{main2} established in \cite{Guo-JAMS-2012}.

\begin{proposition}\label{guo-VPL} Assume that $F^{\infty}_{\pm,0}(x,p)=\mu^{\infty}(p)+\sqrt{\mu^{\infty}(p)}f^{\infty}_{\pm,0}(x,p)\geq0$ and $[f^{\infty}_{\pm,0}(x,p),E^{\infty}_0(x)]$ satisfies conservation laws \eqref{cons-VPL0}.
There exists a sufficiently small constant $\overline{M} > 0$ such that if $\mathcal{E}^{\infty}_{2,2}(0)\leq\overline{M} $,
then the Cauchy problem \eqref{main2} admits a unique global solution $[f^{\infty}(t,x,p), E^{\infty}(t,x)]$ with $F^{\infty}_{\pm}(x,p)=\mu^{\infty}(p)+\sqrt{\mu^{\infty}(p)}f^{\infty}_{\pm}(t,x,p)\geq0$.

Moreover,  for $l\geq m\geq2$, if $\mathcal{E}^{\infty}_{m,l}(0)<\infty$, there exists an increasing continuous function $P_{m,l}(\cdot)$ with $P_{m,l}(0)=0$ and $P_{2,2}\left(\mathcal{E}^{\infty}_{2,2}(0)\right)=\widetilde{C}_2\mathcal{E}^{\infty}_{2,2}(0)$ such that
\begin{align}\label{tdecay-infty}
\begin{aligned}
\mathcal{E}^{\infty}_{m,l}(t)+\int_0^t\mathcal{D}^{\infty}_{m,l}(\tau) \,\mathrm{d}\tau\leq P_{m,l}\left(\mathcal{E}^{\infty}_{m,l}(0)\right),
\qquad\mbox{and} \\
\big\|\partial_t\phi^{\infty}(t)\big\|_{\infty}+\big\|\nabla_x\phi^{\infty}(t)\big\|_{\infty}
+\|f^{\infty}(t)\|_{H^1}\leq C_{l}(1+t)^{-2l+2}\sqrt{\mathcal{E}^{\infty}_{2,l}(0)}.
\end{aligned}
\end{align}
\end{proposition}
 
\begin{remark} Although $E^{\infty}(t,x)=-\nabla_x\phi^{\infty}(t,x)$ can be expressed by $f^{\infty}(t,x,p)$ via the Poisson equation
\begin{align*}
-\Delta_x \phi^{\infty}= \int_{\mathbb R^3}  \sqrt{\mu^{\infty}}\left(f^{\infty}_+-f^{\infty}_-\right)\, \mathrm{d}p,
\end{align*}
we denote $[f^{\infty}(t,x,p), E^{\infty}(t,x)]$ as the solution to the VPL system \eqref{main2}. On the other hand, different from \cite{Guo-JAMS-2012}, the norms of $E^{\infty}(t,x)$ are included in the  definition
 of the energy functionals $\mathcal{E}^{\infty}_{m,\ell}(t)$ and the dissipation rates $\mathcal{D}^{\infty}_{m,l}(t)$ for convenience.
\end{remark}

Now we state our first main result on the global existence of unique solution to the RVML system \eqref{main1-00} for the light speed $c\geq1$.

\begin{theorem}[Global existence of the RVML system]\label{result1RVML}
 Assume that  
\begin{itemize}
\item $c\geq \max_{0\leq i\leq 3}c_i$ and $\ell_i=\ell_3+2(3-i)\geq8$ for $0\leq i\leq 3$;
\item The initial data $[F^c_0(x,p), E^c_0(x), B^c_0(x)]$ satisfy the expansions \eqref{expanion00}, and
$$F^{c}_{\pm,0}(x,p)=\mu^{c}(p)+\sqrt{\mu^{c}(p)}\sum_{i=0}^3\frac{f^{i,c}_{\pm,0}(x,p)}{c^i}\geq0;$$
\item
$[f^{i,c}_{0}(x,p),E^{i,c}_0(x)$, $B^{i,c}_0(x)]$ with $i=0, 1, 2, 3$ satisfy the conservation laws \eqref{cons-RVML0}.
\end{itemize}
There exists a sufficiently small constant $M > 0$, independent of $c$, such that if
$$\sum_{i=0}^3\mathcal{E}^{i,c}(0)<M,$$
 then the Cauchy problem \eqref{main1-00} admits a unique global solution $[F^c(t,x,p), E^c(t,x), B^c(t,x)]$ satisfying:  
\begin{itemize}
\item The conservation laws \eqref{cons-RVML0};
\item $F^c_{\pm}(t,x,p) = \mu^c(p) + \sqrt{\mu^c(p)} \sum_{i=0}^3 c^{-i} f^{i,c}_{\pm}(t,x,p) \geq 0$;
\item The global energy estimate:
 \begin{equation}\label{global}
\sum_{i=0}^3\mathcal{E}^{i,c}(t)+\int_0^t\sum_{i=0}^3\mathcal{D}^{i,c}(\tau) \,\mathrm{d}\tau\leq C_{\ell_0}\sum_{i=0}^3\mathcal{E}^{i,c}(0),\qquad t\in [0, \infty).
\end{equation}
\end{itemize}
In particular, the following results hold:
 \begin{itemize}
\item For any $0 \leq i \leq 2$, $[f^{i,c}(t,x,p), E^{i,c}(t,x), B^{i,c}(t,x)]$ is the unique solution to the Cauchy problem \eqref{mainF0} or \eqref{mainFi} satisfying:
    \begin{align}\label{decay-f012}
\begin{aligned}
\mathcal{E}^{i,c}(t)+\int_0^t\mathcal{D}^{i,c}(\tau) \,\mathrm{d}\tau\leq& C_{\ell_0}\sum_{j=0}^i\mathcal{E}^{j, c}(0)= C_{\ell_i}M_i,  \\
\|f^{i,c}(t)\|^2_{H^5}+\|E^{i,c}(t)\|^2_{H^5}\leq& \mathrm{e}^{-C_{\ell_i}t^{1/3}}\mathcal{E}^{i,c}(0).
\end{aligned}
\end{align}
\item The third-order solution satisfies the decay estimate:
\begin{align}\label{H2-decay}
\|f^{3,c}(t)\|^2_{H^2}+\|E^{3,c}(t)\|^2_{H^2}+\|B^{3,c}(t)-\bar{B}^{3}\|^2_{H^2}\leq C_{\ell_0}c^6(1+t)^{-3}\sum_{i=0}^3\mathcal{E}^{i,c}(0).
\end{align}
\end{itemize}
Here the constant $C_{\ell_0}$ is independent of $c$.
\end{theorem}

\noindent Several remarks are in order:

\begin{remark}[Smallness and Uniformity]
The smallness assumption on the initial data is essential for closing the energy estimates. Notably, the threshold $M > 0$ is independent of the speed of light $c$, which is crucial for obtaining results that hold uniformly in the relativistic parameter.
\end{remark}

\begin{remark}[Decay Properties]
The theorem establishes not only global existence but also provides detailed decay information. The solutions for the lower-order expansions ($i = 0,1,2$) enjoy exponential decay in Sobolev norms, while the third-order term exhibits algebraic decay of order $(1+t)^{-3}$, reflecting its more delicate structure.
\end{remark}

\begin{remark}[Conservation Laws]
The preservation of conservation laws throughout the evolution is a fundamental feature of the construction, ensuring the physical consistency of the solutions and playing a key role in the energy analysis.
\end{remark}

\begin{remark}[Function Space Differences]
We emphasize a key structural difference between the relativistic and classical regimes: in the relativistic setting, the function spaces contain no momentum derivatives, whereas they do appear in the classical case. This distinction arises primarily from the different forms of the equilibrium states in the two regimes.
\end{remark}

\begin{remark}[Global Well-posedness Framework]
The global well-posedness of the Cauchy problems \eqref{mainF0} and \eqref{mainFi} is established in the appendix (see Proposition \ref{RVPL0} and Proposition \ref{RVPLi}), leveraging the key estimates in Lemma \ref{coerci} and Lemma \ref{w-nonlin}. Consequently, our main focus in the proof will be on establishing the solvability of the more challenging system \eqref{mainF3}.
\end{remark}

\begin{remark}[Improved Decay and Physical Interpretation]
We note that the decay estimate \eqref{H2-decay} can be refined as follows: for any $\theta \in [0,1]$,
\begin{align}\label{H2-decayimproved}
c^{-6}\left[\|f^{3,c}(t)\|^2_{H^2} + \|E^{3,c}(t)\|^2_{H^2} + \|B^{3,c}(t) - \bar{B}^{3}\|^2_{H^2}\right] \lesssim c^{-6(1-\theta)}(1 + t)^{-3\theta}.
\end{align}
Combining \eqref{decay-f012} and \eqref{H2-decayimproved}, we observe that the time decay of the normalized fluctuation $\frac{F^c_{\pm}(t,x,p) - \mu^c(p)}{\sqrt{\mu^c(p)}}$ transitions from polynomial to sub-exponential rates in the non-relativistic limit. This transition represents a fundamental distinction between the (VML) and (VPL) systems, highlighting the different long-time behaviors in relativistic versus classical regimes.
\end{remark}

Our second main result concerns the classical limit from the relativistic to the non-relativistic system:

\begin{theorem}[Classical limit of the RVML system]\label{result2RVML}
Assume that:
\begin{itemize}
\item The solution $[F^c(t,x,p), E^c(t,x), B^c(t,x)]$ of the RVML system from Theorem \ref{result1RVML};
\item The solution $[F^{\infty}(t,x,p), E^{\infty}(t,x)]$ of the VPL system \eqref{main2-00} from Proposition \ref{guo-VPL} with $l \geq 16$;
\item The initial data satisfy:
  \begin{align*}
  \|F^{c}_0-F^{\infty}_0\|+\|E^{c}-E^{\infty}_0\|+\big\|B^{c}-\overline{B}^0\big\|\leq \frac{C}{c}.
  \end{align*}
  \end{itemize}

Then the following convergence estimate holds
   \begin{align}\label{c-limit}
  \|F^{c}(t)-F^{\infty}(t)\|+\|E^{c}_0(t)-E^{\infty}_0(t)\|+\big\|B^{c}_0(t)-\overline{B}^0\big\|\leq \frac{C_{\ell_0,l}}{c}.
  \end{align}
\end{theorem}

\begin{remark}[Optimal Convergence Rate]
The estimate provides an explicit convergence rate of $O(1/c)$ for the classical limit from the relativistic Vlasov-Maxwell-Landau system to the Vlasov-Poisson-Landau system. This rate is optimal, as evidenced by the term $\zeta_1 \big(\frac{p}{p^0}\times \sum_{j=0}^2\frac{\bar{B}^{j}}{c^j}\big)\cdot\nabla_p f^{0,c}$ in \eqref{mainf-ci0} for the equation governing $f^{0,c}(t,x,p)-f^{\infty}(t,x,p)$. The presence of these constant magnetic fields $\bar{B}^{j}$ forces a first-order error in $1/c$. This aligns with the optimal rate observed for the relativistic Vlasov-Maxwell to Vlasov-Poisson limit \cite{Schaeffer-CMP-1986}. We note that in the absence of a magnetic field, the convergence of the relativistic Landau operator to its classical counterpart would be of order $O(1/c^2)$.
\end{remark}

\begin{remark}[Uniform Validity and Magnetic Field]
The convergence is uniform in time, a non-trivial result given the different long-time behaviors of the systems. Furthermore, the self-consistent relativistic magnetic field $B^c(t)$ converges to the fixed background field $\overline{B}^0$, confirming the electrostatic nature of the limit.
 \end{remark}

\begin{remark}[Analytical Consistency] The dissipation norm for the relativistic Landau operator \eqref{Sigma-norm}, characterized by \eqref{eqiv-norm}, formally converges to its classical counterpart from \cite{Guo-CMP-2002} as $c \to \infty$. This bridges the analytical structures of the two theories.
\end{remark}

\begin{remark}[Initial Data]
The assumption of initial data closeness is a natural and necessary condition to initialize the two systems in a comparable state, making the subsequent convergence meaningful.
\end{remark}


\section{Basic estimates}

\setcounter{equation}{0}
In this section, we mainly focus on basic estimates on the relativistic Maxwellian $\mu^{c}= \frac{\exp\left\{-cp^0 \right\} }{4 \pi
c  K_2(c^2)}$, the relativistic Landau collision operators $\mathcal{L}h$ \eqref{linearAK-c} and $\varGamma(h,\tilde{h})$ \eqref{Nonlin-c},  and the difference of  collision operators $\mathcal{L}h-\mathbf{L}h$ and $\varGamma(h,\tilde{h})-\Gamma(h,\tilde{h})$. We will also recall some estimates about the collision kernel $\Phi^c(p,q)$ \eqref{cker} and the difference  of the collision kernels $\Phi^c(p,q)-\Phi^{\infty}(p,q)$ from \cite{Strain-Guo-CMP-2004, Lemou-MMAS-2000}, and classical estimates of the Landau collision operators $\mathbf{L}h$ and $\Gamma(h,\tilde{h})$ from \cite{Guo-CMP-2002, Guo-JAMS-2012}.

\subsection{Estimates of the relativistic Maxwellian}
In this part, we estimate  the lower bound and upper bound of the relativistic Maxwellian $\mu^{c}$, and the difference between Maxwellians $\mu^{c}$ and $\mu^{\infty}$.

\begin{lemma}\label{mu-ci} For the relativistic Maxwellian $\mu^{c}$, we have the following bound estimate
\begin{align}\label{upplow-mu}
\Big(1-\frac{4}{c^2}\Big)(2\pi)^{-\frac{3}{2}}\exp\left\{-\frac{|p|^2}{2} \right\}\leq\mu^{c}(p)\leq C\exp\left\{-\frac{3|p|}{8}\min\{|p|,\frac{4c}{3}\} \right\},
\end{align}
where $C>0$ is some constant independent of $c\geq1$.

For the difference between Maxwellians $\mu^{c}$ and $\mu^{\infty}$, we have
\begin{align}
\big|\mu^{c}-\mu^{\infty}\big| \leq& \frac{C}{c^2}\exp\left\{-\frac{|p|}{3}\min\{|p|,\frac{4c}{3}\} \right\},\label{diff-mu0}\\
\big|\sqrt{\mu^{c}}-\sqrt{\mu^{\infty}}\big| \leq& \frac{C}{c^2}\exp\left\{-\frac{|p|}{6}\min\{|p|,\frac{4c}{3}\} \right\}\label{diff-mu1},
\end{align}
where $C>0$ is some constant independent of $c\geq1$.
\end{lemma}
\begin{proof} We first show \eqref{upplow-mu}.
By \eqref{remainder}, we can obtain
\begin{align*}
1+\frac{15}{8c^2}+\frac{7}{128c^4}-\frac{18\exp\left\{\frac{15}{4c^2} \right\}}{3!8^3}\leq \sqrt{\frac{2}{\pi}}c \mathrm{e}^{c^2}K_2(c^2)\leq 1+\frac{15}{8c^2}+\frac{7}{128c^4}+\frac{18\exp\left\{\frac{15}{4c^2} \right\}}{3!8^3},
\end{align*}
which implies that for $c\geq 1$,
\begin{align*}
1\leq\sqrt{\frac{2}{\pi}}c \mathrm{e}^{c^2}K_2(c^2)\leq 1+\frac{4}{c^2}.
\end{align*}
Then we have
\begin{align}\label{muc-lu0}
\Big(1-\frac{4}{c^2}\Big)(2\pi)^{-\frac{3}{2}}\exp\left\{c^2-cp^0 \right\}\leq \frac{\exp\left\{-cp^0 \right\} }{4 \pi
c  K_2(c^2)}\leq (2\pi)^{-\frac{3}{2}}\exp\left\{c^2-cp^0 \right\}.
\end{align}
Note that for $|p|\leq \frac{4c}{3}$,  $p^0\leq \frac{5c}{3}$  and $c+p^0\leq \frac{8c}{3}$. We can further obtain
\begin{align}\label{muc-lu1}
-\frac{|p|^2}{2}\leq c^2-cp^0=\frac{-c|p|^2}{c+p^0}\leq \begin{cases} -\frac{3|p|^2}{8},\quad &|p|\leq \frac{4c}{3},\\
-\frac{c|p|}{2},\quad &|p|\geq \frac{4c}{3} \end{cases} .
\end{align}
Combining \eqref{muc-lu0} and \eqref{muc-lu1} yields \eqref{upplow-mu}.

Now we turn to prove \eqref{diff-mu0}. Note that $1-\mathrm{e}^{-x}\leq x$ for $x\geq0$, and
\begin{align*}
\frac{c-p^0}{2(c+p^0)}|p|^2=\frac{-|p|^4}{2(c+p^0)^2}\leq -\frac{|p|^4}{2c^2}.
\end{align*}
We use  \eqref{muc-lu0} and \eqref{muc-lu1} to have
\begin{align*}
\big|\mu^{c}-\mu^{\infty}\big| \leq& (2\pi)^{-\frac{3}{2}}\exp\left\{c^2-cp^0 \right\}\Big(\frac{4}{c^2}+1-\exp\left\{\frac{c-p^0}{2(c+p^0)}|p|^2 \right\}\Big)\nonumber\\
\leq& (2\pi)^{-\frac{3}{2}}\exp\left\{c^2-cp^0 \right\}\Big(\frac{4}{c^2}+\frac{|p|^4}{8c^2} \Big)\\
\lesssim& \frac{1}{c^2}\exp\left\{-\frac{|p|}{3}\min\{|p|,\frac{4c}{3}\} \right\}.\nonumber
\end{align*}
Then \eqref{diff-mu0} follows. We can obtain \eqref{diff-mu1} in the same way and omit the details for brevity.

\end{proof}

\subsection{Estimates of the collision kernel}
In this subsection, we will derive some estimates about the collision kernel $\Phi^c(p,q)$ \eqref{cker} and the difference  of the collision kernels $\Phi^c(p,q)-\Phi^{\infty}(p,q)$. Before this derivation, we first establish useful estimates in the following three lemmas. 

\begin{lemma}
For any $c \geq 1$, we have
\begin{equation}
\label{p^0q^0 p q lower bound}
2p^0q^0 -2|p||q| \geq \frac {c^2(c^2+|p|^2+|q|^2)}{p^0q^0}\geq c^2  \max \left \{   \frac {p^0}{q^0}, \frac {q^0}{p^0}    \right\} \ge \max \left \{   \frac {\langle p \rangle}{ \langle q \rangle}, \frac {\langle q \rangle}{ \langle p \rangle}  \right \},
\end{equation}
which implies that 
\begin{equation}
\label{p^0q^0 p q lower bound 2}
\frac 1  {2(p^0q^0-p\cdot q +c^2)} \le\frac 1  {2(p^0q^0- |p| | q|) }  \le \frac 1 {c^2}  \min  \left \{   \frac {p^0}{q^0}, \frac {q^0}{p^0}    \right\}   \le  \min \left \{   \frac {\langle p \rangle}{ \langle q \rangle}, \frac {\langle q \rangle}{ \langle p \rangle}  \right \}. 
\end{equation}
\end{lemma}
\begin{proof}
For brevity, we only prove the first inequality in \eqref{p^0q^0 p q lower bound}, which implies \eqref{p^0q^0 p q lower bound 2}. In fact, we have
\[
2p^0q^0 -2|p||q|=\frac {2c^2(c^2+|p|^2+|q|^2)}{p^0q^0 +|p||q|} \frac {c^2(c^2+|p|^2+|q|^2)}{p^0q^0}.
\]
\end{proof}

Now 
we define $s$ the square of the energy in the center-of-momentum system $p+q=0$, and $g$ the relative momentum  by
\begin{equation}
\label{definition g s}
g = \sqrt{2 (p^0q^0 - p \cdot q - c^2)}, \quad s =2 (p^0q^0 - p \cdot q + c^2), \quad  s =g^2+4c^2.
\end{equation}

\begin{lemma}For $g,s$ defined in \eqref{definition g s}, we have
\begin{eqnarray}
c\frac {|p-q|}  {\sqrt{p^0 q^0}} \le& g \leq |p-q|, \quad   g \leq& 5 \min \left \{ \frac {\langle p \rangle \langle q \rangle \sqrt{q^0}} {\sqrt{p^0}} ,  \frac {\langle p \rangle \langle q \rangle \sqrt{p^0}} {\sqrt{q^0}} \right \}, \label{upper and lower bound g}\\
4c^2 \leq & s \leq 4p^0q^0,
\quad  s\geq & c^2  \max \left \{   \frac {p^0}{q^0}, \frac {q^0}{p^0}   \right\}.
\label{upper and lower bound s}
\end{eqnarray}
Thus we can further obtain
\begin{eqnarray}
\frac 1 4    |p-q|^2  \max\{ \langle p \rangle^{-2}, \langle q \rangle^{-2}  \}  \leq&   \frac {(p^0q^0 - p \cdot q)^2-c^4 } {c^2}  \leq & 25\langle p  \rangle^{2}   \langle q \rangle^{2}  \min \left\{   \langle p  \rangle^{2}   ,\langle q \rangle^{2} \right\},\label{upper and lower bound mid term Lambda c p q}\\
\frac {c^2 |p-q|^2}  {p^0 q^0} \le&    \frac {(p^0q^0 - p \cdot q)^2-c^4 } {c^2}    \le &   \frac {p^0 q^0} {c^2}|p-q|^2, 
\label{upper and lower bound mid term Lambda c p q 2} 
\end{eqnarray}
and  
\begin{eqnarray}
\left|\frac {c^2} {p^0q^0}  \frac {(p^0q^0- p \cdot q)^2} {c^4}   -1 \right | &\le&4 \frac {\langle p \rangle^2 \langle q \rangle^2} {c^2}, \label{difference first part of Lambda c}\\
\left  |\left[ \frac {(p^0q^0 - p \cdot q)^2-c^4 } {c^2}  \right]^{-\frac{3}{ 2}}  - |p-q|^{-3}  \right| &\lesssim&\frac {\langle p \rangle^2 \langle q \rangle^2} {c^2}   |p-q|^{-3} . \label{difference second part of Lambda c}
\end{eqnarray}
\end{lemma}
\begin{proof}
The first  inequality in \eqref{upper and lower bound g} for $g$  is proved in \cite{Glassey-Strauss-ARMA-1986} and \cite[Lemma 3.2]{Strain-SIAM-2010}, we only prove the second one here. By the definition of $g$ in \eqref{definition g s}, we have
\begin{align*}
\nonumber
g^2 =& 2(p^0q^0 -p \cdot q -c^2) = 2 \frac {(p^0)^2(q^0)^2 - (p \cdot q + c^2)^2}  {p^0q^0+p \cdot q +c^2}\\
 = &2 \frac {c^2|p|^2 +c^2 |q|^2 + |p|^2 |q|^2   -|p \cdot q|^2 -2 c^2 p \cdot q}  {p^0q^0+p \cdot q +c^2} 
=  2\frac {c^2|p-q|^2 + |p \times q|^2}  {p^0q^0+p \cdot q +c^2}.
\end{align*}
Then we use \eqref{p^0q^0 p q lower bound 2} to further obtain
\begin{align*}
g^2 \le& 2\frac {c^2|p-q|^2 + |p|^2 |q|^2   }  {p^0q^0 - |p || q| +c^2}\\
 \le & 4 \min \left \{ \frac {(c^2|p-q|^2 + |p|^2 |q|^2  )q^0 }  {c^2 p^0} ,  \frac {(c^2|p-q|^2 + |p|^2 |q|^2  )p^0 }  {c^2 q^0}  \right \} 
\\
\le& 24\min \left \{ \frac {\langle p \rangle^2 \langle q \rangle^2  q^0 }  { p^0} ,  \frac {\langle p \rangle^2 \langle q \rangle^2 p^0 }  {q^0}  \right \}.
\end{align*}
The estimates in \eqref{upper and lower bound s} for $s$  follow form \eqref{p^0q^0 p q lower bound} and \eqref{definition g s}. 

For \eqref{upper and lower bound mid term Lambda c p q} and \eqref{upper and lower bound mid term Lambda c p q 2}, we compute that 
\[
(p^0q^0 - p \cdot q)^2-c^4  = (p^0q^0 - p \cdot q-c^2  )(p^0q^0 - p \cdot q+ c^2 ) =g^2 s/4.
\]
Then \eqref{upper and lower bound mid term Lambda c p q} and \eqref{upper and lower bound mid term Lambda c p q 2} follow  by combing above estimates for $g, s$ and the fact that $q^0/c \le \langle q \rangle, p^0/c \le \langle p \rangle$. 
For \eqref{difference first part of Lambda c}, we use 
\begin{align*}
c^2 \le p^0q^0- p \cdot q =c^2+ g^2/2 \le c^2+ |p-q|^2,
\end{align*}
to have
\[
\frac {c^2} {p^0q^0}  \le \frac {c^2} {p^0q^0}  \frac {(p^0q^0- p \cdot q)^2} {c^4}  \le \frac  {(c^2 + |p-q|^2)^2 } {c^4} \le 1 + 4 \frac {\langle p \rangle^2 \langle q \rangle^2} {c^2},
\]
which implies
\[
-\frac {p^2+q^2}  {c^2}\le\frac {c^2} {p^0q^0} -1 \le \frac {c^2} {p^0q^0}  \frac {(p^0q^0- p \cdot q)^2} {c^4}   -1 \le   4 \frac {\langle p \rangle^2 \langle q \rangle^2} {c^2}.
\]
Finally, we come to estimate \eqref{difference second part of Lambda c}.
From \eqref{upper and lower bound mid term Lambda c p q 2},  we have
\[
\frac {c^3} {(p^0q^0)^{3/2}  }   |p-q|^{-3}   \le \left[ \frac {(p^0q^0 - p \cdot q)^2-c^4 } {c^2}  \right]^{-\frac{3}{ 2}}\le\frac{(p^0q^0)^{3/2}  } {c^3} |p-q|^{-3} 
\]
and thus obtain
\begin{align*}
\left  |\left[ \frac {(p^0q^0 - p \cdot q)^2-c^4 } {c^2}  \right]^{-\frac{3}{ 2}}  - |p-q|^{-3}  \right| \le & \max \left \{ \frac{(p^0q^0)^{3/2}  } {c^3} -1,  1- \frac{c^3}    {(p^0q^0)^{3/2} }  \right\} |p-q|^{-3} 
\\
\le & \max \left \{ \frac{(p^0q^0)^{2}  } {c^4} -1,  1- \frac{c^4}    {(p^0q^0)^{4} }   \right\} |p-q|^{-3} 
\\
\le & \left(\frac{(p^0q^0)^{2}  } {c^4} -1  \right)   |p-q|^{-3}  \lesssim \frac {\langle p \rangle^2 \langle q \rangle^2} {c^2}   |p-q|^{-3}. 
\end{align*}
\end{proof}

\begin{lemma}
For any $p, q$ we have
\begin{equation}
\label{p^0 q^0 p 2 lower bound}
\left[(p^0q^0 - p \cdot q)^2-c^4  \right]  |p|^2  \ge  c^2  |p \times q|^2,
\end{equation}
as well as
\begin{equation}
\label{equality p^0 q^0 - p q 2 - c 4 -c 2 p - q 2}
(p^0q^0 - p \cdot q)^2-c^4   -c^2|p-q|^2   = - p\cdot q g^2 +|p \times q|^2.
\end{equation}
For the $\mathcal{S}^c(p, q)$ defined in \eqref{definition S c p q} and any $\bar{p} \in p^{\perp}$, we have
\begin{eqnarray}
\label{equality for p S p}
p_i\mathcal{S}^{c, ij}(p, q)   p_j &=&\frac { (p^0)^2 } {c^2}    |p \times q|^2 =   \frac { (p^0)^2 } {c^2}  |q \times ( p-q)|^2,\\
\label{equality for p S p 1}
p_i \mathcal{S}^{c, ij} (p, q)  \bar{p}_j  &=&  - (q \cdot p)(q \cdot \bar{p} )  + \frac 1 {c^2}( p^0q^0- p \cdot q)   |p|^2( q \cdot \bar{p} ), \\
\label{equality for p 1 S p 1}
\bar{p}_i \mathcal{S}^{c, ij} (p, q)  \bar{p}_j   &=& \frac 1 {c^2} \left[(p^0q^0 - p \cdot q)^2-c^4  \right] |\bar{p} |^2    -   |q \cdot\bar{p} |^2.
\end{eqnarray}
\end{lemma}
\begin{proof}
For \eqref{p^0 q^0 p 2 lower bound}, we have
\begin{align*}
&\left[(p^0q^0 - p \cdot q)^2-c^4  \right]  |p|^2 -   c^2  (p \times q)^2  
\\
&\qquad= [c^2 (|p|^2+ |q|^2) +|p|^2|q|^2 - 2p^0q^0 p \cdot q + (p \cdot q)^2    ]|p|^2 -  c^2  (p \times q)^2  
\\
&\qquad=c^2 |p|^4 + c^2 |q|^2|p|^2 +  |p|^4|q|^2 - 2p^0q^0 p \cdot q |p|^2 + (p \cdot q)^2   | p  | ^2  -c^2 (p \times q)^2
\\
&\qquad= (c^2+ |q|^2) |p|^4 + c^2 (p \cdot q)^2 - 2p^0q^0 p \cdot q |p|^2 + (p \cdot q)^2    |p|^2 
\\
&\qquad= (q^0)^2  |p|^4 + (p^0)^2(p \cdot q)^2 - 2p^0q^0 p \cdot q |p|^2  = (q^0p^2 - p^0 p \cdot q)^2 \ge 0.
\end{align*}
For \eqref{equality p^0 q^0 - p q 2 - c 4 -c 2 p - q 2}, we obtain
\begin{align*}
&(p^0q^0 - p \cdot q)^2-c^4   -c^2|p-q|^2 \\
&\qquad= (c^2+|p|^2)(c^2+|q|^2)    - 2p^0q^0 p \cdot q +|p \cdot q|^2  - c^4 -c^2|p-q|^2 
\\
&\qquad= c^2(|p|^2+|q|^2)  + |p|^2|q|^2 - 2p^0q^0 p \cdot q +|p \cdot q|^2  -c^2|p-q|^2 
\\
&\qquad=   2c^2 p \cdot q - 2p^0q^0 p \cdot q  +  2|p \cdot q|^2  + |p \times q |^2
\\
&\qquad= 2 p \cdot q (c^2+p\cdot q -p^0q^0) + |p \times q |^2 = - p\cdot q g^2 +|p \times q|^2.
\end{align*}
Now we come to the estimates of $\mathcal{S}^c(p, q)$ in \eqref{definition S c p q}. For the estimate \eqref{equality for p S p}, from the definition of $\mathcal{S}^c(p, q)$ and $a_i (b \otimes c)^{i j} d_j =(a \cdot b) (c \cdot d)$, we can easily compute
\[
p_i  S^{c, ij}(p, q) p_j = \frac 1 {c^2} \left[(p^0q^0 - p \cdot q)^2-c^4  \right] |p|^2 - |p|^4 - (p \cdot q)^2 + \frac 1 {c^2}( p^0q^0- p \cdot q ) 2|p|^2 (p \cdot q).
\]
Then we can further obtain
\begin{align*}
p_i  S^{c, ij}(p, q) p_j =&\frac 1 {c^2} \left\{[(p^0)^2 (q^0)^2 -c^4] |p|^2  -c^2|p|^4 -c^2 (p \cdot q)^2   - |p|^2  (p \cdot q)^2\right\}
\\
=& \frac 1 {c^2}  \left\{[(p^0)^2 (q^0)^2 -c^4-c^2 |p|^2 ] |p|^2  -(p^0)^2    (p \cdot q)^2\right\}
\\
=&\frac{ |p|^2|q|^2(p^0)^2   -(p^0)^2   (p \cdot q)^2}{c^2} =  \frac {(p^0)^2  |p \times q|^2} {c^2}  =  \frac { (p^0)^2  |q \times ( p-q)|^2 } {c^2}. 
\end{align*}
For  \eqref{equality for p S p 1}, we  use $\bar{p}  \cdot p =0$ to have
\begin{align*}
p_i S^{c, ij} (p, q)  \bar{p}_j   =& \frac 1 {c^2} \left[(p^0q^0 - p \cdot q)^2-c^4  \right] (p \cdot \bar{p} ) -|p|^2 (p \cdot \bar{p} )  - (q \cdot p) (q \cdot \bar{p} ) 
\\
&+  \frac 1 {c^2}( p^0q^0- p \cdot q)   ( |p|^2( q \cdot \bar{p} ) + (q \cdot p )  (p \cdot \bar{p} ))
\\
 =& -(q \cdot p)(q \cdot \bar{p} )  +  \frac 1 {c^2}( p^0q^0- p \cdot q)   |p|^2( q \cdot \bar{p} )  
\end{align*}
For \eqref{equality for p 1 S p 1}, we compute 
\begin{align*}
\bar{p}_i  S^{c, ij} (p, q)  \bar{p}_j  & = \frac {(p^0q^0 - p \cdot q)^2-c^4} {c^2} |\bar{p} |^2  - |p \cdot \bar{p} |^2 -  |q \cdot \bar{p} |^2 +  \frac {( p^0q^0- p \cdot q)   2(p \cdot \bar{p}  ) (q \cdot \bar{p} )} {c^2} 
\\
&  = \frac  {(p^0q^0 - p \cdot q)^2-c^4} {c^2}  |\bar{p} |^2     -   |q \cdot \bar{p} |^2.
\end{align*}
\end{proof}

With the above preparations, now we estimate the collision kernel $\Phi^c(p,q)$ \eqref{cker} and the difference  of the collision kernels $\Phi^c(p,q)-\Phi^{\infty}(p,q)$.

\begin{lemma}
 For the collision kernel $\Phi^c(p,q)$ given in \eqref{cker}, it holds that
\begin{align}
\left|\Phi^{c,ij}(p,q)\right|\leq&C\min\Big\{\frac{p^0}{c} \langle q\rangle^4, \frac{q^0}{c} \langle p\rangle^4\Big\}|p-q|^{-1},
\label{cker-der01}\\
\Big|\sum_{i}\partial_{p_i}\Phi^{c,ij}(p,q)\Big|\leq & \frac{Cp^0}{c}\langle q\rangle^7|p-q|^{-2},\label{cker-der11}\\
\Big|\partial_{p_i}\partial_{q_j}\Phi^{c,ij}(p,q)\Big|\leq & \frac{C\langle q\rangle}{c^2}|p-q|^{-1},\label{cker-der21}\\
\Phi^{c,ij}(p,q)p_ip_j\leq&\frac{C\big(p^{0}\big)^3}{c^3} \langle q\rangle^6|p-q|^{-1},
\label{cker-der000}\\
\Big|\Phi^{c,ij}(p,q)-\Phi^{\infty,ij}(p,q)\Big|\leq & \frac{C}{c^2}\langle p\rangle^5\langle q\rangle^5|p-q|^{-1}.\label{ker-cinf}
\end{align}
\end{lemma}
\begin{proof} Before the proof of \eqref{cker-der01}, we first claim that
 \begin{equation}
\label{upper and lower bound Lambda c p q}
\frac 1 {125\langle p \rangle^{3} \langle q \rangle^{3}}  \max \left \{   \frac {p^0}{c\langle q \rangle^{6}}      , \frac { q^0}{c \langle p \rangle^{6} }    \right\}     \le   \frac{c^2\Lambda^c(p, q) }{p^0q^0}\le  32    |p-q|^{-3}   \min \left \{\frac {q^0} c \langle p \rangle^4, \frac {p^0} c \langle q \rangle^4  \right\}
\end{equation}
for $\Lambda^c(p, q)$ in \eqref{definition Phi c p q}. In fact, we use \eqref{p^0q^0 p q lower bound} to have 
\[
\frac 1 4 c^4  \max \left \{   \frac {(p^0)^2}{(q^0)^2}, \frac {(q^0)^2}{(p^0)^2}    \right\}  \le (p^0q^0 -|p| |q|)^2\le  (p^0q^0 -p \cdot q)^2 \le 4(p^0)^2(q^0)^2,
\]
which implies that 
\begin{align}
\label{upper and lower bound first term Lambda p q}
\max \left \{   \frac {p^0}{4c\langle q \rangle^{3}}      , \frac { q^0}{4c\langle p \rangle^{3}}      \right\}    \le& \max \left \{   \frac {c^2  p^0}{4(q^0)^3}, \frac {c^2 q^0}{4(p^0)^3}    \right\}  \nonumber\\
 \le&  \frac {(p^0q^0- p \cdot q)^2} {c^2p^0q^0}\le 4 \min \left \{\frac {q^0} c \langle p \rangle, \frac {p^0} c \langle q \rangle  \right\}.
\end{align}

This together with  \eqref{upper and lower bound mid term Lambda c p q} and the fact that $q^0/c \le \langle q \rangle, p^0/c \le \langle p \rangle$ gives 
\eqref{upper and lower bound Lambda c p q}.

Now we come to prove \eqref{cker-der01}. From \eqref{p^0q^0 p q lower bound 2} and \eqref{upper and lower bound mid term Lambda c p q}, we obtain
\begin{align*}
\Phi^{c, ij}(p, q) =& \Lambda^c(p, q) S^{c, ij}(p, q)= \frac {p^0q^0} {c^2} \left[ \frac {\big(q^{\mu}p_{\mu}\big)^2-c^4 } {c^2}  \right]^{-\frac{3}{ 2}}
\\
     &\times   \left[\frac 1 {c^2} \left[\big(q^{\mu}p_{\mu}\big)^2-c^4  \right] \delta_{ij}  - (p -q)_i (p- q)_j - \frac {q^{\mu}p_{\mu} +c^2} {c^2}   (p_i q_j + q_i p_j)  \right]  
\\
\lesssim &\frac {p^0q^0} {c^2}     \left[ \frac {\big(q^{\mu}p_{\mu}\big)^2-c^4 } {c^2}  \right]^{-\frac 1 2} + \frac {p^0q^0} {c^2}   \left[ \frac {\big(q^{\mu}p_{\mu}\big)^2-c^4 } {c^2}  \right]^{-\frac{3}{ 2}} |p-q|^{2} 
\\
&+ \frac {p^0q^0} {c^2}  \left[ \frac {\big(q^{\mu}p_{\mu}\big)^2-c^4 } {c^2}  \right]^{-\frac 1 2} \frac 1 { p^0q^0- p \cdot q + c^2}   |p| |q|
\\
\lesssim &\frac {p^0q^0} {c^2}   \left[  |p-q|^{-1}  \langle q \rangle+ \langle q \rangle^3  |p-q|^{-3} |p-q|^{2} +   |p-q|^{-1} \langle q \rangle   \frac{\langle q \rangle} { \langle p \rangle} |p|  |q|  \right] 
\\
\lesssim  &  \frac {p^0q^0} {c^2} \langle q \rangle^3 |p-q|^{-1}  \lesssim \frac {p^0} {c} \langle q \rangle^4 |p-q|^{-1}  
\end{align*}
and thus \eqref{cker-der01} is  proved by interchanging $p$ and $q$.

By the same computation as in \cite[Lemma 3]{Strain-Guo-CMP-2004} for $c=1$, we have the following expressions of derivatives of $\Phi^{c}(p,q)$ for $c\geq1$:
\begin{align}
 \sum_{i}\partial_{p_i}\Phi^{c,ij}(p,q)=&\frac{-2}{p^0q^0} \Lambda^{c,ij}\left[c^2p_j+\big(q^{\mu}p_{\mu}\big)q_j\right],
\label{cker-der0}\\
\partial_{p_i}\partial_{q_j}\Phi^{c,ij}(p,q)=&\frac{-4q^{\mu}p_{\mu}}{c^2p^0q^0} \left(\frac{1}{c^2}\left(q^{\mu}p_{\mu}\right)^2-c^2\right)^{-\frac{1}{2}}.
\label{cker-der10}
\end{align}
Note that from \eqref{p^0q^0 p q lower bound 2},
\begin{align*}
 \left|c^2p_j-\big(q^{\mu}p_{\mu}\big)q_j\right|= &\left| c^2(q_j-p_j)+g^2q_j\right|\leq c^2|p-q|+2|q|\frac {c^2|p-q|^2 + |(p-q) \times q|^2}  {p^0q^0+p \cdot q +c^2}\\
  \leq & c^2|p-q|+4\frac{\langle q\rangle|p-q|^2(q^0)^2}{\langle p\rangle}\leq 9c^2|p-q|\langle q\rangle^4.
\end{align*}
This together with \eqref{upper and lower bound Lambda c p q} and \eqref{cker-der0} gives \eqref{cker-der11}. For \eqref{cker-der21}, employing the fact $\frac{-4q^{\mu}p_{\mu}}{c^2p^0q^0}\leq \frac{4}{c^2}$, we derive it from \eqref{cker-der10} and \eqref{upper and lower bound mid term Lambda c p q}.

By \eqref{equality for p S p} and \eqref{upper and lower bound Lambda c p q}, we obtain \eqref{cker-der000} as follows:
\begin{align*}
  p_i\Phi^{c, ij}(p, q)  p_j =& \frac {(p^0)^2} {c^2}    |q \times ( p-q)|^2  \frac{c^2\Lambda^c(p, q)}{p^0q^0}\\
   \lesssim&   \frac {(p^0)^2} {c^2}    |q|^2 |p-q|^2 \frac{p^0}{c} \langle q \rangle^{4} |p-q|^{-3} \lesssim 
   \frac {(p^0)^3} {c^3}  \langle q \rangle^6 |p-q|^{-1}.
\end{align*}

Finally, we come to prove \eqref{ker-cinf}. To this end, we first estimate the difference between $\mathcal{S}^{c, ij}(p, q)$ and $\mathcal{S}^{\infty, i j} (p, q)$. We use \eqref{equality p^0 q^0 - p q 2 - c 4 -c 2 p - q 2} to have
\begin{align*}
&|\mathcal{S}^{c, ij}(p, q) -\mathcal{S}^{\infty, ij} (p, q)|\\
 &\qquad =  \Big|\frac {\delta_{ij}} {c^2} \left[(p^0q^0 - p \cdot q)^2-c^4  \right]  -p_i p_j -q_iq_j 
\\
&\qquad+  \frac {1} {c^2}( p^0q^0- p \cdot q)   (p_i q_j + q_i p_j) -        |p-q|^2 \delta_{ij} +    (p-q )_i (p-q)_j \Big|
\\
&\qquad= \left|\frac {\delta_{ij}} {c^2} \left[(p^0q^0 - p \cdot q)^2-c^4   -c^2|p-q|^2   \right] + \frac 1 {c^2}( p^0q^0- p \cdot q -c^2)   (p_iq_j + q_ip_j)  \right|
\\
&\qquad\lesssim \frac 1 {c^2}[ | p | |q|   g^2 + |p \times (p-q)|^2] + \frac 1 {c^2}  |p| |q| g^2
\lesssim  \frac 1 {c^2}  |p-q|^{2}   \langle p \rangle^2 \langle q \rangle^2.
\end{align*}
For the difference between $\Lambda^c(p, q)$ and $\Lambda^\infty (p, q)$ by \eqref{difference first part of Lambda c} and \eqref{difference second part of Lambda c} we first have that 
\begin{align*}
&|\Lambda^c (p, q) -\Lambda^\infty (p, q) |\\
&\qquad=\Big|  \left [\frac {c^2} {p^0q^0}  \frac {(p^0q^0- p \cdot q)^2} {c^4}   -1  \right]  |p-q|^{-3} 
\\
&\qquad+ \frac {c^2} {p^0q^0}  \frac {(p^0q^0- p \cdot q)^2} {c^4}   \left[ \left[ \frac {(p^0q^0 - p \cdot q)^2-c^4 } {c^2}  \right]^{-\frac{3}{ 2}}  - |p-q|^{-3}    \right] \Big|
\\
&\qquad\lesssim \frac {\langle p \rangle^4 \langle q \rangle^4} {c^2}   |p-q|^{-3} .
\end{align*}
Also for the term $S^c(p, q)$ defined in \eqref{definition S c p q}, we have
\begin{align*}
|S^{c, ij}(p, q)| \le&  \frac 1 {c^2} [(p^0q^0 - p \cdot q)^2-c^4 ] +|p-q|^2+  \frac {2g^2} {c^2}|p||q| \\
\lesssim& \frac {g^2 s} {c^2} +|p-q|^2 + |p-q|^2 |p| |q| \lesssim |p-q|^2 \langle p \rangle \langle q \rangle.
\end{align*}
Collecting the above two estimates, we have
\begin{align*}
&|\Phi^{c, ij}(p, q) - \Phi^{\infty, ij} (p, q)| \\
&\qquad\lesssim | \Lambda^c (p, q) -\Lambda^\infty (p, q)||S^{c, ij}(p, q)| +\Lambda^\infty (p, q)  |S^{c, ij}(p, q) - S^{\infty, ij}(p, q)|\\
&\qquad\lesssim\frac {\langle p \rangle^5 \langle q \rangle^5} {c^2}   |p-q|^{-1}. 
\end{align*}
 \end{proof}

Next we recall eigenvalues of $\sigma^c(p)$ and their properties.
\begin{lemma} \label{spec-sigm} The spectrum of $\sigma^c$ consists of the following a simple eigenvalue $\lambda_1(p)$ and a double eigenvalue $\lambda_2(p)$:
\begin{align}\label{upper lower bound lambda 1 lambda 2}
\lambda_1(p)=\sigma^{c,ij}(p)\frac{p_ip_j}{|p|^2}\backsimeq\frac{(p^0)^3}{c^3\langle p\rangle^3},\qquad \lambda_2(p)=\sigma^{c,ij}(p)\frac{\bar{p}_i\bar{p}_j}{|\bar{p}|^2}\backsimeq\frac{p^0}{c\langle p\rangle},
\end{align}
where $\bar{p}\in p^{\bot}$, the eigenspace  perpendicular to $p$. Moreover, for any vector $\xi\in\mathbb{R}^3$, it holds that
\begin{align}\label{sigma i j p equality}
\xi_i \sigma^{c,ij}(p)\xi_j=\lambda_1(p)\frac{(p\cdot\xi)^2}{|p|^2}+\lambda_2(p)\frac{(p\times\xi)^2}{|p|^2}.
\end{align}
Here and in the sequel, Einstein's summation convention is used.
\end{lemma}
\begin{proof} The proof of this lemma for $c=1$ was given in \cite[Proposition 4.3, 4.4]{ Lemou-MMAS-2000}. Here we derive \eqref{upper lower bound lambda 1 lambda 2} and \eqref{sigma i j p equality} for general $c\geq1$.

We first prove \eqref{sigma i j p equality}. To this end, we split $\xi$ into two parts
\[
\xi = \xi_1+ \xi_2, \quad \xi_1= \frac {p(p\cdot  \xi)} {|p|^2}   ,\quad \xi_1 \parallel p, \quad \xi_2 \perp p,\quad \xi_1 \perp \xi_2
\]
and thus we have
\[
\xi_i \sigma^{c, ij} (p)  \xi_j  = \xi_{1, i} \sigma^{c, ij} (p)\xi_{1, j}+  2 \xi_{1, i}    \sigma^{c, ij}  (p) \xi_{2, j} + \xi_{2, i}\sigma^{c, ij}  (p)\xi_{2, j}.
\]
From the definition of $\lambda_{1}(p)$ and $\lambda_{2}(p)$, we have
\[
 \xi_{1, i} \sigma^{c, ij}  (p)  \xi_{1, j} = \lambda_{1}(p )|\xi_{1}|^2 = \lambda_{1}(p)  \left  |\frac {p} {|p|} \cdot  \xi   \right|^2 , \quad  \xi_{2, i} \sigma^{c, ij}  \xi_{2, j} =\lambda_{2}(p) |\xi_{2}|^2 =  \lambda_{2}(p) \left   |\frac {p} {|p|} \times  \xi   \right|^2 .
\]
By symmetry, to prove \eqref{sigma i j p equality}, we only need to prove that for any $p, \bar{p}\in \R$ satisfying  $p \cdot \bar{p}=0$,
\[
p \sigma_{c}(p)   \bar{p}  = 0.
\]
For arbitrarily given $p$ and $\bar{p}$, we can always change coordinates such that $p=(0, 0, |p|),  \bar{p}=(0, | \bar{p} |, 0)$. Under such coordinate, using the spherical coordinate system for $q$, we have 
\begin{equation}
\label{coordinate relation p q p 1}
p \cdot q  = |p| |q| \cos \theta,\quad \bar{p} \cdot q =  | \bar{p} | |q| \sin \theta \sin \phi ,\quad |p \times q|^2 = |p|^2 |q|^2 \sin^2 \theta .
\end{equation}
 Noting that $\Lambda^c(p, q) $  is of the form $|p|, |q|$, and $p \cdot q$,  we use  \eqref{equality for p S p 1} to have
\begin{align*}
p_i \sigma^{c, ij}  \bar{p}_j = &\int_{\R^3} p_i \mathcal{S}^{c, ij} (p, q) \bar{p}_j   \Lambda^c (p, q) \mu^c   (q)\, \mathrm{d}q 
\\
=  &  \int_{\R^3}    \left[- (q \cdot p)(q \cdot \bar{p})  +  \frac 1 {c^2}( p_0q^0- p \cdot q)    |p|^2( q \cdot \bar{p})   \right]    \Lambda^c (p, q)  \mu^c   (q)\, \mathrm{d}q
\\
=  & \int_{0}^\infty\int_0^\pi \int_0^{2 \pi} \left [- |p| |q| \cos \theta | \bar{p} | |q| \sin \theta \sin \phi  +  \frac 1 {c^2}( p^0q^0- |p|  |q| \cos \theta)    |p|^2   | \bar{p}  | |q|  \sin \theta \sin \phi     \right]  
\\
&\times \Lambda^c (|p|, |q|, |p| |q|\cos \theta)  \mu^c(|q|)  \sin \theta\, \mathrm{d }\phi \mathrm{d} \theta   \mathrm{ d}|q|. 
\end{align*}
Integrate with respect to $\phi$ gives $0$ for the above integral and thus \eqref{sigma i j p equality} follows. 

Now we come to prove \eqref{upper lower bound lambda 1 lambda 2}. To this end, we need to compute the lower and upper bound for $\lambda_{1}(p)$ and $\lambda_{2} (p)$. We only prove the case $|p| \ge 1$ as $|p| \le 1$ can be proved easily. 

For the lower bound of $\lambda_{1} (p)$, choosing coordinates as \eqref{coordinate relation p q p 1}, we use \eqref{upplow-mu}, \eqref{equality for p S p}, \eqref{upper and lower bound Lambda c p q}, and the fact $c/q^0 \ge \langle q \rangle^{-1}$ to have
\begin{align*}
\lambda_{1}(p) =& \frac 1  {|p^2|}\int_{\R^3}    p_i S^{c, ij}(p, q)    p_j  \Lambda^c(p, q)  \mu^c   (q)dq
\\
\ge & A \int_{\R^3}   \frac { (p^0)^2 } {c^2}    \frac{ |p \times q  |^2}  {|p|^2}      
\frac {p^0}{c}  \langle p \rangle^{-3} \langle q \rangle^{-9}      \mu^c   (q)\,
  \mathrm{d}q\\
\ge&   A   \frac {(p^0)^3} {c^3}  \langle  p  \rangle^{-3}      \int_{\R^3}    \frac{ |p \times   q  |^2}  {|p|^2}  \exp\left\{-\frac{|p|^2}{2}\right\}\, \mathrm{d}q
\\
\ge &  A   \frac {(p^0)^3} {c^3}  \langle  p  \rangle^{-3}        \int_{0}^\infty\int_0^\pi \int_0^{2 \pi}  \sin^2  \theta   \exp\left\{-\frac{|p|^2}{2}\right\} \, \mathrm{d} \phi \mathrm{d} \theta    \mathrm{d}|q| \\
\ge&    A \frac {(p^0)^3} {c^3}  \langle  p  \rangle^{-3}    
\end{align*}
for some constant $A>0$ independent of $c$. While for the upper bound, we use \eqref{upplow-mu}, \eqref{equality for p S p}, \eqref{upper and lower bound Lambda c p q},  the facts $p \times q = (p-q) \times q$ and $q^0/c \le \langle  q \rangle$ to have
\begin{align*}
\lambda_1 (p)
=& \frac 1  {|p^2|}\int_{\R^3}    p_i S^{c, ij}(p, q)    p_j  \Lambda^c(p, q) \mu^c   (q)\, \mathrm{d}q
\\
\lesssim & \int_{\R^3}    \frac {(p^0)^2} {c^2}   \frac{ |p \times q  |^2}  {|p|^2}   |p-q|^{-3}   \frac {p^0} {c}   \langle q \rangle^4      e^{- \frac 1 2|q| } \, \mathrm{d}q
\\
\lesssim &       \frac {(p^0)^3} {c^3}  \langle  p  \rangle^{-2}      \int_{\R^3}      |p-q|^{-3}   \langle q \rangle^4  { |(p-q) \times q  |^2}    \exp\left\{- \frac 3 8|q|\right\} \, \mathrm{d}q
\\
\lesssim &       \frac {(p^0)^3} {c^3}  \langle  p  \rangle^{-2}      \int_{\R^3}      |p-q|^{-1}     \langle q \rangle^6   \exp\left\{- \frac 3 8|q|\right\} \, \mathrm{d}q
\lesssim    \frac {(p^0)^3} {c^3}  \langle  p  \rangle^{-3}.
\end{align*}

 Next we come to treat the case of $\lambda_{2} (p, \bar{p})$. We first prove that $\lambda_{2} (p, \bar{p})=\lambda_{2} (p)$, i.e. $\lambda_{2}$ is independent of the choice of $\bar{p} \perp p$. We claim that for any $\bar{p} \in  p^{\perp}$, 
\begin{equation}
\label{claim p 1 p integral 1 2}
\int_{\R^3}   \frac1 {|\bar{p}|^2}  (\bar{p} \cdot q)^2  \Lambda^c(p, q) \mu^c(q) \,\mathrm{d}q  =  \frac 1 2  \int_{\R^3}   \frac1 {|p|^2}  (p \times q)^2  \Lambda^c (p, q)   \mu^c   (q) \,\mathrm{d}q  ,\quad \forall p_1 \perp p.
\end{equation}
In fact,  this is due to the fact that the tangent space of $p$ is of dimension 2. Similarly as before, choosing coordinates as \eqref{coordinate relation p q p 1},  we have
\begin{align*}
&\int_{\R^3}   \frac   1 {|\bar{p}|^2}  (\bar{p} \cdot q)^2  \Lambda^c(p, q) \mu^c(q)\,\mathrm{d}q \\
&\qquad=    \int_{0}^\infty\int_0^\pi \int_0^{2 \pi}   |q|^2 \sin^2  \theta \sin^2  \phi \Lambda^c (|p| ,|q|, |p| |q|\cos \theta)  \mu^c(|q|)   \sin \theta\, \mathrm{d} \phi \mathrm{d} \theta    \mathrm{d}|q|.  
\end{align*}
Obviously the right hand side is independent of the choice of $\bar{p} \perp p$. Also we have
\begin{align*}
&\int_{\R^3}   \frac1 {|p|^2}  (p \times q)^2  \Lambda^c(p, q) \mu^c(q) \,\mathrm{d}q    \\
&\qquad=    \int_{0}^\infty\int_0^\pi \int_0^{2 \pi}   |q|^2 \sin^2  \theta \sin^2  \phi \Lambda^c (|p| ,|q|, |p| |q|\cos \theta)  \mu^c(|q|)   \sin \theta\, \mathrm{d} \phi \mathrm{d} \theta    \mathrm{d}|q|.    
\end{align*}
The claim \eqref{claim p 1 p integral 1 2}  is verified by the above two equalities together with the following fact
\[
\int_0^{2\pi} \sin^2 \phi \,\mathrm{d}\phi =\frac 1 2\int_0^{2\pi} \,\mathrm{d} \phi.
\]
Thus we further use \eqref{equality for p 1 S p 1} to have 
\[
\lambda_{2}(p, \bar{p})  = \int_{\R^3}  \left ( \frac 1 {c^2} \left[(p^0q^0 - p \cdot q)^2-c^4  \right]  -   \frac1 {2 |p|^2}  (p \times q)^2 \right   )  \Lambda^c (p, q) \mu^c   (q)\, \mathrm{d}q,\quad \forall \bar{p} \perp p,
 \]
which implies that  $\lambda_{2}(p, \bar{p})=\lambda_{2} (p)$. 

Now we come to compute the upper and lower bound for $\lambda_{2}(p)$. For the lower bound, we use \eqref{upplow-mu}, \eqref{upper and lower bound mid term Lambda c p q}, \eqref{p^0 q^0 p 2 lower bound}, and \eqref{upper and lower bound Lambda c p q} to have
\begin{align*}
\lambda_{2} (p)  
 \ge  & \int_{\R^3}    \frac 1 {2c^2} \left[(p^0q^0 - p \cdot q)^2-c^4  \right]    \Lambda^c (p, q)\mu^c   (q) \, \mathrm{d}q
 \\
\ge & \bar{A}\int_{\R^3}     |p-q|^2   \langle q \rangle^{-2}   \frac {p^0}{c}  \langle p \rangle^{-3} \langle q \rangle^{-9}     \mu^c   (q)\, \mathrm{d}q
\\
\ge &  \bar{A}    \frac {p^0} {c}  \langle  p  \rangle^{-3}      \int_{\R^3}      |p-q|^2   \exp\left\{-\frac{|q|^2}{2} \right\} \, \mathrm{d}q\ge A \frac {p^0} {c}  \langle  p  \rangle^{-1}    
\end{align*}
for some constant $\bar{A} >0$. While for the upper bound,  we use \eqref{upplow-mu}, \eqref{upper and lower bound mid term Lambda c p q},  and \eqref{upper and lower bound first term Lambda p q} to have
\begin{align*}
\lambda_{2}(p) 
\le  & \int_{\R^3}    \frac 1 {c^2} \left[(p^0q^0 - p \cdot q)^2-c^4  \right]    \Lambda^c (p, q) \mu^c   (q)\, \mathrm{d}q
\\
\le &    \int_{\R^3}    \frac {c^2} {p^0q^0}  \frac {(p^0q^0- p \cdot q)^2} {c^4} \left[(p^0q^0 - p \cdot q)^2-c^2  \right]^{-\frac 1 2} \mu^c   (q)\, \mathrm{d}q
\\
\lesssim & \int_{\R^3}    \frac {p^0} {c}  \langle q \rangle^2   |p-q|^{-1}        \exp\left\{-\frac{3|q|}{8} \right\} \, \mathrm{d}q
\lesssim \frac {p^0} {c}  \langle  p  \rangle^{-1}.    
\end{align*}
\end{proof}

From Lemma \ref{spec-sigm}, one can easily deduce that
\begin{corollary}\label{corollary-1}
For the $|\cdot|_{\sigma}$ norm defined in \eqref{Sigma-norm}, we have
 \begin{align}\label{eqiv-norm}
    |h(t,x)|_{\sigma}^2\backsimeq& \Big|\Big(\frac{p^0}{c\langle p\rangle}\Big)^{\frac{1}{2}} h(t,x)\Big|^2_{L^2}+\Big|\Big(\frac{p^0}{c\langle p\rangle}\Big)^{\frac{3}{2}}\frac{p}{|p|}\cdot\nabla_{p} h(t,x)\Big|^2_{L^2}\nonumber\\
    &+ \Big|\Big(\frac{p^0}{c\langle p\rangle}\Big)^{\frac{1}{2}}\frac{p}{|p|}\times\nabla_{p} h(t,x)\Big|^2_{L^2}.
 \end{align}
\end{corollary}

\subsection{Estimates of the  collision operators}
In this subsection, we first establish estimates of the relativistic Landau collision operators $\mathcal{L}h$ and $\varGamma(h,\tilde{h})$,  then provide the difference estimates of  collision operators $\mathcal{L}h-\mathbf{L}h$ and $\varGamma(h,\tilde{h})-\Gamma(h,\tilde{h})$.

\subsubsection{Estimates of  $\mathcal{L}h$ }
We will prove the uniform coercivity of  $\mathcal{L}h$ with respect to the light speed $c\geq1$ and its weighted estimate successively.
\begin{lemma}\label{coerci} When $c\geq1$ is sufficiently large, there is a uniform constant $\delta_0>0$ independent of $c$ such that
\begin{align}\label{coerci0}
(\mathcal{L}h, h)\geq \delta_0 |\{I-\mathcal{P}\}h|_{\sigma}^2
\end{align}
\end{lemma}
\begin{proof} For any $c \ge 1$ fixed, we can easily prove \eqref{coerci0} with some coefficient $\delta_0(c) > 0$ by using arguments as in \cite[Lemma 8]{Strain-Guo-CMP-2004}. However we need a uniform coefficient $\delta_0$ independent of $c$ and have to prove the result for all $c$ such that $ c$ properly large. 

 For convenience, assume that $h=\{I-\mathcal{P}\}h$. Denote $\chi(\cdot)$ as a smooth function such that $0 \le \chi(\cdot)\le 1$, $\chi(x) =1$ if $0 \le x \le 1$ and $\chi(x) =1$ if $x \ge 2$. For any constant $R>1$, we also denote that $\chi_R(x) :=\chi(|x|/R)$.
 
 Recall that $\mathcal{L}_{\pm}h=\mathcal{A}_{\pm}h+\mathcal{K}_{\pm}h$ given in \eqref{linearAK-c}. As in \cite[Lemma 6]{Strain-Guo-CMP-2004}, one has
\begin{align}
\mathcal{A}_{\pm}h=&-2\partial_{p_i}\left(\sigma^{c,ij}\partial_{p_j}h_{\pm}\right)
+\frac12\sigma^{c,ij}\frac{cp_i}{p^0}\frac{cp_j}{p^0}h_{\pm}
-\partial_{p_i}\Big(\sigma^{c,ij}\frac{cp_j}{p^0}\Big)h_{\pm},\label{express-A}\\
\mathcal{K}_{\pm}h=&\left(\mu^c(p)\right)^{-\frac12}\partial_{p_i}\Big[\mu^c(p)\int_{\mathbb{R}^3}
\Phi^{c,ij}\sqrt{\mu^c(q)}\partial_{q_j}\left(h_{+}+h_-\right)\,\mathrm{d} q\Big]\nonumber\\
&+\left(\mu^c(p)\right)^{-\frac12}\partial_{p_i}\Big[\mu^c(p)\int_{\mathbb{R}^3}
\Phi^{c,ij}\sqrt{\mu^c(q)}\frac{cq_j}{2q^0}\left(h_{+}+h_-\right)\,\mathrm{d} q\Big].\label{express-K}
\end{align}
Then we have
\begin{align}\label{coerci-hh}
(\mathcal{L}h, h)= & |h|_{\sigma}^2 -\Big(\partial_{p_i}\Big(\sigma^{c,ij}\frac{cp_j}{p^0}\Big)h, h\Big)+\sum_{\pm} \left(\mathcal{K}_{\pm}h, h\right).
\end{align}
Now we claim that for any $m>1$, there is a positive constant $C(m)$ independent of $c$ such that
\begin{align}\label{coerci-hK}
\begin{aligned}
&\Big|\Big(\partial_{p_i}\Big(\sigma^{c,ij}\frac{cp_j}{p^0}\Big)h, h\Big)\Big|+\sum_{\pm} |\left(\mathcal{K}_{\pm}h, h\right)|\leq \left(\frac{1}{2}+\frac{C}{m}\right)|h|_{\sigma}^2+C(m)\int_{|p|\leq 2m} \langle p\rangle^{-3}|h|^2\,\mathrm{d} p.
\end{aligned}
\end{align}
Once the claim holds true, we will further derive a lower bound estimate of $|h|_{L^2_{-3/2}}$ and combine it with \eqref{coerci-hK} to obtain the uniform coercivity estimate \eqref{coerci0}. We first verify \eqref{coerci-hK}.

\noindent\underline{{\it Step 1. Verification of \eqref{coerci-hK}:}}   
 Denote the first term and the second term in the left hand side of \eqref{coerci-hK} as $\mathcal{I}_1$ and $\mathcal{I}_2$, respectively. Now we use the definition of $|\cdot|_{\sigma}$ in \eqref{Sigma-norm} to have
\begin{align*}
\mathcal{I}_1=& \Big|\Big(\sigma^{c,ij}\frac{cp_j}{p^0}h, \partial_{p_i}h\Big)\Big|\\
\leq&\Big(\int_{\mathbb{R}^3}
\sigma^{c,ij}\frac{cp_i}{p^0}\frac{cp_j}{p^0}|h(p)|^2\,\mathrm{d} p\Big)^{\frac12}\Big(\int_{\mathbb{R}^3}
\sigma^{c,ij}\partial_{p_i}h(p)\partial_{p_j}h(p)\,\mathrm{d} p\Big)^{\frac12}\\
\leq& \frac{1}{2}\Big(2\int_{\mathbb{R}^3}
\sigma^{c,ij}\partial_{p_i}h(p)\partial_{p_j}h(p)\,\mathrm{d} p+\frac{1}{2}\int_{\mathbb{R}^3}
\sigma^{c,ij}\frac{cp_i}{p^0}\frac{cp_j}{p^0}|h(p)|^2\,\mathrm{d} p\Big)=\frac12|h|^2_{\sigma}.
\end{align*}
For $\mathcal{I}_2$, we have
\begin{align*}
\mathcal{I}_2=&\Big|\int_{\mathbb{R}^6}\sqrt{\mu^c(p)}\sqrt{\mu^c(q)}
\Phi^{c,ij}\sum_{\pm}\Big[\frac{cp_i}{2p^0}h_{\pm}(p)+\partial_{p_i}h_{\pm}(p)\Big]
\sum_{\pm}\Big[\frac{cq_j}{2q^0}h_{\pm}(q)+\partial_{q_j}h_{\pm}(q)\Big]\,\mathrm{d} p\mathrm{d} q\Big|
\end{align*}
We decompose $\mathcal{I}_2$ into two parts
\begin{align*}
  \mathcal{I}_2=& \Big|\int_{\mathbb{R}^6}\chi_m(p)\chi_m(q)\Big| +\Big| \int_{\mathbb{R}^6}\left(1-\chi_m(p)\chi_m(q)\right)\Big|  :=\mathcal{I}_2^1+\mathcal{I}_2^2
\end{align*}
We integrate by parts with respect to $p$ and $q$ to delete  derivatives of $h_{\pm}$ in $\mathcal{I}_2^1$, and
use \eqref{cker-der01}, \eqref{cker-der11}, \eqref{cker-der21} to have
\begin{align*}
|\mathcal{I}_2^1|\leq& C(m) \int_{\mathbb{R}^6}\chi_m(p)\chi_m(q)|h(p)||h(q)|\sum_{i,j}\Big(|\Phi^{c,ij}|+|\partial_{p_i}\Phi^{c,ij}|+|\partial_{q_j}\Phi^{c,ij}|
+|\partial_{p_i}\partial_{q_j}\Phi^{c,ij}|\Big)\,\mathrm{d} p\mathrm{d} q\\
\leq& C(m) \int_{\mathbb{R}^6}\chi_m(p)\chi_m(q)|h(p)||h(q)|\Big(|p-q|^{-1}+|p-q|^{-2}
\Big)\,\mathrm{d} p\mathrm{d} q\\
\leq& C(m) \Big(\int_{\mathbb{R}^6}\chi_m(p)\chi_m(q)|h(p)|^2\Big(|p-q|^{-1}+|p-q|^{-2}
\Big)\,\mathrm{d} p\mathrm{d} q\Big)^{\frac12}\\
&\times\Big(\int_{\mathbb{R}^6}\chi_m(p)\chi_m(q)|h(q)|^2\Big(|p-q|^{-1}+|p-q|^{-2}
\Big)\,\mathrm{d} p\mathrm{d} q\Big)^{\frac12}\\
\leq& C(m)\int_{\mathbb{R}^3}\chi_m(p)|h(p)|^2\,\mathrm{d} p=C(m)\int_{|p|\leq 2m}\langle p\rangle^{-3}|h|^2\,\mathrm{d} p.
\end{align*}
Noting that
\begin{align*}
\left(1-\chi_m(p)\chi_m(q)\right)\sqrt{\mu^c(p)}\sqrt{\mu^c(q)}\leq \frac{C}{m}\big[\mu^c(p)\mu^c(q)\big]^{\frac{1}{4}},
\end{align*}
we can bound $\mathcal{I}_2^2$ by

\begin{align*}
\mathcal{I}_2^2\leq&\frac{C}{m}\Big(\int_{\mathbb{R}^6}
\big[\mu^c(p)\mu^c(q)\big]^{\frac{1}{4}}
\Phi^{c,ij}\sum_{\pm}\Big[\frac{cp_i}{2p^0}h_{\pm}(p)+\partial_{p_i}h_{\pm}(p)\Big]
\sum_{\pm}\Big[\frac{cp_j}{2p^0}h_{\pm}(p)+\partial_{p_j}h_{\pm}(p)\Big]\,\mathrm{d} p\mathrm{d} q\Big)^{\frac12}\\
& \times\Big(\int_{\mathbb{R}^6}
\big[\mu^c(p)\mu^c(q)\big]^{\frac{1}{4}}
\Phi^{c,ij}\sum_{\pm}\Big[\frac{cq_i}{2q^0}h_{\pm}(q)+\partial_{q_i}h_{\pm}(q)\Big]
\sum_{\pm}\Big[\frac{cq_j}{q2^0}h_{\pm}(q)+\partial_{q_j}h_{\pm}(q)\Big]\,\mathrm{d} p\mathrm{d} q\Big)^{\frac12}\\
\leq&\frac{C}{m}\Big(\int_{\mathbb{R}^6}
\big[\mu^c(p)\mu^c(q)\big]^{\frac{1}{4}}
\Phi^{c,ij}\sum_{\pm}\Big[\partial_{p_i}h_{\pm}(p)\partial_{p_j}h_{\pm}(p)
+\frac{1}{4}\frac{cp_i}{p^0}\frac{cp_j}{p^0}|h_{\pm}|^2(p)\Big]\,\mathrm{d} p\mathrm{d} q\Big)^{\frac12}\\
&\times\Big(\int_{\mathbb{R}^6}
\big[\mu^c(p)\mu^c(q)\big]^{\frac{1}{4}}
\Phi^{c,ij}\sum_{\pm}\Big[\partial_{q_i}h_{\pm}(q)\partial_{q_j}h_{\pm}(q)
+\frac{1}{4}\frac{cq_i}{q^0}\frac{cq_j}{2q^0}|h_{\pm}|^2(q)\Big]\,\mathrm{d} p\mathrm{d} q\Big)^{\frac12}.
\end{align*}
Note that
\begin{align}\label{eqiv-sigma}
\int_{\mathbb{R}^3} \Phi^{c,ij}(p,q)
\big[\mu^c(p)\big]^{\frac{1}{4}}\,\mathrm{d} q\simeq \sigma^{c,ij}(p),\qquad\int_{\mathbb{R}^3} \Phi^{c,ij}(p,q)
\big[\mu^c(q)\big]^{\frac{1}{4}}\,\mathrm{d} p\simeq \sigma^{c,ij}(q).
\end{align}
We can use the definition of $|\cdot|_{\sigma}$ to further estimate $\mathcal{I}_2^2$ as
\begin{align*}
\mathcal{I}_2^2\leq& \frac{C}{m}
\Big(\int_{\mathbb{R}^3}
\sigma^{c,ij}(p)\sum_{\pm}\Big[\partial_{p_i}h_{\pm}(p)\partial_{p_j}h_{\pm}(p)
+\frac{1}{4}\frac{cp_i}{p^0}\frac{cp_j}{p^0}|h_{\pm}|^2(p)\Big]\,\mathrm{d} p\Big)^{\frac12}\\
&\times\Big(\int_{\mathbb{R}^3}
\sigma^{c,ij}(q)\sum_{\pm}\Big[\partial_{q_i}h_{\pm}(q)\partial_{q_j}h_{\pm}(q)
+\frac{1}{4}\frac{cq_i}{q^0}\frac{cq_j}{q^0}|h_{\pm}|^2(q)\Big]\,\mathrm{d} q\Big)^{\frac12}
\leq\frac{C}{m}|h|_{\sigma}^2.
\end{align*}
Combing above estimates of $\mathcal{I}_2^1$ and $\mathcal{I}_2^2$ gives
\begin{align*}
\left|\mathcal{I}_2\right|\leq \frac{C}{m}|h|_{\sigma}^2+C(m)\int_{|p|\leq m}|h|^2\,\mathrm{d} p,
\end{align*}
which together with the estimate of  $\mathcal{I}_1$ yields \eqref{coerci-hK}.

 \noindent\underline{{\it Step 2. Estimation of $|h|_{L^2_{-3/2}}$:}}  As in \cite[Lemma 1]{Strain-Guo-CMP-2004}, we can obtain that 
 \begin{align}
\label{representation for L c}
\nonumber
(\mathcal{L}  h , h )= &  \int_{\R^3} \int_{\R^3}    \mu^c (p)  \mu^c  (q)   [ \partial_{p_i} ( [\mu^c(p)]^{-\frac 1 2 } h   (p) )  - \partial_{q_i}([\mu^c(q)]^{-\frac 1 2 }  h(q) )    ] \\
&\times \Phi^{c, ij}(p, q)     \partial_{p_j} ( [\mu^c(p)]^{-\frac 1 2 }   h(p) )  \,  \mathrm{d}p \mathrm{d}q
\\ \nonumber
= &\frac 1  2    \int_{\R^3} \int_{\R^3}   \mu^c  (p)  \mu^c  (q)   [  \partial_{p_i}  ([\mu^c(p)]^{-\frac 1 2 }h   (p) )  -\partial_{q_i}  ([\mu^c(q)]^{-\frac 1 2 }   h(q) )    ]  \Phi^{c, ij} (p, q)   
\\ 
&\times [ \partial_{p_j}  ( [\mu^c(p)]^{-\frac 1 2 }h(p) )  - \partial_{q_j}  ( [\mu^c(q)]^{-\frac 1 2 } h(q) )    ] \,  \mathrm{d}p \mathrm{d}q. \nonumber
\end{align}
By the positivity of $\Phi^c (p, q)$, we have  
\begin{align*}
(\mathcal{L} h , h )
\ge   &\frac 1  2    \int_{\R^3} \int_{\R^3}   \chi_R^2(p) \chi_R^2(q)   \mu^c   (p)  \mu^c   (q)   [ \partial_{p_i}  ( [\mu^c(p)]^{-\frac 1 2 }h(p) )  - \partial_{q_i}  ( [\mu^c(q)]^{-\frac 1 2 } h(q) )    ]  
\\
& \times \Phi^{c, ij}(p, q) [ \partial_{p_j}  ( [\mu^c(p)]^{-\frac 1 2 }h(p) )  - \partial_{q_j}  ( [\mu^c(q)]^{-\frac 1 2 } h(q) )    ] \,  \mathrm{d}p \mathrm{d}q.
\end{align*}
We easily compute that 
\begin{align*}
&\chi_R(p) \chi_R(q)\sqrt{\mu^c(p)}  \sqrt{\mu^c(q)} [ \partial_{p_i}  ( [\mu^c(p)]^{-\frac 1 2 }h(p) )  - \partial_{q_i}  ( [\mu^c(q)]^{-\frac 1 2 } h(q) )    ]
\\
=&     \sqrt{\mu^c(p)}  \sqrt{\mu^c(q)} [  \chi_R(q) \partial_{p_i}  ( [\mu^c(p)]^{-\frac 1 2 } \chi_R(p) h(p) )  -  \chi_R(p) \partial_{q_i}  (   [\mu^c(q)]^{-\frac 1 2 }\chi_R(q)  h(q) )    ]
\\
&-   [  \chi_R(q)\sqrt{\mu^c(q)}    \partial_{p_i}  (\chi_R(p) ) h   (p) - \chi_R(p) \sqrt{\mu^c(p)} \partial_{q_i} (  \chi_R(q)  ) h  (q) ].
\end{align*}
Then since $\Phi_c(p, q)$ is positive symmetric definite matrix, using the following fact
\[
a_i\Phi^{c, ij}   (p, q) a_j + 4b_i\Phi^{c, ij}(p, q) b_j  \ge 4 b_i\Phi^{c, ij}   (p, q) a_j   ,\quad \forall a, b \in \R^3,
\]
we have
\begin{equation}
\label{positivity of Phi c}
(a_i -b_i )\Phi^{c, ij}(p, q) (a_j -    b_j)    \ge \frac 1 2   a_i    \Phi^{c, ij}   (p, q) a_j  -  \frac 3 2  b_i  \Phi^{c, ij} (p, q) b_j  ,\quad \forall a, b \in \R^3.
\end{equation}
This implies
\begin{align}\label{coer-lower}
(\mathcal{L} h , h )
\ge &\frac 1  4     \int_{\R^3} \int_{\R^3}   \mu^c  (p)  \mu^c    (q)   [  \chi_R(q) \partial_{p_i}  ( [\mu^c(p)]^{-\frac 1 2 } \chi_R(p) h(p) )  -  \chi_R(p) \partial_{q_i}  (   [\mu^c(q)]^{-\frac 1 2 }\chi_R(q)  h(q) )    ]
\nonumber\\
&\times     \Phi^{c, ij}(p, q)[  \chi_R(q) \partial_{p_j}  ( [\mu^c(p)]^{-\frac 1 2 } \chi_R(p) h(p) )  -  \chi_R(p) \partial_{q_j}  (   [\mu^c(q)]^{-\frac 1 2 }\chi_R(q)  h(q) )    ]  \,  \mathrm{d}p \mathrm{d}q
\nonumber\\
- &     \int_{\R^3} \int_{\R^3}     [  \chi_R(q)\sqrt{\mu^c(q)}    \partial_{p_i}  (\chi_R(p) ) h   (p) - \chi_R(p) \sqrt{\mu^c(p)} \partial_{q_i} (  \chi_R(q)  ) h   (q) ]
    \Phi^{c, ij}(p, q)
\\
&\times    [  \chi_R(q)\sqrt{\mu^c(q)}    \partial_{p_j}  (\chi_R(p) ) h  (p) - \chi_R(p) \sqrt{\mu^c(p)} \partial_{q_j} (  \chi_R(q)  ) h   (q) ]    \,  \mathrm{d}p \mathrm{d}q :=T_1 +T_2.\nonumber
\end{align}
Using the symmetry of $p$ and $q$, \eqref{cker-der000},  and the fact $  \partial_{p_i}  (\chi_R (p)) = \chi_R'(|p|)  \frac {p_i} {|p|}$ with $|\chi_R' | \lesssim 1/R$, we have
\begin{align}\label{coer-lower2}
|T_2| \le & 4 \int_{\R^3} \int_{\R^3}  |  \chi_R(q) |^2    \mu^c   (q)  |h   (p) |^2      \partial_{p_i}  (\chi_R(p) )   \Phi^{c, ij}(p, q)    \partial_{p_j} (\chi_R(p) )   \,  \mathrm{d}p \mathrm{d}q\nonumber
\\
\lesssim& \frac 1 {R^2}  \int_{\R^3} \int_{\R^3}    \mu^c  (q)  |h   (p) |^2   \frac {p_i} {\langle   p   \rangle} \Phi^{c, ij}(p, q)   \frac {p_j} {\langle p \rangle}    dp dq \lesssim \frac 1 {R^2}   |  h |_{\sigma}^2.
\end{align}
We also compute that
\begin{align*}
&  \chi_R(q) \partial_{p_i}  \big( [\mu^c(p)]^{-\frac 1 2 } \chi_R(p) h(p) \big)  -  \chi_R(p) \partial_{q_i}  \big(   [\mu^c(q)]^{-\frac 1 2 }\chi_R(q)  h(q) \big)
\\
 &\qquad=    \big[  \partial_{p_i}  \big( [\mu^c(p)]^{-\frac 1 2 } \chi_R(p) h(p) \big)  - \partial_{q_i}  \big(   [\mu^c(q)]^{-\frac 1 2 }\chi_R(q)  h(q) \big)\big]    
\\
&\qquad-  \big[    (1-  \chi_R(q))  \partial_{p_i}  \big( [\mu^c(p)]^{-\frac 1 2 } \chi_R(p) h(p) \big)  - (1-  \chi_R(p))    \partial_{q_i}  \big(   [\mu^c(q)]^{-\frac 1 2 }\chi_R(q)  h(q) \big)   \big ]. 
\end{align*}
For $T_1$, we use \eqref{positivity of Phi c} again to have
\begin{align}\label{coer-lower1}
T_1\ge &  \frac 1 8\int_{\R^3} \int_{\R^3}   \mu^c  (p)  \mu^c    (q)    [  \partial_{p_i}  ( [\mu^c(p)]^{-\frac 1 2 } \chi_R(p) h(p) )  - \partial_{q_i}  (   [\mu^c(q)]^{-\frac 1 2 }\chi_R(q)  h(q) )]     \nonumber
\\
&\times     \Phi^{c, ij} (p, q) \big[  \partial_{p_j} \big ( [\mu^c(p)]^{-\frac 1 2 } \chi_R(p) h(p)\big)  - \partial_{q_j}  \big(   [\mu^c(q)]^{-\frac 1 2 }\chi_R(q)  h(q) \big)\big]      \,  \mathrm{d}p \mathrm{d}q\nonumber
\\
&-  \int_{\R^3} \int_{\R^3}    \mu^c   (p)  \mu^c   (q)    \big[    (1-  \chi_R(q))  \partial_{p_i}  \big( [\mu^c(p)]^{-\frac 1 2 } \chi_R(p) h(p) \big) \\
& - (1-  \chi_R(p))    \partial_{q_i} \big (   [\mu^c(q)]^{-\frac 1 2 }\chi_R(q)  h(q) \big)    \big]   
     \Phi^{c, ij}(p, q)    \big[    (1-  \chi_R(q))  \partial_{p_i}  \big( [\mu^c(p)]^{-\frac 1 2 } \chi_R(p) h(p) \big) \nonumber\\
      &- (1-  \chi_R(p))    \partial_{q_i}  \big(   [\mu^c(q)]^{-\frac 1 2 }\chi_R(q)  h(q) \big)    \big]  \,  \mathrm{d}p \mathrm{d}q \nonumber\\
:=&  T_{11} +T_{12}.\nonumber
\end{align}
For  $T_{11}$, we use \eqref{representation for L c} to have
\[
T_{11}=  \frac 1 4(\mathcal{L} (\chi_R h), \chi_R h).
\]
While for  $T_{12}$, by symmetry of $p, q$ and \eqref{positivity of Phi c}, we have
\begin{align*}
|T_{12} |
\lesssim &\int_{\R^3} \int_{\R^3} \mu^c  (p)  \mu^c   (q)     [ (1-  \chi_R(q))  \partial_{p_i}  ( [\mu^c(p)]^{-\frac 1 2 } \chi_R(p) h(p) )    ]   \Phi^{c, ij}(p, q) 
\\
&\times    [  (1-  \chi_R(q))  \partial_{p_j}  ( [\mu^c(p)]^{-\frac 1 2 } \chi_R(p) h(p) )  ]  \,  \mathrm{d}p \mathrm{d}q
\end{align*}
Using the fact that 
\[
(1- \chi_R(q))  \sqrt{\mu^c(q)} \lesssim \frac 1 R ( \mu^c(q) )^{\frac 1 4}
\]
As the estimation of $\mathcal{I}_5$ in Lemma \ref{w-nonlin}, we have 
\[
|T_{12}  | \lesssim \frac 1 R  | \chi_R h  |_{\sigma}^2  \lesssim \frac 1 R  | h |_{\sigma}^2.
\]
For the last term $T_{11}$, we compute that 
\begin{align*}
(\mathcal{L} (\chi_R h), \chi_R h)  =  (\mathbf{L} (\chi_R h), \chi_R h)   - ( (\mathbf{L}-\mathcal{L}) (\chi_R h), \chi_R h)  
\end{align*}
where $\mathbf{L} $ is the classical Landau operator. For the second term,  by Lemma \ref{diff-LL}, we have 
\begin{align*}
| ( (\mathbf{L}-\mathcal{L}) (\chi_R h), \chi_R h)  | \lesssim \frac 1 {c^2} |\chi_R h  |_{H^1_{8}}^2 \lesssim \frac {R^{10}} {c^2}  |h  |_{\sigma}. 
\end{align*}
For the first term, since $\mathbf{L}$ is the classical Landau operator,  using Lemma 5 in \cite{Guo-CMP-2002}  we have
\begin{align*}
(\mathbf{L} (\chi_R h), \chi_R h) \ge & \delta | (I-\mathbf{P}) ( \chi_R h) |_{L^2_{-3/ 2}}^2  \nonumber
\\
\ge & \delta  |  \chi_R h |_{L^2_{-3/ 2}}^2  - \delta | \mathbf{P} (\chi_R h) |_{L^2_{-3/ 2}}^2 
\\
\ge & \delta  |     h    |_{L^2_{-3/ 2}}^2  -  \delta |   (1-\chi_R )  h  |_{L^2_{-3/2}}^2  - \delta |  \mathbf{P} (\chi_R h )   |_{L^2_{-3/2}}^2 \nonumber
\end{align*} 
for some constant $\delta   >0$. We easily compute that 
\begin{align*}
| (1-\chi_R )   h |_{L^2_{-\frac 3 2}}^2 \lesssim \frac 1 R |  (1-\chi_R ) \langle p \rangle  h |_{L^2_{-3/2}}^2 \lesssim   \frac 1 R   | h |_{\sigma}^2.
\end{align*}
Due to our assumption $\mathcal{P} h =0$, we have $\mathcal{P} (\chi_R h ) =\mathcal{P} [(1-\chi_R ) h] $ and further use \eqref{rp-cp} to obtain
\begin{align*}
| \mathbf{P} (\chi_R h ) |_{L^2_{-3/ 2}}^2 \lesssim \left(\frac 1 {c^2}+1\right)| \mathcal{P} (\chi_R h ) |_{L^2_{-3/ 2}}^2  \lesssim  (\frac 1 {c^2} + \frac 1 R )|  h |_{L^2_{-1/ 2}}^2 \lesssim    (\frac 1 {c^2} + \frac 1 R ) |   h  |_{\sigma}^2.
\end{align*}
Collecting the above estimates in \eqref{coer-lower1}, we have 
\begin{align}\label{coer-lower11}
T_1\geq \frac{\delta}{4}  |h|_{L^2_{-3/ 2}}^2-C(\frac 1 {c^2} + \frac 1 R ) |   h  |_{\sigma}^2.
\end{align}
Then we insert \eqref{coer-lower2}  and \eqref{coer-lower11}  in \eqref{coer-lower}, and choose $R=c^{\frac 1 {10}}$ to obtain
\begin{align}\label{coer-lower000}
(\mathcal{L} h , h )
\ge \frac{\delta}{4}  |h|_{L^2_{-3/ 2}}^2-C(\frac {R^{10}} {c^2} + \frac 1 R ) |   h  |_{\sigma}^2\geq  \frac{\delta}{4}  |h|_{L^2_{-3/ 2}}^2-\frac{C}{c}|   h  |_{\sigma}^2.
\end{align}

Finally, we combine \eqref{coerci-hh} and \eqref{coerci-hK} to have
\begin{align*}
\begin{aligned}
&(\mathcal{L} h,  h) \geq \left(\frac{1}{2}-\frac{C}{m}\right)|h|_{\sigma}^2-C(m)\int_{|p|\leq 2m} \langle p\rangle^{-3}|h|^2\,\mathrm{d} p.
\end{aligned}
\end{align*}
Then we choose $m$ suitably large such that $\frac{C}{m}\leq 1/4$. For fixed large $m$ and sufficiently large $c\gg m$, we further make a linear combination of \eqref{coer-lower000} and the above estimate to derive \eqref{coerci0}. 

\end{proof}

Now we turn to the weighted estimate of $\mathcal{L}h$.
\begin{lemma}\label{coerci-w} For any $c\geq1$, $\ell>0$ and $0\leq\vartheta\leq \frac{1}{32}$, there is a large constant  $m>1$ and a uniform constant $C(m,\ell)$ independent of $c$ such that
\begin{align}\label{coerci0-w}
(\mathcal{L}h, w_{\ell}^2 h)\geq \frac14 |h|_{\sigma,w_{\ell}}^2-C(m,\ell) \Big|\big(\mu^c\big)^{1/8}h\Big|_{\sigma}^2.
\end{align}
\end{lemma}
\begin{proof}
As in Lemma \ref{coerci}, we  use \eqref{express-A} and \eqref{express-K} to have
\begin{align}\label{coerci-hhw}
(\mathcal{L}h, w_{\ell}^2 h)= & |h|_{\sigma,w_{\ell}}^2 -\Big(\partial_{p_i}\Big(\sigma^{c,ij}\frac{cp_j}{p^0}\Big)h, w_{\ell}^2 h\Big)+\sum_{\pm} \left(\mathcal{K}_{\pm}h, w_{\ell}^2 h\right).
\end{align}
Denote the last two terms in the right hand side of \eqref{coerci-hhw} as $\mathcal{I}_3$ and $\mathcal{I}_4$ separately. Noting that $\frac{p^0}{c\langle p\rangle}\leq1$,
we use the definition of the weighted dissipation norm $|\cdot|_{\sigma, w_{\ell}}$ in \eqref{ell-norm} to have
\begin{align*}
\mathcal{I}_3=& \Big(\sigma^{c,ij}\frac{cp_j}{p^0}h, w_{\ell}^2 h\partial_{p_i}h\Big)+\Big(\sigma^{c,ij}\frac{cp_j}{p^0}h, \Big(\frac{5\ell p^0}{c\langle p\rangle^2}+\frac{\vartheta p^0}{c\langle p\rangle\ln(\mathrm{e}+t)}\Big) \frac{cp_i}{p^0} w_{\ell}^2 h\Big)\\
\leq&\Big(\int_{\mathbb{R}^3}
\sigma^{c,ij}\frac{cp_i}{p^0}\frac{cp_j}{p^0}w_{\ell}^2(p)|h(p)|^2\,\mathrm{d} p\Big)^{\frac12}\Big(\int_{\mathbb{R}^3}
\sigma^{c,ij}w_{\ell}^2(p)\partial_{p_i}h(p)\partial_{p_j}h(p)\,\mathrm{d} p\Big)^{\frac12}\\
&+2\vartheta |h|_{\sigma,w_{\ell}}^2+\Big(\sigma^{c,ij}\frac{cp_j}{p^0}h, \chi_m(p)\frac{5\ell p^0}{c\langle p\rangle^2} \frac{cp_i}{p^0} w_{\ell}^2 h\Big)\\
& +\Big(\sigma^{c,ij}\frac{cp_j}{p^0}h, (1-\chi_m(p))\frac{5\ell p^0}{c\langle p\rangle^2} \frac{cp_i}{p^0} w_{\ell}^2 h\Big)\\
\leq& \frac{1}{2}\int_{\mathbb{R}^3}
\sigma^{c,ij}w_{\ell}^2(p)\Big(2\partial_{p_i}h(p)\partial_{p_j}h(p)
+\frac{1}{2}\frac{cp_i}{p^0}\frac{cp_j}{p^0}|h(p)|^2\Big)\,\mathrm{d} p\\
&+2\vartheta |h|_{\sigma,w_{\ell}}^2+C(m,\ell)\Big|\big(\mu^c\big)^{1/4}h\Big|_{\sigma}^2+\frac{10\ell}{m}\Big(\int_{\mathbb{R}^3}
\sigma^{c,ij}\frac{cp_i}{p^0}\frac{cp_j}{p^0}w_{\ell}^2(p)|h(p)|^2\,\mathrm{d} p\Big)^{\frac12}\\
&\times\Big(\int_{\mathbb{R}^3}
\sigma^{c,ij}w_{\ell}^2(p)\partial_{p_i}h(p)\partial_{p_j}h(p)\,\mathrm{d} p\Big)^{\frac12}\\
\leq&\Big(\frac12+2\vartheta+\frac{10\ell}{m}\Big)|h|^2_{\sigma,\ell}+C(m,\ell)\Big|\big(\mu^c\big)^{1/4}h\Big|_{\sigma}^2.
\end{align*}
For $\mathcal{I}_4$, we use \eqref{upplow-mu}, \eqref{cker-der01}, and \eqref{express-K} to have
\begin{align*}
\mathcal{I}_4=&\int_{\mathbb{R}^6}\sqrt{\mu^c(p)}\sqrt{\mu^c(q)}
\Phi^{c,ij}\sum_{\pm}\Big[\frac{cq_j}{q^0}h_{\pm}(q)+\partial_{q_j}h_{\pm}(q)\Big]\\
&\times\sum_{\pm}\Big[\Big(\frac{5\ell p^0}{c\langle p\rangle^2}+\frac{\vartheta p^0}{c\langle p\rangle\ln(\mathrm{e}+t)}\Big) \frac{cp_i}{p^0}w_{\ell}^2(p)h_{\pm}(p)+w_{\ell}^2(p)\partial_{p_i}h_{\pm}(p)\Big]
\,\mathrm{d} p\mathrm{d} q\\
\leq& C(\ell)\int_{\mathbb{R}^6} \big[\mu^c(p)\mu^c(q)\big]^{\frac{3}{8}}|p-q|^{-1}\big(|h(p)|+|\nabla_ph(p)|\big)\big(|h(q)|+|\nabla_qh(q)|\big)\,\mathrm{d} p\mathrm{d} q\\
\leq&C(\ell)\Big(\int_{\mathbb{R}^6} \big[\mu^c(p)\mu^c(q)\big]^{\frac{3}{8}}|p-q|^{-1}\big(|h(p)|^2+|\nabla_ph(p)|^2\big)\,\mathrm{d} p\mathrm{d} q\Big)^{\frac12}\\
&\times\Big(\int_{\mathbb{R}^6} \big[\mu^c(p)\mu^c(q)\big]^{\frac{3}{8}}|p-q|^{-1}\big(|h(q)|^2+|\nabla_qh(q)|^2\big)\,\mathrm{d} p\mathrm{d} q\Big)^{\frac12}\\
\leq& C(\ell)\Big|\big(\mu^c\big)^{1/8}h\Big|_{\sigma}^2.
\end{align*}
We collect the estimates of $\mathcal{I}_3$ and $\mathcal{I}_4$ in \eqref{coerci-hhw}, and choose $m>1$ suitably large to obtain \eqref{coerci0-w}.
\end{proof}

\subsubsection{Estimates of  $\varGamma(h,\tilde{h})$ }
In this part, we are devoted to the weighted estimate of the nonlinear relativistic Landau collision operator $\varGamma(h,\tilde{h})$.

\begin{lemma}\label{w-nonlin} Let  $\bar{h}(t,x,p)=[\bar{h}_+(t,x,p), \bar{h}_-(t,x,p)]^{\textit{t}}$, $0\leq\vartheta\leq \frac{1}{32}$, and $\ell\geq0$. There is a uniform constant $C(\ell)>1$ independent of $c$ such that
\begin{align}\label{nonlin-w}
\big|\big(\varGamma(h,\tilde{h}), w_{\ell}^2 \bar{h}\big)\big|\leq C(\ell)  \Big(\Big|\big(\mu^c\big)^{\frac{1}{20}}h\Big|_{L^2}|\tilde{h}|_{\sigma,w_{\ell}} +\Big|\big(\mu^c\big)^{\frac{1}{20}}h\Big|_{\sigma}|\tilde{h}|_{L^2_{w_{\ell}}} \Big)|\bar{h}|_{\sigma,w_{\ell}}.
\end{align}
\end{lemma}
\begin{proof} From \eqref{Nonlin-c}, we have
\begin{align}\label{nonlin-w00}
\big(\varGamma(h,\tilde{h}), w_{\ell}^2 \bar{h}\big)=&-\sum_{\pm}\int_{\mathbb{R}^6} \sqrt{\mu^c(q)}
\Phi^{c,ij}\sum_{\pm}h_{\pm}(q)\Big[\frac{-cp_j}{2p^0}\tilde{h}_{\pm}(p)+\partial_{p_j}\tilde{h}_{\pm}(p)\Big]\nonumber\\
&\times \Big[\Big(\frac{5\ell p^0}{c\langle p\rangle^2}+\frac{\vartheta p^0}{c\langle p\rangle\ln(\mathrm{e}+t)}+\frac{1}{2}\Big) \frac{cp_i}{p^0} w_{\ell}^2 \bar{h}_{\pm}(p)+w_{\ell}^2(p)\partial_{p_i}\bar{h}_{\pm}(p)\Big]
\,\mathrm{d} p\mathrm{d} q\nonumber\\
&+\sum_{\pm}\int_{\mathbb{R}^6} \sqrt{\mu^c(q)}
\Phi^{c,ij}\sum_{\pm}\Big[-\frac{cq_j}{2q^0}h_{\pm}(q)+\partial_{q_j}h_{\pm}(q)\Big]\tilde{h}_{\pm}(p)\\
&\times \Big[\Big(\frac{5\ell p^0}{c\langle p\rangle^2}+\frac{\vartheta p^0}{c\langle p\rangle\ln(\mathrm{e}+t)}+\frac{1}{2}\Big) \frac{cp_i}{p^0} w_{\ell}^2 \bar{h}_{\pm}(p)+w_{\ell}^2(p)\partial_{p_i}\bar{h}_{\pm}(p)\Big]
\,\mathrm{d} p\mathrm{d} q\nonumber\\
:=&\mathcal{I}_5+\mathcal{I}_6.\nonumber
\end{align}

Now we treat $\mathcal{I}_5$.  We will use \eqref{eqiv-norm} in Corollary \ref{corollary-1} to estimate this term. To this end, we make a decomposition of the momentum derivative $\nabla_p$:
$$\nabla_pf(p)=\frac{\nabla_pf(p)\cdot p}{|p|^2}p+\Big(\nabla_pf(p)-\frac{\nabla_pf(p)\cdot p}{|p|^2}p\Big):=\nabla_{\parallel}f(p)+\nabla_{\perp}f(p).$$
Then we use \eqref{cker-der01}, \eqref{cker-der000}, and \eqref{eqiv-norm} to have
\begin{align*}
|\mathcal{I}_5|\leq& 2\Big(\int_{\mathbb{R}^6}\sqrt{\mu^c(q)}|h(q)|
w_{\ell}^2(p)\Big[\Phi^{c}\Big(\frac{-cp}{2p^0}\tilde{h}_{\pm}(p)+\nabla_{\parallel}\tilde{h}_{\pm}(p)+\nabla_{\perp}\tilde{h}_{\pm}(p)\Big)\Big]\\
&\times\Big[\frac{-cp}{2p^0}\tilde{h}_{\pm}(p)+\nabla_{\parallel}\tilde{h}_{\pm}(p)+\nabla_{\perp}\tilde{h}_{\pm}(p)\Big]\,\mathrm{d} p\mathrm{d} q\Big)^{\frac12}\Big(\int_{\mathbb{R}^6}\sqrt{\mu^c(q)}|h(q)|w_{\ell}^2(p)
\\
&\times\Big\{\Phi^{c}\Big[\Big(\frac{5\ell p^0}{c\langle p\rangle^2}+\frac{\vartheta p^0}{c\langle p\rangle\ln(\mathrm{e}+t)}+\frac{1}{2}\Big) \frac{cp}{p^0} \bar{h}_{\pm}(p)+\nabla_{\parallel}\bar{h}_{\pm}(p)+\nabla_{\perp}\bar{h}_{\pm}(p)\Big]\Big\}\\
&\cdot\Big[\Big(\frac{5\ell p^0}{c\langle p\rangle^2}+\frac{\vartheta p^0}{c\langle p\rangle\ln(\mathrm{e}+t)}+\frac{1}{2}\Big) \frac{cp}{p^0}  \bar{h}_{\pm}(p)+\nabla_{\parallel}\bar{h}_{\pm}(p)+\nabla_{\perp}\bar{h}_{\pm}(p)\Big]
\,\mathrm{d} p\mathrm{d} q\Big)^{\frac12}\\
\leq& C(\ell)\Big(\int_{\mathbb{R}^6}\big[\mu^c(q)\big]^{\frac{1}{4}}|h(q)||p-q|^{-1}
w_{\ell}^2(p)\\
&\times\Big[\frac{\big(p^0\big)^3}{c^3\langle p\rangle^2}\big|\nabla_{\parallel}\tilde{h}_{\pm}(p)\big|^2+\frac{p^0}{c}\left(\big|\nabla_{\perp}\tilde{h}_{\pm}(p)\big|^2+\big|\tilde{h}_{\pm}(p)\big|^2\right)\Big]\,\mathrm{d} p\mathrm{d} q\Big)^{\frac12}\\
&\times\Big(\int_{\mathbb{R}^6}\big[\mu^c(q)\big]^{\frac{1}{4}}|h(q)||p-q|^{-1}
w_{\ell}^2(p)\\
&\times\Big[\frac{\big(p^0\big)^3}{c^3\langle p\rangle^2}\big|\nabla_{\parallel}\bar{h}_{\pm}(p)\big|^2+\frac{p^0}{c}\left(\big|\nabla_{\perp}\bar{h}_{\pm}(p)\big|^2+\big|\bar{h}_{\pm}(p)\big|^2\right)\Big]\,\mathrm{d} p\mathrm{d} q\Big)^{\frac12}.
\end{align*}
Noting that
\begin{align*}
&\int_{\mathbb{R}^3}\big[\mu^c(q)\big]^{\frac{1}{4}}|h(q)||p-q|^{-1}\,\mathrm{d} q\\
&\qquad\leq \Big(\int_{\mathbb{R}^3}\big[\mu^c(q)\big]^{\frac{1}{4}}|p-q|^{-2}\,\mathrm{d} q\Big)^{\frac12}\Big(\int_{\mathbb{R}^3}\big[\mu^c(q)\big]^{\frac{1}{4}}|h(q)|^2\,\mathrm{d} q\Big)^{\frac12}\\
&\qquad\leq C\langle p\rangle^{-1} \Big|\big(\mu^c\big)^{1/8}h\Big|_{L^2}
\end{align*}
by \eqref{upplow-mu}, we can further bound $\mathcal{I}_5$ as
\begin{align*}
|\mathcal{I}_5|\leq&C(\ell)\Big|\big(\mu^c\big)^{1/8}h\Big|_{L^2}|\tilde{h}|_{\sigma,w_{\ell}} |\bar{h}|_{\sigma,w_{\ell}}.
\end{align*}

To estimate $\mathcal{I}_6$, we divide its estimation into two cases: $|p|\leq \frac{10}{9}\langle q\rangle$ and
$|p|> \frac{10}{9}\langle q\rangle$.

\noindent $\bullet$ Case 1 ($|p|\leq \frac{10}{9}\langle q\rangle$): Note that $\mu^c(q)\leq C \big(\mu^c(p)\big)^{3/5}$ in this case. This implies that for $0\leq\vartheta\leq \frac{1}{32}$,
$$\sqrt{\mu^c(q)}w_{\ell}^2(p)\leq w_{\ell}^2(p)\big(\mu^c(q)\big)^{\frac{1}{8}}\big(\mu^c(p)\big)^{1/4} \leq C(\ell)\big(\mu^c(q)\big)^{\frac{1}{8}}\big(\mu^c(p)\big)^{\frac{1}{16}}$$
by \eqref{upplow-mu}.
Then we use \eqref{cker-der01} to have
\begin{align*}
|\mathcal{I}_6|\leq& C(\ell)\int_{\mathbb{R}^6} \big[\mu^c(q)\big]^{\frac{1}{9}}\big[\mu^c(p)\big]^{\frac{1}{18}}|p-q|^{-1}\big(|h(q)|+|\nabla_qh(q)|\big)|\tilde{h}(p)|
\big(|\bar{h}(p)|+|\nabla_p\bar{h}(p)|\big)\,\mathrm{d} p\mathrm{d} q\\
\leq&C(\ell)\int_{\mathbb{R}^3}\Big(\int_{\mathbb{R}^3} \big[\mu^c(q)\big]^{\frac{1}{20}}|p-q|^{-2}\,\mathrm{d} q\Big)^{\frac12}\Big(\int_{\mathbb{R}^3} \big[\mu^c(q)\big]^{\frac{1}{8}}\big(|h(q)|^2+|\nabla_qh(q)|^2\big)\,\mathrm{d} q\Big)^{\frac12}\\
&\times \big[\mu^c(p)\big]^{\frac{1}{18}}\big(|h(q)|+|\nabla_qh(q)|\big)|\tilde{h}(p)|
\big(|\bar{h}(p)|+|\nabla_p\bar{h}(p)|\big)\,\mathrm{d} p\\
\leq& C(\ell)\Big|\big(\mu^c\big)^{\frac{1}{20}}h\Big|_{\sigma}|\tilde{h}|_{L^2} |\bar{h}|_{\sigma} .
\end{align*}

\noindent $\bullet$ Case 2 ($|p|> \frac{10}{9}\langle q\rangle$): In this case, it holds that $|p-q|>\langle p\rangle-\langle q\rangle>\frac{\langle p\rangle}{10}$. By the definition of the norm $|\cdot|_{\sigma, w_{\ell}}$ in \eqref{ell-norm}, \eqref{cker-der01}, and \eqref{eqiv-sigma}, we get
\begin{align*}
|\mathcal{I}_6|\leq& \sum_{\pm}\Big(\int_{\mathbb{R}^6} \sqrt{\mu^c(q)}
\Phi^{c,ij}\sum_{\pm}\Big[-\frac{cq_i}{2q^0}h_{\pm}(q)+\partial_{q_i}h_{\pm}(q)\Big]\\
&\times\sum_{\pm}\Big[-\frac{cq_j}{2q^0}h_{\pm}(q)+\partial_{q_j}h_{\pm}(q)\Big]w_{\ell}^2(p)\big|\tilde{h}_{\pm}(p)\big|^2
\,\mathrm{d} p\mathrm{d} q\Big)^{\frac12}\\
&\times\Big(\int_{\mathbb{R}^6} \sqrt{\mu^c(q)}
\Phi^{c,ij}\sum_{\pm}w_{\ell}^2(p)\Big[\Big(\frac{5\ell p^0}{c\langle p\rangle^2}+\frac{\vartheta p^0}{c\langle p\rangle\ln(\mathrm{e}+t)}+\frac{1}{2}\Big) \frac{cp_i}{p^0} \bar{h}_{\pm}(p)+\partial_{p_i}\bar{h}_{\pm}(p)\Big]\\
&\times \Big[\Big(\frac{5\ell p^0}{c\langle p\rangle^2}+\frac{\vartheta p^0}{c\langle p\rangle\ln(\mathrm{e}+t)}+\frac{1}{2}\Big) \frac{cp_j}{p^0} \bar{h}_{\pm}(p)+\partial_{p_j}\bar{h}_{\pm}(p)\Big]
\,\mathrm{d} p\mathrm{d} q\Big)^{\frac12}\\
\leq& C(\ell)\Big(\int_{\mathbb{R}^6}\big[\mu^c(q)\big]^{\frac{1}{3}}\frac{p^0}{c\langle p\rangle}
\Big[|h(q)|^2+\big|\nabla_{q}h(q)\big|^2\Big]w_{\ell}^2(p)\big|\tilde{h}(p)\big|^2
\,\mathrm{d} p\mathrm{d} q\Big)^{\frac12}|\bar{h}|_{\sigma,w_{\ell}}\\
\leq&C(\ell)\Big|\big(\mu^c\big)^{1/8}h\Big|_{\sigma}|\tilde{h}|_{L^2_{w_{\ell}}} |\bar{h}|_{\sigma,w_{\ell}}.
\end{align*}
We collect the estimates in the above two cases to have
$$|\mathcal{I}_6|\leq C(\ell)\Big|\big(\mu^c\big)^{1/9}h\Big|_{\sigma}|\tilde{h}|_{L^2_{w_{\ell}}} |\bar{h}|_{\sigma,w_{\ell}}.$$
Then we plug the estimates of $\mathcal{I}_5$ and $\mathcal{I}_6$ in \eqref{nonlin-w00} to derive \eqref{nonlin-w}.

\end{proof}

\subsection{Difference estimates of the  collision operators}
In this subsection, we will establish the difference estimates of $\mathcal{L}h-\mathbf{L}h$ and $\varGamma(h,\tilde{h})-\Gamma(h,\tilde{h})$, which will be crucially used in the proof of classical limit from the RVML system \eqref{main1-00} to the VPL system \eqref{main2-00}.

\subsubsection{Estimates of  $\mathcal{L}h-\mathbf{L}h$ }
Before its estimation, we first list the explicit form of $\mathbf{A}_{\pm}h$ and $\mathbf{K}_{\pm}h$ in \eqref{linearAK-inf}.
 As in \cite[Lemma 1]{Guo-CMP-2002}, we can obtain
\begin{align*}
\mathbf{A}_{\pm}h=&-2\partial_{p_i}\left(\sigma^{\infty,ij}\partial_{p_j}h_{\pm}\right)
+\frac12\sigma^{\infty,ij}p_i p_j h_{\pm}
-\partial_{p_i}\big(\sigma^{\infty,ij}p_j\big)h_{\pm},\\
\mathbf{K}_{\pm}h=&\left(\mu^{\infty}(p)\right)^{-\frac12}\partial_{p_i}\Big[\mu^{\infty}(p)\int_{\mathbb{R}^3}
\Phi^{\infty,ij}\sqrt{\mu^{\infty}(q)}\partial_{q_j}\left(h_{+}+h_-\right)\,\mathrm{d} q\Big]\nonumber\\
&+\left(\mu^{\infty}(p)\right)^{-\frac12}\partial_{p_i}\Big[\mu^{\infty}(p)\int_{\mathbb{R}^3}
\Phi^{\infty,ij}\sqrt{\mu^{\infty}(q)}\frac{q_j}{2}\left(h_{+}+h_-\right)\,\mathrm{d} q\Big].
\end{align*}

We have the following estimate of $\mathcal{L}h-\mathbf{L}h$.
\begin{lemma}\label{diff-LL} For any $c\geq1$, there is a uniform constant $C>0$ independent of $c$ such that
\begin{align}\label{diff-LL0}
|(\mathcal{L}h-\mathbf{L}h, \tilde{h})|\leq \frac{C}{c^2}|h|_{H^1_{p,8}}|\tilde{h}|_{H^1_{p,-2}}.
\end{align}
Here and below, for $l\in \mathbb{R}$, we use the following norms
$$|h|^2_{H^1_{p,l}}:=\sum_{i=0}^1\sum_{\pm}\left|\lag p\rag^{l}\nabla_p^i h_{\pm}\right|^2_{L^2},\qquad \|h\|^2_{H^1_{p,l}}:=\sum_{i=0}^1\sum_{\pm}\left\|\lag p\rag^{l}\nabla_p^i h_{\pm}\right\|^2.$$
\end{lemma}
\begin{proof} 

 Note that
\begin{align}\label{LLe}
|(\mathcal{L}h-\mathbf{L}h, \tilde{h})|\leq|(\mathcal{A}h-\mathbf{A}h), \tilde{h})|+|(\mathcal{K}h-\mathbf{K}h), \tilde{h})|:=\mathcal{I}_7+\mathcal{I}_8.
\end{align}
We first estimate $\mathcal{I}_7$.

\noindent $\bullet$ {\it Estimation of $\mathcal{I}_7$:} By the explicit expressions of $\mathcal{A}_{\pm}$ and $\mathbf{A}_{\pm}$, we integrate by parts with respect to $p$ to get
\begin{align*}
\mathcal{I}_7\leq&  \sum_{\pm}\Big|\Big(\big(\sigma^{c,ij}-\sigma^{\infty,ij}\big)\big(2\partial_{p_j}h_{\pm}+\frac{cp_j}{p^0}h_{\pm}\big), \partial_{p_i}\tilde{h}_{\pm}\Big)\Big|\\
&+\sum_{\pm}\Big|\Big(\big(\sigma^{c,ij}
-\sigma^{\infty,ij}\big)\frac{cp_j}{p^0}\partial_{p_i}h_{\pm}, \tilde{h}_{\pm}\Big)\Big|+\sum_{\pm}\Big|\Big(\frac12\big(\sigma^{c,ij}-\sigma^{\infty,ij}\big)\frac{cp_i}{p^0}\frac{cp_j}{p^0}h_{\pm},\tilde{h}_{\pm}\Big)\Big|\\
&+\sum_{\pm}\Big|\Big(\sigma^{\infty,ij}\big(\frac{c}{p^0}-1\big)p_jh_{\pm}, \partial_{p_i}\tilde{h}_{\pm}\Big)\Big|+\sum_{\pm}\Big|\Big(\sigma^{\infty,ij}
\big(\frac{c}{p^0}-1\big)p_j\partial_{p_i}h_{\pm}, \tilde{h}_{\pm}\Big)\Big|\\
&+\sum_{\pm}\Big|\Big(\frac12\sigma^{\infty,ij}\Big(\frac{c^2}{\big(p^0\big)^2}-1\Big)p_ip_jh_{\pm},\tilde{h}_{\pm}\Big)\Big|.
\end{align*}
Note that
\begin{align}\label{diff-pp}
\begin{aligned}
1-\frac{c}{p^0}=&\frac{|p|^2}{p^0(c+p^0)}\leq \frac{|p|^2}{2c^2},\qquad 1-\frac{c^2}{\big(p^0\big)^2}=\big(1-\frac{c}{p^0}\big)\big(1+\frac{c}{p^0}\big)\leq \frac{|p|^2}{c^2},\\
\sigma^{\infty,ij}=&\int_{{\mathbb{R}}^{3}}\Phi^{\infty,ij}(p,q)\mu^{\infty}(q)\, \mathrm{d} q\leq \int_{{\mathbb{R}}^{3}}|p-q|^{-1}\mu^{\infty}(q)\, \mathrm{d} q\leq \frac{C}{\langle p\rangle},
\end{aligned}
\end{align}
and
\begin{align*}
&\big|\sigma^{c,ij}(p)-\sigma^{\infty,ij}(p)\big|\\
&\qquad \leq \int_{{\mathbb{R}}^{3}}\big|\Phi^{c,ij}(p,q)-\Phi^{\infty,ij}(p,q)\big|\mu^{\infty}(q)\, \mathrm{d} q+\int_{{\mathbb{R}}^{3}}\Phi^{\infty,ij}(p,q)\big|\mu^{c}(q)-\mu^{\infty}(q)\big|\, \mathrm{d} q\\
&\qquad \leq  \int_{{\mathbb{R}}^{3}}\frac{C}{c^2}\langle p\rangle^5\langle q\rangle^5|p-q|^{-1}\mu^{\infty}(q)\, \mathrm{d} q+\int_{{\mathbb{R}}^{3}}\frac{C}{c^2}|p-q|^{-1}\exp\left\{-\frac{|p|}{3}\min\{|p|,\frac{4c}{3} \right\}\big|\, \mathrm{d} q\\
&\qquad \leq \frac{C}{c^2}\langle p\rangle^4
\end{align*}
by \eqref{diff-mu0} and \eqref{ker-cinf}. Then we collect the above estimates to arrive at that
\begin{align*}
\mathcal{I}_7\leq \frac{C}{c^2}|h|_{H^1_{p,8}}|\tilde{h}|_{H^1_{p,-2}}.
\end{align*}

\noindent $\bullet$ {\it Estimation of $\mathcal{I}_8$:} By the explicit expressions of $\mathcal{K}_{\pm}$ and $\mathbf{K}_{\pm}$, we integrate by parts with respect to $p$ to get
\begin{align*}
\mathcal{I}_8=&-\int_{\mathbb{R}^6}\sqrt{\mu^c(p)}\sqrt{\mu^c(q)}
\Phi^{c,ij}\sum_{\pm}\Big[\frac{cp_i}{2p^0}\tilde{h}_{\pm}(p)+\partial_{p_i}\tilde{h}_{\pm}(p)\Big]
\sum_{\pm}\Big[\frac{cq_j}{2q^0}h_{\pm}(q)+\partial_{q_j}h_{\pm}(q)\Big]\,\mathrm{d} p\mathrm{d} q\\
&+\int_{\mathbb{R}^6}\sqrt{\mu^{\infty}(p)}\sqrt{\mu^{\infty}(q)}
\Phi^{\infty,ij}\sum_{\pm}\Big[\frac{p_i}{2}\tilde{h}_{\pm}(p)+\partial_{p_i}\tilde{h}_{\pm}(p)\Big]
\sum_{\pm}\Big[\frac{q_j}{2}h_{\pm}(q)+\partial_{q_j}h_{\pm}(q)\Big]\,\mathrm{d} p\mathrm{d} q\\
=&-\int_{\mathbb{R}^6}\left(\sqrt{\mu^c(p)}\sqrt{\mu^c(q)} \Phi^{c,ij}-\sqrt{\mu^{\infty}(p)}\sqrt{\mu^{\infty}(q)}
\Phi^{\infty,ij}\right)\\
&\times\sum_{\pm}\Big[\frac{cp_i}{2p^0}\tilde{h}_{\pm}(p)+\partial_{p_i}\tilde{h}_{\pm}(p)\Big]
\sum_{\pm}\Big[\frac{cq_j}{2q^0}h_{\pm}(q)+\partial_{q_j}h_{\pm}(q)\Big]\,\mathrm{d} p\mathrm{d} q\\
&+\int_{\mathbb{R}^6}\sqrt{\mu^{\infty}(p)}\sqrt{\mu^{\infty}(q)}
\Phi^{\infty,ij}\Big(\sum_{\pm}\Big[\frac{p_i}{2}\tilde{h}_{\pm}(p)+\partial_{p_i}\tilde{h}_{\pm}(p)\Big]
\sum_{\pm}\Big[\frac{q_j}{2}h_{\pm}(q)+\partial_{q_j}h_{\pm}(q)\Big]\\
&-\sum_{\pm}\Big[\frac{cp_i}{2p^0}\tilde{h}_{\pm}(p)+\partial_{p_i}\tilde{h}_{\pm}(p)\Big]
\sum_{\pm}\Big[\frac{cq_j}{2q^0}h_{\pm}(q)+\partial_{q_j}h_{\pm}(q)\Big]\Big)\,\mathrm{d} p\mathrm{d} q\\
:=&\mathcal{I}_8^1+\mathcal{I}_8^2.
\end{align*}
Note that
\begin{align}\label{mumu-cinf}
&\big|\sqrt{\mu^c(p)}\sqrt{\mu^c(q)} \Phi^{c,ij}-\sqrt{\mu^{\infty}(p)}\sqrt{\mu^{\infty}(q)}
\Phi^{\infty,ij}\big|\nonumber\\
&\qquad \leq \big|\big(\sqrt{\mu^c(p)}\sqrt{\mu^c(q)} -\sqrt{\mu^{\infty}(p)}\sqrt{\mu^{\infty}(q)}\big)
\Phi^{c,ij}\big|\nonumber\\
&\qquad \quad+\big|\sqrt{\mu^{\infty}(p)}\sqrt{\mu^{\infty}(q)}\big(\Phi^{c,ij}-
\Phi^{\infty,ij}\big)\big|\\
&\qquad \leq \frac{C}{c^2}\exp\left\{-\frac{|p|}{6}\min\{|p|,\frac{4c}{3}-\frac{|q|}{6}\min\{|q|,\frac{4c}{3}\} \right\}\times\min\Big\{\frac{p^0}{c} \langle q\rangle^4, \frac{q^0}{c} \langle p\rangle^4\Big\}|p-q|^{-1}\nonumber\\
&\qquad \quad+\frac{C}{c^2}\exp\left\{-\frac{|p|^2+|q|^2}{4} \right\}\langle p\rangle^5\langle q\rangle^5|p-q|^{-1}\nonumber\\
&\qquad \leq \frac{C}{c^2}\exp\left\{-\frac{|p|+|q|}{8} \right\}|p-q|^{-1}\nonumber
\end{align}
by \eqref{upplow-mu}, \eqref{diff-mu1}, \eqref{cker-der01}, and \eqref{ker-cinf}. Then we can bound $\mathcal{I}_8^1$ as
\begin{align*}
|\mathcal{I}_8^1|\leq&\frac{C}{c^2}\Big(\int_{\mathbb{R}^6}\exp\left\{-\frac{|p|+|q|}{10} \right\}|p-q|^{-2}\big(|h(q)|^2+\big|\nabla_{q}h(q)\big|^2\big)
\,\mathrm{d} p\mathrm{d} q\Big)^{\frac12}\\
&\times\Big(\int_{\mathbb{R}^6}\exp\left\{-\frac{|p|+|q|}{10} \right\}\big(|\tilde{h}(p)|^2+\big|\nabla_{p}\tilde{h}(p)\big|^2\big)
\,\mathrm{d} p\mathrm{d} q\Big)^{\frac12}\\
\leq &\frac{C}{c^2}|h|_{H^1_p}|\tilde{h}|_{H^1_{p,-2}}.
\end{align*}
For $\mathcal{I}_8^2$, we use \eqref{diff-pp} to have
\begin{align*}
|\mathcal{I}_8^2|\leq&\int_{\mathbb{R}^6}\exp\left\{-\frac{|p|^2+|q|^2}{4} \right\}|p-q|^{-1}\Big[\frac{|p_iq_j|}{4}\big(1-\frac{c^2}{p^0q^0}\big)|\tilde{h}(p)||h(q)|\\
&+\frac{|p_i|}{2}\big(1-\frac{c}{q^0}\big)|\tilde{h}(p)|\big|\partial_{q_j}h(q)\big|+\frac{|q_j|}{2}\big(1-\frac{c}{p^0}\big)\big|\partial_{p_i}\tilde{h}(p)\big||h(q)|\Big]\,\mathrm{d} p\mathrm{d} q\\
\leq&\frac{C}{c^2}\Big(\int_{\mathbb{R}^6}\exp\left\{-\frac{|p|^2+|q|^2}{8} \right\}|p-q|^{-2}\big(|h(q)|^2+\big|\nabla_{q}h(q)\big|^2\big)
\,\mathrm{d} p\mathrm{d} q\Big)^{\frac12}\\
&\times\Big(\int_{\mathbb{R}^6}\exp\left\{-\frac{|p|^2+|q|^2}{8} \right\}\big(|\tilde{h}(p)|^2+\big|\nabla_{p}\tilde{h}(p)\big|^2\big)
\,\mathrm{d} p\mathrm{d} q\Big)^{\frac12}\\
\leq &\frac{C}{c^2}|h|_{H^1_p}|\tilde{h}|_{H^1_{p,-2}}.
\end{align*}
Then we combine the estimate of $\mathcal{I}_8^1$ and $\mathcal{I}_8^2$ to get
\begin{align*}
|\mathcal{I}_8|\leq &\frac{C}{c^2}|h|_{H^1_p}|\tilde{h}|_{H^1_{p,-2}}.
\end{align*}
Finally, we collect the estimate of $\mathcal{I}_7$ and $\mathcal{I}_8$ in \eqref{LLe} to deduce \eqref{diff-LL0}.
\end{proof}

\subsubsection{Estimates of  $\varGamma(h,\tilde{h})-\Gamma(h,\tilde{h})$ }
In this part, we deal with the difference of the nonlinear collision operators: $\varGamma(h,\tilde{h})-\Gamma(h,\tilde{h})$ .

\begin{lemma}\label{diff-GG} Let  $\bar{h}(t,x,p)=[\bar{h}_+(t,x,p), \bar{h}_-(t,x,p)]^{\textit{t}}$. There is a uniform constant $C>1$ independent of $c$ such that
\begin{align}\label{diff-GG0}
\big|\big(\varGamma(h,\tilde{h})-\Gamma(h,\tilde{h}),  \bar{h}\big)\big|\leq \frac{C}{c^2}|h|_{H^1_p}|\tilde{h}|_{H^1_{p,8}}|\bar{h}|_{H^1_{p,-2}}.
\end{align}
\end{lemma}
\begin{proof} From \eqref{Nonlin-c}, we have
\begin{align*}
&\big(\varGamma(h,\tilde{h})-\Gamma(h,\tilde{h}),  \bar{h}\big)\\
&\qquad=-\sum_{\pm}\int_{\mathbb{R}^6} \sqrt{\mu^c(q)}
\Phi^{c,ij}\sum_{\pm}h_{\pm}(q)\Big[\frac{-cp_j}{2p^0}\tilde{h}_{\pm}(p)+\partial_{p_j}\tilde{h}_{\pm}(p)\Big]\nonumber\\
&\qquad\quad\times \Big[\frac{cp_i}{2p^0}\bar{h}_{\pm}(p)+\partial_{p_i}\bar{h}_{\pm}(p)\Big]
\,\mathrm{d} p\mathrm{d} q\nonumber\\
&\qquad\quad+\sum_{\pm}\int_{\mathbb{R}^6} \sqrt{\mu^{\infty}(q)}
\Phi^{\infty,ij}\sum_{\pm}h_{\pm}(q)\Big[\frac{-p_j}{2}\tilde{h}_{\pm}(p)+\partial_{p_j}\tilde{h}_{\pm}(p)\Big]\nonumber\\
&\qquad\quad\times \Big[\frac{p_i}{2}\bar{h}_{\pm}(p)+\partial_{p_i}\bar{h}_{\pm}(p)\Big]
\,\mathrm{d} p\mathrm{d} q\nonumber\\
&\qquad\quad+\sum_{\pm}\int_{\mathbb{R}^6} \sqrt{\mu^c(q)}
\Phi^{c,ij}\sum_{\pm}\Big[-\frac{cq_j}{2q^0}h_{\pm}(q)+\partial_{q_j}h_{\pm}(q)\Big]\tilde{h}_{\pm}(p)\\
&\qquad\quad\times \Big[\frac{cp_i}{p^0}\bar{h}_{\pm}(p)+\partial_{p_i}\bar{h}_{\pm}(p)\Big]
\,\mathrm{d} p\mathrm{d} q\nonumber\\
&\qquad\quad-\sum_{\pm}\int_{\mathbb{R}^6} \sqrt{\mu^{\infty}(q)}
\Phi^{\infty,ij}\sum_{\pm}\Big[-\frac{q_j}{2}h_{\pm}(q)+\partial_{q_j}h_{\pm}(q)\Big]\tilde{h}_{\pm}(p)\\
&\qquad\quad\times \Big[\frac{p_i}{2}\bar{h}_{\pm}(p)+\partial_{p_i}\bar{h}_{\pm}(p)\Big]
\,\mathrm{d} p\mathrm{d} q\nonumber.
\end{align*}
Then we further rewrite the above equality as
\begin{align}\label{nonlin-w0}
&\big(\varGamma(h,\tilde{h})-\Gamma(h,\tilde{h}),  \bar{h}\big)\nonumber\\
&\qquad=-\sum_{\pm}\int_{\mathbb{R}^6}\big( \sqrt{\mu^c(q)}
\Phi^{c,ij}-\sqrt{\mu^{\infty}(q)}
\Phi^{\infty,ij}\big)\sum_{\pm}h_{\pm}(q)\nonumber\\
&\qquad\quad\times \Big[\frac{-cp_j}{2p^0}\tilde{h}_{\pm}(p)+\partial_{p_j}\tilde{h}_{\pm}(p)\Big]
\Big[\frac{cp_i}{2p^0}\bar{h}_{\pm}(p)+\partial_{p_i}\bar{h}_{\pm}(p)\Big]
\,\mathrm{d} p\mathrm{d} q\nonumber\\
&\qquad\quad+\sum_{\pm}\int_{\mathbb{R}^6} \sqrt{\mu^{\infty}(q)}
\Phi^{\infty,ij}\sum_{\pm}h_{\pm}(q)\Big[\Big(\frac{-p_j}{2}\tilde{h}_{\pm}(p)
+\partial_{p_j}\tilde{h}_{\pm}(p)\Big)\Big(\frac{p_i}{2}\bar{h}_{\pm}(p)+\partial_{p_i}\bar{h}_{\pm}(p)\Big)\nonumber\\
&\qquad\quad- \Big(\frac{cp_j}{2p^0}\tilde{h}_{\pm}(p)+\partial_{p_j}\tilde{h}_{\pm}(p)\Big)
\Big(\frac{cp_i}{2p^0}\bar{h}_{\pm}(p)+\partial_{p_i}\bar{h}_{\pm}(p)\Big)\Big]
\,\mathrm{d} p\mathrm{d} q\nonumber\\
&\qquad\quad+\sum_{\pm}\int_{\mathbb{R}^6} \big( \sqrt{\mu^c(q)}
\Phi^{c,ij}-\sqrt{\mu^{\infty}(q)}
\Phi^{\infty,ij}\big)\tilde{h}_{\pm}(p)\\
&\qquad\quad\times\sum_{\pm}\Big[-\frac{cq_j}{2q^0}h_{\pm}(q)+\partial_{q_j}h_{\pm}(q)\Big] \Big[\frac{cp_i}{p^0}\bar{h}_{\pm}(p)+\partial_{p_i}\bar{h}_{\pm}(p)\Big]
\,\mathrm{d} p\mathrm{d} q\nonumber\\
&\qquad\quad-\sum_{\pm}\int_{\mathbb{R}^6} \sqrt{\mu^{\infty}(q)}
\Phi^{\infty,ij}\tilde{h}_{\pm}(p)\Big[\Big(\frac{-p_i}{2}\bar{h}_{\pm}(p)+\partial_{p_i}\bar{h}_{\pm}(p)\Big)
\sum_{\pm}\Big(-\frac{q_j}{2}h_{\pm}(q)+\partial_{q_j}h_{\pm}(q)\Big)\nonumber\\
&\qquad\quad- \Big(\frac{-cp_i}{p^0}\bar{h}_{\pm}(p)+\partial_{p_i}\bar{h}_{\pm}(p)\Big)
\sum_{\pm}\Big(-\frac{cq_j}{2q^0}h_{\pm}(q)+\partial_{q_j}h_{\pm}(q)\Big)\Big]
\,\mathrm{d} p\mathrm{d} q\nonumber\\
&\qquad:=\mathcal{I}_9+\mathcal{I}_{10}+\mathcal{I}_{11}+\mathcal{I}_{12}.\nonumber
\end{align}

\noindent $\bullet$ {\it Estimation of $\mathcal{I}_9$ and $\mathcal{I}_{11}$:} Via similar estimation in \eqref{mumu-cinf}, one has
\begin{align*}\label{mumu-cinf}
&\big|\sqrt{\mu^c(q)} \Phi^{c,ij}-\sqrt{\mu^{\infty}(q)}
\Phi^{\infty,ij}\big|\nonumber\\
&\qquad \leq \big|\big(\sqrt{\mu^c(q)} -\sqrt{\mu^{\infty}(q)}\big)
\Phi^{c,ij}\big|+\big|\sqrt{\mu^{\infty}(q)}\big(\Phi^{c,ij}-
\Phi^{\infty,ij}\big)\big|\\
&\qquad \leq \frac{C}{c^2}\exp\left\{-\frac{|q|}{6}\min\{|q|,c\} \right\}\times\min\Big\{\frac{p^0}{c} \langle q\rangle^4, \frac{q^0}{c} \langle p\rangle^4\Big\}|p-q|^{-1}\nonumber\\
&\qquad \quad+\frac{C}{c^2}\exp\left\{-\frac{|q|^2}{4} \right\}\langle p\rangle^5\langle q\rangle^5|p-q|^{-1}\nonumber\\
&\qquad \leq \frac{C}{c^2}\exp\left\{-\frac{|q|}{8} \right\}\langle p\rangle^5|p-q|^{-1}.\nonumber
\end{align*}
Then we have
\begin{align*}
|\mathcal{I}_9|\leq&\frac{C}{c^2}\Big(\int_{\mathbb{R}^3}\exp\left\{-\frac{|q|}{10} \right\}|p-q|^{-2}\langle p\rangle^2\,\mathrm{d} q\Big)^{\frac12}\Big(\int_{\mathbb{R}^3}|h(q)|^2
\,\mathrm{d} q\Big)^{\frac12}\\
&\times\Big(\int_{\mathbb{R}^3}\langle p\rangle^{16}\big(|\tilde{h}(p)|^2+\big|\nabla_{p}\tilde{h}(p)\big|^2\big)
\,\mathrm{d} p\Big)^{\frac12}\Big(\int_{\mathbb{R}^3}\langle p\rangle^{-4}\big(|\bar{h}(p)|^2+\big|\nabla_{p}\bar{h}(p)\big|^2\big)
\,\mathrm{d} p\Big)^{\frac12}\\
\leq &\frac{C}{c^2}|h|_{L^2}|\tilde{h}|_{H^1_{p,8}}|\bar{h}|_{H^1_{p,-2}},
\end{align*}
and
\begin{align*}
|\mathcal{I}_{11}|\leq&\frac{C}{c^2}\Big(\int_{\mathbb{R}^3}\exp\left\{-\frac{|q|}{10} \right\}|p-q|^{-2}\langle p\rangle^2\,\mathrm{d} q\Big)^{\frac12}\Big(\int_{\mathbb{R}^3}|h(q)|^2+|\nabla_qh(q)|^2
\,\mathrm{d} q\Big)^{\frac12}\\
&\times\Big(\int_{\mathbb{R}^3}\langle p\rangle^{14}|\tilde{h}(p)|^2
\,\mathrm{d} p\Big)^{\frac12}\Big(\int_{\mathbb{R}^3}\langle p\rangle^{-4}\big(\bar{h}(p)|^2+\big|\nabla_{p}\bar{h}(p)\big|^2\big)
\,\mathrm{d} p\Big)^{\frac12}\\
\leq &\frac{C}{c^2}|h|_{H^1_p}|\tilde{h}|_{L^2_7}|\bar{h}|_{H^1_{p,-2}}.
\end{align*}

\noindent $\bullet$ {\it Estimation of $\mathcal{I}_{10}$ and $\mathcal{I}_{12}$:}
By \eqref{diff-pp}, we get
\begin{align*}
|\mathcal{I}_{10}|\leq&\int_{\mathbb{R}^6}\exp\left\{-\frac{|q|^2}{4} \right\}|p-q|^{-1}|h(q)|\Big[\frac{|p_ip_j|}{4}\big(1-\frac{c^2}{\big(p^0)^2}\big)|\tilde{h}(p)||\bar{h}(p)|\\
&+\frac{|p_j|}{2}\big(1-\frac{c}{p^0}\big)|\tilde{h}(p)|\big|\partial_{p_i}\bar{h}(q)\big|
+\frac{|p_i|}{2}\big(1-\frac{c}{p^0}\big)\big|\partial_{p_j}\tilde{h}(p)\big||\bar{h}(p)|\Big]\,\mathrm{d} p\mathrm{d} q\\
\leq&\frac{C}{c^2}\Big(\int_{\mathbb{R}^3}\exp\left\{-\frac{|q|^2}{2} \right\}|p-q|^{-2}\langle p\rangle^{2}
\,\mathrm{d} q\Big)^{\frac12}\Big(\int_{\mathbb{R}^3} |h(q)|^2
\,\mathrm{d} q\Big)^{\frac12}\\
&\times\Big(\int_{\mathbb{R}^3}\langle p\rangle^{10}\big(|\tilde{h}(p)|^2+\big|\nabla_{p}\tilde{h}(p)\big|^2\big)
\,\mathrm{d} p\Big)^{\frac12}\Big(\int_{\mathbb{R}^3}\langle p\rangle^{-4}\big(|p|^2||\bar{h}(p)|^2+\big|\nabla_{p}\bar{h}(p)\big|^2\big)
\,\mathrm{d} p\Big)^{\frac12}\\
\leq &\frac{C}{c^2}|h|_{L^2}|\tilde{h}|_{H^1_{p,5}}|h|_{H^1_{p,-2}},
\end{align*}
and
\begin{align*}
|\mathcal{I}_{12}|\leq&\int_{\mathbb{R}^6}\exp\left\{-\frac{|q|^2}{4} \right\}|p-q|^{-1}|\tilde{h}(p)|\Big[\frac{|p_iq_j|}{4}\big(1-\frac{c^2}{p^0q^0}\big)|h(q)||\bar{h}(p)|\\
&+\frac{|q_j|}{2}\big(1-\frac{c}{q^0}\big)|h(q)|\big|\partial_{p_i}\bar{h}(p)\big|
+\frac{|p_i|}{2}\big(1-\frac{c}{p^0}\big)\big|\partial_{q_j}h(q)\big||\bar{h}(p)|\Big]\,\mathrm{d} p\mathrm{d} q\\
\leq&\frac{C}{c^2}\Big(\int_{\mathbb{R}^3}\exp\left\{-\frac{|q|^2}{4} \right\}|p-q|^{-2}\langle p\rangle^2
\,\mathrm{d} q\Big)^{\frac12}\Big(\int_{\mathbb{R}^3} |h(q)|^2+\big|\nabla_{q}h(q)\big|^2
\,\mathrm{d} q\Big)^{\frac12}\\
&\times\Big(\int_{\mathbb{R}^3}\langle p\rangle^{8}|\tilde{h}(p)|^2
\,\mathrm{d} p\Big)^{\frac12}\Big(\int_{\mathbb{R}^3}\langle p\rangle^{-4}\big(|\bar{h}(p)|^2+\big|\nabla_{p}\bar{h}(p)\big|^2\big)
\,\mathrm{d} p\Big)^{\frac12}\\
\leq &\frac{C}{c^2}|h|_{L^2}|\tilde{h}|_{H^1_{p,4}}|h|_{H^1_{p,-2}},
\end{align*}

\end{proof}

\section{Global existence of the RVML system}

\setcounter{equation}{0}
In this section, we are devoted to the global existence of a unique solution to the Cauchy problem of the RVML system \eqref{main1-00}. This goal will be realized through the local-in-time existence, the {\it a priori} energy estimate, and the continuation argument.

\subsection{The local-in-time existence}
In this subsection, we construct a local-in-time existence of a unique solution to the RVML system \eqref{main1-00}.

\begin{proposition}\label{local-exi}
Assume that  
\begin{itemize}
\item $c\geq \max_{0\leq i\leq 3}c_i$ and $\ell_i=\ell_3+2(3-i)\geq8$ for $0\leq i\leq 3$;
\item The initial datum $\big[F^c_{0}(x,p), E^c_{0}(x), B^c_0(x)\big]$ satisfies expansions \eqref{expanion00}, and 
$$F^{3,c}_{\pm,0}(x,p)=\mu^{c}(p)+\sqrt{\mu^{c}(p)}\sum_{i=0}^3\frac{f^{i,c}_{\pm,0}(x,p)}{c^i}\geq0;$$
\item
$[f^{i,c}_{0}(x,p),E^{i,c}_0(x)$, $B^{i,c}_0(x)]$  for $0\leq i\leq 3$ satisfy conservation laws \eqref{cons-RVML0};
\item  For small enough constants $M_i$ with $i=0, 1, 2$, $[f^{i,c}_{0}(x,p),E^{i,c}_0(x)$, $B^{i,c}_0(x)]$ is the unique solution to the Cauchy problem \eqref{mainF0} or \eqref{mainFi} satisfying
    \begin{align*}
\begin{aligned}
\mathcal{E}^{i,c}(t)+\int_0^t\mathcal{D}^{i,c}(\tau) \,\mathrm{d}\tau\leq& C_{\ell_i}\mathcal{E}^{i,c}(0), \qquad \mbox{and}\\
\|f^{i,c}(t)\|^2_{H^5}+\|E^{i,c}(t)\|^2_{H^5}\leq& \mathrm{e}^{-C_{\ell_i}t^{1/3}}\mathcal{E}^{i,c}(0).
\end{aligned}
\end{align*}
\end{itemize}
There exist sufficiently small constants  $M$ with $M>\sum_{i=0}^3\mathcal{E}^{i,c}(0)$ and $\overline{T}>0$, which are independent of $c$,  such that if $\sum_{i=0}^3\mathcal{E}^{i,c}(0)\leq \frac{M}{2} $,
then the Cauchy problem \eqref{main1-00} admits a unique global solution $[F^{c}(t,x,p), E^{c}(t,x), B^{c}(t,x)]$ satisfying conservation laws \eqref{cons-RVML0}, $F^{c}_{\pm}(x,p)=\mu^{c}(p)+\sqrt{\mu^{c}(p)}\sum_{i=0}^3c^{-i}f^{i,c}_{\pm}(t,x,p)\geq0$, and
\begin{equation*}
\sum_{i=0}^3\mathcal{E}^{i,c}(t)+\int_0^t\sum_{i=0}^3\mathcal{D}^{i,c}(\tau) \,\mathrm{d}\tau\leq M,\qquad t\in [0, \overline{T}].
\end{equation*}
\end{proposition}
\begin{proof}
The proof can be proceeded by modifying that in \cite[Theorem 6]{Strain-Guo-CMP-2004} and we omit the details for brevity.
\end{proof}

\subsection{{\it A priori} estimates}
In this subsection, we establish {\it a priori} estimates of the unique solution $[F^{c}(t,x,p), E^{c}(t,x), B^{c}(t,x)]$ to the Cauchy problem of RVML system \eqref{main1-00} constructed in Proposition \ref{local-exi}. Make the following {\it a priori} assumption
\begin{equation}\label{apriori-assu}
\frac{(1+t)^3}{c^6}\mathcal{E}^{3,c}_2(t)+\sum_{i=0}^3\mathcal{E}^{i,c}(t)\leq M,\qquad t\in [0, \overline{T}].
\end{equation}
Based on the {\it a priori} assumption \eqref{apriori-assu}, we will derive  microscopic dissipation estimates, weighted energy estimates,  macroscopic dissipation and electromagnetic dissipation  estimates of the solution $[f^{3,c}(t,x,p)$, $E^{3,c}(t,x), B^{3,c}(t,x)]$, and the time decay of $\mathcal{E}^{3,c}_2(t)$ in a successive way.

\subsubsection{Microscopic dissipation estimates}

In this part, we  derive microscopic dissipation estimates of the unique solution $[f^{3,c}(t,x,p), E^{3,c}(t,x), B^{3,c}(t,x)]$ to the Cauchy problem of RVML system \eqref{mainF3}. These estimates include lower order microscopic dissipation estimates, which enjoy rapid time-decay, and  higher order microscopic dissipation estimates. We will first derive the lower order microscopic dissipation estimates.  
\vskip 0.2cm

\noindent\underline{{\it Step 1. Lower order microscopic dissipation estimates.}} 
Before deriving the estimates for the remainders $[f^{3,c}(t,x,p), E^{3,c}(t,x), B^{3,c}(t,x)]$, we first establish crucial estimates in the following lemmas. The first lemma deals with the momentum growth terms involving the electric field without momentum derivatives.
\begin{lemma}\label{grow-low}
Let  $[F^{c}(t,x,p), E^{c}(t,x), B^{c}(t,x)]$ be the unique solution to the Cauchy problem of RVML system \eqref{main1-00} constructed in Proposition \ref{local-exi}.
Under the assumption \eqref{apriori-assu}, there is a constant $C>0$, which is independent of $c$, such that for $0\leq t\leq \overline{T}$,
\begin{align}\label{growEfj-low}
&\sum_{|\alpha|=0}^2\sum_{\alpha_1\leq\alpha}\Big|\Big\langle \frac{1}{2}\zeta_1\frac{cp}{p^0}\cdot \partial^{\alpha_1}E^{3,c} \sum_{j=0}^2\frac{\partial^{\alpha-\alpha_1}f^{j,c}}{c^j}, \partial^{\alpha}f^{3,c}\Big\rangle\Big|\nonumber\\
&\hspace{1cm}\leq o(1) \big\|f^{3,c}\big\|_{H^2_{\sigma}}^2+CM\min\Big\{\frac{ \mathcal{E}^{3,c}_2(t)}{(1+t)^{9}},\sum_{j=0}^2\mathcal{D}^{j,c}(t)\Big\},
   \end{align}
\begin{align}\label{growEjf-low}
&\sum_{|\alpha|=0}^2\sum_{\alpha_1\leq\alpha}\Big|\Big\langle \frac{1}{2}\zeta_1\frac{cp}{p^0}\cdot \sum_{j=0}^2\frac{\partial^{\alpha_1}E^{j,c} }{c^j} \partial^{\alpha-\alpha_1}f^{3,c}, \partial^{\alpha}f^{3,c}\Big\rangle\Big|\leq o(1) \big\|f^{3,c}\big\|_{H^2_{\sigma}}^2+CM \sum_{j=0}^2 \mathcal{D}^{j,c}_2(t),
   \end{align}   
\begin{align}\label{growEf-low}
&\sum_{|\alpha|=0}^2\sum_{\alpha_1\leq\alpha}\Big|\Big\langle \frac{1}{2}\zeta_1\frac{cp}{p^0}\cdot \sum_{j=0}^2\frac{\partial^{\alpha_1}E^{3,c} }{c^3} \partial^{\alpha-\alpha_1}f^{3,c}, \partial^{\alpha}f^{3,c}\Big\rangle\Big|
    \leq\,o(1) \big\|f^{3,c}\big\|_{H^2_{\sigma}}^2+CM \mathcal{D}^{3,c}_{2}(t),
   \end{align}    
 and
 \begin{align}\label{growEjfj-low}
&\sum_{|\alpha|=0}^2\sum_{\alpha_1\leq\alpha}\Big|\Big\langle \frac{1}{2}\zeta_1\frac{cp}{p^0}\cdot \sum_{\substack{j_1+j_2\geq3\\1\leq j_1,j_2< 3}}\frac{\partial^{\alpha_1}E^{j_1,c} \partial^{\alpha-\alpha_1}f^{j_2,c}}{c^{j_1+j_2-3}}, \partial^{\alpha}f^{3,c}\Big\rangle\Big|\nonumber\\
&\qquad\leq o(1) \big\|f^{3,c}\big\|_{H^2_{\sigma}}^2+CM \sum_{j=0}^2 \mathcal{D}^{j,c}_{2}(t).
   \end{align}   
\end{lemma}
\begin{proof} For brevity, we only prove \eqref{growEfj-low} and \eqref{growEf-low} since \eqref{growEjf-low} and \eqref{growEjfj-low} can be proved similarly. We first show \eqref{growEfj-low}.

Denote the L.H.S. of \eqref{growEfj-low} as $\mathcal{I}_{01}$. Note that $\frac{c}{p^0}\leq1$. 
We apply the Gagliard-Nirenberg inequality in bounded domains, and  use 
 \eqref{decay-f0}, \eqref{decay-fi}, and the {\it a priori} assumption \eqref{apriori-assu} to have
\begin{align*}
\mathcal{I}_{01}\leq& \frac{1}{2}\sum_{|\alpha|=0}^2\sum_{j=0}^2\Big(\big\|E^{3,c}\big\|_{L^{\infty}} \big\|\langle p\rangle^{\frac{3}{2}}\frac{\partial^{\alpha}f^{j,c}}{c^j}\big\|+\sum_{|\alpha_1|=1}\big\|\nabla_xE^{3,c}\big\|_{L^{3}} \big\|\langle p\rangle^{\frac{3}{2}}\frac{\nabla_x^{|\alpha|-1}f^{j,c}}{c^j}\big\|_{L^{6}_xL^2_p} \\
&+\big\|\nabla^2_xE^{3,c}\big\| \big\|\langle p\rangle^{\frac{3}{2}}\frac{f^{j,c}}{c^j}\big\|_{L^{\infty}_xL^2_p}\Big)\big\|\langle p\rangle^{-\frac{1}{2}}\partial^{\alpha}f^{3,c}\big\|\\
\leq&C \sum_{j=0}^2\sum_{|\alpha|=0}^2\big\|E^{3,c}\big\|_{H^2} \big\|\langle p\rangle^{\frac{3}{2}}\frac{f^{j,c}}{c^j}\big\|_{H^2} \big\|\partial^{\alpha}f^{3,c}\big\|_{\sigma} \\
\leq& o(1) \big\|f^{3,c}\big\|_{H^2_{\sigma}}^2+ C\sum_{j=0}^2\big\|\langle p\rangle^{\frac{3}{2}}\frac{f^{j,c}}{c^j}\big\|^2_{H^2}
\big\|E^{3,c}\big\|^2_{H^2}
\\
\leq& o(1) \big\|f^{3,c}\big\|_{H^2_{\sigma}}^2+CM \sum_{j=0}^2\min\left\{ \mathrm{e}^{-C_{\ell_j}(1+t)^{1/3}}\mathcal{E}^{3,c}(t), \mathcal{D}^{j,c}(t)\right\}\\
\leq & o(1) \big\|f^{3,c}\big\|_{H^2_{\sigma}}^2+CM\min\Big\{\frac{ \mathcal{E}^{3,c}_2(t)}{(1+t)^{9}},\sum_{j=0}^2\mathcal{D}^{j,c}(t)\Big\}.
\end{align*}

Now we turn to prove \eqref{growEf-low}. Denote the L.H.S. of \eqref{growEf-low} as $\mathcal{I}_{02}$. We also apply the Gagliard-Nirenberg inequality in bounded domains, and  use the definition of $\mathcal{D}^{3,c}_{m}(t)$ in \eqref{nowei-ed3} and the {\it a priori} assumption \eqref{apriori-assu} to have
\begin{align*}
\mathcal{I}_{02}\leq& \frac{1}{2c^3}\sum_{|\alpha|=0}^2\Big(\big\|E^{3,c}\big\|_{L^{\infty}} \big\|\langle p\rangle^{\frac{3}{2}}\partial^{\alpha}f^{3,c}\big\|+\sum_{|\alpha_1|=1}\big\|\nabla_xE^{3,c}\big\|_{L^{3}} \big\|\langle p\rangle^{\frac{3}{2}}\nabla_x^{|\alpha|-1}f^{3,c}\big\|_{L^{6}_xL^2_p} \\
&+\big\|\nabla^2_xE^{3,c}\big\| \big\|\langle p\rangle^{\frac{3}{2}}f^{3,c}\big\|_{L^{\infty}_xL^2_p}\Big)\big\|\langle p\rangle^{-\frac{1}{2}}\partial^{\alpha}f^{3,c}\big\|\\
\leq&\frac{C}{c^3} \sum_{|\alpha|=0}^2 \big\|E^{3,c}\big\|_{H^1}^{\frac{1}{2}}\big\|E^{3,c}\big\|_{H^3}^{\frac{1}{2}} \big\|\langle p\rangle^{\frac{7}{2}}f^{3,c}\big\|^{\frac{1}{2}}_{H^2}\big\|\langle p\rangle^{-\frac{1}{2}}f^{3,c}\big\|_{H^2}^{\frac{1}{2}} \|\partial^{\alpha}f^{3,c}\|_{\sigma} \\
\leq& o(1) \big\|f^{3,c}\big\|_{H^2_{\sigma}}^2+ \frac{C}{c^6}\Big(\big\|E^{3,c}\big\|_{H^3}^2+\big\|\langle p\rangle^{\frac{7}{2}}f^{3,c}\big\|^2_{H^2}\Big)\Big(\big\|E^{3,c}\big\|_{H^1}^2+\|f^{3,c}\big\|^2_{H^2_{\sigma}}\Big)\\
\leq& o(1) \big\|f^{3,c}\big\|_{H^2_{\sigma}}^2+CM \mathcal{D}^{3,c}_{2}(t).
\end{align*}

\end{proof}   

The second lemma is to control the momentum derivative  terms associated with the Lorentz force.
\begin{lemma}\label{deriv-low}
Under the assumptions in Lemma \ref{grow-low}, there is a constant $C>0$, which is independent of $c$, such that for $0\leq t\leq \overline{T}$,
\begin{align}\label{derivEfj-low}
&\sum_{|\alpha|=0}^2\sum_{\alpha_1\leq\alpha}\Big|\Big\langle \zeta_1\partial^{\alpha_1}E^{3,c}\cdot\sum_{j=0}^2\frac{\partial^{\alpha-\alpha_1}\nabla_pf^{j,c}}{c^j}, \partial^{\alpha}f^{3,c}\Big\rangle\Big|\\
&\qquad\leq o(1) \big\|f^{3,c}\big\|_{H^2_{\sigma}}^2+CM\min\Big\{\frac{ \mathcal{E}^{3,c}_2(t)}{(1+t)^{9}},\sum_{j=0}^2\mathcal{D}^{j,c}(t)\Big\},
   \end{align}
\begin{align*}
&\sum_{|\alpha|=0}^2\sum_{\alpha_1\leq\alpha}\Big|\Big\langle \zeta_1\sum_{j=0}^2\frac{\partial^{\alpha_1}E^{j,c} }{c^j} \cdot\partial^{\alpha-\alpha_1}\nabla_pf^{3,c}, \partial^{\alpha}f^{3,c}\Big\rangle\Big|\leq o(1) \big\|f^{3,c}\big\|_{H^2_{\sigma}}^2+CM \sum_{j=0}^2 \mathcal{D}^{j,c}_{2}(t),
   \end{align*}   

\begin{align}\label{derivEf-low}
&\sum_{|\alpha|=0}^2\sum_{\alpha_1\leq\alpha}\Big|\Big\langle \zeta_1\partial^{\alpha_1}\big(E^{3,c}+\frac{p}{p^0}\times B^{3,c}\big)\cdot\frac{\partial^{\alpha-\alpha_1}\nabla_pf^{3,c}}{c^3} , \partial^{\alpha}f^{3,c}\Big\rangle\Big|\nonumber\\
 &\qquad   \leq\,o(1) \big\|f^{3,c}\big\|_{H^2_{\sigma}}^2+CM \mathcal{D}^{3,c}_{2}(t),
   \end{align}    
 and
 \begin{align*}
&\sum_{|\alpha|=0}^2\sum_{\alpha_1\leq\alpha}\Big|\Big\langle \zeta_1 \sum_{\substack{j_1+j_2\geq3\\1\leq j_1,j_2< 3}}\frac{\partial^{\alpha_1}E^{j_1,c} \cdot \partial^{\alpha-\alpha_1}\nabla_pf^{j_2,c}}{c^{j_1+j_2-3}}, \partial^{\alpha}f^{3,c}\Big\rangle\Big|\nonumber\\
 &\qquad    \leq\,o(1) \big\|f^{3,c}\big\|_{H^2_{\sigma}}^2+CM \sum_{j=0}^2 \mathcal{D}^{j,c}_{2}(t).
   \end{align*}   
\end{lemma}
\begin{proof} As in Lemma \ref{grow-low}, we only prove \eqref{derivEfj-low} and \eqref{derivEf-low}. We first show \eqref{derivEfj-low}.

Denote the L.H.S. of \eqref{derivEfj-low} as $\mathcal{I}_{03}$.  
We integrate by part w.r.t. $p$, apply the Gagliard-Nirenberg inequality in bounded domains, and  use 
 \eqref{decay-f0}, \eqref{decay-fi}, and the {\it a priori} assumption \eqref{apriori-assu} to have
\begin{align*}
\mathcal{I}_{03}\leq& \sum_{|\alpha|=0}^2\sum_{j=0}^2\Big(\Big\|E^{3,c}\big\|_{L^{\infty}} \big\|\langle p\rangle^{\frac{3}{2}}\frac{\partial^{\alpha}f^{j,c}}{c^j}\big\|+\sum_{|\alpha_1|=1}\|\nabla_xE^{3,c}\|_{L^{3}} \big\|\langle p\rangle^{\frac{3}{2}}\frac{\nabla_x^{|\alpha|-1}f^{j,c}}{c^j}\big\|_{L^{6}_xL^2_p} \\
&+\|\nabla^2_xE^{3,c}\| \big\|\langle p\rangle^{\frac{3}{2}}\frac{f^{j,c}}{c^j}\big\|_{L^{\infty}_xL^2_p}\Big)\big\|\langle p\rangle^{-\frac{3}{2}}\partial^{\alpha}\nabla_pf^{3,c}\big\|\\
\leq&C \sum_{j=0}^2\sum_{|\alpha|=0}^2\big\|E^{3,c}\big\|_{H^2} \big\|\langle p\rangle^{\frac{3}{2}}\frac{f^{j,c}}{c^j}\big\|_{H^2} \big\|\partial^{\alpha}f^{3,c}\big\|_{\sigma} \\
\leq& o(1) \big\|f^{3,c}\big\|_{H^2_{\sigma}}^2+CM\min\Big\{\frac{ \mathcal{E}^{3,c}_2(t)}{(1+t)^{9}},\sum_{j=0}^2\mathcal{D}^{j,c}(t)\Big\}.
\end{align*}

Now we turn to prove \eqref{derivEf-low}. Denote the L.H.S. of \eqref{derivEf-low} as $\mathcal{I}_{04}$. Note that $\frac{|p|}{p^0}\leq1$. We apply integration by parts w.r.t. $p$, and use similar arguments as the derivation of  \eqref{growEf-low} to have
\begin{align*}
\mathcal{I}_{04}\leq& \frac{1}{c^3}\sum_{|\alpha|=1}^2\Big[\big(\|\nabla_xE^{3,c}\|_{L^{3}} +\|\nabla_xB^{3,c}\|_{L^{3}} \big)\big\|\langle p\rangle^{\frac{3}{2}}\nabla_x^{|\alpha|-1}f^{3,c}\big\|_{L^{6}_xL^2_p} \\
&+\big(\|\nabla^2_xE^{3,c}\|+\|\nabla^2_xB^{3,c}\|\big) \big\|\langle p\rangle^{\frac{3}{2}}f^{3,c}\big\|_{L^{\infty}_xL^2_p}\Big]\big\|\langle p\rangle^{-\frac{3}{2}}\partial^{\alpha}\nabla_pf^{3,c}\big\|\\
\leq&\frac{C}{c^3} \sum_{|\alpha|=0}^2 \big(\Big\|E^{3,c}\big\|_{H^1}+\|B^{3,c}-\bar{B}^{3}\|_{H^1}\big)^{\frac{1}{2}}\big(\Big\|E^{3,c}\big\|_{H^3}
+\|B^{3,c}-\bar{B}^{3}\|_{H^3}\big)^{\frac{1}{2}} \\
&\times\big\|\langle p\rangle^{\frac{7}{2}}f^{3,c}\big\|^{\frac{1}{2}}_{H^2}\big\|\langle p\rangle^{-\frac{1}{2}}f^{3,c}\big\|_{H^2}^{\frac{1}{2}} \|\partial^{\alpha}f^{3,c}\|_{\sigma} \\
\leq& o(1) \big\|f^{3,c}\big\|_{H^2_{\sigma}}^2+ \frac{C}{c^6}\Big(\Big\|E^{3,c}\big\|_{H^3}^2+\|B^{3,c}-\bar{B}^{3}\|_{H^3}^2+\big\|\langle p\rangle^{\frac{7}{2}}f^{3,c}\big\|^2_{H^2}\Big)\\
&\times\Big(\Big\|E^{3,c}\big\|_{H^1}^2+\|B^{3,c}-\bar{B}^{3}\|_{H^1}^2+\|f^{3,c}\big\|^2_{H^2_{\sigma}}\Big)\\
\leq& o(1) \big\|f^{3,c}\big\|_{H^2_{\sigma}}^2+CM \mathcal{D}^{3,c}_{2}(t).
\end{align*}

\end{proof}  
 
Based on the above estimates, we can obtain the following lower order estimates of $[f^{3,c}(t, x, p)$, $E^ {3,c}(t, x)$, $B^{3,c}(t, x)]$.
\begin{proposition}\label{2norm00} Under the assumptions in Lemma \ref{grow-low}, there are $c-$independent constants $\delta_0>0$ and $C>0$ such that for $0\leq t\leq \overline{T}$,
\begin{align} \label{2norm}
&\frac{\mathrm{d}}{\mathrm{d}t}\mathcal{E}^{3,c}_{2}(t)+
\delta_0 \|\{I-\mathcal{P}\}f^{3,c}\|^2_{H^2_{\sigma}}\nonumber\\
&\hspace{1cm}\leq \big(o(1)+\sum_{j=0}^3\frac{C\sqrt{M}}{c^j}\big)\big\|f^{3,c}\big\|_{H^2_{\sigma}}^2+\frac{CM \mathcal{E}^{3,c}_2(t)}{(1+t)^{9}}
+CM\sum_{j=0}^3 \mathcal{D}^{j,c}_{2}(t).
\end{align}
\end{proposition}
\begin{proof} Taking $\vartheta=\ell=0$ in \eqref{nonlin-w}, we apply Gagliard-Nirenberg inequality in bounded domains, and  use the {\it a priori} assumption \eqref{apriori-assu} to have 
\begin{align}\label{nonlin-2j3}
&\sum_{|\alpha|=0}^2\sum_{\alpha_1\leq\alpha}\Big|\Big\langle \sum_{j=0}^2\frac{1}{c^j}\left[\varGamma\left(\partial^{\alpha_1}f^{j,c}, \partial^{\alpha-\alpha_1}f^{3,c}\right)+\varGamma\left(\partial^{\alpha_1}f^{3,c}, \partial^{\alpha-\alpha_1}f^{j,c}\right)\right]\nonumber\\
&\hspace{1cm}+\frac{1}{c^3}\varGamma\left(\partial^{\alpha_1}f^{3,c}, \partial^{\alpha-\alpha_1}f^{3,c}\right), \partial^{\alpha}f^{3,c}\Big\rangle\Big|\nonumber\\
&\qquad\leq C\Big[\sum_{j=0}^2
\big(\frac{\|f^{j,c}\|_{H^2}}{c^j}\big\|f^{3,c}\big\|_{H^2_{\sigma}}+\frac{\|f^{j,c}\|_{H^2_{\sigma}}}{c^j}\big\|f^{3,c}\big\|_{H^2}\big)
+\frac{\big\|f^{3,c}\big\|_{H^2}}{c^3}\big\|f^{3,c}\big\|_{H^2_{\sigma}}\Big]\big\|f^{3,c}\big\|_{H^2_{\sigma}}\\
&\qquad\leq \Big(o(1)+\sum_{j=0}^3\frac{C\sqrt{M}}{c^{j}}\Big)\big\|f^{3,c}\big\|_{H^2_{\sigma}}^2
+\sum_{j=0}^2\frac{CM}{c^{2j}}\|f^{j,c}\|_{H^2_{\sigma}}^2,
\nonumber
\end{align}
and
\begin{align}\label{nonlin-2jj}
&\sum_{|\alpha|=0}^2\sum_{\alpha_1\leq\alpha}\Big|\Big\langle \sum_{\substack{j_1+j_2\geq3   \\1\leq j_1,j_2<3}}\frac{\varGamma\left(\partial^{\alpha_1}f^{j_1,c}, \partial^{\alpha-\alpha_1}f^{j_2,c}\right)}{c^{j_1+j_2-3}}, \partial^{\alpha}f^{3,c}\Big\rangle\Big|\nonumber\\
&\qquad\leq C\sum_{\substack{j_1+j_2\geq3   \\1\leq j_1,j_2<3}}\sum_{j=0}^2
\frac{\|f^{j_1,c}\|_{H^2}\|f^{j_2,c}\|_{H^2_{\sigma}}
+\|f^{j_1,c}\|_{H^2_{\sigma}}\|f^{j_2,c}\|_{H^2}}{c^{j_1+j_2-3}}\big\|f^{3,c}\big\|_{H^2_{\sigma}}\\
&\qquad\leq o(1)\big\|f^{3,c}\big\|_{H^2_{\sigma}}^2+\sum_{j=0}^2CM\|f^{j,c}\|_{H^2_{\sigma}}^2.
\nonumber
\end{align}
Now we apply $\partial^{\alpha}$ with $|\alpha|\leq2$ to the equation of $f^{3,c}(t, x, p)$ in \eqref{mainF3}, take the inner product with $\partial^{\alpha}f^{3,c}(t, x, p)$ in $L^2_xL^2_p$, and use Lemma \ref{coerci}, Lemma \ref{grow-low}, Lemma \ref{deriv-low}, \eqref{nonlin-2j3} and \eqref{nonlin-2jj} to have
\begin{align} \label{2normf3}
&\frac{\mathrm{d}}{\mathrm{d}t}\|f^{3,c}(t)\|^2_{H^2}+
\delta_0 \|\{I-\mathcal{P}\}f^{3,c}\|^2_{H^2_{\sigma}}+\sum_{|\alpha|=0}^2\Big\langle \partial^{\alpha}E^{3,c}\cdot\frac{cp}{p^0}\sqrt{\mu^c}\zeta_0,\partial^{\alpha}f^{3,c}\Big\rangle\nonumber\\
&\qquad\leq 
\big(o(1)+\sum_{j=0}^3\frac{CM}{c^j}\big)\big\|f^{3,c}\big\|_{H^2_{\sigma}}^2+\frac{CM \mathcal{E}^{3,c}_2(t)}{(1+t)^{9}}
+CM\sum_{j=0}^2\|f^{j,c}\|_{H^2_{\sigma}}^2.
\end{align}
On the other hand, from the Maxwell system in \eqref{mainF3}, we can obtain
\begin{align} \label{2normEB}
\frac{\mathrm{d}}{\mathrm{d}t}\left(\|E^{3,c}(t)\|^2_{H^2}+\|B^{3,c}(t)-\bar{B}^{3}\|^2_{H^2}\right)=\sum_{|\alpha|=0}^2\Big\langle \partial^{\alpha}E^{3,c}\cdot\frac{cp}{p^0}\sqrt{\mu^c}\zeta_0,\partial^{\alpha}f^{3,c}\Big\rangle.
\end{align}
Then we combine \eqref{2normf3} and \eqref{2normEB} to obtain \eqref{2norm}.
\end{proof}
\vskip 0.2cm

\noindent\underline{{\it Step 2. Higher order microscopic dissipation estimates.}} 
As in Step $1$, we first establish crucial estimates in the following two lemmas, which correspond to Lemma \ref{grow-low} and \ref{deriv-low}. 
\begin{lemma}\label{grow-high}
Under the assumptions in Lemma \ref{grow-low}, there is a constant $C>0$, which is independent of $c$, such that for $0\leq t\leq \overline{T}$,
\begin{align}\label{growEfj-high}
&\sum_{|\alpha|=3}^5\sum_{\alpha_1\leq\alpha}\Big|\Big\langle \frac{1}{2}\zeta_1\frac{cp}{p^0}\cdot \partial^{\alpha_1}E^{3,c} \sum_{j=0}^2\frac{\partial^{\alpha-\alpha_1}f^{j,c}}{c^j}, \partial^{\alpha}f^{3,c}\Big\rangle\Big|
    \leq\,o(1) \big\|f^{3,c}\big\|_{H^5_{\sigma}}^2+CM \sum_{j=0}^2 \frac{\mathcal{D}^{j,c}(t)}{c^{2j}},
   \end{align}
\begin{align*}
&\sum_{|\alpha|=3}^5\sum_{\alpha_1\leq\alpha}\Big|\Big\langle \frac{1}{2}\zeta_1\frac{cp}{p^0}\cdot \sum_{j=0}^2\frac{\partial^{\alpha_1}E^{j,c} }{c^j} \partial^{\alpha-\alpha_1}f^{3,c}, \partial^{\alpha}f^{3,c}\Big\rangle\Big|
    \leq\,o(1) \big\|f^{3,c}\big\|_{H^5_{\sigma}}^2+CM  \mathcal{D}^{3,c}(t),
   \end{align*}   
\begin{align}\label{growEf-high}
&\sum_{|\alpha|=3}^5\sum_{\alpha_1\leq\alpha}\Big|\Big\langle \frac{1}{2}\zeta_1\frac{cp}{p^0}\cdot \sum_{j=0}^2\frac{\partial^{\alpha_1}E^{3,c} }{c^3} \partial^{\alpha-\alpha_1}f^{3,c}, \partial^{\alpha}f^{3,c}\Big\rangle\Big| \leq o(1) \big\|f^{3,c}\big\|_{H^5_{\sigma}}^2+CM \mathcal{D}^{3,c}(t),
   \end{align}    
 and
 \begin{align*}
&\sum_{|\alpha|=3}^5\sum_{\alpha_1\leq\alpha}\Big|\Big\langle \frac{1}{2}\zeta_1\frac{cp}{p^0}\cdot \sum_{\substack{j_1+j_2\geq3\\1\leq j_1,j_2< 3}}\frac{\partial^{\alpha_1}E^{j_1,c} \partial^{\alpha-\alpha_1}f^{j_2,c}}{c^{j_1+j_2-3}}, \partial^{\alpha}f^{3,c}\Big\rangle\Big|\nonumber\\
 &\qquad    \leq\,o(1) \big\|f^{3,c}\big\|_{H^5_{\sigma}}^2+CM \sum_{j=0}^2 \mathcal{D}^{j,c}(t).
   \end{align*}   
\end{lemma}
\begin{proof} As in Lemma \ref{grow-low}, we only prove \eqref{growEfj-high} and \eqref{growEf-high}. We first show \eqref{growEfj-high}.

As the estimation of $\mathcal{I}_{01}$, we take $L^{\infty}-L^2-L^2$ for $|\alpha_1|=0$, $L^{3}-L^6-L^2$ for $|\alpha_1|=1$, and $L^2-L^{\infty}-L^2$ for $|\alpha_1|\geq2$, and bound the L.H.S. of \eqref{growEfj-high} by
\begin{align*}
C \sum_{|\alpha|=3}^5\sum_{j=0}^2 \Big\|E^{3,c}\big\|_{H^5} \big\|\langle p\rangle^{\frac{3}{2}}\frac{f^{j,c}}{c^j}\big\|_{H^5} \|\partial^{\alpha}f^{3,c}\|_{\sigma} 
\leq& o(1) \big\|f^{3,c}\big\|_{H^5_{\sigma}}^2+CM \sum_{j=0}^2 \frac{\mathcal{D}^{j,c}(t)}{c^{2j}}.
\end{align*}

For \eqref{growEf-high}, we denote its L.H.S. of \eqref{growEf-high} as $\mathcal{I}_{05}$. Then we take $L^{\infty}-L^2-L^2$ for $|\alpha_1|=0$, $L^{3}-L^6-L^2$ for $|\alpha_1|=1$, and $L^2-L^{\infty}-L^2$ for $|\alpha_1|\geq2$, and use the definition of $\mathcal{D}^{3,c}(t)$ in \eqref{dissirat-3c} and the {\it a priori} assumption \eqref{apriori-assu} to have
\begin{align*}
\mathcal{I}_{05}\leq& \frac{C}{c^3}\sum_{|\alpha|=3}^5\Big( \Big\|E^{3,c}\big\|_{H^2}\big\|\langle p\rangle^{\frac{3}{2}}f^{3,c}\big\|_{H^5}+\sum_{|\alpha_1|\geq3}\|\nabla_x^{|\alpha_1|}E^{3,c}\big\|\langle p\rangle^{\frac{3}{2}}\partial^{\alpha-\alpha_1}f^{3,c}\big\|_{H^2}\Big)\big\|\langle p\rangle^{-\frac{1}{2}}\partial^{\alpha}f^{3,c}\big\|\\
\leq& o(1) \big\|f^{3,c}\big\|_{H^5_{\sigma}}^2+ \frac{C}{c^6}\Big(\frac{\Big\|E^{3,c}\big\|^2_{H^2}\mathcal{D}^{3,c}(t)}{(1+t)^{(1+\epsilon_{0})/2}}
+\Big\|E^{3,c}\big\|_{H^5}^2\mathcal{D}^{3,c}_{4,w}(t)\Big)\\
\leq& o(1) \big\|f^{3,c}\big\|_{H^2_{\sigma}}^2+CM \mathcal{D}^{3,c}(t).
\end{align*}

\end{proof}   

The second lemma is to control the momentum derivative  terms associated with the Lorentz force.
\begin{lemma}\label{deriv-high}
Under the assumptions in Lemma \ref{grow-low}, there is a constant $C>0$, which is independent of $c$, such that for $0\leq t\leq \overline{T}$,
\begin{align}\label{derivEfj-high}
&\sum_{|\alpha|=3}^5\sum_{\alpha_1\leq\alpha}\Big|\Big\langle \zeta_1\partial^{\alpha_1}E^{3,c}\cdot\sum_{j=0}^2\frac{\partial^{\alpha-\alpha_1}\nabla_pf^{j,c}}{c^j}, \partial^{\alpha}f^{3,c}\Big\rangle\Big|
    \leq\,o(1) \big\|f^{3,c}\big\|_{H^5_{\sigma}}^2+CM \sum_{j=0}^2 \frac{\mathcal{D}^{j,c}(t)}{c^{2j}},
   \end{align}
\begin{align*}
&\sum_{|\alpha|=3}^5\sum_{\alpha_1\leq\alpha}\Big|\Big\langle \zeta_1\sum_{j=0}^2\frac{\partial^{\alpha_1}E^{j,c} }{c^j} \cdot\partial^{\alpha-\alpha_1}\nabla_pf^{3,c}, \partial^{\alpha}f^{3,c}\Big\rangle\Big|
    \leq\,o(1) \big\|f^{3,c}\big\|_{H^5_{\sigma}}^2+CM \mathcal{D}^{3,c}(t),
   \end{align*}   

\begin{align}\label{derivEf-high}
&\sum_{|\alpha|=3}^5\sum_{\alpha_1\leq\alpha}\Big|\Big\langle \zeta_1\partial^{\alpha_1}\big(E^{3,c}+\frac{p}{p^0}\times B^{3,c}\big)\cdot\frac{\partial^{\alpha-\alpha_1}\nabla_pf^{3,c}}{c^3} , \partial^{\alpha}f^{3,c}\Big\rangle\Big|\nonumber\\
&\qquad\leq o(1) \big\|f^{3,c}\big\|_{H^5_{\sigma}}^2+CM \mathcal{D}^{3,c}(t),
   \end{align}    
 and
 \begin{align*}
&\sum_{|\alpha|=3}^5\sum_{\alpha_1\leq\alpha}\Big|\Big\langle \zeta_1 \sum_{\substack{j_1+j_2\geq3\\1\leq j_1,j_2< 3}}\frac{\partial^{\alpha_1}E^{j_1,c} \cdot \partial^{\alpha-\alpha_1}\nabla_pf^{j_2,c}}{c^{j_1+j_2-3}}, \partial^{\alpha}f^{3,c}\Big\rangle\Big|\nonumber\\
   &\qquad \leq\,o(1) \big\|f^{3,c}\big\|_{H^5_{\sigma}}^2+CM \sum_{j=0}^2 \mathcal{D}^{j,c}(t).
   \end{align*}   
\end{lemma}
\begin{proof} As before, we only prove \eqref{derivEfj-high} and \eqref{derivEf-high}.

For the L.H.S. of \eqref{derivEfj-high}, we take $L^{\infty}-L^2-L^2$ for $|\alpha_1|=0$, $L^{3}-L^6-L^2$ for $|\alpha_1|=1$, and $L^2-L^{\infty}-L^2$ for $|\alpha_1|\geq2$, and use the definition of $\mathcal{D}^{i,c}(t)$ in \eqref{dissirat-ic} with $i=0, 1, 2, 3$ to obtain its upper bound
\begin{align*}
C \sum_{|\alpha|=3}^5\sum_{j=0}^2 \Big\|E^{3,c}\big\|_{H^5} \big\|\langle p\rangle^{\frac{1}{2}}\frac{\nabla_pf^{j,c}}{c^j}\big\|_{H^5} \|\partial^{\alpha}f^{3,c}\|_{\sigma} 
\leq o(1) \big\|f^{3,c}\big\|_{H^5_{\sigma}}^2+CM \sum_{j=0}^2 \frac{\mathcal{D}^{j,c}(t)}{c^{2j}}.
\end{align*}

Now we turn to prove \eqref{derivEf-high}. Denote the L.H.S. of \eqref{derivEf-high} as $\mathcal{I}_{06}$. We first apply integration by parts w.r.t. $p$, and use similar arguments as the derivation of  \eqref{growEf-high} to have
\begin{align*}
\mathcal{I}_{06}\leq& \frac{C}{c^3}\sum_{|\alpha|=3}^5\Big[ \big(\Big\|E^{3,c}\big\|_{H^2}+\|B^{3,c}-\bar{B}^{3}\|_{H^2}\big)\big\|\langle p\rangle^{\frac{3}{2}}f^{3,c}\big\|_{H^5}\\
&+\sum_{|\alpha_1|\geq3}\big(\|\nabla_x^{|\alpha_1|}E^{3,c}+\|\nabla_x^{|\alpha_1|}B^{3,c}\|\big)\big\|\langle p\rangle^{\frac{3}{2}}\partial^{\alpha-\alpha_1}f^{3,c}\big\|_{H^2}\Big]\big\|\langle p\rangle^{-\frac{3}{2}}\partial^{\alpha}\nabla_pf^{3,c}\big\|\\
\leq& o(1) \big\|f^{3,c}\big\|_{H^5_{\sigma}}^2+CM \mathcal{D}^{3,c}(t).
\end{align*}

\end{proof}   

\begin{proposition}\label{6norm00} Under the assumptions in Lemma \ref{grow-low}, there are $c-$independent constants $\delta_0>0$ and $C>0$ such that for $0\leq t\leq \overline{T}$,
\begin{align} \label{6norm}
&\frac{\mathrm{d}}{\mathrm{d}t}\left(\|f^{3,c}(t)\|^2_{H^5}+\|E^{3,c}(t)\|^2_{H^5}
+\|B^{3,c}(t)-\bar{B}^{3}\|^2_{H^5}\right)+
\delta_0 \|\{I-\mathcal{P}\}f^{3,c}\|^2_{H^5_{\sigma}}\nonumber\\
&\qquad\leq o(1) \big\|f^{3,c}\big\|_{H^5_{\sigma}}^2
+CM\sum_{j=0}^3\mathcal{D}^{j,c}(t).
\end{align}
\end{proposition}
\begin{proof} As the estimation of \eqref{nonlin-2j3} and \eqref{nonlin-2jj}, we can obtain 
\begin{align}\label{nonlin-6j3}
&\sum_{|\alpha|=0}^5\sum_{\alpha_1\leq\alpha}\Big|\Big\langle \sum_{j=0}^2\frac{1}{c^j}\left[\varGamma\left(\partial^{\alpha_1}f^{j,c}, \partial^{\alpha-\alpha_1}f^{3,c}\right)+\varGamma\left(\partial^{\alpha_1}f^{3,c}, \partial^{\alpha-\alpha_1}f^{j,c}\right)\right]\nonumber\\
&\hspace{1cm}+\frac{1}{c^3}\varGamma\left(\partial^{\alpha_1}f^{3,c}, \partial^{\alpha-\alpha_1}f^{3,c}\right), \partial^{\alpha}f^{3,c}\Big\rangle\Big|\nonumber\\
&\qquad\leq \Big(o(1)+\sum_{j=0}^3\frac{C\sqrt{M}}{c^{j}}\Big)\big\|f^{3,c}\big\|_{H^5_{\sigma}}^2
+\sum_{j=0}^2\frac{CM}{c^{2j}}\|f^{j,c}\|_{H^5_{\sigma}}^2,
\end{align}
and
\begin{align}\label{nonlin-6jj}
&\sum_{|\alpha|=0}^5\sum_{\alpha_1\leq\alpha}\Big|\Big\langle \sum_{\substack{j_1+j_2\geq3   \\1\leq j_1,j_2<3}}\frac{\varGamma\left(\partial^{\alpha_1}f^{j_1,c}, \partial^{\alpha-\alpha_1}f^{j_2,c}\right)}{c^{j_1+j_2-3}}, \partial^{\alpha}f^{3,c}\Big\rangle\Big|\nonumber\\
&\qquad\leq o(1)\big\|f^{3,c}\big\|_{H^5_{\sigma}}^2+\sum_{j=0}^2CM\|f^{j,c}\|_{H^5_{\sigma}}^2.
\end{align}
Now we apply $\partial^{\alpha}$ with $|\alpha|\leq5$ to the equations of $[f^{3,c}(t, x, p), E^{3,c}(t, x), B^{3,c}(t, x)]$ in \eqref{mainF3}, take the corresponding inner products  in $H^{5}$, and use Lemma \ref{coerci}, Lemma \ref{grow-low}, Lemma \ref{deriv-low}, Lemma \ref{grow-high}, Lemma \ref{deriv-high}, \eqref{nonlin-6j3} and \eqref{nonlin-6jj} to deduce \eqref{6norm}.
\end{proof}

\subsubsection{Weighted energy estimates}
In this part, we  derive weighted energy estimates of the unique solution $[f^{3,c}(t,x,p), E^{3,c}(t,x), B^{3,c}(t,x)]$ to the Cauchy problem of RVML system \eqref{mainF3}. These estimates include non-highest order weighted energy estimates of $f^{3,c}(t,x,p)$, whose weighted norms don't increase with time, and the highest order weighted energy estimates of $f^{3,c}(t,x,p)$, whose weighted norms increase with time. We will first derive the lower order microscopic dissipation estimates.  
\vskip 0.2cm

\noindent\underline{{\it Step 1. Non-highest order weighted energy estimates.}} 
 In this step, as preparations for the weighted energy estimates of $f^{3,c}(t,x,p)$,  we first treat the weighted energy estimate of the linear term $- E^{3,c}\cdot\frac{cp}{p^0}\sqrt{\mu^c}\zeta_0$ which is difficult since $\partial^{\alpha}E^{3,c}(t, x)$ with $|\alpha|\leq 5$ is not dissipative due to the absence of coefficient $\frac{1}{c}$. Then we will deal with the weighted energy estimates of nonlinear terms involving the electromagnetic field.

In the following lemma, by similar arguments as in \cite{Jiang-Lei-Zhao-JFA-2024}, we  proceed the weighted energy estimate of the linear term $- E^{3,c}\cdot\frac{cp}{p^0}\sqrt{\mu^c}\zeta_0$.

\begin{lemma}\label{lin-wei}
Under the assumptions in Lemma \ref{grow-low}, there is a constant $C>0$, which is independent of $c$, such that for $0\leq t\leq \overline{T}$,
\begin{align}\label{lin-wei0}
&\sum_{|\alpha|=0}^4\Big\langle -\partial^{\alpha}E^{3,c}\cdot\frac{cp}{p^0}\sqrt{\mu^c}\zeta_0, w_{\ell_3-|\alpha|}^2\partial^{\alpha}f^{3,c}\Big\rangle\nonumber\\
&\qquad   \leq\sum_{|\alpha|=0}^4\frac{-1}{4C_{\alpha}(t)}\frac{\mathrm{d}}{\mathrm{d}t}\Big\langle \Big(\frac{cp}{p^0}\sqrt{\mu^c}\zeta_0, w_{\ell_3-|\alpha|}^2\partial^{\alpha}f^{3,c}\Big),\Big(\frac{cp}{p^0}\sqrt{\mu^c}\zeta_0, w_{\ell_3-|\alpha|}^2\partial^{\alpha}f^{3,c}\Big)\Big\rangle\nonumber\\
&\qquad    +  C\big\|f^{3,c}\big\|_{H^5_{\sigma}}^2+CM \sum_{j=0}^2\mathcal{D}^{j,c}(t).
   \end{align}

\end{lemma}
\begin{proof} Denote the L.H.S. of \eqref{lin-wei0} as $\mathcal{I}_{07}$. To estimate $\mathcal{I}_{07}$, we first obtain an expression of $\partial^{\alpha}E^{3,c}(t, x)$ with $|\alpha|\leq 4$. 
 For $i=1, 2, 3$, we apply $\partial^{\alpha}$ to the equation of $f^{3,c}(t,x,p)$ in \eqref{mainF3} and multiply the resultant by $w_{\ell_3-|\alpha|}^2(p)\frac{cp_i}{p^0}\sqrt{\mu^c}$ to have
 \begin{align*}
    & \Big( \partial_t\big[\partial^{\alpha}f^{3,c}_{+}-\partial^{\alpha}f^{3,c}_{-}\big],w_{\ell_3-|\alpha|}^2\frac{cp_i}{p^0}\sqrt{\mu^c}\Big)
    + \Big( \frac{cp}{p^0}\cdot \big[\partial^{\alpha}\nabla_xf^{3,c}_{+}-\partial^{\alpha}\nabla_xf^{3,c}_{-}\big],w_{\ell_3-|\alpha|}^2
    \frac{cp_i}{p^0}\sqrt{\mu^c}\Big)\nonumber\\
    &\qquad=2\partial^{\alpha}E^{3,c}_i\Big( \frac{cp_i}{p^0}\sqrt{\mu^c},w_{\ell_3-|\alpha|}^2\frac{cp_i}{p^0}\sqrt{\mu^c}\Big)+\Big( \partial^{\alpha}\mathcal{I}_{R,+}-\partial^{\alpha}\mathcal{I}_{R,-},w_{\ell_3-|\alpha|}^2\frac{cp_i}{p^0}\sqrt{\mu^c}\Big),
 \end{align*}
 where $\mathcal{I}_{R}$ is used to denote the remaining terms of the equation of $f^{3,c}(t, x, p)$ except for the first three terms in \eqref{mainF3}. Denote 
 $$C_{\alpha}(t):=\Big( \frac{cp_i}{p^0}\sqrt{\mu^c},w_{\ell_3-|\alpha|}^2\frac{cp_i}{p^0}\sqrt{\mu^c}\Big)>0,$$ 
 which has a uniform lower bound. Then we get
 \begin{align*}
    \partial^{\alpha}E^{3,c}_i=&\frac{1}{2C_{\alpha}(t)} \partial_t\Big( \big[\partial^{\alpha}f^{3,c}_{+}-\partial^{\alpha}f^{3,c}_{-}\big], w_{\ell_3-|\alpha|}^2\frac{cp_i}{p^0}\sqrt{\mu^c}\Big)\\
    &-\frac{1}{2C_{\alpha}(t)} \Big( \big[\partial^{\alpha}f^{3,c}_{+}-\partial^{\alpha}f^{3,c}_{-}\big], \partial_t\big(w_{\ell_3-|\alpha|}^2\big)\frac{cp_i}{p^0}\sqrt{\mu^c}\Big)\\
    &+ \frac{1}{2C_{\alpha}(t)}\Big( \frac{cp}{p^0}\cdot \big[\partial^{\alpha}\nabla_xf^{3,c}_{+}-\partial^{\alpha}\nabla_xf^{3,c}_{-}\big],w_{\ell_3-|\alpha|}^2
    \frac{cp_i}{p^0}\sqrt{\mu^c}\Big)\nonumber\\
    &-\frac{1}{2C_{\alpha}(t)}\Big( \partial^{\alpha}\mathcal{I}_{R,+}-\partial^{\alpha}\mathcal{I}_{R,-}, w_{\ell_3-|\alpha|}^2\frac{cp_i}{p^0}\sqrt{\mu^c}\Big).
 \end{align*}
 Now we plug the above expression of $\partial^{\alpha}E^{3,c}_i$ in $\mathcal{I}_{07}$ to have
 \begin{align*}
   \mathcal{I}_{07}= & \sum_{|\alpha|=0}^4\frac{-1}{4C_{\alpha}(t)}\frac{\mathrm{d}}{\mathrm{d}t}\Big\langle \Big(\frac{cp}{p^0}\sqrt{\mu^c}\zeta_0, w_{\ell_3-|\alpha|}^2\partial^{\alpha}f^{3,c}\Big),\Big(\frac{cp}{p^0}\sqrt{\mu^c}\zeta_0, w_{\ell_3-|\alpha|}^2\partial^{\alpha}f^{3,c}\Big)\Big\rangle \\
      &+\sum_{|\alpha|=0}^4\Big\langle \frac{1}{2C_{\alpha}(t)}\Big(\frac{cp}{p^0}\sqrt{\mu^c}\zeta_0, \partial_t\big(w_{\ell_3-|\alpha|}^2\big)\partial^{\alpha}f^{3,c}\Big),\Big(\frac{cp}{p^0}\sqrt{\mu^c}\zeta_0, w_{\ell_3-|\alpha|}^2\partial^{\alpha}f^{3,c}\Big)\Big\rangle  \\
   &-\sum_{|\alpha|=0}^4\Big\langle \frac{1}{2C_{\alpha}(t)}\Big(\frac{cp}{p^0}\sqrt{\mu^c}\zeta_0, w_{\ell_3-|\alpha|}^2\frac{cp}{p^0}\cdot\partial^{\alpha}\nabla_xf^{3,c}\Big),\Big(\frac{cp}{p^0}\sqrt{\mu^c}\zeta_0, w_{\ell_3-|\alpha|}^2\partial^{\alpha}f^{3,c}\Big)\Big\rangle \\
   &+\sum_{|\alpha|=0}^4\Big\langle \frac{1}{2C_{\alpha}(t)}\Big(\frac{cp}{p^0}\sqrt{\mu^c}\zeta_0, w_{\ell_3-|\alpha|}^2\partial^{\alpha}\mathcal{I}_{R}\Big),\Big(\frac{cp}{p^0}\sqrt{\mu^c}\zeta_0, w_{\ell_3-|\alpha|}^2\partial^{\alpha}f^{3,c}\Big)\Big\rangle.  
 \end{align*}
 Then we can use the {\it a priori} assumption
\eqref{apriori-assu} to further estimate $\mathcal{I}_{07}$ as
 \begin{align*}
   \mathcal{I}_{07}\leq & \sum_{|\alpha|=0}^4\frac{-1}{4C_{\alpha}(t)}\frac{\mathrm{d}}{\mathrm{d}t}\Big\langle \Big(\frac{cp}{p^0}\sqrt{\mu^c}\zeta_0, w_{\ell_3-|\alpha|}^2\partial^{\alpha}f^{3,c}\Big),\Big(\frac{cp}{p^0}\sqrt{\mu^c}\zeta_0, w_{\ell_3-|\alpha|}^2\partial^{\alpha}f^{3,c}\Big)\Big\rangle \\
   &+C\| f^{3,c}\|_{H^5_{\sigma}}^2 +CM \sum_{j=0}^2 \mathcal{D}^{j,c}(t).  
 \end{align*}
 \end{proof}
 
Now we come to establish weighted energy estimates of the nonlinear terms involving the Lorenz force. In the following two lemmas, we will treat the momentum growth terms involving the electric field without momentum derivatives and the momentum derivative terms involving the Lorenz force successively. 
\begin{lemma}\label{grow-wei}
Under the assumptions in Lemma \ref{grow-low}, there is a constant $C>0$, which is independent of $c$, such that for $0\leq t\leq \overline{T}$,
\begin{align}\label{growEfj-wei}
&\sum_{|\alpha|=0}^4\sum_{\alpha_1\leq\alpha}\Big|\Big\langle \frac{1}{2}\zeta_1\frac{cp}{p^0}\cdot \partial^{\alpha_1}E^{3,c} \sum_{j=0}^2\frac{\partial^{\alpha-\alpha_1}f^{j,c}}{c^j}, w_{\ell_3-|\alpha|}^2\partial^{\alpha}f^{3,c}\Big\rangle\Big|\nonumber\\
  &\qquad  \leq\,o(1)  \big\|f^{3,c}\big\|_{H^{4,3}_{\sigma,w}}^2+CM \sum_{j=0}^2 \frac{\mathcal{D}^{j,c}(t)}{c^{2j}},
   \end{align}
   
\begin{align*}
&\sum_{|\alpha|=0}^4\sum_{\alpha_1\leq\alpha}\Big|\Big\langle \frac{1}{2}\zeta_1\frac{cp}{p^0}\cdot \sum_{j=0}^2\frac{\partial^{\alpha_1}E^{j,c} }{c^j} \partial^{\alpha-\alpha_1}f^{3,c}, w_{\ell_3-|\alpha|}^2\partial^{\alpha}f^{3,c}\Big\rangle\Big|\nonumber\\
   &\qquad \leq\,o(1) \big\|f^{3,c}\big\|_{H^{4,3}_{\sigma,w}}^2+CM \mathcal{D}^{3,c}(t),
   \end{align*}   
\begin{align}\label{growEf-wei}
&\sum_{|\alpha|=0}^4\sum_{\alpha_1\leq\alpha}\Big|\Big\langle \frac{1}{2}\zeta_1\frac{cp}{p^0}\cdot \sum_{j=0}^2\frac{\partial^{\alpha_1}E^{3,c} }{c^3} \partial^{\alpha-\alpha_1}f^{3,c}, w_{\ell_3-|\alpha|}^2\partial^{\alpha}f^{3,c}\Big\rangle\Big|\nonumber\\
&\qquad \leq o(1) \big\|f^{3,c}\big\|_{H^{4,3}_{\sigma,w}}^2+CM \mathcal{D}^{3,c}(t),
   \end{align}    
 and
 \begin{align*}
&\sum_{|\alpha|=0}^4\sum_{\alpha_1\leq\alpha}\Big|\Big\langle \frac{1}{2}\zeta_1\frac{cp}{p^0}\cdot \sum_{\substack{j_1+j_2\geq3\\1\leq j_1,j_2< 3}}\frac{\partial^{\alpha_1}E^{j_1,c} \partial^{\alpha-\alpha_1}f^{j_2,c}}{c^{j_1+j_2-3}}, w_{\ell_3-|\alpha|}^2\partial^{\alpha}f^{3,c}\Big\rangle\Big|\nonumber\\
  &\qquad  \leq\,o(1) \big\|f^{3,c}\big\|_{H^{4,3}_{\sigma,w}}^2+CM \sum_{j=0}^2 \mathcal{D}^{j,c}(t).
   \end{align*}   
\end{lemma}
\begin{proof} As in Lemma \ref{grow-low}, we only prove \eqref{growEfj-wei} and \eqref{growEf-wei}. We first show \eqref{growEfj-wei}.

We take $L^{\infty}-L^2-L^2$ for $|\alpha_1|=0$, $L^{3}-L^6-L^2$ for $|\alpha_1|=1$, and $L^2-L^{\infty}-L^2$ for $|\alpha_1|\geq2$, and   bound the L.H.S. of \eqref{growEfj-wei} by
\begin{align*}
&C \sum_{|\alpha|=0}^4\sum_{j=0}^2\Big( \Big\|E^{3,c}\big\|_{H^2} \big\|\langle p\rangle^{\frac{3}{2}}\frac{f^{j,c}}{c^j}\big\|_{H^{4,3}_{w}} +\sum_{|\alpha_1|\geq3}\|\nabla_x^{|\alpha_1|}E^{3,c}\| \big\|\langle p\rangle^{\frac{3}{2}}\frac{f^{j,c}}{c^j}\big\|_{H^{3,3}_{w}} \Big)\big\|f^{3,c}\big\|_{H^{4,3}_{\sigma,w}} \\
&\qquad\leq o(1) \big\|f^{3,c}\big\|_{H^{4,3}_{\sigma,w}}^2+CM \sum_{j=0}^2 \frac{\mathcal{D}^{j,c}(t)}{c^{2j}}.
\end{align*}
Here we used the definition of the norm $\mathcal{D}^{j,c}(t)$ given in \eqref{dissirat-ic} for  $j=0, 1, 2$ and the fact $\ell_j=\ell_3+2(3-j)$.

Now we turn to prove \eqref{growEf-wei}.  We denote its L.H.S. as $\mathcal{I}_{08}$, take $L^{\infty}-L^2-L^2$ for $|\alpha_1|=0$, $L^{3}-L^6-L^2$ for $|\alpha_1|=1$, and $L^2-L^{\infty}-L^2$ for $|\alpha_1|\geq2$, and use the definition of $\mathcal{D}^{3,c}(t)$ in \eqref{dissirat-3c} and the {\it a priori} assumption \eqref{apriori-assu} to have
\begin{align*}
\mathcal{I}_{08}\leq& \frac{C}{c^3}\sum_{|\alpha|=0}^4\Big( \Big\|E^{3,c}\big\|_{H^2}\big\|\langle p\rangle^{\frac{1}{2}}f^{3,c}\big\|_{H^{4,3}_{w}}^2+\sum_{|\alpha_1|\geq3}\|\nabla_x^{|\alpha_1|}E^{3,c}\|\\
&\times\big\|\langle p\rangle^{\frac{3}{2}}w_{\ell_3-|\alpha|}\partial^{\alpha-\alpha_1}f^{3,c}\big\|_{H^2}\big\|\langle p\rangle^{-\frac{1}{2}}w_{\ell_3-|\alpha|}\partial^{\alpha}f^{3,c}\big\|\Big)\\
\leq& o(1) \big\|f^{3,c}\big\|_{H^{4,3}_{\sigma,w}}^2+ \frac{C}{c^6}\Big(c^3\Big\|E^{3,c}\big\|_{H^2}Y(t)^{-1}\mathcal{D}^{3,c}(t)
+\Big\|E^{3,c}\big\|_{H^4}^2\mathcal{D}^{3,c}_{4,w}(t)\Big)\\
\leq& o(1) \big\|f^{3,c}\big\|_{H^{4,3}_{\sigma,w}}^2+CM \mathcal{D}^{3,c}(t).
\end{align*}

\end{proof}   

\begin{lemma}\label{deriv-wei}
Under the assumptions in Lemma \ref{grow-low}, there is a constant $C>0$, which is independent of $c$, such that for $0\leq t\leq \overline{T}$,
\begin{align}\label{derivEfj-wei}
&\sum_{|\alpha|=0}^4\sum_{\alpha_1\leq\alpha}\Big|\Big\langle \zeta_1\partial^{\alpha_1}E^{3,c}\cdot\sum_{j=0}^2\frac{\partial^{\alpha-\alpha_1}\nabla_pf^{j,c}}{c^j}, w_{\ell_3-|\alpha|}^2\partial^{\alpha}f^{3,c}\Big\rangle\Big|\nonumber\\
  &\qquad  \leq\,o(1)  \big\|f^{3,c}\big\|_{H^{4,3}_{\sigma,w}}^2+CM \sum_{j=0}^2 \frac{\mathcal{D}^{j,c}(t)}{c^{2j}},
   \end{align}
\begin{align*}
&\sum_{|\alpha|=0}^4\sum_{\alpha_1\leq\alpha}\Big|\Big\langle \zeta_1\sum_{j=0}^2\frac{\partial^{\alpha_1}E^{j,c} }{c^j} \cdot\partial^{\alpha-\alpha_1}\nabla_pf^{3,c}, w_{\ell_3-|\alpha|}^2\partial^{\alpha}f^{3,c}\Big\rangle\Big|\nonumber\\
  &\qquad  \leq\,o(1)  \big\|f^{3,c}\big\|_{H^{4,3}_{\sigma,w}}^2+CM \mathcal{D}^{3,c}(t),
   \end{align*}   

\begin{align}\label{derivEf-wei}
&\sum_{|\alpha|=0}^4\sum_{\alpha_1\leq\alpha}\Big|\Big\langle \zeta_1\partial^{\alpha_1}\big(E^{3,c}+\frac{p}{p^0}\times B^{3,c}\big)\cdot\frac{\partial^{\alpha-\alpha_1}\nabla_pf^{3,c}}{c^3} , w_{\ell_3-|\alpha|}^2\partial^{\alpha}f^{3,c}\Big\rangle\Big|\nonumber\\
&\qquad\leq \Big(o(1)+\frac{C\sqrt{M}}{c^3}\Big) \big\|f^{3,c}\big\|_{H^{4,3}_{\sigma,w}}^2+CM \mathcal{D}^{3,c}(t),
   \end{align}    
 and
 \begin{align*}
&\sum_{|\alpha|=0}^4\sum_{\alpha_1\leq\alpha}\Big|\Big\langle \zeta_1 \sum_{\substack{j_1+j_2\geq3\\1\leq j_1,j_2< 3}}\frac{\partial^{\alpha_1}E^{j_1,c} \cdot \partial^{\alpha-\alpha_1}\nabla_pf^{j_2,c}}{c^{j_1+j_2-3}}, w_{\ell_3-|\alpha|}^2\partial^{\alpha}f^{3,c}\Big\rangle\Big|\nonumber\\
  &\qquad  \leq\,o(1)  \big\|f^{3,c}\big\|_{H^{4,3}_{\sigma,w}}^2+CM \sum_{j=0}^2 \frac{\mathcal{D}^{j,c}(t)}{c^{2j}}.
   \end{align*}   
\end{lemma}
\begin{proof} As before, we only prove \eqref{derivEfj-wei} and \eqref{derivEf-wei}.

As the estimation of \eqref{growEfj-wei}, we can bound the L.H.S. of \eqref{derivEfj-wei} by
\begin{align*}
&C \sum_{|\alpha|=0}^4\sum_{j=0}^2\Big( \big(\|E^{3,c}\big\|_{H^2}+\|B^{3,c}-\overline{B}^{3,c}\big\|_{H^2}\big) \big\|\langle p\rangle^{\frac{1}{2}}\frac{\nabla_pf^{j,c}}{c^j}\big\|_{H^{4,3}_{w}} \\
&\qquad+\sum_{|\alpha_1|\geq3}\left(\|\nabla_x^{|\alpha_1|}E^{3,c}\| +\|\nabla_x^{|\alpha_1|}B^{3,c}\| \right)\big\|\langle p\rangle^{\frac{1}{2}}\frac{\nabla_pf^{j,c}}{c^j}\big\|_{H^{3,3}_{w}} \Big)\big\|f^{3,c}\big\|_{H^{4,3}_{\sigma,w}} \\
&\qquad\leq o(1) \big\|f^{3,c}\big\|_{H^{4,3}_{\sigma,w}}^2+CM \sum_{j=0}^2 \frac{\mathcal{D}^{j,c}(t)}{c^{2j}}.
\end{align*}

Now we turn to prove \eqref{derivEf-wei}. Denote its L.H.S. as $\mathcal{I}_{09}$. We first apply integration by parts w.r.t. $p$, then use $L^{\infty}-L^2-L^2$ for $|\alpha_1|\leq2$ and $L^2-L^{\infty}-L^2$ for $|\alpha_1|\geq3$ to have
\begin{align*}
\mathcal{I}_{09}\leq& \frac{C}{c^3}\sum_{|\alpha|=0}^4 \|E^{3,c}\|_{H^2}\big\|f^{3,c}\big\|^2_{H^{4,3}_{w}}
+\big\|\langle p\rangle^{-\frac{3}{2}}w_{\ell_3-|\alpha|}\partial^{\alpha}\nabla_pf^{3,c}\big\|\\ &\times\Big(\sum_{1\leq|\alpha_1|\leq2}\big(\|\nabla_x^{|\alpha_1|}E^{3,c}\|_{H^2}
+\|\nabla_x^{|\alpha_1|}B^{3,c}\|_{H^2}\big) \big\|\langle p\rangle^{\frac{3}{2}}w_{\ell_3-|\alpha|}\partial^{\alpha-\alpha_1}f^{3,c}\big\|\\
&+\sum_{|\alpha_1|\geq3}\big(\|\nabla_x^{|\alpha_1|}E^{3,c}+\|\nabla_x^{|\alpha_1|}B^{3,c}\|\big) \big\|\langle p\rangle^{\frac{3}{2}}w_{\ell_3-|\alpha|}\partial^{\alpha-\alpha_1}f^{3,c}\big\|_{H^2}\Big)\Big]\\
\leq& \Big(o(1)+\frac{C\sqrt{M}}{c^3}\Big) \big\|f^{3,c}\big\|_{H^{4,3}_{\sigma,w}}^2+CM \mathcal{D}^{3,c}(t).
\end{align*}

\end{proof}   

Now we are at the position to prove the weighted energy estimates of  $f^{3,c}(t, x, p)$ except for the non-highest order derivatives. 

\begin{proposition}\label{5norm-wei00} Under the assumptions in Lemma \ref{grow-low}, there is a $c-$independent constant $C>0$ such that for $0\leq t\leq \overline{T}$,
\begin{align} \label{5norm-wei}
&\frac{\mathrm{d}}{\mathrm{d}t}\Big(\|f^{3,c}(t)\|^2_{H^{4,3}_{w}}
+\sum_{|\alpha|=0}^4\frac{1}{4C_{\alpha}(t)}\Big\langle \Big(\frac{cp}{p^0}\sqrt{\mu^c}\zeta_0, w_{\ell_3-|\alpha|}^2\partial^{\alpha}f^{3,c}\Big),\Big(\frac{cp}{p^0}\sqrt{\mu^c}\zeta_0, w_{\ell_3-|\alpha|}^2\partial^{\alpha}f^{3,c}\Big)\Big\rangle\Big)\nonumber\\
&+Y(t)\big\|\langle p\rangle^{\frac{1}{2}}f^{3,c}\big\|^2_{H^{4,3}_{w}}+\frac{1}{5} \|f^{3,c}(t)\|^2_{H^{4,3}_{w,\sigma}}
\leq C \big\|f^{3,c}\big\|_{H^4_{\sigma}}^2
+CM\sum_{j=0}^3\mathcal{D}^{j,c}(t). 
\end{align}
\end{proposition}
\begin{proof} By \eqref{nonlin-w} in Lemma \ref{w-nonlin}, we take $L^{\infty}-L^2-L^2$ for $|\alpha_1|\leq3$ and $L^2-L^{\infty}-L^2$ for $|\alpha_1|\geq3$ to obtain 
\begin{align}\label{nonlin-w5j3}
&\sum_{|\alpha|=0}^4\sum_{\alpha_1\leq\alpha}\Big|\Big\langle \sum_{j=0}^2\frac{1}{c^j}\left[\varGamma\left(\partial^{\alpha_1}f^{j,c}, \partial^{\alpha-\alpha_1}f^{3,c}\right)+\varGamma\left(\partial^{\alpha_1}f^{3,c}, \partial^{\alpha-\alpha_1}f^{j,c}\right)\right]\nonumber\\
&\qquad+\frac{1}{c^3}\varGamma\left(\partial^{\alpha_1}f^{3,c}, \partial^{\alpha-\alpha_1}f^{3,c}\right), w_{\ell_3-|\alpha|}^2\partial^{\alpha}f^{3,c}\Big\rangle\Big|\nonumber\\
&\qquad \leq C\sum_{|\alpha|=0}^4\bigg[\sum_{|\alpha_1|\leq2}\bigg(
\sum_{j=0}^2\frac{1}{c^j}\Big(\big\|\big(\mu^c\big)^{\frac{1}{9}}\partial^{\alpha_1}f^{j,c}\big\|_{H^2}
\|w_{\ell_3-|\alpha|}\partial^{\alpha-\alpha_1}f^{3,c}\|_{\sigma}\nonumber\\
&\qquad+
\big\|\big(\mu^c\big)^{\frac{1}{20}}\partial^{\alpha_1}f^{j,c}\big\|_{H^2_{\sigma}}
\|w_{\ell_3-|\alpha|}\partial^{\alpha-\alpha_1}f^{3,c}\|+
\big\|\big(\mu^c\big)^{\frac{1}{20}}\partial^{\alpha_1}f^{3,c}\big\|_{H^2}
\|w_{\ell_3-|\alpha|}\partial^{\alpha-\alpha_1}f^{j,c}\|_{\sigma}\nonumber\\
&\qquad+
\big\|\big(\mu^c\big)^{\frac{1}{20}}\partial^{\alpha_1}f^{3,c}\big\|_{H^2_{\sigma}}
\|w_{\ell_3-|\alpha|}\partial^{\alpha-\alpha_1}f^{j,c}\|\Big)
+\frac{1}{c^3}\big\|\big(\mu^c\big)^{\frac{1}{20}}\partial^{\alpha_1}f^{3,c}\big\|_{H^2}
\|w_{\ell_3-|\alpha|}\partial^{\alpha-\alpha_1}f^{3,c}\|_{\sigma}\nonumber\\
&\qquad+\frac{1}{c^3}
\big\|\big(\mu^c\big)^{\frac{1}{20}}\partial^{\alpha_1}f^{3,c}\big\|_{H^2_{\sigma}}
\|w_{\ell_3-|\alpha|}\partial^{\alpha-\alpha_1}f^{3,c}\|\bigg)\nonumber\\
&\qquad+\sum_{|\alpha_1|=3}^4\bigg(
\sum_{j=0}^2\frac{1}{c^j}\Big(\big\|\big(\mu^c\big)^{\frac{1}{20}}\partial^{\alpha_1}f^{j,c}\big\|
\|w_{\ell_3-|\alpha|}\partial^{\alpha-\alpha_1}f^{3,c}\|_{H^2_{\sigma}}\nonumber\\
&\qquad+
\big\|\big(\mu^c\big)^{\frac{1}{20}}\partial^{\alpha_1}f^{j,c}\big\|_{\sigma}
\|w_{\ell_3-|\alpha|}\partial^{\alpha-\alpha_1}f^{3,c}\|_{H^2}+
\big\|\big(\mu^c\big)^{\frac{1}{20}}\partial^{\alpha_1}f^{3,c}\big\|
\|w_{\ell_3-|\alpha|}\partial^{\alpha-\alpha_1}f^{j,c}\|_{H^2_{\sigma}}\nonumber\\
&\qquad+
\big\|\big(\mu^c\big)^{\frac{1}{20}}\partial^{\alpha_1}f^{3,c}\big\|_{\sigma}
\|w_{\ell_3-|\alpha|}\partial^{\alpha-\alpha_1}f^{j,c}\|_{H^2}\Big)
+\frac{1}{c^3}\big\|\big(\mu^c\big)^{\frac{1}{20}}\partial^{\alpha_1}f^{3,c}\big\|
\|w_{\ell_3-|\alpha|}\partial^{\alpha-\alpha_1}f^{3,c}\|_{H^2_{\sigma}}\nonumber\\
&\qquad+\frac{1}{c^3}
\big\|\big(\mu^c\big)^{\frac{1}{20}}\partial^{\alpha_1}f^{3,c}\big\|_{\sigma}
\|w_{\ell_3-|\alpha|}\partial^{\alpha-\alpha_1}f^{3,c}\|_{H^2}\bigg)\bigg]
\|w_{\ell_3-|\alpha|}\partial^{\alpha}f^{3,c}\|_{H^2_{\sigma}}\nonumber\\
&\qquad\leq o(1)\big\|f^{3,c}\big\|_{H^{4,3}_{w,\sigma}}^2
+\sum_{j=0}^3\frac{1}{c^{2j}}\Big(\|f^{j,c}\|_{H^4}^2\big\|f^{3,c}\big\|_{H^{4,3}_{w,\sigma}}^2
+\big\|f^{3,c}\big\|_{H^{4,3}_{w}}^2\|f^{j,c}\|_{H^{4}_{\sigma}}^2\nonumber\\
&\qquad+\big\|f^{3,c}\big\|_{H^4}^2\|f^{j,c}\|_{H^{4,3}_{w,\sigma}}^2
+\|f^{j,c}\|_{H^{4,3}_{w}}^2\big\|f^{3,c}\big\|_{H^{4}_{\sigma}}^2\Big)\nonumber\\
&\qquad\leq \Big(o(1)+\sum_{j=0}^3\frac{CM}{c^{2j}}\Big)\big\|f^{3,c}\big\|_{H^{4,3}_{w,\sigma}}^2
+\sum_{j=0}^2\frac{CM\mathcal{D}^{j,c}(t)}{c^{2j}},
\end{align}
and
\begin{align}\label{nonlin-w5jj}
&\sum_{|\alpha|=0}^4\sum_{\alpha_1\leq\alpha}\Big|\Big\langle \sum_{\substack{j_1+j_2\geq3   \\1\leq j_1,j_2<3}}\frac{\varGamma\left(\partial^{\alpha_1}f^{j_1,c}, \partial^{\alpha-\alpha_1}f^{j_2,c}\right)}{c^{j_1+j_2-3}}, w_{\ell_3-|\alpha|}^2\partial^{\alpha}f^{3,c}\Big\rangle\Big|\nonumber\\
&\qquad\leq o(1)\big\|f^{3,c}\big\|_{H^{4,3}_{w,\sigma}}^2
+\sum_{\substack{j_1+j_2\geq3   \\1\leq j_1,j_2<3}}\frac{1}{c^{2(j_1+j_2-3)}}\Big(\|f^{j_1,c}\|_{H^4}^2\|f^{j_2,c}\|_{H^{4,3}_{w,\sigma}}^2
+\|f^{j_1,c}\|_{H^{4,3}_{w}}^2\|f^{j_2,c}\|_{H^{4}_{\sigma}}^2\Big)\nonumber\\
&\qquad\leq o(1)\big\|f^{3,c}\big\|_{H^{4,3}_{w,\sigma}}^2
+\sum_{j=0}^2CM\mathcal{D}^{j,c}(t).
\end{align}
Now we apply $\partial^{\alpha}$ with $|\alpha|\leq4$ to the equations of $f^{3,c}(t, x, p)$ in \eqref{mainF3}, take the  inner products with $w_{\ell_3-|\alpha|}^2\partial^{\alpha}f^{3,c}$ in $H^{4}$, and use Lemma \ref{coerci-w}, Lemma \ref{lin-wei}, Lemma \ref{grow-wei}, Lemma \ref{deriv-wei},  \eqref{nonlin-w5j3} and \eqref{nonlin-w5jj} to deduce \eqref{5norm-wei}.
\end{proof}

\noindent\underline{{\it Step 2. Highest order weighted energy estimates.}} 
 To proceed the highest order weighted energy estimates of $f^{3,c}(t,x,p)$, as in Lemma \ref{grow-wei} and Lemma \ref{deriv-wei}, we first derive the weighted energy estimates of nonlinear terms involving the electromagnetic field in the following two lemmas.

\begin{lemma}\label{grow-hwei}
Under the assumptions in Lemma \ref{grow-low}, there is a constant $C>0$, which is independent of $c$, such that for $0\leq t\leq \overline{T}$,
\begin{align}\label{growEfj-hwei}
&\sum_{|\alpha|=5}\sum_{\alpha_1\leq\alpha}\Big|\Big\langle \frac{1}{2}\zeta_1\frac{cp}{p^0}\cdot \partial^{\alpha_1}E^{3,c} \sum_{j=0}^2\frac{\partial^{\alpha-\alpha_1}f^{j,c}}{c^j}, (1+t)^{-\frac{1+\epsilon_{0}}{2}}w_{\ell_3-5}^2\partial^{\alpha}f^{3,c}\Big\rangle\Big|\nonumber\\
  &\qquad  \leq\,o(1)  (1+t)^{-\frac{1+\epsilon_{0}}{2}}\|w_{\ell_3-5}\nabla_x^{5}f^{3,c}\|_{\sigma}^2+CM \sum_{j=0}^2 \frac{\mathcal{D}^{j,c}(t)}{c^{2j}},
   \end{align}
   
\begin{align*}
&\sum_{|\alpha|=5}\sum_{\alpha_1\leq\alpha}\Big|\Big\langle \frac{1}{2}\zeta_1\frac{cp}{p^0}\cdot \sum_{j=0}^2\frac{\partial^{\alpha_1}E^{j,c} }{c^j} \partial^{\alpha-\alpha_1}f^{3,c}, (1+t)^{-\frac{1+\epsilon_{0}}{2}}w_{\ell_3-5}^2\partial^{\alpha}f^{3,c}\Big\rangle\Big|\nonumber\\
   &\qquad \leq\,o(1) (1+t)^{-\frac{1+\epsilon_{0}}{2}}\|w_{\ell_3-5}\nabla_x^{5}f^{3,c}\|_{\sigma}^2+CM \mathcal{D}^{3,c}(t),
   \end{align*}   
\begin{align}\label{growEf-hwei}
&\sum_{|\alpha|=5}\sum_{\alpha_1\leq\alpha}\Big|\Big\langle \frac{1}{2}\zeta_1\frac{cp}{p^0}\cdot \sum_{j=0}^2\frac{\partial^{\alpha_1}E^{3,c} }{c^3} \partial^{\alpha-\alpha_1}f^{3,c}, (1+t)^{-\frac{1+\epsilon_{0}}{2}}w_{\ell_3-5}^2\partial^{\alpha}f^{3,c}\Big\rangle\Big|\nonumber\\
&\qquad \leq o(1) (1+t)^{-\frac{1+\epsilon_{0}}{2}}\|w_{\ell_3-5}\nabla_x^{5}f^{3,c}\|_{\sigma}^2+CM \mathcal{D}^{3,c}(t),
   \end{align}    
 and
 \begin{align*}
&\sum_{|\alpha|=5}\sum_{\alpha_1\leq\alpha}\Big|\Big\langle \frac{1}{2}\zeta_1\frac{cp}{p^0}\cdot \sum_{\substack{j_1+j_2\geq3\\1\leq j_1,j_2< 3}}\frac{\partial^{\alpha_1}E^{j_1,c} \partial^{\alpha-\alpha_1}f^{j_2,c}}{c^{j_1+j_2-3}}, (1+t)^{-\frac{1+\epsilon_{0}}{2}}w_{\ell_3-5}^2\partial^{\alpha}f^{3,c}\Big\rangle\Big|\nonumber\\
  &\qquad  \leq\,o(1) (1+t)^{-\frac{1+\epsilon_{0}}{2}}\|w_{\ell_3-5}\nabla_x^{5}f^{3,c}\|_{\sigma}^2+CM \sum_{j=0}^2 \mathcal{D}^{j,c}(t).
   \end{align*}   
\end{lemma}
\begin{proof} As in in Lemma \ref{grow-low}, we only prove \eqref{growEfj-hwei} and \eqref{growEf-hwei}. We first show \eqref{growEfj-hwei}.

As the estimation of \eqref{growEfj-wei}, we take $L^{\infty}-L^2-L^2$ for $|\alpha_1|=0$, $L^{3}-L^6-L^2$ for $|\alpha_1|=1$, and $L^2-L^{\infty}-L^2$ for $|\alpha_1|\geq2$, and   bound the L.H.S. of \eqref{growEfj-hwei} by
\begin{align*}
&C(1+t)^{-\frac{1+\epsilon_{0}}{2}} \sum_{|\alpha|=5}\sum_{j=0}^2\Big( \Big\|E^{3,c}\big\|_{H^2} \big\|\langle p\rangle^{\frac{3}{2}}\frac{f^{j,c}}{c^j}\big\|_{H^{5,3}_{w}} \\ &+\sum_{|\alpha_1|\geq3}\|\nabla_x^{|\alpha_1|}E^{3,c}\| \big\|\langle p\rangle^{\frac{3}{2}}\frac{f^{j,c}}{c^j}\big\|_{H^{5,3}_{w}} \Big)\|w_{\ell_3-5}\nabla_x^{5}f^{3,c}\|_{\sigma} \\
&\qquad\leq o(1)(1+t)^{-\frac{1+\epsilon_{0}}{2}} \|w_{\ell_3-5}\nabla_x^{5}f^{3,c}\|_{\sigma}^2+CM \sum_{j=0}^2 \frac{\mathcal{D}^{j,c}(t)}{c^{2j}}.
\end{align*}
Here we used the definition of the norm $\mathcal{D}^{j,c}(t)$ given in \eqref{dissirat-ic} for  $j=0, 1, 2$ and the fact $\ell_j=\ell_3+2(3-j)$.

To prove \eqref{growEf-hwei}, we denote its L.H.S. as $\mathcal{I}_{10}$. We further take $L^{\infty}-L^2-L^2$ for $|\alpha_1|=0$, $L^{3}-L^6-L^2$ for $|\alpha_1|=1$, and $L^2-L^{\infty}-L^2$ for $|\alpha_1|\geq2$, and use the definition of $\mathcal{D}^{3,c}(t)$ in \eqref{dissirat-3c} and the {\it a priori} assumption \eqref{apriori-assu} to have
\begin{align*}
\mathcal{I}_{10}\leq& \frac{C(1+t)^{-\frac{1+\epsilon_{0}}{2}}}{c^3}\sum_{|\alpha|=5}m\Big( \Big\|E^{3,c}\big\|_{H^2}\big\|\langle p\rangle^{\frac{1}{2}}f^{3,c}\big\|_{H^{5,3}_{w}}^2+\sum_{|\alpha_1|\geq3}\|\nabla_x^{|\alpha_1|}E^{3,c}\|\\
&\times\big\|\langle p\rangle^{\frac{3}{2}}w_{\ell_3-5}\partial^{\alpha-\alpha_1}f^{3,c}\big\|_{H^2}\big\|\langle p\rangle^{-\frac{1}{2}}w_{\ell_3-5}\partial^{\alpha}f^{3,c}\big\|\Big)\\
\leq& o(1) (1+t)^{-\frac{1+\epsilon_{0}}{2}} \|w_{\ell_3-5}\nabla_x^{5}f^{3,c}\|_{\sigma}^2+ \frac{C}{c^6}\Big(c^3\Big\|E^{3,c}\big\|^2_{H^2}Y(t)^{-1}\mathcal{D}^{3,c}(t)
+\Big\|E^{3,c}\big\|_{H^5}^2\mathcal{D}^{3,c}_{4,w}(t)\Big)\\
\leq& o(1) (1+t)^{-\frac{1+\epsilon_{0}}{2}} \|w_{\ell_3-5}\nabla_x^{5}f^{3,c}\|_{\sigma}^2+CM \mathcal{D}^{3,c}(t).
\end{align*}

\end{proof}   

\begin{lemma}\label{deriv-hwei}
Under the assumptions in Lemma \ref{grow-low}, there is a constant $C>0$, which is independent of $c$, such that for $0\leq t\leq \overline{T}$,
\begin{align}\label{derivEfj-hwei}
&\sum_{|\alpha|=5}\sum_{\alpha_1\leq\alpha}\Big|\Big\langle \zeta_1\partial^{\alpha_1}E^{3,c}\cdot\sum_{j=0}^2\frac{\partial^{\alpha-\alpha_1}\nabla_pf^{j,c}}{c^j}, (1+t)^{-\frac{1+\epsilon_{0}}{2}}w_{\ell_3-5}^2\partial^{\alpha}f^{3,c}\Big\rangle\Big|\nonumber\\
  &\qquad  \leq\,o(1)  (1+t)^{-\frac{1+\epsilon_{0}}{2}}\|w_{\ell_3-5}\nabla_x^{5}f^{3,c}\|_{\sigma}^2+CM \sum_{j=0}^2 \frac{\mathcal{D}^{j,c}(t)}{c^{2j}},
   \end{align}
\begin{align}\label{derivEjf-hwei}
&\sum_{|\alpha|=5}\sum_{\alpha_1\leq\alpha}\Big|\Big\langle \zeta_1\sum_{j=0}^2\frac{\partial^{\alpha_1}E^{j,c} }{c^j} \cdot\partial^{\alpha-\alpha_1}\nabla_pf^{3,c}, (1+t)^{-\frac{1+\epsilon_{0}}{2}}w_{\ell_3-5}^2\partial^{\alpha}f^{3,c}\Big\rangle\Big|\nonumber\\
  &\qquad  \leq\,o(1)  (1+t)^{-\frac{1+\epsilon_{0}}{2}}\|w_{\ell_3-5}\nabla_x^{5}f^{3,c}\|_{\sigma}^2+CM \mathcal{D}^{3,c}(t),
   \end{align}   

\begin{align}\label{derivEf-hwei}
&\sum_{|\alpha|=5}\sum_{\alpha_1\leq\alpha}\Big|\Big\langle \zeta_1\partial^{\alpha_1}\big(E^{3,c}+\frac{p}{p^0}\times B^{3,c}\big)\cdot\frac{\partial^{\alpha-\alpha_1}\nabla_pf^{3,c}}{c^3} , (1+t)^{-\frac{1+\epsilon_{0}}{2}}w_{\ell_3-5}^2\partial^{\alpha}f^{3,c}\Big\rangle\Big|\nonumber\\
  &\qquad  \leq\,o(1)  (1+t)^{-\frac{1+\epsilon_{0}}{2}}\|w_{\ell_3-5}\nabla_x^{5}f^{3,c}\|_{\sigma}^2+CM \mathcal{D}^{3,c}(t),
   \end{align}    
 and
 \begin{align}\label{derivEjfj-hwei}
&\sum_{|\alpha|=5}\sum_{\alpha_1\leq\alpha}\Big|\Big\langle \zeta_1 \sum_{\substack{j_1+j_2\geq3\\1\leq j_1,j_2< 3}}\frac{\partial^{\alpha_1}E^{j_1,c} \cdot \partial^{\alpha-\alpha_1}\nabla_pf^{j_2,c}}{c^{j_1+j_2-3}}, (1+t)^{-\frac{1+\epsilon_{0}}{2}}w_{\ell_3-5}^2\partial^{\alpha}f^{3,c}\Big\rangle\Big|\nonumber\\
  &\qquad  \leq\,o(1)  (1+t)^{-\frac{1+\epsilon_{0}}{2}}\|w_{\ell_3-5}\nabla_x^{5}f^{3,c}\|_{\sigma}^2+CM \sum_{j=0}^2 \frac{\mathcal{D}^{j,c}(t)}{c^{2j}}.
   \end{align}   
\end{lemma}
\begin{proof} As before, we only prove \eqref{derivEfj-hwei} and \eqref{derivEf-hwei}.

As the estimation of \eqref{growEfj-hwei}, we can bound the L.H.S. of \eqref{derivEfj-hwei} by
\begin{align*}
&C (1+t)^{-\frac{1+\epsilon_{0}}{2}}\sum_{|\alpha|=5}\sum_{j=0}^2\Big( \Big\|E^{3,c}\big\|_{H^2} \big\|\langle p\rangle^{\frac{1}{2}}\frac{\nabla_pf^{j,c}}{c^j}\big\|_{H^{5,3}_{w}} \\
&\qquad+\sum_{|\alpha_1|\geq3}\|\nabla_x^{|\alpha_1|}E^{3,c}\| \big\|\langle p\rangle^{\frac{1}{2}}\frac{\nabla_pf^{j,c}}{c^j}\big\|_{H^{5,3}_{w}} \Big)\|w_{\ell_3-5}\nabla_x^{5}f^{3,c}\|_{\sigma} \\
&\qquad\leq o(1)(1+t)^{-\frac{1+\epsilon_{0}}{2}} \|w_{\ell_3-5}\nabla_x^{5}f^{3,c}\|_{\sigma}^2+CM \sum_{j=0}^2 \frac{\mathcal{D}^{j,c}(t)}{c^{2j}}.
\end{align*}

Now we turn to prove \eqref{derivEf-hwei}. Denote the L.H.S. of \eqref{derivEf-hwei} as $\mathcal{I}_{10}$. We first apply integration by parts w.r.t. $p$, and $L^{\infty}-L^2-L^2$ for $|\alpha_1|\leq2$ and $L^2-L^{\infty}-L^2$ for $|\alpha_1|\geq3$ to have
\begin{align*}
\mathcal{I}_{10}\leq& \frac{C(1+t)^{-\frac{1+\epsilon_{0}}{2}}}{c^3}\sum_{|\alpha|=5}\Big[ \big(\Big\|E^{3,c}\big\|_{H^2}+\|B^{3,c}-\bar{B}^{3}\|_{H^2}\big)\|w_{\ell_3-5}\nabla_x^{5}f^{3,c}\|^2+\big\|\langle p\rangle^{-\frac{3}{2}}w_{\ell_3-5}\partial^{\alpha}\nabla_pf^{3,c}\big\|\\ &\times\Big(\sum_{1\leq|\alpha_1|\leq2}\big(\|\nabla_x^{|\alpha_1|}E^{3,c}\|_{H^2}
+\|\nabla_x^{|\alpha_1|}B^{3,c}\|_{H^2}\big) \big\|\langle p\rangle^{\frac{3}{2}}w_{\ell_3-5}\partial^{\alpha-\alpha_1}f^{3,c}\big\|\\
&+\sum_{|\alpha_1|\geq3}\big(\|\nabla_x^{|\alpha_1|}E^{3,c}+\|\nabla_x^{|\alpha_1|}B^{3,c}\|\big) \big\|\langle p\rangle^{\frac{3}{2}}w_{\ell_3-|\alpha|}\partial^{\alpha-\alpha_1}f^{3,c}\big\|_{H^2}\Big)\Big]\\
\leq& \frac{C(1+t)^{-\frac{1+\epsilon_{0}}{2}}}{c^3} \big(\Big\|E^{3,c}\big\|_{H^2}+\|B^{3,c}-\bar{B}^{3}\|_{H^2}\big)\big\|\langle p\rangle^{\frac{1}{2}}w_{\ell_3-5}\nabla_x^{5}f^{3,c}\|\big\|\langle p\rangle^{-\frac{1}{2}}w_{\ell_3-5}\nabla_x^{5}f^{3,c}\|\\
 &+C(1+t)^{-\frac{1+\epsilon_{0}}{2}}\Big[ o(1)\|w_{\ell_3-5}\nabla_x^{5}f^{3,c}\|_{\sigma}^2 +\frac{C}{c^6}\Big(\Big\|E^{3,c}\big\|^2_{H^5}
+\|B^{3,c}-\bar{B}^{3}\|^2_{H^5}\Big)\big\|f^{3,c}\big\|_{H^{4,3}_{w,\sigma}}^2\Big]  \\
\leq& o(1) (1+t)^{-\frac{1+\epsilon_{0}}{2}} \|w_{\ell_3-5}\nabla_x^{5}f^{3,c}\|_{\sigma}^2+CM \mathcal{D}^{3,c}(t).
\end{align*}

\end{proof}   

Based on the above weighted estimates, we can finally derive the highest order weighted energy estimates of  $f^{3,c}(t, x, p)$. 

\begin{proposition}\label{6norm-wei00} Under the assumptions in Lemma \ref{grow-low}, there is a $c-$independent constant $C>0$ such that for $0\leq t\leq \overline{T}$,
\begin{align} \label{6norm-wei}
&\frac{\mathrm{d}}{\mathrm{d}t}\Big[(1+t)^{-\frac{1+\epsilon_{0}}{2}} \|w_{\ell_3-5}\nabla_x^{5}f^{3,c}\|^2\Big]+(1+t)^{-\frac{1+\epsilon_{0}}{2}}\Big( \frac{1+\epsilon_{0}}{2(1+t)}\|w_{\ell_3-5}\nabla_x^{5}f^{3,c}\|^2\nonumber\\
&\qquad+Y(t)\big\|\langle p\rangle^{\frac{1}{2}}w_{\ell_3-5}\nabla_x^{5}f^{3,c}\big\|^2+\frac{1}{5} \|w_{\ell_3-5}\nabla_x^{5}f^{3,c}(t)\|^2_{\sigma}\Big)
\nonumber\\
&\qquad\leq C(1+t)^{-(1+\epsilon_{0})}\|\nabla_x^5E^{3,c}\|^2+C\|\nabla_x^5f^{3,c}\|^2_{\sigma}
+CM\sum_{j=0}^3\mathcal{D}^{j,c}(t). 
\end{align}
\end{proposition}
\begin{proof} Note that
\begin{align}\label{lin-hwei0}
&\sum_{|\alpha|=5}\Big\langle \partial^{\alpha}E^{3,c}\cdot\frac{cp}{p^0}\sqrt{\mu^c}\zeta_0, (1+t)^{-\frac{1+\epsilon_{0}}{2}} w_{\ell_3-5}^2\partial^{\alpha}f^{3,c}\Big\rangle\nonumber\\
&\qquad   \leq C(1+t)^{-(1+\epsilon_{0})}\|\nabla_x^5E^{3,c}\|^2+C\|\nabla_x^5f^{3,c}\|^2_{\sigma}.
   \end{align}

 As the estimation in \eqref{nonlin-w5j3}, we  can obtain 
\begin{align}\label{nonlin-w6j3}
&\sum_{|\alpha|=5}\sum_{\alpha_1\leq\alpha}\Big|\Big\langle \sum_{j=0}^2\frac{1}{c^j}\left[\varGamma\left(\partial^{\alpha_1}f^{j,c}, \partial^{\alpha-\alpha_1}f^{3,c}\right)+\varGamma\left(\partial^{\alpha_1}f^{3,c}, \partial^{\alpha-\alpha_1}f^{j,c}\right)\right]\nonumber\\
&\qquad+\frac{1}{c^3}\varGamma\left(\partial^{\alpha_1}f^{3,c}, \partial^{\alpha-\alpha_1}f^{3,c}\right), (1+t)^{-\frac{1+\epsilon_{0}}{2}}w_{\ell_3-5}^2\partial^{\alpha}f^{3,c}\Big\rangle\Big|\nonumber\\
&\qquad\leq o(1)(1+t)^{-\frac{1+\epsilon_{0}}{2}}\|w_{\ell_3-5}\nabla_x^{5}f^{3,c}(t)\|^2_{\sigma}
+\sum_{j=0}^3\frac{(1+t)^{-\frac{1+\epsilon_{0}}{2}}}{c^{2j}}
\nonumber\\
&\qquad\times\Big(\|f^{j,c}\|_{H^5}^2\big\|f^{3,c}\big\|_{H^{5,3}_{w,\sigma}}^2
+\big\|f^{3,c}\big\|_{H^{5,3}_{w}}^2\|f^{j,c}\|_{H^{5}_{\sigma}}^2
+\big\|f^{3,c}\big\|_{H^5}^2\|f^{j,c}\|_{H^{5,3}_{w,\sigma}}^2\nonumber\\
&\qquad
+\|f^{j,c}\|_{H^{5,3}_{w}}^2\big\|f^{3,c}\big\|_{H^{5}_{\sigma}}^2
+\big\|f^{3,c}\big\|_{H^5}^2\|f^{3,c}\|_{H^{5,3}_{w,\sigma}}^2
+\|f^{3,c}\|_{H^{5,3}_{w}}^2\big\|f^{3,c}\big\|_{H^{5}_{\sigma}}^2\Big)\nonumber\\
&\qquad\leq o(1)(1+t)^{-\frac{1+\epsilon_{0}}{2}}\|w_{\ell_3-5}\nabla_x^{5}f^{3,c}(t)\|^2_{\sigma}
+\sum_{j=0}^3\frac{CM\mathcal{D}^{j,c}(t)}{c^{2j}},
\end{align}
and
\begin{align}\label{nonlin-w6jj}
&\sum_{|\alpha|=0}^5\sum_{\alpha_1\leq\alpha}\Big|\Big\langle \sum_{\substack{j_1+j_2\geq3   \\1\leq j_1,j_2<3}}\frac{\varGamma\left(\partial^{\alpha_1}f^{j_1,c}, \partial^{\alpha-\alpha_1}f^{j_2,c}\right)}{c^{j_1+j_2-3}}, (1+t)^{-\frac{1+\epsilon_{0}}{2}}w_{\ell_3-5}^2\partial^{\alpha}f^{3,c}\Big\rangle\Big|\nonumber\\
&\qquad\leq o(1)(1+t)^{-\frac{1+\epsilon_{0}}{2}}\|w_{\ell_3-5}\nabla_x^{5}f^{3,c}(t)\|^2_{\sigma}
+\sum_{\substack{j_1+j_2\geq3   \\1\leq j_1,j_2<3}}\frac{(1+t)^{-\frac{1+\epsilon_{0}}{2}}}{c^{2(j_1+j_2-3)}}\nonumber\\
&\qquad \times\Big(\|f^{j_1,c}\|_{H^5}^2\|f^{j_2,c}\|_{H^{5,3}_{w,\sigma}}^2
+\|f^{j_1,c}\|_{H^{5,3}_{w}}^2\|f^{j_2,c}\|_{H^{5}_{\sigma}}^2\Big)\nonumber\\
&\qquad\leq o(1)(1+t)^{-\frac{1+\epsilon_{0}}{2}}\|w_{\ell_3-5}\nabla_x^{5}f^{3,c}(t)\|^2_{\sigma}
+\sum_{j=0}^2CM\mathcal{D}^{j,c}(t).
\end{align}
Now we apply $\partial^{\alpha}$ with $|\alpha|=5$ to the equation of $f^{3,c}(t, x, p)$ in \eqref{mainF3}, take the  inner products with $(1+t)^{-\frac{1+\epsilon_{0}}{2}}w_{\ell_3-5}^2\partial^{\alpha}f^{3,c}$ in $H^{5}$, and use Lemma \ref{coerci-w}, Lemma \ref{grow-hwei}, Lemma \ref{deriv-hwei}, \eqref{lin-hwei0}, \eqref{nonlin-w6j3}, and \eqref{nonlin-w6jj} to deduce \eqref{6norm-wei}.
\end{proof}

\subsubsection{Macroscopic dissipation and electromagnetic dissipation  estimates}
In this part, we will establish the macroscopic dissipation and electromagnetic dissipation  estimates of
\eqref{mainF3}. To start, we derive the macroscopic dissipation of \eqref{mainF3}.

\noindent\underline{{\it Step 1. Macroscopic dissipation estimates.}}
We first rewrite \eqref{mainF3} as follows
\begin{align}\label{mainF3-mac}
\begin{aligned}
&\partial_t f^{3,c} + \frac{cp}{p^0}\cdot \nabla_x f^{3,c}+ E^{3,c}\cdot\frac{cp}{p^0}\sqrt{\mu^c}\zeta_0+\Big(\zeta_1 \frac{p}{p^0}\times  \sum_{j=0}^3\frac{\bar{B}^{j}}{c^j}\Big)\cdot\nabla_p f^{3,c}+\mathcal{L}f^{3,c} =\widetilde{\mathcal{Q}},
\end{aligned}
\end{align}
where $\widetilde{\mathcal{Q}}$ denote nonlinear terms in \eqref{mainF3}.
Before deducing the macroscopic dissipation estimates, we first obtain some integration estimates of $p$ for later use.

\begin{lemma}\label{p-integr}
For integers $k, m\geq0$, it holds that
\begin{align}
\int_{\mathbb R^3} |p|^{2m}\mu^c(p)\, \mathrm{d} p&=\frac{(2m+2)!}{2^{m+1}(m+1)!}\frac{K_{m+2}(c^2)}{K_{2}(c^2)},\label{pm-integr}\\
\int_{\mathbb R^3} |p|^{2m}\left(p^0\right)^{2k-1}\mu^c(p)\, \mathrm{d} p&=\sum_{i=0}^k\binom{k}{i}\frac{2[(m+i+1)]!}{2^{m+i+1}(m+i+1)!}c^{2(k-i)-1}\frac{K_{m+i+1}(c^2)}{K_{2}(c^2)}.\label{pmk-integr}
\end{align}
\end{lemma}
\begin{proof}
Recall the expression of $\mu^c(p)= \frac{\exp\left\{-cp^0 \right\} }{4 \pi
c  K_2(c^2)}$. Let $\lambda=cp^0=c\sqrt{|p|^2+c^2}$. One has
$$|p|=\frac{\sqrt{\lambda^2-c^4}}{c},\qquad |p|\,\mathrm{d}|p|=\frac{p^0}{c}\,\mathrm{d}\lambda=\frac{\lambda}{c^2}\,\mathrm{d}\lambda.$$
Then we use the definition of $K_j(\cdot)$ in \eqref{defini-kj0} to compute the integration in \eqref{pm-integr} as
\begin{align*}
\int_{\mathbb R^3} |p|^{2m}\mu^c(p)\, \mathrm{d} p&
=\int_{0}^{\infty} |p|^{2m+2}\frac{\exp\left\{-cp^0 \right\} }{
c  K_2(c^2)}\, \mathrm{d} |p|=\int_{c^2}^{\infty} \frac{\lambda(\lambda^2-c^4)^{m+1/2}\mathrm{e}^{-\lambda} }{
c^{2m+4}  K_2(c^2)}\, \mathrm{d} \lambda\\
&=\frac{(2m+2)!}{2^{m+1}(m+1)!}\frac{K_{m+2}(c^2)}{K_{2}(c^2)}.
\end{align*}
Similarly, for the integration in \eqref{pmk-integr}, we use the definition of $K_j(\cdot)$ in \eqref{defini-kj}
 to have
\begin{align*}
\int_{\mathbb R^3} |p|^{2m}\left(p^0\right)^{2k-1}\mu^c(p)\, \mathrm{d} p
&=\int_{c^2}^{\infty} \frac{(\lambda^2-c^4)^{m+1/2}\lambda^{2k}\mathrm{e}^{-\lambda} }{
c^{2m+2k+3}  K_2(c^2)}\, \mathrm{d} \lambda\\
&=\sum_{i=0}^k\binom{k}{i}\int_{c^2}^{\infty} \frac{(\lambda^2-c^4)^{m+i+1/2}c^{4(k-i)}\mathrm{e}^{-\lambda} }{
c^{2m+2k+3}  K_2(c^2)}\, \mathrm{d} \lambda\\
&=\sum_{i=0}^k\binom{k}{i}\frac{2[(m+i+1)]!}{2^{m+i+1}(m+i+1)!}c^{2(k-i)-1}\frac{K_{m+i+1}(c^2)}{K_{2}(c^2)}.
\end{align*}

\end{proof}

Now we derive local conservation laws of \eqref{mainF3-mac}.
Taking inner product of \eqref{mainF3-mac} with $[\sqrt{\mu^c},\sqrt{\mu^c}]^{\textit{t}}$, $\big[\frac{p}{C_{\mathbf{b}}}\sqrt{\mu^c}, \frac{p}{C_{\mathbf{b}}}\sqrt{\mu^c}\big]^{\textit{t}},$ and $\big[\frac{p^{0}-C_0}{C_{\mathbf{c}}}\sqrt{\mu^c}, \frac{p^{0}-C_0}{C_{\mathbf{c}}}\sqrt{\mu^c}]^{\textit{t}}$ in $L^2_p$, we can obtain
\begin{align}\label{conslaw-a}
&\partial_t \mathbf{a}^3_{\pm}+ \frac{\nabla_x\cdot \mathbf{b}^3}{C_{\mathbf{b}}}=(\widetilde{\mathcal{Q}}, \sqrt{\mu^c}),
\end{align}
\begin{align}\label{conslaw-b}
&\partial_t \mathbf{b}^3+ \frac{\nabla_x(\mathbf{a}^3_{+}+\mathbf{a}^3_{-})}{2C_{\mathbf{b}}}+\frac{\nabla_x \mathbf{c}^3}{cC_{\mathbf{b}}C_{\mathbf{c}}}+\Big(\frac{cp}{p^0}\cdot \nabla_x \big[\{I-\mathcal{P}\}_+f^{3,c}+\{I-\mathcal{P}\}_-f^{3,c}\big],\frac{p}{C_{\mathbf{b}}}\sqrt{\mu^c}\Big)\nonumber\\
&\hspace{2cm}=\big(\widetilde{\mathcal{Q}}_{+}+\widetilde{\mathcal{Q}}_{-}, \frac{p}{C_{\mathbf{b}}}\sqrt{\mu^c}\big),
\end{align}
\begin{align}\label{conslaw-c}
&\partial_t \mathbf{c}^3+ \frac{\nabla_x\cdot\mathbf{b}^3}{cC_{\mathbf{b}}C_{\mathbf{c}}}+\Big(\frac{cp}{p^0}\cdot \nabla_x \big[\{I-\mathcal{P}\}_+f^{3,c}+\{I-\mathcal{P}\}_-f^{3,c}\big],\frac{p^{0}-C_0}{C_{\mathbf{c}}}\sqrt{\mu^c}\Big)\nonumber\\
&\hspace{2cm}=\big(\widetilde{\mathcal{Q}}_{+}+\widetilde{\mathcal{Q}}_{-}, \frac{p^{0}-C_0}{C_{\mathbf{c}}}\sqrt{\mu^c}\big).
\end{align}
We also need equations of  $\nabla_x\mathbf{a}_{\pm}, \nabla_x\mathbf{b} $, and $\nabla_x\mathbf{c}$.
For $1\leq i, j\leq 3$, we take inner product of \eqref{mainF3-mac} with $\big[\frac{p^0-cC_{\mathbf{b}}^2}{C_{\mathbf{c}}}p_i\sqrt{\mu^c}, \frac{p^0-cC_{\mathbf{b}}^2}{C_{\mathbf{c}}}p_i\sqrt{\mu^c}\big]^{\textit{t}}$, $\big[\frac{p^2_i-C_{\mathbf{b}}^2}{C_{\mathbf{b}}}\sqrt{\mu^c}, \frac{p^2_i-C_{\mathbf{b}}^2}{C_{\mathbf{b}}}\sqrt{\mu^c}\big]^{\textit{t}}$, $\big[\frac{p_ip_j}{C_{\mathbf{b}}}\sqrt{\mu^c}, \frac{p_ip_j}{C_{\mathbf{b}}}\sqrt{\mu^c}\big]^{\textit{t}}$ with $i\neq j$,   and $[c(p^{0}-C_{\mathbf{a}})p_i\sqrt{\mu^c}, c(p^{0}-C_{\mathbf{a}})p_i\sqrt{\mu^c}]^{\textit{t}}$ in $L^2_p$ to have 
\begin{align}\label{macro-c}
&\partial_t\Big(f^{3,c}_{+}+f^{3,c}_{-}, \frac{p^0-cC_{\mathbf{b}}^2}{C_{\mathbf{c}}}p_i\sqrt{\mu^c}\Big) +2\rho_{\mathbf{c}}\partial_i\mathbf{c}^3+\Big(\mathcal{L}_{+}f^{3,c}+\mathcal{L}_{-}f^{3,c}, \frac{p^0-cC_{\mathbf{b}}^2}{C_{\mathbf{c}}}p_i\sqrt{\mu^c}\Big)\nonumber\\
& +\Big(\frac{cp}{p^0}\cdot \nabla_x \big[\{I-\mathcal{P}\}_+f^{3,c}+\{I-\mathcal{P}\}_-f^{3,c}\big],
\frac{p^0-cC_{\mathbf{b}}^2}{C_{\mathbf{c}}}p_i\sqrt{\mu^c}\Big)
\nonumber\\
&\hspace{2cm}=\Big(\widetilde{\mathcal{Q}}_{+}+\widetilde{\mathcal{Q}}_{-}, \frac{p^{0}-cC_{\mathbf{b}}^2}{C_{\mathbf{c}}}p_i\sqrt{\mu^c}\Big),
\end{align}
\begin{align}\label{macro-bi}
&\partial_t\Big(f^{3,c}_{+}+f^{3,c}_{-}, \frac{p^2_i-C_{\mathbf{b}}^2}{C_{\mathbf{b}}}\sqrt{\mu^c}\Big) +4\partial_i\mathbf{b}^3_i+\Big(\mathcal{L}_{+}f^{3,c}+\mathcal{L}_{-}f^{3,c}, \frac{p^2_i-C_{\mathbf{b}}^2}{C_{\mathbf{b}}}\sqrt{\mu^c}\Big)\nonumber\\
& +\Big(\frac{cp}{p^0}\cdot \nabla_x \big[\{I-\mathcal{P}\}_+f^{3,c}+\{I-\mathcal{P}\}_-f^{3,c}\big],
\frac{p^2_i-C_{\mathbf{b}}^2}{C_{\mathbf{b}}}\sqrt{\mu^c}\Big)
\nonumber\\
&\hspace{2cm}=\Big(\widetilde{\mathcal{Q}}_{+}+\widetilde{\mathcal{Q}}_{-}, \frac{p^2_i-C_{\mathbf{b}}^2}{C_{\mathbf{b}}}\sqrt{\mu^c}\Big),
\end{align}
\begin{align}\label{macro-bij}
&\partial_t\Big(f^{3,c}_{+}+f^{3,c}_{-}, \frac{p_ip_j}{C_{\mathbf{b}}}\sqrt{\mu^c}\Big) +2\big(\partial_i\mathbf{b}^3_j+\partial_j\mathbf{b}^3_i\big)+\Big(\mathcal{L}_{+}f^{3,c}+\mathcal{L}_{-}f^{3,c}, \frac{p_ip_j}{C_{\mathbf{b}}}\sqrt{\mu^c}\Big)\nonumber\\
& +\Big(\frac{cp}{p^0}\cdot \nabla_x \big[\{I-\mathcal{P}\}_+f^{3,c}+\{I-\mathcal{P}\}_-f^{3,c}\big],
\frac{p_ip_j}{C_{\mathbf{b}}}\sqrt{\mu^c}\Big)
\nonumber\\
&\hspace{2cm}=\Big(\widetilde{\mathcal{Q}}_{+}+\widetilde{\mathcal{Q}}_{-}, \frac{p_ip_j}{C_{\mathbf{b}}}\sqrt{\mu^c}\Big),
\end{align}
and
\begin{align}\label{macro-a}
&\partial_t\big(f^{3,c}_{\pm}, c(p^0-C_{\mathbf{a}})p_i\sqrt{\mu^c}\big) -\rho_{\mathbf{a}}\partial_i\mathbf{a}^3_{\pm}+\Big(\frac{cp}{p^0}\cdot \nabla_x \{I-\mathcal{P}\}_{\pm}f^{3,c}],
c(p^0-C_{\mathbf{a}})p_i\sqrt{\mu^c}\Big)
\nonumber\\
&\pm\rho_{\mathbf{a}}E^{3,c}_i \pm\Big( \Big(\frac{p}{p^0}\times  \sum_{j=0}^3\frac{\bar{B}^{j}}{c^j}\Big)\cdot \frac{\mathbf{b}}{C_{\mathbf{b}}}, c(p^0-C_{\mathbf{a}})p_i\mu^c\Big)\nonumber\\
&\mp\Big( \Big(\frac{p}{p^0}\times  \sum_{j=0}^3\frac{\bar{B}^{j}}{c^j}\Big)\cdot\nabla_p\big[c(p^0-C_{\mathbf{a}})p_i\sqrt{\mu^c}\big] , \{I-\mathcal{P}\}_{\pm}f^{3,c}\Big)
\nonumber\\
&+\big(\mathcal{L}_{\pm}f^{3,c}, c(p^0-C_{\mathbf{a}})p_i\sqrt{\mu^c}\big)=\big(\widetilde{\mathcal{Q}}_{\pm}, c(p^0-C_{\mathbf{a}})p_i\sqrt{\mu^c}\big).
\end{align}
Here the constants $\rho_{\mathbf{c}}, C_{\mathbf{a}}, \rho_{\mathbf{a}}$ are 
\begin{align*}
  \rho_{\mathbf{c}}&=\frac{C_{\mathbf{c}}^2+1/c^2}{C_{\mathbf{c}}^2},\qquad \rho_{\mathbf{a}}=-c\Big(c\frac{K_3(c^2)}{K_2(c^2)}-C_{\mathbf{a}}\Big), \qquad
 C_{\mathbf{a}}= c\Big(-c^2\frac{K^2_3(c^2)}{K^2_2(c^2)}+6\frac{K_3(c^2)}{K_2(c^2)}+c^2\Big).
\end{align*}
From \eqref{transform} and \eqref{acurate} in Lemma \ref{K01p}, it holds that
\begin{align*}
  \rho_{\mathbf{c}}&= \frac{5}{3}+\frac{O(1)}{c^{2}}, \qquad
  \rho_{\mathbf{a}}=\frac{5}{2}+\frac{O(1)}{c^{2}},\qquad  C_{\mathbf{a}}=c+\frac{5}{c}+\frac{O(1)}{c^{3}},\\
  \frac{|p^0-cC_{\mathbf{b}}^2|}{C_{\mathbf{c}}}&=\frac{|p^0-c-c\big(K_3(c^2)/K_2(c^2)-1\big)|}
  {\sqrt{3/(2c^2)+O(1)/c^4}}\leq\frac{|p|^2+5}{2}+\frac{O(1)}{c},\\
  |c(p^0-C_{\mathbf{a}})|&=c\Big|p^0-c-\frac{5}{c}-\frac{O(1)}{c^2}\Big|\leq\frac{|p|^2}{2}+5+\frac{O(1)}{c}.
\end{align*}
Now we come to the explicit macroscopic dissipation estimation of $f^{3,c}$.

\noindent $\bullet$ Case 1 (Estimation of $\big\|\nabla_x\mathbf{c}^3\big\|$):
We multiply \eqref{macro-c} by $\partial_i\mathbf{c}^3$ and use Lemma \ref{w-nonlin} with $\vartheta=\ell=0$ to have
\begin{align*}
&\sum_{i=1}^3\Big\langle\partial_t\Big(f^{3,c}_{+}+f^{3,c}_{-}, \frac{p^0-cC_{\mathbf{b}}^2}{C_{\mathbf{c}}}p_i\sqrt{\mu^c}\Big), \partial_i\mathbf{c}^3\Big\rangle +2\rho_{\mathbf{c}}\big\|\nabla_x\mathbf{c}^3\big\|^2\\
&\qquad=-\sum_{i=1}^3\Big\langle\Big(\mathcal{L}_{+}f^{3,c}+\mathcal{L}_{-}f^{3,c}, \frac{p^0-cC_{\mathbf{b}}^2}{C_{\mathbf{c}}}p_i\sqrt{\mu^c}\Big),\partial_i\mathbf{c}^3\Big\rangle \nonumber\\
&\qquad -\sum_{i=1}^3\Big\langle\Big(\frac{cp}{p^0}\cdot \nabla_x \big[\{I-\mathcal{P}\}_+f^{3,c}+\{I-\mathcal{P}\}_-f^{3,c}\big],
\frac{p^0-cC_{\mathbf{b}}^2}{C_{\mathbf{c}}}p_i\sqrt{\mu^c}\Big), \partial_i\mathbf{c}^3\Big\rangle
\nonumber\\
&\qquad+\sum_{i=1}^3\Big\langle\Big(\widetilde{\mathcal{Q}}_{+}+\widetilde{\mathcal{Q}}_{-}, \frac{p^{0}-cC_{\mathbf{b}}^2}{C_{\mathbf{c}}}p_i\sqrt{\mu^c}\Big), \partial_i\mathbf{c}^3\Big\rangle\\
&\qquad \leq o(1) \big\|\nabla_x\mathbf{c}^3\big\|^2+\left|\{I-\mathcal{P}\}f^{3,c}\right|_{H^1_{\sigma}}^2+CM\sum_{j=0}^3
\mathcal{D}^{j,c}_2,
\end{align*}
where $C$ is a uniform constant independent of $c$.  
Note that 
\begin{align*}
   & \Big(f^{3,c}_{+}+f^{3,c}_{-}, \frac{p^0-cC_{\mathbf{b}}^2}{C_{\mathbf{c}}}p_i\sqrt{\mu^c}\Big)
   =\Big(\frac{5}{\sqrt{6}c^2}+\frac{O(1)}{c^3}\Big)\mathbf{b}_i^3\\
 &\qquad+\Big(\big[\{I-\mathcal{P}\}_+f^{3,c}+\{I-\mathcal{P}\}_-f^{3,c}\big],
\frac{p^0-cC_{\mathbf{b}}^2}{C_{\mathbf{c}}}p_i\sqrt{\mu^c}\Big),  
    \end{align*}
and
\begin{align*}
   &\sum_{i=1}^3\Big\langle\partial_t\Big(f^{3,c}_{+}+f^{3,c}_{-}, \frac{p^0-cC_{\mathbf{b}}^2}{C_{\mathbf{c}}}p_i\sqrt{\mu^c}\Big), \partial_i\mathbf{c}^3\Big\rangle \\
&\qquad=\sum_{i=1}^3\frac{\mathrm{d}}{\mathrm{d}t}\Big(f^{3,c}_{+}+f^{3,c}_{-}, \frac{p^0-cC_{\mathbf{b}}^2}{C_{\mathbf{c}}}p_i\sqrt{\mu^c}\Big), \partial_i\mathbf{c}^3\Big\rangle+\sum_{i=1}^3\Big(\partial_i[f^{3,c}_{+}+f^{3,c}_{-}], \frac{p^0-cC_{\mathbf{b}}^2}{C_{\mathbf{c}}}p_i\sqrt{\mu^c}\Big), \partial_t\mathbf{c}^3\Big\rangle\\
&\qquad\leq\sum_{i=1}^3\frac{\mathrm{d}}{\mathrm{d}t}\Big(f^{3,c}_{+}+f^{3,c}_{-}, \frac{p^0-cC_{\mathbf{b}}^2}{C_{\mathbf{c}}}p_i\sqrt{\mu^c}\Big), \partial_i\mathbf{c}^3\Big\rangle
+C\Big(\|\nabla_x\{I-\mathcal{P}\}f^{3,c}\|_{\sigma}+\frac{\|\nabla_x \mathbf{b}^3\|}{c^2}\Big)\\
&\qquad \times \Big(\|\nabla_x\{I-\mathcal{P}\}f^{3,c}\|_{\sigma}+\|\nabla_x \mathbf{b}^3\|+\sqrt{M}\sum_{j=0}^3
\sqrt{\mathcal{D}^{j,c}}_2\Big)  
    \end{align*}
by \eqref{conslaw-c} for some $c-$independent constant $C>0$. Then we can further obtain
\begin{align}\label{dissi-c}
   &\sum_{i=1}^3\frac{\mathrm{d}}{\mathrm{d}t}\Big\langle\Big(f^{3,c}_{+}+f^{3,c}_{-}, \frac{p^0-cC_{\mathbf{b}}^2}{C_{\mathbf{c}}}p_i\sqrt{\mu^c}\Big), \partial_i\mathbf{c}^3\Big\rangle+\rho_{\mathbf{c}}\big\|\nabla_x\mathbf{c}^3\big\|^2\nonumber\\
&\qquad\leq C\Big(\|\{I-\mathcal{P}\}f^{3,c}\|_{H^1_{\sigma}}^2+\frac{\|\nabla_x \mathbf{b}^3\|^2}{c^2}+M\sum_{j=0}^3
\mathcal{D}^{j,c}_2\Big) . 
    \end{align}

\noindent $\bullet$ Case 2 (Estimation of $\big\|\nabla_x\mathbf{b}^3\big\|$): 
Note that
\begin{align*}
&4\sum_{i=1}^3\big\|\partial_i\mathbf{b}^3_i\big\|^2+\sum_{i\neq j}2\Big(\big\|\partial_i\mathbf{b}^3_j\big\|^2+\int_{{\mathbb T}^3} \partial_i\mathbf{b}^3_i \partial_j\mathbf{b}^3_j\, \mathrm{d} x\Big)\\
  &\qquad =4\sum_{i=1}^3\big\|\partial_i\mathbf{b}^3_i\big\|^2+\sum_{i\neq j}2\int_{{\mathbb T}^3} \partial_i\mathbf{b}^3_j\big(\partial_i\mathbf{b}^3_j+ \partial_j\mathbf{b}^3_i\big)\, \mathrm{d} x\Big) \\
  &\qquad=2\left(\big\|\nabla_x\mathbf{b}^3\big\|^2+\big\|\nabla_x\cdot\mathbf{b}^3\big\|^2\right)  
\end{align*}
Multiplying \eqref{macro-bi} and \eqref{macro-bij} by $\partial_i\mathbf{b}^3_i$ and $\partial_i\mathbf{b}^3_j$ respectively, we can use Lemma \ref{w-nonlin} with $\vartheta=\ell=0$ to obtain
\begin{align*}
&\Big\langle\partial_t[f^{3,c}_{+}+f^{3,c}_{-}], \sum_{i=1}^3\partial_i\mathbf{b}^3_i\frac{p^2_i-C_{\mathbf{b}}^2}{C_{\mathbf{b}}}\sqrt{\mu^c}
+\sum_{i\neq j}\partial_i\mathbf{b}^3_j\frac{p_ip_j}{C_{\mathbf{b}}}\sqrt{\mu^c}\Big\rangle\\
  &\qquad+2\left(\big\|\nabla_x\mathbf{b}^3\big\|^2+\big\|\nabla_x\cdot\mathbf{b}^3\big\|^2\right)  
\nonumber\\
&\qquad \leq o(1) \big\|\nabla_x\mathbf{b}^3\big\|^2+C\Big(\left\|\{I-\mathcal{P}\}f^{3,c}\right\|_{H^1_{\sigma}}^2+M\sum_{j=0}^3
\mathcal{D}^{j,c}_2\Big).
\end{align*}
where $C$ is a uniform constant independent of $c$.  
Note that 
\begin{align*}
   & \int_{\mathbb{R}^3} \frac{(p_i^2-C_{\mathbf{b}}^2)(p^0-C_{0})}{C_{\mathbf{c}}}\mu^c(p)\,\mathrm{d} p
   =\int_{\mathbb{R}^3} \frac{(p_i^2-C_{\mathbf{b}}^2)p^0}{C_{\mathbf{c}}}\mu^c(p)\,\mathrm{d} p
   =\frac{C_{\mathbf{a}}}{c^2C_{\mathbf{c}}}=\frac{\sqrt{6}}{3}+\frac{O(1)}{c^2}.  
    \end{align*}
As the derivation of \eqref{dissi-c}, we use \eqref{conslaw-b} to have
\begin{align}\label{dissi-b}
   &\frac{\mathrm{d}}{\mathrm{d}t}\Big\langle[f^{3,c}_{+}+f^{3,c}_{-}], \sum_{i=1}^3\partial_i\mathbf{b}^3_i\frac{p^2_i-C_{\mathbf{b}}^2}{C_{\mathbf{b}}}\sqrt{\mu^c}
+\sum_{i\neq j}\partial_i\mathbf{b}^3_j\frac{p_ip_j}{C_{\mathbf{b}}}\sqrt{\mu^c}\Big\rangle+\big\|\nabla_x\mathbf{b}^3\big\|^2+\big\|\nabla_x\cdot\mathbf{b}^3\big\|^2  
\nonumber\\
&\qquad\leq o(1) \big\|\nabla_x(\mathbf{a}^3_{+}+\mathbf{a}^3_{-})\big\|^2+ C\Big(\|\{I-\mathcal{P}\}f^{3,c}\|_{H^1_{\sigma}}^2+\|\nabla_x \mathbf{c}^3\|^2+M\sum_{j=0}^3
\mathcal{D}^{j,c}_2\Big), 
    \end{align}
where  $C>0$ is a $c-$independent constant.

\noindent $\bullet$ Case 3 (Estimation of $\big\|\nabla_x\mathbf{a}^3_{\pm}\big\|$): 
From \eqref{macro-a}, one has
\begin{align}\label{macro-aa}
&\partial_t\big([f^{3,c}_{+}+f^{3,c}_{-}], c(p^0-C_{\mathbf{a}})p_i\sqrt{\mu^c}\big) -\rho_{\mathbf{a}}\partial_i\big(\mathbf{a}^3_{+}+\mathbf{a}^3_{-}\big)\nonumber\\
&+\Big(\frac{cp}{p^0}\cdot \nabla_x \big[\{I-\mathcal{P}\}_+f^{3,c}+\{I-\mathcal{P}\}_-f^{3,c}\big],
c(p^0-C_{\mathbf{a}})p_i\sqrt{\mu^c}\Big)
\nonumber\\
&+\big(\mathcal{L}_{+}f^{3,c}+\mathcal{L}_{-}f^{3,c}, c(p^0-C_{\mathbf{a}})p_i\sqrt{\mu^c}\big)=\big(\widetilde{\mathcal{Q}}_{+}+\widetilde{\mathcal{Q}}_{-}, c(p^0-C_{\mathbf{a}})p_i\sqrt{\mu^c}\big),
\end{align}
and
\begin{align}\label{macro-as}
&\partial_t\big([f^{3,c}_{+}-f^{3,c}_{-}], c(p^0-C_{\mathbf{a}})p_i\sqrt{\mu^c}\big) -\rho_{\mathbf{a}}\partial_i\big(\mathbf{a}^3_{+}-\mathbf{a}^3_{-}\big)\nonumber\\
&+\Big(\frac{cp}{p^0}\cdot \nabla_x \big[\{I-\mathcal{P}\}_+f^{3,c}-\{I-\mathcal{P}\}_-f^{3,c}\big],
c(p^0-C_{\mathbf{a}})p_i\sqrt{\mu^c}\Big)
\nonumber\\
&+2\rho_{\mathbf{a}}E^{3,c}_i +2\Big( \Big(\frac{p}{p^0}\times  \sum_{j=0}^3\frac{\bar{B}^{j}}{c^j}\Big)\cdot \frac{\mathbf{b}^3}{C_{\mathbf{b}}}, c(p^0-C_{\mathbf{a}})p_i\mu^c\Big)\nonumber\\
&+2\Big( \Big(\frac{p}{p^0}\times  \sum_{j=0}^3\frac{\bar{B}^{j}}{c^j}\Big)\cdot\nabla_p\big[c(p^0-C_{\mathbf{a}})p_i\sqrt{\mu^c}\big] , \{I-\mathcal{P}\}_{\pm}f^{3,c}\Big)
\nonumber\\
&+\big(\mathcal{L}_{+}f^{3,c}-\mathcal{L}_{-}f^{3,c}, c(p^0-C_{\mathbf{a}})p_i\sqrt{\mu^c}\big)=\big(\widetilde{\mathcal{Q}}_{+}-\widetilde{\mathcal{Q}}_{-}, c(p^0-C_{\mathbf{a}})p_i\sqrt{\mu^c}\big).
\end{align}
We take the inner product of \eqref{macro-aa} and $\partial_i\big(\mathbf{a}^3_{+}+\mathbf{a}^3_{-}\big)$ in $L^2_x$, and use Lemma \ref{w-nonlin} with $\vartheta=\ell=0$ to obtain
\begin{align*}
&-\sum_{i=1}^3\big\langle\big(\partial_t[f^{3,c}_{+}+f^{3,c}_{-}], c(p^0-C_{\mathbf{a}})p_i\sqrt{\mu^c}\big),\partial_i\big(\mathbf{a}^3_{+}+\mathbf{a}^3_{-}\big)\big\rangle 
+\rho_{\mathbf{a}}\big\|\nabla_x\big(\mathbf{a}^3_{+}+\mathbf{a}^3_{-}\big)\big\|^2\\
  &\qquad \leq o(1) \big\|\nabla_x\big(\mathbf{a}^3_{+}+\mathbf{a}^3_{-}\big)\big\|^2
+C\Big(\left\|\{I-\mathcal{P}\}f^{3,c}\right\|_{H^1_{\sigma}}^2+M\sum_{j=0}^3
\mathcal{D}^{j,c}_2\Big).
\end{align*}
where $C$ is a uniform constant independent of $c$.  
Note that from Lemma \ref{p-integr},
\begin{align*}
   & \int_{\mathbb{R}^3} c(p^0-C_{\mathbf{a}})p_i^2\mu^c(p)\,\mathrm{d} p
   =c\Big[c+\frac{5K_3(c^2)}{cK_2(c^2)}-\frac{K_3(c^2)}{K_2(c^2)}C_{\mathbf{a}}\Big] =-\frac{5}{2}+\frac{O(1)}{c^2}.  
    \end{align*}
As the derivation of \eqref{dissi-c}, we use \eqref{conslaw-a} to have
\begin{align}\label{dissi-aa}
   &-\sum_{i=1}^3\frac{\mathrm{d}}{\mathrm{d}t}
   \big\langle\big(f^{3,c}_{+}+f^{3,c}_{-}, c(p^0-C_{\mathbf{a}})p_i\sqrt{\mu^c}\big),\partial_i\big(\mathbf{a}^3_{+}+\mathbf{a}^3_{-}\big)\big\rangle 
+\frac{\rho_{\mathbf{a}}}{2}\big\|\nabla_x\big(\mathbf{a}^3_{+}+\mathbf{a}^3_{-}\big)\big\|^2
\nonumber\\
&\qquad\leq C\Big(\|\{I-\mathcal{P}\}f^{3,c}\|_{H^1_{\sigma}}^2+\|\nabla_x \cdot\mathbf{b}^3\|^2+M\sum_{j=0}^3
\mathcal{D}^{j,c}_2\Big), 
    \end{align}
where  $C>0$ is a $c-$independent constant.

On the other hand, note the facts 
\begin{align*}
  &\nabla_x \cdot E^{3,c}=\mathbf{a}^3_{+}-\mathbf{a}^3_{-},\qquad -\int_{\mathbb{T}^3}\nabla_x\big(\mathbf{a}^3_{+}-\mathbf{a}^3_{-}\big)\cdot E^{3,c}\,\mathrm{d}x=\|\mathbf{a}^3_{+}-\mathbf{a}^3_{-}\|^2, \\
  & \Big( \Big|\frac{p}{p^0}\times  \sum_{j=0}^3\frac{\bar{B}^{j}}{c^j}\Big| , c(p^0-C_{\mathbf{a}})|p|\mu^c\Big)=\frac{O(1)}{c}.
\end{align*}
From \eqref{macro-as}, as the derivation of \eqref{dissi-aa}, we can obtain 
\begin{align}\label{dissi-as}
   &-\sum_{i=1}^3\frac{\mathrm{d}}{\mathrm{d}t}
   \big\langle\big(f^{3,c}_{+}-f^{3,c}_{-}, c(p^0-C_{\mathbf{a}})p_i\sqrt{\mu^c}\big),\partial_i\big(\mathbf{a}^3_{+}-\mathbf{a}^3_{-}\big)\big\rangle\nonumber\\
&\qquad+\frac{\rho_{\mathbf{a}}}{2}\big\|\nabla_x\big(\mathbf{a}^3_{+}-\mathbf{a}^3_{-}\big)\big\|^2
+\rho_{\mathbf{a}}\big\|\mathbf{a}^3_{+}-\mathbf{a}^3_{-}\big\|^2
\nonumber\\
&\qquad\leq C\Big(\|\{I-\mathcal{P}\}f^{3,c}\|_{H^1_{\sigma}}^2+\|\nabla_x \mathbf{b}^3\|^2+M\sum_{j=0}^3
\mathcal{D}^{j,c}_2\Big), 
    \end{align}
where  $C>0$ is a $c-$independent constant.

Now we choose constants $\kappa_1$, $o(1)$ in \eqref{dissi-b} small enough such that 
$$C\kappa_1\leq \frac{1}{2},\qquad o(1)\leq \frac{\kappa_1\rho_{\mathbf{a}}}{4}.$$
Then we make a linear combination of \eqref{dissi-aa} with \eqref{dissi-b} to have
\begin{align}\label{dissicom-ab}
   &\frac{\mathrm{d}}{\mathrm{d}t}\Big(\Big\langle[f^{3,c}_{+}+f^{3,c}_{-}], \sum_{i=1}^3\partial_i\mathbf{b}^3_i\frac{p^2_i-C_{\mathbf{b}}^2}{C_{\mathbf{b}}}\sqrt{\mu^c}
+\sum_{i\neq j}\partial_i\mathbf{b}^3_j\frac{p_ip_j}{C_{\mathbf{b}}}\sqrt{\mu^c}\Big\rangle\nonumber\\
&-\kappa_1\sum_{i=1}^3
   \big\langle\big(f^{3,c}_{+}+f^{3,c}_{-}, c(p^0-C_{\mathbf{a}})p_i\sqrt{\mu^c}\big),\partial_i\big(\mathbf{a}^3_{+}
   +\mathbf{a}^3_{-}\big)\big\rangle\Big)\nonumber\\ 
&+\frac{\kappa_1\rho_{\mathbf{a}}}{4}\big\|\nabla_x\big(\mathbf{a}^3_{+}+\mathbf{a}^3_{-}\big)\big\|^2
+\big\|\nabla_x\mathbf{b}^3\big\|^2+\frac{1}{2}\big\|\nabla_x\cdot\mathbf{b}^3\big\|^2  
\nonumber\\
&\qquad\leq C_1\Big(\|\{I-\mathcal{P}\}f^{3,c}\|_{H^1_{\sigma}}^2+\|\nabla_x \mathbf{c}^3\|^2+M\sum_{j=0}^3
\mathcal{D}^{j,c}\Big), 
    \end{align}
where $C_1$ is a large $c-$independent constant. We further choose $\kappa_2$ small and $c_3$ large such that for $c\geq c_3$,
$$\kappa_2 C_1\leq \frac{\rho_{\mathbf{c}}}{2}, \qquad\frac{C}{c^2}\leq \frac{\kappa_2}{4}.$$
Then we make a linear combination of \eqref{dissicom-ab} with \eqref{dissi-c} to get
 \begin{align}\label{dissicom-abc}
   &\frac{\mathrm{d}}{\mathrm{d}t}\Big[\kappa_2\Big(\Big\langle[f^{3,c}_{+}+f^{3,c}_{-}], \sum_{i=1}^3\partial_i\mathbf{b}^3_i\frac{p^2_i-C_{\mathbf{b}}^2}{C_{\mathbf{b}}}\sqrt{\mu^c}
+\sum_{i\neq j}\partial_i\mathbf{b}^3_j\frac{p_ip_j}{C_{\mathbf{b}}}\sqrt{\mu^c}\Big\rangle\nonumber\\
&-\kappa_1\sum_{i=1}^3
   \big\langle\big(\partial_t[f^{3,c}_{+}+f^{3,c}_{-}], c(p^0-C_{\mathbf{a}})p_i\sqrt{\mu^c}\big),\partial_i\big(\mathbf{a}^3_{+}
   +\mathbf{a}^3_{-}\big)\big\rangle\Big)\nonumber\\
   &+\sum_{i=1}^3\Big\langle\Big(f^{3,c}_{+}+f^{3,c}_{-}, \frac{p^0-cC_{\mathbf{b}}^2}{C_{\mathbf{c}}}p_i\sqrt{\mu^c}\Big), \partial_i\mathbf{c}^3\Big\rangle\Big]\nonumber\\
   & +\frac{\kappa_1\kappa_2\rho_{\mathbf{a}}}{4}\big\|\nabla_x\big(\mathbf{a}^3_{+}+\mathbf{a}^3_{-}\big)\big\|^2
+\frac{3\kappa_2}{4}\big\|\nabla_x\mathbf{b}^3\big\|^2+\frac{\kappa_2}{2}\big\|\nabla_x\cdot\mathbf{b}^3\big\|^2
+\frac{\rho_{\mathbf{c}}}{2}\big\|\nabla_x\mathbf{c}^3\big\|^2  
\nonumber\\
&\qquad\leq C_2\Big(\|\{I-\mathcal{P}\}f^{3,c}\|_{H^1_{\sigma}}+M\sum_{j=0}^3
\mathcal{D}^{j,c}_2\Big), 
    \end{align}
where $C_2$ is a large $c-$independent constant.
Finally, we multiply \eqref{dissi-as} by $\frac{\kappa_2}{4C}$ and add the resultant to \eqref{dissicom-abc} to have
\begin{align}\label{dissicom-aabc0}
   &\frac{\mathrm{d}}{\mathrm{d}t}\mathcal{E}^{3,c,1}_{\mathbf{abc}}(t)+
   \frac{\kappa_2}{4C}\Big[\frac{\rho_{\mathbf{a}}}{2}\big\|\nabla_x\big(\mathbf{a}^3_{+}-\mathbf{a}^3_{-}\big)\big\|^2
+\rho_{\mathbf{a}}\big\|\mathbf{a}^3_{+}-\mathbf{a}^3_{-}\big\|^2\Big]\nonumber\\
   & +\frac{\kappa_1\kappa_2\rho_{\mathbf{a}}}{4}\big\|\nabla_x\big(\mathbf{a}^3_{+}+\mathbf{a}^3_{-}\big)\big\|^2
+\frac{\kappa_2}{2}\big\|\nabla_x\mathbf{b}^3\big\|^2+\frac{\kappa_2}{2\widehat{}}\big\|\nabla_x\cdot\mathbf{b}^3\big\|^2
+\frac{\rho_{\mathbf{c}}}{2}\big\|\nabla_x\mathbf{c}^3\big\|^2  
\nonumber\\
&\qquad\leq C_{3}\Big(\|\{I-\mathcal{P}\}f^{3,c}\|_{H^1_{\sigma}}^2+M\sum_{j=0}^3
\mathcal{D}^{j,c}_2\Big), 
    \end{align}
where $C_3$ is a $c-$independent constant, and 
\begin{align*}
\mathcal{E}^{3,c,1}_{\mathbf{abc}}(t)=&\kappa_2\Big(\Big\langle[f^{3,c}_{+}+f^{3,c}_{-}], \sum_{i=1}^3\partial_i\mathbf{b}^3_i\frac{p^2_i-C_{\mathbf{b}}^2}{C_{\mathbf{b}}}\sqrt{\mu^c}
+\sum_{i\neq j}\partial_i\mathbf{b}^3_j\frac{p_ip_j}{C_{\mathbf{b}}}\sqrt{\mu^c}\Big\rangle\nonumber\\
&-\kappa_1\sum_{i=1}^3
   \big\langle\big(f^{3,c}_{+}+f^{3,c}_{-}, c(p^0-C_{\mathbf{a}})p_i\sqrt{\mu^c}\big),\partial_i\big(\mathbf{a}^3_{+}
   +\mathbf{a}^3_{-}\big)\big\rangle\Big)\nonumber\\
   &+\sum_{i=1}^3\Big\langle\Big(f^{3,c}_{+}+f^{3,c}_{-}, \frac{p^0-cC_{\mathbf{b}}^2}{C_{\mathbf{c}}}p_i\sqrt{\mu^c}\Big), \partial_i\mathbf{c}^3\Big\rangle\nonumber\\
   &+\frac{1}{4C}\sum_{i=1}^3
   \big\langle\big(f^{3,c}_{+}-f^{3,c}_{-}, c(p^0-C_{\mathbf{a}})p_i\sqrt{\mu^c}\big),\partial_i\big(\mathbf{a}^3_{+}-\mathbf{a}^3_{-}\big)\big\rangle.
\end{align*}

In the same way as the derivation of \eqref{dissicom-aabc0}, we can obtain that for $2\leq m\leq 5$,
\begin{align}\label{dissicom-aabc}
   &\frac{\mathrm{d}}{\mathrm{d}t}\mathcal{E}^{3,c,m}_{\mathbf{abc}}(t)
   +\frac{\kappa_2}{4C}\Big[\frac{\rho_{\mathbf{a}}}{2}\big\|\nabla_x\big(\mathbf{a}^3_{+}-\mathbf{a}^3_{-}\big)\big\|^2_{H^{m-1}}
+\rho_{\mathbf{a}}\big\|\mathbf{a}^3_{+}-\mathbf{a}^3_{-}\big\|^2_{H^{m-1}}\Big]\nonumber\\
   & +\frac{\kappa_1\kappa_2\rho_{\mathbf{a}}}{4}\big\|\nabla_x\big(\mathbf{a}^3_{+}+\mathbf{a}^3_{-}\big)\big\|^2_{H^{m-1}}
+\frac{\kappa_2}{2}\big\|\nabla_x\mathbf{b}^3\big\|^2_{H^4}+\frac{\kappa_2}{2}\big\|\nabla_x\cdot\mathbf{b}^3\big\|^2_{H^{m-1}}
+\frac{\rho_{\mathbf{c}}}{2}\big\|\nabla_x\mathbf{c}^3\big\|^2_{H^{m-1}}  
\nonumber\\
&\qquad\leq C_{m+2}
\Big(\|\{I-\mathcal{P}\}f^{3,c}\|_{H^m_{\sigma}}^2+M\sum_{j=0}^3 \mathcal{D}^{j,c}_m\Big), 
    \end{align}
where $C_{m+2}$ are $c-$independent constants, and 
\begin{align*}
\mathcal{E}^{3,c,m}_{\mathbf{abc}}(t)=&\sum_{|\alpha|=0}^{m-1}\Big[\kappa_2\Big(\Big\langle\partial^{\alpha}[f^{3,c}_{+}+f^{3,c}_{-}], \sum_{i=1}^3\partial_i\partial^{\alpha}\mathbf{b}^3_i\frac{p^2_i-C_{\mathbf{b}}^2}{C_{\mathbf{b}}}\sqrt{\mu^c}
+\sum_{i\neq j}\partial_i\partial^{\alpha}\mathbf{b}^3_j\frac{p_ip_j}{C_{\mathbf{b}}}\sqrt{\mu^c}\Big\rangle\nonumber\\
&-\kappa_1\sum_{i=1}^3
   \big\langle\big(\partial^{\alpha}[f^{3,c}_{+}+f^{3,c}_{-}], c(p^0-C_{\mathbf{a}})p_i\sqrt{\mu^c}\big),\partial_i\partial^{\alpha}\big(\mathbf{a}^3_{+}
   +\mathbf{a}^3_{-}\big)\big\rangle\Big)\nonumber\\
   &+\sum_{i=1}^3\Big\langle\Big(\partial^{\alpha}[f^{3,c}_{+}+f^{3,c}_{-}], \frac{p^0-cC_{\mathbf{b}}^2}{C_{\mathbf{c}}}p_i\sqrt{\mu^c}\Big), \partial_i\partial^{\alpha}\mathbf{c}^3\Big\rangle\nonumber\\
   &+\frac{1}{4C}\sum_{i=1}^3
   \big\langle\big(\partial^{\alpha}[f^{3,c}_{+}-f^{3,c}_{-}], c(p^0-C_{\mathbf{a}})p_i\sqrt{\mu^c}\big),\partial_i\partial^{\alpha}\big(\mathbf{a}^3_{+}-\mathbf{a}^3_{-}\big)\big\rangle\Big].
\end{align*}
 
Then in Step 1, we have proved the following proposition.
\begin{proposition}\label{abc-26} Under the assumptions in Lemma \ref{grow-low}, we further assume that $c\geq c_3$. Then for $2\leq m\leq 5$, there exist $c-$independent constant $C_{0m}>0$ and a functional $\mathcal{E}^{3,c,m}_{\mathbf{abc}}(t)$,  which satisfy 
$$|\mathcal{E}^{3,c,m}_{\mathbf{abc}}(t)|\leq C_{0m} \|f^{3,c}(t)\|_{H^m}^2 $$
 such that \eqref{dissicom-aabc} holds.
\end{proposition}

\noindent\underline{{\it Step 2. Electromagnetic dissipation estimates.}} We will derive the dissipation estimates of the electric field $E^{3,c}(t,x)$ and the magnetic field $B^{3,c}(t,x)$ separately.

\noindent $\bullet$ Case 1 (Dissipation estimates of the electric field): For $1\leq m\leq 5$, we apply $\partial^{\alpha}$ to \eqref{macro-a}, make a difference between the corresponding two equations, and take inner product of the resultant with $c^{-2}\partial^{\alpha}E^{3,c}(t,x)$ to obtain
\begin{align*}
&-\sum_{|\alpha|=0}^{m}\sum_{i=1}^3
   \big\langle\partial^{\alpha}\partial_t\big(f^{3,c}_{+}-f^{3,c}_{-}, c(p^0-C_{\mathbf{a}})p_i\sqrt{\mu^c}\big),\frac{\partial^{\alpha}E^{3,c}}{c^2}\big\rangle
   +\frac{\rho_{\mathbf{a}}}{c^2}\big\|\mathbf{a}^3_{+}-\mathbf{a}^3_{-}\big\|^2_{H^m}
   +\frac{2\rho_{\mathbf{a}}}{c^2}\Big\|E^{3,c}\big\|^2_{H^m}
\nonumber\\
&\qquad\leq\sum_{|\alpha|=0}^{m}\sum_{i=1}^3 2\Big|\Big\langle\Big( \Big(\frac{p}{p^0}\times  \sum_{j=0}^3\frac{\bar{B}^{j}}{c^j}\Big)\cdot \frac{\partial^{\alpha}\mathbf{b}}{C_{\mathbf{b}}}, c(p^0-C_{\mathbf{a}})p_i\mu^c\Big),\frac{\partial^{\alpha}E^{3,c}}{c^2}\Big\rangle\Big|\nonumber\\
&\qquad+\sum_{|\alpha|=0}^{m}\sum_{i=1}^3 2\Big|\Big\langle\Big( \Big(\frac{p}{p^0}\times  \sum_{j=0}^3\frac{\bar{B}^{j}}{c^j}\Big)\cdot\nabla_p\big[c(p^0-C_{\mathbf{a}})p_i\sqrt{\mu^c}\big] , \\ &\hspace{3cm}\partial^{\alpha}\big[\{I-\mathcal{P}\}_+f^{3,c}-\{I-\mathcal{P}\}_-f^{3,c}\big]\Big),
\frac{\partial^{\alpha}E^{3,c}}{c^2}\Big\rangle\Big|
\nonumber\\
&\qquad+\sum_{|\alpha|=0}^{m}\sum_{i=1}^3\big|\big\langle\big(\mathcal{L}_{+}\partial^{\alpha}f^{3,c}
-\mathcal{L}_{-}\partial^{\alpha}f^{3,c}, c(p^0-C_{\mathbf{a}})p_i\sqrt{\mu^c}\big), \frac{\partial^{\alpha}E^{3,c}}{c^2}\big\rangle\big|\\
&\qquad+\sum_{|\alpha|=0}^{m}\sum_{i=1}^3\big|\big\langle
\big(\partial^{\alpha}\widetilde{\mathcal{Q}}_{+}-\partial^{\alpha}\widetilde{\mathcal{Q}}_{-}, c(p^0-C_{\mathbf{a}})p_i\sqrt{\mu^c}\big), \frac{\partial^{\alpha}E^{3,c}}{c^2}\big\rangle\big|\\
&\qquad\leq \frac{o(1)}{c^2}\Big\|E^{3,c}\big\|^2_{H^m}+
\frac{C}{c^2}\Big(\|\{I-\mathcal{P}\}f^{3,c}\|_{H^{m+1}_{\sigma}}^2+M\sum_{j=0}^3 \mathcal{D}^{j,c}_{\max\{2,m\}}\Big).
\end{align*}
On the other hand, from the equation of $E^{3,c}(t,x)$ in \eqref{mainF3}, we get
\begin{align*}
&\sum_{|\alpha|=0}^{m}\sum_{i=1}^3
   \big\langle\partial^{\alpha}\partial_t\big(f^{3,c}_{+}-f^{3,c}_{-}, c(p^0-C_{\mathbf{a}})p_i\sqrt{\mu^c}\big),\frac{\partial^{\alpha}E^{3,c}}{c^2}\big\rangle\\
   &\qquad =\sum_{|\alpha|=0}^{m}\sum_{i=1}^3\frac{\mathrm{d}}{\mathrm{d}t}
   \big\langle\partial^{\alpha}\big(f^{3,c}_{+}-f^{3,c}_{-}, c(p^0-C_{\mathbf{a}})p_i\sqrt{\mu^c}\big),\frac{\partial^{\alpha}E^{3,c}}{c^2}\big\rangle\\
   &\qquad -\sum_{|\alpha|=0}^{m}\sum_{i=1}^3
   \big\langle\partial^{\alpha}\big(f^{3,c}_{+}-f^{3,c}_{-}, c(p^0-C_{\mathbf{a}})p_i\sqrt{\mu^c}\big),\frac{\nabla_x \times \partial^{\alpha}B^{3,c}}{c}\big\rangle\\
   &\qquad +\sum_{|\alpha|=0}^{m}\sum_{i=1}^3
   \big\langle\partial^{\alpha}\big(f^{3,c}_{+}-f^{3,c}_{-}, c(p^0-C_{\mathbf{a}})p_i\sqrt{\mu^c}\big),\int_{\mathbb R^3}\frac{p}{cp^0}\sqrt{\mu^c}\partial^{\alpha}\left(f^{3,c}_{+}-f^{3,c}_-\right)\, \mathrm{d} p\big\rangle\\
   &\qquad \leq\sum_{|\alpha|=0}^{m}\sum_{i=1}^3\frac{\mathrm{d}}{\mathrm{d}t}
   \big\langle\partial^{\alpha}\big(f^{3,c}_{+}-f^{3,c}_{-}, c(p^0-C_{\mathbf{a}})p_i\sqrt{\mu^c}\big),\frac{\partial^{\alpha}E^{3,c}}{c^2}\big\rangle\\
   &\qquad +\frac{o(1)}{c^2}\|\nabla_xB^{3,c}\|^2_{H^{\max\{0,m-1\}}}+C\Big(\|\{I-\mathcal{P}\}f^{3,c}\|^2_{H^{m+1}_{\sigma}}+M\sum_{j=0}^3
\mathcal{D}^{j,c}_{\max\{2,m\}}\Big).
 \end{align*}  
Then, for $|\alpha|\leq m\leq 4$, we combine the above two estimates to have
\begin{align}\label{dissi-Em}
&\sum_{|\alpha|=0}^{m}\frac{\mathrm{d}}{\mathrm{d}t}
   \big\langle\partial^{\alpha}\big(f^{3,c}_{+}-f^{3,c}_{-}, c(p^0-C_{\mathbf{a}})p_i\sqrt{\mu^c}\big),\frac{\partial^{\alpha}E^{3,c}}{c^2}\big\rangle
   +\frac{\rho_{\mathbf{a}}}{c^2}\big\|\mathbf{a}^3_{+}-\mathbf{a}^3_{-}\big\|^2_{H^m}
   +\frac{3\rho_{\mathbf{a}}}{2c^2}\Big\|E^{3,c}\big\|^2_{H^m}
\nonumber\\
&\qquad\leq \frac{o(1)}{c^2}\|\nabla_xB^{3,c}\|^2_{H^{m-1}}+\frac{C}{c^2}\Big(\|\{I-\mathcal{P}\}f^{3,c}\|_{H^{m+1}_{\sigma}}^2+M\sum_{j=0}^3
\mathcal{D}^{j,c}_{\max\{2,m\}}\Big).
\end{align}

\noindent $\bullet$ Case 2 (Dissipation estimates of the magnetic field): 
Recall that
\begin{align*}
  &\frac{1}{c}\partial_t E^{3,c}-  \nabla_x \times B^{3,c} =- \frac{1}{c}\int_{\mathbb R^3}\frac{cp}{p^0}\sqrt{\mu^c}\left(f^{3,c}_{+}-f^{3,c}_-\right)\, \mathrm{d} p.
\end{align*}
 For $0\leq m\leq 3$, we apply $\partial^{\alpha}$ to the above equation, and take the inner product of the resultant with $\frac{1}{c^2}\nabla_x \times\partial^{\alpha}B^{3,c}(t, x)$ in $L^2_x$ to have
 \begin{align*}
  &-\frac{1}{c^3}\sum_{|\alpha|=0}^m\langle\partial_t \partial^{\alpha}E^{3,c}, \nabla_x \times\partial^{\alpha}B^{3,c}\rangle+\frac{1}{c^2} \|\nabla_x \times \partial^{\alpha}B^{3,c}\|_{H^m}^2 \\
  &\qquad= \frac{1}{c^3}\sum_{|\alpha|=0}^m\Big\langle\int_{\mathbb R^3}\frac{cp}{p^0}\sqrt{\mu^c}\partial^{\alpha}\left(f^{3,c}_{+}-f^{3,c}_-\right)\, \mathrm{d} p, \nabla_x \times\partial^{\alpha}B^{3,c}\Big\rangle\\
 &\qquad\leq \frac{1}{4c^2} \|\nabla_x \times B^{3,c}\|_{H^m}^2+\frac{C}{c^2}\|\{I-\mathcal{P}\}f^{3,c}\|_{H^{m}_{\sigma}}^2.
\end{align*}
Note that $\nabla_x \cdot B^{3,c}=0$, $\partial_t B^{3,c}+ c\nabla_x \times E^{3,c}=0$, and
\begin{align*}
  -c\|\nabla_x \times\partial^{\alpha}E^{3,c}\|^2&=\langle \partial_t\partial^{\alpha}B^{3,c}, \nabla_x \times\partial^{\alpha}E^{3,c}\rangle=
  \langle \partial^{\alpha}E^{3,c}, \nabla_x \times\partial^{\alpha}\partial_tB^{3,c}\rangle\\
  &=\frac{\mathrm{d}}{\mathrm{d}t}\langle \partial^{\alpha}E^{3,c}, \nabla_x \times\partial^{\alpha}B^{3,c}\rangle-\langle \partial^{\alpha}\partial_tE^{3,c}, \nabla_x \times\partial^{\alpha}B^{3,c}\rangle.
  \end{align*}
 We can further obtain 
 \begin{align}\label{dissi-Bm}
  &-\frac{1}{c^3}\sum_{|\alpha|=0}^m\frac{\mathrm{d}}{\mathrm{d}t}\langle \partial^{\alpha}E^{3,c}, \nabla_x \times\partial^{\alpha}B^{3,c}\rangle+\frac{3}{4c^2} \|\nabla_x \times B^{3,c}\|_{H^m}^2 \nonumber\\
   &\qquad\leq \frac{1}{c^2} \|\nabla_x \times E^{3,c}\|_{H^m}^2+\frac{C}{c^2}\|\{I-\mathcal{P}\}f^{3,c}\|_{H^{m}_{\sigma}}^2.
\end{align}
  Now for $1\leq m\leq 4$ in \eqref{dissi-Em}, we multiply \eqref{dissi-Bm} by $\frac{\rho_{\mathbf{a}}}{2}$, and add the resultant to \eqref{dissi-Em} to obtain
 \begin{align}\label{dissi-EBm}
&\frac{\mathrm{d}}{\mathrm{d}t} \mathcal{E}^{3,c,m}_{EB}(t)
   +\frac{\rho_{\mathbf{a}}}{c^2}\big\|\mathbf{a}^3_{+}-\mathbf{a}^3_{-}\big\|^2_{H^m}
   +\frac{\rho_{\mathbf{a}}}{c^2}\Big\|E^{3,c}\big\|^2_{H^m}+\frac{1}{2c^2} \|\nabla_x  B^{3,c}\|_{H^{m-1}}^2 
\nonumber\\
&\qquad\leq \frac{\overline{C}_{0m}}{c^2}\Big(\|\{I-\mathcal{P}\}f^{3,c}\|_{H^{m+1}_{\sigma}}^2+M\sum_{j=0}^3
\mathcal{D}^{j,c}_{\max\{2,m\}}\Big),
\end{align} 
  where 
  \begin{align*}
    \mathcal{E}^{3,c,m}_{EB}(t)= & -\sum_{|\alpha|=0}^{m}
   \big\langle\partial^{\alpha}\big(f^{3,c}_{+}-f^{3,c}_{-}, c(p^0-C_{\mathbf{a}})p_i\sqrt{\mu^c}\big),\frac{\partial^{\alpha}E^{3,c}}{c^2}\big\rangle\\
   &-\frac{\rho_{\mathbf{a}}}{2c^3}\sum_{|\alpha|=0}^{m-1}\langle \partial^{\alpha}E^{3,c}, \nabla_x \times\partial^{\alpha}B^{3,c}\rangle.
  \end{align*}
  
\begin{proposition}\label{EB-15} Under the assumptions in Lemma \ref{grow-low}, we further assume that $c\geq c_3$. Then for $1\leq m\leq 4$, there exist $c-$independent constant $\overline{C}_{0m}>0$ and a functional $\mathcal{E}^{3,c,m}_{EB}(t)$,  which satisfy 
$$|\mathcal{E}^{3,c,m}_{EB}(t)|\leq \frac{\overline{C}_{0m}}{c^2} \big(\|f^{3,c}(t)\|_{H^m}^2+\|E^{3,c}(t)\|_{H^m}^2+\|\nabla_xB^{3,c}(t)\|_{H^{m-1}}^2\big) $$
 such that \eqref{dissi-EBm} holds.
\end{proposition}  

Now we are at the position to deduce closed {\it a priori} estimates. 
We first show the uniform boundness of $\mathcal{E}^{3,c}_{2}(t)$ under the smallness assumption of $M$.

We first choose $\kappa_3$ small enough such that
$$ \kappa_3\Big(C_{4}+\frac{\overline{C}_{01}}{c^2_3}\Big)\leq \frac{1}{4}\delta_0.$$ 
Letting $m=2$ in \eqref{dissicom-aabc} and $m=1$ in \eqref{dissi-EBm}, we make a linear combination of \eqref{2norm}, \eqref{dissicom-aabc} and  \eqref{dissi-EBm} to have
 \begin{align}\label{2norm-bound}
&\frac{\mathrm{d}}{\mathrm{d}t} \Big[\mathcal{E}^{3,c}_{2}(t)+\kappa_3\Big(\mathcal{E}^{3,c,2}_{\mathbf{abc}}(t)
+\mathcal{E}^{3,c,1}_{EB}(t)\Big)\Big]+
\frac{3\delta_0 }{4}\|\{I-\mathcal{P}\}f^{3,c}\|^2_{H^2_{\sigma}}\nonumber\\
&+\kappa_3\Big(\frac{\kappa_2\rho_{\mathbf{a}}}{8C}\big\|\nabla_x\big(\mathbf{a}^3_{+}-\mathbf{a}^3_{-}\big)\big\|^2_{H^{1}}
+\frac{\kappa_2\rho_{\mathbf{a}}}{4C}\big\|\mathbf{a}^3_{+}-\mathbf{a}^3_{-}\big\|^2_{H^{1}}
+\frac{\kappa_1\kappa_2\rho_{\mathbf{a}}}{4}\big\|\nabla_x\big(\mathbf{a}^3_{+}+\mathbf{a}^3_{-}\big)\big\|^2_{H^{1}}\nonumber\\
   & +\frac{\kappa_2}{2}\big\|\nabla_x\mathbf{b}^3\big\|^2_{H^1}+\frac{\kappa_2}{2}\big\|\nabla_x\cdot\mathbf{b}^3\big\|^2_{H^{1}}
+\frac{\rho_{\mathbf{c}}}{2}\big\|\nabla_x\mathbf{c}^3\big\|^2_{H^{1}}+\frac{\rho_{\mathbf{a}}}{c^2}\Big\|E^{3,c}\big\|^2_{H^1}+\frac{1}{2c^2} \|\nabla_x  B^{3,c}\|^2\Big) \nonumber\\
&\qquad\leq \big(o(1)+\sum_{j=0}^3\frac{C\sqrt{M}}{c^j}\big)\big\|f^{3,c}\big\|_{H^2_{\sigma}}^2
+CM\sum_{j=0}^3 \mathcal{D}^{j,c}_{2}(t).  
\end{align}
The choice of $\kappa_3$ implies that
$$\widetilde{\mathcal{E}}^{3,c}_{2}(t):=\mathcal{E}^{3,c}_{2}(t)+\kappa_3\Big(\mathcal{E}^{3,c,2}_{\mathbf{abc}}(t)
+\mathcal{E}^{3,c,1}_{EB}(t)\Big)\backsimeq\mathcal{E}^{3,c}_{2}(t).$$
On the other hand, from \eqref{cons-RVMLi}, one applies Poincaré's inequality to have
\begin{align}\label{L2-abc2}
  &\big\|B^{3,c}-\overline{B}^3\big\|^2\leq C_p \big\|\nabla_x  B^{3,c}\big\|^2,\qquad\|\mathbf{a}^{3}_{\pm}\|^2 \leq C_p\big\|\nabla_x \mathbf{a}^{3}_{\pm}\|^2,\nonumber  \\
  &\|\mathbf{b}^{3}\|^2\leq \big|\int_{\mathbb{T}^3} \mathbf{b}^{3}_+\,\mathrm{d} x\big|^2+C_p\big|\nabla_x \mathbf{b}^{3}\|^2\leq \frac{C}{c^2}\sum_{j=0}^3\big\|E^{j,c}\big\|^2+C_p\big\|\nabla_x \mathbf{b}^{3}\|^2,\nonumber\\
   &\|\mathbf{c}^{3}\|^2\leq \big|\int_{\mathbb{T}^3} \mathbf{c}^{3}_+\,\mathrm{d} x\big|^2+C_p\big\|\nabla_x \mathbf{c}^{3}\big\|^2\leq C\Big(\sum_{j=0}^3\big\|E^{j,c}\big\|^4+C_p^2 \big\|\nabla_x  B^{3,c}\big\|^4\Big)+C_p\big\|\nabla_x \mathbf{c}^{3}\big\|^2,
\end{align} 
  where $C_p$ is the Poincaré constant. \eqref{2norm-bound} and \eqref{L2-abc2} implies that, for $M$ sufficiently small,
\begin{align}\label{nowei-2ed3}
  \frac{\mathrm{d}}{\mathrm{d}t}\widetilde{\mathcal{E}}^{3,c}_{2}(t)+ \mathcal{D}^{3,c}_{2}(t)\leq \frac{CM \mathcal{E}^{3,c}_2(t)}{(1+t)^{9}}+ CM\sum_{j=0}^2 \mathcal{D}^{j,c}_{2}(t), 
    \end{align}
where $ \mathcal{D}^{3,c}_{2}(t)$ is defined in \eqref{nowei-ed3} with $m=2$.

Letting $i=m=2$ in \eqref{dissip-ied} and assuming $M$ small enough such that $CM\leq \frac12$ in \eqref{nowei-2ed3}, we combine \eqref{nowei-2ed3} and \eqref{dissip-ied} to have
\begin{align}\label{nowei-2edi3}
  \frac{\mathrm{d}}{\mathrm{d}t}\sum_{j=0}^3\widetilde{\mathcal{E}}^{j,c}_{2}(t)+ \frac12\sum_{j=0}^3\mathcal{D}^{j,c}_{2}(t)\leq  \frac{CM \mathcal{E}^{3,c}_2(t)}{(1+t)^{9}}. 
    \end{align}
     
Similarly, we choose $\kappa_4$ small enough such that
$$ \kappa_4\Big(C_{7}+\frac{\overline{C}_{04}}{c^2_3}\Big)\leq \frac{1}{4}\delta_0.$$ 
Letting $m=5$ in \eqref{dissicom-aabc} and $m=4$ in \eqref{dissi-EBm}, we make a linear combination of \eqref{6norm} in Proposition \ref{6norm00}, \eqref{dissicom-aabc} and  \eqref{dissi-EBm} to have
 \begin{align}\label{6norm-bound}
&\frac{\mathrm{d}}{\mathrm{d}t} \widetilde{\mathcal{E}}^{3,c}_{5}(t)+
\frac{3\delta_0 }{4}\|\{I-\mathcal{P}\}f^{3,c}\|^2_{H^2_{\sigma}}\nonumber\\
&+\kappa_4\Big(\frac{\kappa_2\rho_{\mathbf{a}}}{8C}\big\|\nabla_x\big(\mathbf{a}^3_{+}-\mathbf{a}^3_{-}\big)\big\|^2_{H^{4}}
+\frac{\kappa_2\rho_{\mathbf{a}}}{4C}\big\|\mathbf{a}^3_{+}-\mathbf{a}^3_{-}\big\|^2_{H^{4}}
+\frac{\kappa_1\kappa_2\rho_{\mathbf{a}}}{4}\big\|\nabla_x\big(\mathbf{a}^3_{+}+\mathbf{a}^3_{-}\big)\big\|^2_{H^{4}}\nonumber\\
   & +\frac{\kappa_2}{2}\big\|\nabla_x\mathbf{b}^3\big\|^2_{H^4}+\frac{\kappa_2}{2}\big\|\nabla_x\cdot\mathbf{b}^3\big\|^2_{H^{4}}
+\frac{\rho_{\mathbf{c}}}{2}\big\|\nabla_x\mathbf{c}^3\big\|^2_{H^{4}}
+\frac{\rho_{\mathbf{a}}}{c^2}\big\|E^{3,c}\big\|^2_{H^4}+\frac{1}{2c^2} \|\nabla_x  B^{3,c}\|^2_4\Big) \nonumber\\
&\qquad\leq \big(o(1)+\sum_{j=0}^3\frac{C\sqrt{M}}{c^j}\big)\big\|f^{3,c}\big\|_{H^5_{\sigma}}^2
+CM\sum_{j=0}^3 \mathcal{D}^{j,c}(t),  
\end{align}
where
$$\widetilde{\mathcal{E}}^{3,c}_{5}(t):=\mathcal{E}^{3,c}_{5}(t)+\kappa_4\Big(\mathcal{E}^{3,c,5}_{\mathbf{abc}}(t)
+\mathcal{E}^{3,c,4}_{EB}(t)\Big)\backsimeq\mathcal{E}^{3,c}_{5}(t).$$

Then, for small enough $M$, we collect \eqref{L2-abc2} and \eqref{6norm-bound} to obtain that
there exists a dissipation norm $\mathcal{D}^{3,c}_{5}(t)$, which    is defined in \eqref{nowei-ed3} with $m=5$, such that
\begin{align}\label{nowei-6ed3}
  \frac{\mathrm{d}}{\mathrm{d}t}\widetilde{\mathcal{E}}^{3,c}_{5}(t)+ \mathcal{D}^{3,c}_{5}(t)\leq  CM\sum_{j=0}^3 \mathcal{D}^{j,c}(t). 
    \end{align}
For $M$ small enough such that $CM\leq \frac12$ in \eqref{nowei-2ed3}, we combine \eqref{nowei-6ed3} and \eqref{dissip-ied}  with $m=5$ to have
\begin{align}\label{nowei-6edi3}
  \frac{\mathrm{d}}{\mathrm{d}t}\sum_{j=0}^3\widetilde{\mathcal{E}}^{j,c}_{5}(t)+ \frac12\sum_{j=0}^3\mathcal{D}^{j,c}_{5}(t)\leq  0. 
    \end{align}
  Now we choose another constant $N_0$ large enough such that
 $$\frac{N_0 }{2}\geq C,$$
 where $C$ is the coefficient of $(1+t)^{-(1+\epsilon_{0})}\|\nabla_x^5E^{3,c}\|^2$ in \eqref{6norm-wei}.
 We also choose a constant $N_1$ large enough such that 
$$\frac{N_1}{8}\mathcal{D}^{3,c}_{5}(t) \geq C,$$
 where $C$ is the sum of the coefficients of  $\big\|f^{3,c}\big\|^2_{H^4
 _{\sigma}}$ in \eqref{5norm-wei} and  $\|\nabla_x^5f^{3,c}\|^2_{\sigma}$ in \eqref{6norm-wei}.
Then, for large constants $N_0, N_1$ chosen above and 
 $M$ sufficiently small, we make a linear combination of \eqref{5norm-wei} in Proposition \ref{5norm-wei00}, \eqref{6norm-wei} in Proposition \ref{6norm-wei00}, and $ \big[N_1+N_0(1+t)^{-\epsilon_{0}}\big]\eqref{nowei-6edi3}$ to have
\begin{align} \label{sumnorm-wei}
&\frac{\mathrm{d}}{\mathrm{d}t}\widetilde{\mathcal{E}}^{3,c}(t)+\widetilde{\mathcal{D}}^{3,c}(t)
\leq  CM\sum_{j=0}^2\mathcal{D}^{j,c}(t), 
\end{align}
where
\begin{align*}
  \widetilde{\mathcal{E}}^{3,c}(t):=&\sum_{|\alpha|=0}^4\frac{1}{4C_{\alpha}(t)}\Big\langle \Big(\frac{cp}{p^0}\sqrt{\mu^c}\zeta_0, w_{\ell_3-|\alpha|}^2\partial^{\alpha}f^{3,c}\Big),\Big(\frac{cp}{p^0}\sqrt{\mu^c}\zeta_0, w_{\ell_3-|\alpha|}^2\partial^{\alpha}f^{3,c}\Big)\Big\rangle\nonumber\\
&+\|f^{3,c}(t)\|^2_{H^{4,3}_{w}}+(1+t)^{-\frac{1+\epsilon_{0}}{2}} \|w_{\ell_3-5}\nabla_x^{5}f^{3,c}\|^2+\big[N_1+N_0(1+t)^{-\epsilon_{0}}\big]
\sum_{j=0}^3\widetilde{\mathcal{E}}^{j,c}_{5}(t),
\end{align*}
\begin{align*}
  \widetilde{\mathcal{D}}^{3,c}(t):=&\frac{1}{2}\Big[Y(t)\big\|\langle p\rangle^{\frac{1}{2}}f^{3,c}\big\|^2_{H^{4,3}_{w}}+\frac{1}{5} \|f^{3,c}(t)\|^2_{H^{4,3}_{w,\sigma}}+(1+t)^{-\frac{1+\epsilon_{0}}{2}}\Big( \frac{1+\epsilon_{0}}{2(1+t)}\|w_{\ell_3-5}\nabla_x^{5}f^{3,c}\|^2\nonumber\\
&+Y(t)\big\|\langle p\rangle^{\frac{1}{2}}w_{\ell_3-5}\nabla_x^{5}f^{3,c}\big\|^2+\frac{1}{5} \|w_{\ell_3-5}\nabla_x^{5}f^{3,c}(t)\|^2_{\sigma}\Big)
\Big]\nonumber\\
&+\frac{N_1+N_0 \epsilon_{0}(1+t)^{-(1+\epsilon_{0})}}{2}\sum_{j=0}^3\mathcal{D}^{3,c}_{5}(t)+\frac{N_0 \epsilon_{0}}{2}(1+t)^{-(1+\epsilon_{0})}\widetilde{\mathcal{E}}^{3,c}_{5}(t).
\end{align*}  

Here we used the definitions of $\mathcal{E}^{3,c}(t)$ in \eqref{enerfun-3c} and $\mathcal{D}^{3,c}(t)$ in \eqref{dissirat-3c}, and the facts
\begin{align*}
  \widetilde{\mathcal{E}}^{3,c}(t)\backsimeq \mathcal{E}^{3,c}(t),\qquad \widetilde{\mathcal{D}}^{3,c}(t)\backsimeq \mathcal{D}^{3,c}(t).
\end{align*}  

Letting $i=2$ in the second inequality of \eqref{dissip-ied}, for $M$ sufficient small, we collect \eqref{sumnorm-wei} and \eqref{dissip-ied} together to have
\begin{align} \label{sumnorm-whole}
&\frac{\mathrm{d}}{\mathrm{d}t}\sum_{j=0}^3\widetilde{\mathcal{E}}^{j,c}(t)+\frac{1}{2}\widetilde{\mathcal{D}}^{3,c}(t)
+\frac{1}{2}\sum_{j=0}^2\mathcal{D}^{j,c}(t)
\leq  0.
\end{align}
\eqref{sumnorm-whole} implies that
\begin{align} \label{sumnorm-uni}
&\sum_{j=0}^3\widetilde{\mathcal{E}}^{j,c}(t)\leq  \sum_{j=0}^3\widetilde{\mathcal{E}}^{j,c}(0).
\end{align}

\subsubsection{Time decay of $\mathcal{E}^{3,c}_2(t)$}
In this part, we will derive the time decay of $\mathcal{E}^{3,c}_2(t)$.

 Note that 
 \begin{align*}
   \big\|f^{j,c}\big\|_{H^2}^2\leq& \big\|\langle p\rangle^{\frac{3}{2}}f^{j,c}\big\|_{H^2}^{\frac{1}{2}}\big\|\langle p\rangle^{-\frac{1}{2}}f^{j,c}\big\|_{H^2}^{\frac{3}{2}}\leq C\big\|\langle p\rangle^{\frac{3}{2}}f^{j,c}\big\|_{H^2}^{\frac{1}{2}}\big\|f^{j,c}\big\|_{H^2_{\sigma}}^{\frac{3}{2}}, \qquad 0\leq j\leq 3,\\
   \big\|[E^{3,c}, B^{3,c}]\big\|_{H^2}^2\leq& \Big(c^6\big\|[E^{3,c}, B^{3,c}]\big\|_{H^5}^2\Big)^{\frac{1}{4}}\Big(\frac{1}{c^2}\big\|[E^{3,c}, B^{3,c}]\big\|_{H^1}^2\Big)^{\frac{3}{4}},\\
   \big\|E^{j,c}\big\|_{H^2}^2\leq& C \big\|E^{j,c}\big\|_{H^5}^{\frac{1}{2}}\big\|E^{j,c}\big\|_{H^1}^{\frac{3}{2}}, \qquad 0\leq j\leq 2.
 \end{align*}
Then, from \eqref{sumnorm-uni}, we have
\begin{align*}
  \mathcal{E}^{3,c}_2(t)\leq& C\big[\mathcal{D}^{3,c}_2(t)\big]^{\frac{3}{4}}\Big(\big\|\langle p\rangle^{\frac{3}{2}}f^{j,c}\big\|_{H^2}^2+c^6\big\|[E^{3,c}, B^{3,c}]\big\|_{H^5}^2\Big)^{\frac{1}{4}}
  \leq  Cc^{\frac{3}{2}}\big[\mathcal{D}^{3,c}_2(t)\big]^{\frac{3}{4}}\big[\mathcal{E}^{3,c}(0)\big]^{\frac{1}{4}}.
\end{align*}
This implies
\begin{align}\label{dissi-3ed}
  \mathcal{D}^{3,c}_2(t)\geq &\frac{1}{C } \big[\mathcal{E}^{3,c}_2(t)\big]^{\frac{4}{3}}\big[c^6\mathcal{E}^{3,c}(0)\big]^{-\frac{1}{3}}.
\end{align}
Similarly, for $0\leq j\leq 2$, we can obtain
\begin{align}\label{dissi-jed}
  \mathcal{D}^{j,c}_2(t)\geq &\frac{1}{C } \big[\mathcal{E}^{j,c}_2(t)\big]^{\frac{4}{3}}\big[\mathcal{E}^{j,c}(0)\big]^{-\frac{1}{3}}.
\end{align}
Putting \eqref{nowei-2edi3}, \eqref{dissi-3ed}, and \eqref{dissi-jed} together, we can obtain
\begin{align*}
  \frac{\mathrm{d}}{\mathrm{d}t}\sum_{j=0}^3\widetilde{\mathcal{E}}^{j,c}_{2}(t)+ \frac{1}{C }\Big[\sum_{j=0}^3\widetilde{\mathcal{E}}^{j,c}_{2}(t)\Big]^{\frac{4}{3}}
  \Big[c^6\sum_{j=0}^3\mathcal{E}^{j,c}(0)\Big]^{-\frac{1}{3}}\leq  \frac{CM \mathcal{\widetilde{E}}^{3,c}_2(t)}{(1+t)^{9}}. 
    \end{align*}
This implies
\begin{align}\label{decay-2E}
  \sum_{j=0}^3\mathcal{E}^{j,c}_2(t)\leq & C c^6\sum_{j=0}^3\mathcal{E}^{j,c}(0)(1+t)^{-3}.
\end{align}

\noindent {\it Proof of Theorem \ref{result1RVML}:} Under the {\it a priori} assumption \eqref{apriori-assu}, we have shown \eqref{sumnorm-whole}, which implies \eqref{global} for $t\in [0, \overline{T}]$ since 
\begin{align*}
  \sum_{j=0}^3\widetilde{\mathcal{E}}^{j,c}(t)\backsimeq \sum_{j=0}^3\mathcal{E}^{3,c}(t),\qquad \widetilde{\mathcal{D}}^{3,c}(t)\backsimeq \mathcal{D}^{3,c}(t).
\end{align*}  
While \eqref{decay-2E} implies \eqref{H2-decay} for $t\in [0, \overline{T}]$. With the estimates \eqref{sumnorm-whole} and  \eqref{decay-2E}, we can show that the {\it a priori} assumption \eqref{apriori-assu} will always holds true if $\sum_{j=0}^3\mathcal{E}^{j,c}(0)$ is chosen sufficiently small. 
Then, we can combine the local existence result Proposition \ref{local-exi} and the {\it a priori} estimates \eqref{sumnorm-whole} and  \eqref{decay-2E} to obtain the global existence a unique solution to the system \eqref{mainF3} via the continuation argument. Then, \eqref{global} and \eqref{H2-decay} hold.  Moreover, by similar arguments in \cite{Strain-Guo-CMP-2004}, we can show that $F^{c}_{\pm}(t,x,p)=\mu^{c}(p)+\sqrt{\mu^{c}(p)}\sum_{i=0}^3c^{-i}f^{i,c}_{\pm}(t,x,p)\geq0$.

\section{Classical limit of the RVML system}

\setcounter{equation}{0}

In this section, we are devoted to the proof of our second main result Theorem \ref{result2RVML},
the classical limit from the RVML system \eqref{main1-00} to the VPL system \eqref{main2-00}. Due to the expansions \eqref{expanion0t}, this limit can be reduced to the classical limit from the RVML system \eqref{mainF0}, the reformulated form of \eqref{mainF0-00},  to the reformulated VPL system \eqref{main2}, the reformulated form of \eqref{main2-00}. To start, we define the difference between the solution $[f^{0,c}(t,x,p), E^{0,c}(t,x)]$ to the Cauchy problem \eqref{mainF0} and the solution $[f^{\infty}(t,x,p), E^{\infty}(t,x)]$ to the Cauchy problem \eqref{main2}:
\begin{align*}
  f(t,x,p):=f^{0,c}(t,x,p)-f^{\infty}(t,x,p),\qquad E_f(t,x):=E^{0,c}(t,x)-E^{\infty}(t,x).
\end{align*}
Then $[f(t,x,p), E_f(t,x)]$ satisfies the following Cauchy problem:
\begin{align}\label{mainf-ci0}
&\partial_t f + \frac{cp}{p^0}\cdot \nabla_x f+ E_f\cdot\frac{cp}{p^0}\sqrt{\mu^c}\zeta_0-\frac{1}{2}\zeta_1\frac{cp}{p^0}\cdot \big(E_f f^{0,c}+E^{\infty}f\big)\nonumber\\
&\hspace{0.2cm}+\zeta_1 \big(E_f \cdot \nabla_p f^{0,c}+ E^{0,c}\cdot\nabla_p f\big)+\mathcal{L}f =\varGamma\left(f^{0,c}, f\right)+\varGamma\left(f, f^{\infty}\right)+\mathcal{ Q}_f,\nonumber\\
&\partial_t E_f =- \int_{\mathbb R^3}\frac{cp}{p^0}\sqrt{\mu^c}\left(f_{+}-f_-\right)\, \mathrm{d} p
+\int_{\mathbb R^3}\big(\frac{cp}{p^0}\sqrt{\mu^c}-p\sqrt{\mu^{\infty}}\big)\left(f^{\infty}_{+}-f^{\infty}_-\right)\, \mathrm{d} p, \\\rule{0ex}{1.5em}
& \nabla_x\cdot E_f= \int_{\mathbb R^3}  \sqrt{\mu^c}\left(f_+-f_-\right)\, \mathrm{d}p+
 \int_{\mathbb R^3}  \big(\sqrt{\mu^c}-\sqrt{\mu^{\infty}}\big)\left(f^{{\infty}}_+-f^{{\infty}}_-\right)\, \mathrm{d}p, \qquad
 \nabla_x \times E_f=0,\nonumber\\\rule{0ex}{1.0em}
&f(0,x,p)=f^{0,c}_{0}(x,p),\qquad E_f(0,x)=E_{f,0}(x), \nonumber
\end{align}
where
\begin{align*}
  \mathcal{ Q}_f= & -\big(\frac{cp}{p^0}-p\big)\cdot \nabla_x f^{\infty}+E^{\infty}\cdot\big(\frac{cp}{p^0}\sqrt{\mu^c}-p\sqrt{\mu^{\infty}}\big)\zeta_0
  +\frac{1}{2}\zeta_1\big(\frac{cp}{p^0}-p\big)\cdot E^{\infty} f^{\infty}\\
  &+\zeta_1 \big(\frac{p}{p^0}\times  \sum_{j=0}^2\frac{\bar{B}^{j}}{c^j}\big)\cdot\nabla_p f^{0,c} -\big(\mathcal{L} f^{\infty}-\mathbf{L} f^{\infty}\big)+\big[\varGamma\left(f^{\infty}, f^{\infty}\right)-\Gamma\left(f^{\infty}, f^{\infty}\right)\big].
\end{align*}

\noindent {\it Proof of Theorem \ref{result2RVML}:} 
Note that for $E^{\infty}=-\nabla_x\phi^{\infty}$, as in \cite{Guo-JAMS-2012}, 
\begin{align*}
  \mathrm{e}^{\pm\phi^{\infty}} f_{\pm}\big(\frac{cp}{p^0}\cdot \nabla_x f_{\pm} \pm \frac{1}{2} \frac{cp}{p^0}\cdot \nabla_x\phi^{\infty}f_{\pm}\big)=\frac{1}{2} \frac{cp}{p^0}\cdot \nabla_x\big(\phi^{\pm\infty}f^2_{\pm}\big),
\end{align*}
and 
\begin{align*}
  &\sum_{\pm}\big|\big\langle \frac{1}{2}\frac{cp}{p^0}\cdot E_{f,\pm}f^{0,c}, \mathrm{e}^{\pm\phi^{\infty}} f_{\pm}\big\rangle\big|\leq \big\|\langle p\rangle f^{0,c}\big\|_{\infty}\big(\|f\|^2+\|E_f\|^2\big)\\
  &\qquad \leq \big\| f^{0,c}\big\|^{\frac{1}{2}}_{\infty}\big\|\langle p\rangle^2 f^{0,c}\big\|^{\frac{1}{2}}_{\infty}\big(\|f\|^2+\|E_f\|^2\big)\leq \mathrm{e}^{- C_{\ell_0}t^{1/3}}\sqrt{\mathcal{E}^{0,c}(0)}\big(\|f\|^2+\|E_f\|^2\big),
   \end{align*} 
  \begin{align*}
  &\sum_{\pm}\big|\big\langle \zeta_1 E_{f,\pm} \cdot \nabla_p f^{0,c}_{\pm}, \mathrm{e}^{\pm\phi^{\infty}} f_{\pm}\big\rangle\big|\leq\mathrm{e}^{\|\phi^{\infty}\|_{\infty}}\big|\big\langle \zeta_1 E_f f^{0,c}, \nabla_pf\big\rangle\big|\\
  &\qquad \leq \mathrm{e}^{\|\phi^{\infty}\|_{\infty}}\big|\big\langle \zeta_1 E_f f^{0,c}, \nabla_p\mathcal{P}f\big\rangle\big|+\big|\big\langle \zeta_1 E_f f^{0,c}, \nabla_p\{I-\mathcal{P}\}f\big\rangle\big|\\
  &\qquad \leq o(1) \|\{I-\mathcal{P}\}f\|_{\sigma}^2+C\Big[\big\|f^{0,c}\big\|_{\infty}\big(\|f\|^2+\|E_f\|^2\big)
  +\big\|\langle p\rangle^{\frac{3}{2}}f^{0,c}\big\|^2_{\infty}\|E_f\|^2\Big]\\
  &\qquad \leq o(1) \|\{I-\mathcal{P}\}f\|_{\sigma}^2+C\big[\mathrm{e}^{-C_{\ell_0}t^{1/3}}\sqrt{\mathcal{E}^{0,c}(0)}
  +\mathcal{D}^{0,c}(t)\big]\big(\|f\|^2+\|E_f\|^2\big),
\end{align*}
\begin{align*}
&\sum_{\pm}\big|\big\langle  \varGamma_{\pm}\left(f^{0,c}, f\right)+\varGamma_{\pm}\left(f, f^{\infty}\right), \mathrm{e}^{\pm\phi^{\infty}} f_{\pm}\big\rangle\big|\\
  &\qquad \leq o(1) \|\{I-\mathcal{P}\}f\|_{\sigma}^2+ C\Big[\Big(\big\|f^{0,c}\big\|_{H^2}^2+\Big\|\big|f^{\infty}\big|_{L^2}\Big\|_{L^{\infty}_x}\Big)
  \|f\|^2_{\sigma}\\
  &\qquad\quad+\Big(\big\|f^{0,c}\big\|_{H^2_{\sigma}}^2
  +\Big\|\big|f^{\infty}\big|_{\sigma}\Big\|_{L^{\infty}_x}\Big)\|f\|^2\Big]\\
  &\qquad \leq \Big[o(1) +C\big(\big\|f^{0,c}\big\|_{H^2}^2+\big\|f^{\infty}\big\|_{H^2}^2\big)\Big]
  \|\{I-\mathcal{P}\}f\|_{\sigma}^2\\
  &\qquad\quad+ C\Big(\big\|f^{\infty}\big\|_{H^1}\big\|f^{\infty}\big\|_{H^2}
  +\big\|\langle p\rangle^{\frac{3}{2}}\{I-{\bf{P}}\}f^{\infty}\big\|_{H^2_{\bm{\sigma}}}^2\big)\|f\|^2
  \\&\qquad \leq \Big[o(1) +C\big(\mathcal{E}^{0,c}(0)+\mathcal{E}^{\infty}_{2,2}(0)\big)\Big]\|\{I-\mathcal{P}\}f\|_{\sigma}^2\\
 &\qquad\quad +C\big[\mathrm{e}^{-C_{\ell_0}t^{1/3}}\mathcal{E}^{0,c}(0)
  +(1+t)^{-2}\mathcal{E}^{\infty}_{2,2}(0)+\mathcal{D}^{0,c}(t)+\mathcal{D}^{\infty}_{2,2}(t)\big]\|f\|^2
\end{align*}
from \eqref{tdecay-infty}, \eqref{nonlin-w} with $\vartheta=\ell=0$, and \eqref{decay-f0}.

On the other hand, from \eqref{tdecay-infty} in Proposition \ref{guo-VPL} with $l\geq 16$, \eqref{diff-mu1} in Lemma \ref{mu-ci}, Lemma \ref{diff-LL}, and Lemma \ref{diff-GG}, we have
\begin{align*}
  &\sum_{\pm}\big|\big\langle  \mathcal{ Q}_{f,\pm}, \mathrm{e}^{\pm\phi^{\infty}} f_{\pm}\big\rangle\big|\\
  &\qquad\leq \frac{C\|f\|}{c^2}\Big(\big\|\langle p\rangle^{3}\nabla_xf^{\infty}\big\|+\big\|E^{\infty}\big\|\Big)+\frac{C}{c}\big\|\langle p\rangle^{\frac{5}{2}} f^{0,c}\big\|\|\langle p\rangle^{-\frac{3}{2}} \nabla_pf\|\\
&\qquad+\frac{C}{c^2}\big\|\|f^{\infty}\|_{H^1_{p,8}}\big\|_{L^2_x}\|f\|_{\sigma}+\frac{C}{c^2}
  \big\|\|f^{\infty}\|^2_{H^1_{p,8}}\big\|_{H^1_x}\|f\|_{\sigma}\\
 &\qquad \leq\frac{C_l\|f\|}{c^2}(1+t)^{-5}\sqrt{\mathcal{E}^{\infty}_{2,6}(0)}
  +\frac{C_{\ell_0}}{c}\mathrm{e}^{-C_{\ell_0}t^{1/3}}\sqrt{\mathcal{E}^{0,c}(0)}\big(\|f\|+\|\{I-\mathcal{P}\}f\|_{\sigma}\big)\\
  &\qquad+\frac{C_l}{c^2}(1+t)^{-15}\sqrt{\mathcal{E}^{\infty}_{2,16}(0)}
  \big(\|f\|+\|\{I-\mathcal{P}\}f\|_{\sigma}\big)\\
  &\qquad +\frac{C_l}{c^2}\big(\|f^{\infty}\|\|f^{\infty}\|_{H^2}+\|\langle p\rangle^{10}\{I-\mathbf{P}\}f^{\infty}\|^2_{H^1_{\bm{\sigma}}}\big)\|f\|\\
   &\qquad +\frac{C_l}{c^2}\big(\|f^{\infty}\|\|f^{\infty}\|_{H^2}+\|\{I-\mathbf{P}\}f^{\infty}\|_{\tilde{H}^2_8}\|\langle p\rangle^{10}\{I-\mathbf{P}\}f^{\infty}\|_{H^1_{\bm{\sigma}}}\big)\|\{I-\mathcal{P}\}f\|_{\sigma}\\
 &\qquad \leq o(1)\|\{I-\mathcal{P}\}f\|_{\sigma}^2+C_{\ell_0,l}
 \Big[(1+t)^{-5}+\mathcal{D}^{\infty}_{2,10}(t)\Big]\|f\|^2\\
 &\qquad+\frac{C_{\ell_0,l}}{c^2}
 \Big[(1+t)^{-5}\big[\mathcal{E}^{\infty}_{2,16}(0)+\mathcal{E}^{0,c}(0)\big]
 +\mathcal{D}^{\infty}_{2,10}(t)\Big].
\end{align*}  
Then, we take the inner products of the equation of $f(t, x, p)$ with $\mathrm{e}^{\pm\phi^{\infty}} f_{\pm}(t, x, p)$ and the equation of $E_f(t, x)$ with $E_f(t, x)$ in \eqref{mainf-ci0} respectively, and use \eqref{coerci0} and the smallness of $\mathcal{E}^{0,c}(0)$ to obtain
\begin{align*}
  &\frac{1}{2}\frac{\mathrm{d}}{\mathrm{d}t}\sum_{\pm}\big(\big\|\mathrm{e}^{\pm\phi^{\infty}/2}f_{\pm}\big\|^2
  +\big\|\mathrm{e}^{\pm\phi^{\infty}/2}E_{f,\pm}\big\|^2\big) +\frac{\delta_0}{2}\|\{I-\mathcal{P}\}f\|_{\sigma}^2 \nonumber\\
  &\qquad \leq  C_{\ell_0,l}\Big[\mathrm{e}^{-C_{\ell_0}t^{1/3}}\sqrt{\mathcal{E}^{0,c}(0)}
  +\mathcal{D}^{0,c}(t)+(1+t)^{-2}+\mathcal{D}^{\infty}_{2,10}(t)\Big]\big(\|f\|^2+\|E_f\|^2\big)\nonumber\\
    &\qquad +\frac{C_{\ell_0,l}}{c^2}
 \Big[(1+t)^{-2}\big[\mathcal{E}^{\infty}_{2,16}(0)+\mathcal{E}^{0,c}(0)\big]+\mathcal{D}^{\infty}_{2,10}(t)\Big].\nonumber
\end{align*}
From \eqref{tdecay-infty}, \eqref{decay-f0}, and the Gronwall Lemma as in \cite[Lemma 5]{Guo-JAMS-2012}, we can further obtain
\begin{align}\label{fe-ci}
&\|f(t)\|^2+\|E_f(t)\|^2+\int_{0}^{t}\|\{I-\mathcal{P}\}f(s)\|_{\sigma}^2\,\mathrm{d}s\\
&\qquad\leq
C_{\ell_0,l}\big(\|f_0\|^2+\|E_{f,0}\|^2\big)+
\frac{C_{\ell_0,l}}{c^2}
 \big[\mathcal{E}^{\infty}_{2,16}(0)+\mathcal{E}^{0,c}(0)\big].\nonumber
\end{align}
Combing \eqref{fe-ci}, the expansions in \eqref{expanion0t} and the estimates \eqref{diff-mu0}, and \eqref{diff-mu1} for the difference between Maxwellians $\mu^{c}$ and $\mu^{\infty}$, we get \eqref{c-limit}. 

\section{Appendix}

\setcounter{equation}{0}
This section is composed of two parts. In the first part, we list important properties and estimates of the modified Bessel functions. In the second part, we construct coefficients $[f^{i,c}(t,x, p), E^{i,c}(t,x)]$  with $i=0, 1, 2$ in the expansions \eqref{expanion0t} for the RVML system \eqref{main1-00}.

\subsection{Modified Bessel functions and properties}

In this subsection, we recall expressions of the modified Bessel functions, their basic properties, and several important estimates with respect to these functions, which are proved in \cite{Ruggeri-Xiao-Zhao-ARMA-2021}. Now we give the modified Bessel functions and collect their basic properties:
\begin{lemma}\cite{Groot-Leeuwen-Weert-1980,Oliver-1974,Wa}\label{def-pro}
	Let
	$K_j(\gamma)$ be the Bessel function defined by
	\begin{equation}\label{defini-kj} K_j(\gamma)=\frac{(2^j)j!}{(2j)!}\frac{1}{\gamma^j}\int_{\lambda=\gamma}^{\lambda=\infty}\mathrm{e}^{-\lambda}(\lambda^2-\gamma^2)^{j-1/2}d\lambda,\quad(j\geq0).
	\end{equation}
	Then the following identities hold:
	\begin{eqnarray} &K_j(\gamma)=\frac{2^{j-1}(j-1)!}{(2j-2)!}\frac{1}{\gamma^j}\int_{\lambda=\gamma}^{\lambda
=\infty}\mathrm{e}^{-\lambda}\lambda(\lambda^2-\gamma^2)^{j-3/2}d\lambda,\quad(j>0),\label{defini-kj0}\\
	&K_{j+1}(\gamma)=2j\frac{K_j(\gamma)}{\gamma}+K_{j-1}(\gamma),\quad(j\geq1), \label{transform}\\
	&K_j(\gamma)<K_{j+1}(\gamma),\quad (j\geq0), \nonumber
	\end{eqnarray}
	and
	\begin{eqnarray}
	&\frac{d}{d\gamma}\left(\frac{K_j(\gamma)}{\gamma^j}\right)=-\frac{K_{j+1}(\gamma)}{\gamma^j},
\quad(j\geq0),\nonumber\\
	&\displaystyle K_{j}(\gamma)=\sqrt{\frac{\pi}{2\gamma}}\mathrm{e}^{-\gamma}\left(\gamma_{j,n}(\gamma)\gamma^{-n}
+\sum_{m=0}^{n-1}A_{j,m}\gamma^{-m}\right),\quad(j\geq0,~n\geq1),\label{remainder}
	\end{eqnarray}
	where expressions of the coefficients in (\ref{remainder}) are
	
	\begin{align*} \label{coefficient}
	\begin{aligned}
	&A_{j,0}=1\\
	&A_{j,m}=\frac{(4j^2-1)(4j^2-3^2)\cdots(4j^2-(2m-1)^2)}{m!8^m},\quad(j\geq0,~m\geq1), \\
	&|\gamma_{j,n}(\gamma)|\leq2\mathrm{e}^{[j^2-1/4]\gamma^{-1}}|A_{j,n}|,\quad(j\geq0,~n\geq1).
	\end{aligned}
	\end{align*}
	%
	
\end{lemma}

For later use, we also need two different estimates:
\begin{lemma}\cite{Ruggeri-Xiao-Zhao-ARMA-2021}\label{K01p} Let $\gamma\in(\sqrt{2}, \infty)$. Then $\frac{K_0(\gamma)}{K_1(\gamma)}$ satisfies:
	\begin{equation*}
	1-\frac{1}{2\gamma}\leq\frac{K_0(\gamma)}{K_1(\gamma)}\leq 1-\frac{1}{2\gamma}+\frac{3}{8\gamma^2}+\frac{3}{16\gamma^3}.
	\end{equation*}
	Moreover, for $\gamma\in(2, \infty)$, it holds that
	\begin{align} \label{acurate}
	\begin{aligned}
	&\frac{K_0(\gamma)}{K_1(\gamma)}\geq 1-\frac{1}{2\gamma}+\frac{3}{8\gamma^2}-\frac{3}{8\gamma^3}+\frac{63}{128\gamma^4}-\frac{31}{20\gamma^5}, \\
	&\frac{K_0(\gamma)}{K_1(\gamma)}\leq 1-\frac{1}{2\gamma}+\frac{3}{8\gamma^2}-\frac{3}{8\gamma^3}+\frac{63}{128\gamma^4}+\frac{7}{8\gamma^5}.
	\end{aligned}
	\end{align}
\end{lemma}

\subsection{Construction of the coefficients }
In this part, we are devoted to the construction of coefficients $[f^{i,c}(t,x, p), E^{i,c}(t,x)]$  with $i=0, 1, 2$ in the expansions \eqref{expanion0t} for the RVML system \eqref{main1-00}.

To start, we provide the construction of $[f^{0,c}(t,x, p), E^{0,c}(t,x)]$,  the unique solution to the Cauchy problem \eqref{mainF0}.


\begin{proposition}\label{RVPL0}
Assume that $\ell_0\geq 14$, $[f^{0,c}_{0}(x,p),E^{0,c}_0(x),B^{0,c}_0(x)]$  satisfies conservation laws \eqref{cons-RVML00}.
 There exist a large constant $c_0>0$ and a sufficiently small constant $M_0 > 0$ such that if  $c\geq c_0$,  $\mathcal{E}^{0, c}(0)\leq M_0 $, 
then the Cauchy problem \eqref{mainF0} admits a unique global solution $[f^{0,c}(t,x,p), E^{0,c}(t,x)]$ satisfying 
\begin{align}\label{dissip-0ed}
\begin{aligned}
\frac{\mathrm{d}}{\mathrm{d}t}\widetilde{\mathcal{E}}^{0,c}_{m}(t)
+\mathcal{D}^{0,c}_{m}(t)\leq& 0,\qquad \widetilde{\mathcal{E}}^{0,c}_{m}(t)\backsimeq \mathcal{E}^{0,c}_{m}(t),\qquad 2\leq m\leq 5,\\
\frac{\mathrm{d}}{\mathrm{d}t}\widetilde{\mathcal{E}}^{0,c}(t)
+\mathcal{D}^{0,c}(t)\leq& 0,\qquad \widetilde{\mathcal{E}}^{0,c}_{2,2}(t)\backsimeq \mathcal{E}^{0,c}_{2,2}(t).
\end{aligned}
\end{align}
 Furthermore, it holds that
\begin{align}\label{decay-f0}
\begin{aligned}
\mathcal{E}^{0,c}(t)+\int_0^t\mathcal{D}^{0,c}(\tau) \,\mathrm{d}\tau\leq& C_{\ell_0}\mathcal{E}^{0,c}(0), \qquad \mbox{and}\\
\mathcal{E}^{0,c}(t)\leq& \mathrm{e}^{-C_{\ell_0}(1+t)^{1/3}}\mathcal{E}^{0,c}(0).
\end{aligned}
\end{align}
Here $\mathcal{E}^{0,c}_{m}(t), \mathcal{E}^{0,c}(t)$ and $\mathcal{D}^{0,c}_{m}(t), \mathcal{D}^{0,c}(t)$ are given in \eqref{enerfun-ic} and \eqref{dissirat-ic} with $i=0$ respectively, $\mathcal{E}^{i,c}_m(t)$ is defined in \eqref{nowei-edi}, $C_{\ell_0}$ is some constant independent of $c\geq c_0$.

\end{proposition}

\begin{proof}
Note that $[f^{0,c}(t,x, p), E^{0,c}(t,x)]$ satisfies a relativistic version VPL system with an additional linear term
$$\zeta_1 \frac{p}{p^0}\times  \sum_{j=0}^2\frac{\bar{B}^{j}}{c^j}\cdot\nabla_p f^{0,c},$$
where $\bar{B}^0, \bar{B}^1$, and $\bar{B}^2$ are constant magnetic fields. Since we will construct a solution to the Cauchy problem \eqref{mainF0} in a norm space $\mathcal{E}^{0,c}(t)$, which is given in \eqref{enerfun-ic} with $i=0$, this additional term always disappears in our energy estimates and doesn't make essential changes to the  dissipation estimates  of the macroscopic part $\mathcal{P}f^{0,c}$ and the electric field $E^{0,c}$. 
%
The energy estimates in \eqref{dissip-0ed}, \eqref{decay-f0} and the proof of global existence of the Cauchy problem \eqref{mainF0} can be established similarly as in \cite{Guo-JAMS-2012}. For brevity, we omit the details  and only prove the time decay estimate in \eqref{decay-f0} via the interpolation techniques introduced in \cite{XXZ-JFA-2017}.

To this end, we first denote
\begin{eqnarray*}
\mathcal{E}^{0,c}(t)=\int_{\mathbb{R}^3}g^2(t,p)\,\mathrm{d}p.    
\end{eqnarray*}
From the definition of $\mathcal{D}^{0,c}(t)$, we have
\begin{align*}
  C_{\ell_0}\int_{\mathbb{R}^3}\langle p\rangle^{-1}g^2(t,p) \,\mathrm{d}p+C_{\ell_0}(1+t)^{-1}\ln ^{-2}(\mathrm{e}+t)\int_{\mathbb{R}^3}\langle p\rangle  g^2(t,p) \,\mathrm{d}p\lesssim \mathcal{D}^{0,c}(t).
  \end{align*}
We further use simple interpolation inequality w.r.t. the power of $p$ to have 
\begin{align*}
  (1+t)^{-\frac{1}{2}}\ln ^{-1}(\mathrm{e}+t)\int_{\mathbb{R}^3}  g^2(t,p) \,\mathrm{d}p\lesssim \mathcal{D}^{0,c}(t). 
\end{align*}
Then we use \eqref{dissip-0ed} to have
\begin{align*}
\frac{\mathrm{d}}{\mathrm{d}t}\int_{\mathbb{R}^3}g^2(t,p)\,\mathrm{d}p+C_{\ell_0}(1+t)^{-\frac{1}{2}}\ln ^{-1}(\mathrm{e}+t)\int_{\mathbb{R}^3}  g^2(t,p) \,\mathrm{d}p\lesssim 0
\end{align*}
which implies 
\begin{align*}
\int_{\mathbb{R}^3}g^2(t,p)\,\mathrm{d}p\leq&  \exp\left\{-C_{\ell_0}\int_0^t(1+s)^{-\frac{1}{2}}\ln ^{-1}(\mathrm{e}+s)\,\mathrm{d}s\right\}\int_{\mathbb{R}^3}g^2(0,p)\,\mathrm{d}p\\
\leq& \exp\left\{-C_{\ell_0}(1+t)^{\frac{1}{2}}\ln ^{-1}(\mathrm{e}+t)\right\}\int_{\mathbb{R}^3}g^2(0,p)\,\mathrm{d}p.
\end{align*}  
This gives $\eqref{decay-f0}_2$.

\end{proof}

Based on Proposition \ref{RVPL0}, we can successively solve the Cauchy problem \eqref{mainFi} for $i=1, 2$.

\begin{proposition}\label{RVPLi}
For $i=1, 2$, assume that $\ell_i\geq 14-2i$, $M_j$ with $0\leq j\leq i-1$ are small enough, $[f^{i,c}_{0}(x,p),E^{i,c}_0(x),B^{i,c}_0(x)]$  satisfies conservation laws \eqref{cons-RVML0i}.
 There exist a large constant $c_i>0$ and a sufficiently small constant $M_i > 0$ such that if $c\geq c_i$, 
 $$\sum_{j=0}^i\mathcal{E}^{j, c}(0)\leq M_i, $$ 
then the Cauchy problem \eqref{mainFi} admits a unique global solution $[f^{i,c}(t,x,p), E^{i,c}(t,x)]$ satisfying 
\begin{align}\label{dissip-ied}
\begin{aligned}
\frac{\mathrm{d}}{\mathrm{d}t}\sum_{j=0}^i\widetilde{\mathcal{E}}^{j,c}_{m}(t)
+\sum_{j=0}^i\mathcal{D}^{j,c}_{m}(t)\leq& 0,\qquad \widetilde{\mathcal{E}}^{j,c}_{m}(t)\backsimeq \mathcal{E}^{j,c}_{m}(t),\qquad 2\leq m\leq 5,\\
\frac{\mathrm{d}}{\mathrm{d}t}\sum_{j=0}^i\widetilde{\mathcal{E}}^{j,c}(t)
+\sum_{j=0}^i\mathcal{D}^{j,c}(t)\leq& 0,\qquad \widetilde{\mathcal{E}}^{j,c}(t)\backsimeq \mathcal{E}^{j,c}(t).
\end{aligned}
\end{align}
Furthermore, it holds that
\begin{align}\label{decay-fi}
\begin{aligned}
\sum_{j=0}^i\mathcal{E}^{j,c}(t)+\int_0^t\sum_{j=0}^i\mathcal{D}^{j,c}(\tau) \,\mathrm{d}\tau\leq& C_{\ell_j}\sum_{j=0}^i\mathcal{E}^{j,c}(0), \qquad \mbox{and}\\
\mathcal{E}^{i,c}_5(t)\leq& \mathrm{e}^{-C_{\ell_j}(1+t)^{1/3}}\sum_{j=0}^i\mathcal{E}^{j,c}(0).
\end{aligned}
\end{align}
Here $ \mathcal{E}^{i,c}(t)$ and $\mathcal{D}^{i,c}(t)$ are given in \eqref{enerfun-ic} and \eqref{dissirat-ic} respectively, $\mathcal{E}^{i,c}_m(t)$ and $\mathcal{D}^{i,c}_m(t)$ are defined in \eqref{nowei-edi}, and $C_{\ell_j}$ is some constant independent of $c\geq c_i$.

\end{proposition}

%
%

\bigbreak
\begin{center}
{\bf Acknowledgment}
\end{center}
The research of LBH was supported by National Natural Science Foundation of China  (Project No.~11771236) and New Cornerstone Investigator Program (Project No.~100001127). The research of YJL was supported by the National Natural Science Foundation of China (Project No.~12171176). The research of QHX was supported by the National Natural Science Foundation of China (Project No.~12271506) and the National Key Research and Development Program of China (Project No.~2020YFA0714200).

\medskip
\noindent\textbf{Data Availability Statement:}
Data sharing is not applicable to this article as no datasets were generated or analysed during the current study.

\noindent\textbf{Conflict of Interest:}
The authors declare that they have no conflict of interest.

\end{document}